\documentclass{aims}
\usepackage{amsmath}
  \usepackage{paralist}
  \usepackage{algorithm}
  \usepackage{amssymb}
  \usepackage{pifont}
  \usepackage{amsfonts,amssymb,listings,txfonts,multirow}
  \usepackage{fancyhdr}
  \usepackage{soul,color,xcolor}
  \usepackage{epstopdf}
  \usepackage{graphics} 
  \usepackage{epsfig} 
\usepackage{graphicx}  \usepackage{epstopdf}
 \usepackage[colorlinks=true]{hyperref}
 \usepackage{algorithmic}
   \usepackage[mathscr]{eucal}
   \usepackage[capitalize]{cleveref}
   \usepackage{numcompress}
\usepackage[caption=false]{subfig}
\usepackage{multirow}
\usepackage{morefloats}
\hypersetup{urlcolor=blue, citecolor=red}

\def\currentvolume{X}
 \def\currentissue{X}
  \def\currentyear{200X}
   \def\currentmonth{XX}
    \def\ppages{X--XX}

 \crefformat{equation}{\textup{#2(#1)#3}}
\crefrangeformat{equation}{\textup{#3(#1)#4--#5(#2)#6}}
\crefmultiformat{equation}{\textup{#2(#1)#3}}{ and \textup{#2(#1)#3}}
{, \textup{#2(#1)#3}}{, and \textup{#2(#1)#3}}
\crefrangemultiformat{equation}{\textup{#3(#1)#4--#5(#2)#6}}%
{ and \textup{#3(#1)#4--#5(#2)#6}}{, \textup{#3(#1)#4--#5(#2)#6}}{, and \textup{#3(#1)#4--#5(#2)#6}}

\Crefformat{equation}{#2Equation~\textup{(#1)}#3}
\Crefrangeformat{equation}{Equations~\textup{#3(#1)#4--#5(#2)#6}}
\Crefmultiformat{equation}{Equations~\textup{#2(#1)#3}}{ and \textup{#2(#1)#3}}
{, \textup{#2(#1)#3}}{, and \textup{#2(#1)#3}}
\Crefrangemultiformat{equation}{Equations~\textup{#3(#1)#4--#5(#2)#6}}%
{ and \textup{#3(#1)#4--#5(#2)#6}}{, \textup{#3(#1)#4--#5(#2)#6}}{, and \textup{#3(#1)#4--#5(#2)#6}}

\crefdefaultlabelformat{#2\textup{#1}#3}

\theoremstyle{definition}

\newtheorem{remark}{Remark}

\title[Two Stage Continuous Domain Regularization for Piecewise Constant Image Restoration]{Two Stage Continuous Domain Regularization for Piecewise Constant Image Restoration}

\author[Cai, Choi, and Zhou]{ }

 \email{jfcai@ust.hk}
 \email{jaycjk@tongji.edu.cn}
 \email{zhouchenchuan@tongji.edu.cn}

\thanks{The research of the first author is supported by the Hong Kong Research Grants Council (HKRGC) GRF 16306319 and 16309220. The researches of the second and third authors are supported by the National Natural Science Foundation of China (NSFC) Youth Program 11901436, the Shanghai Science and Technology Committee 20JC1413500, and the Fundamental Research Funds for the Central Universities.}

\thanks{$^*$ Corresponding author: Jae Kyu Choi}

\usepackage{amsopn}

\newcommand{\argmin}{\operatornamewithlimits{argmin}}

\DeclareMathAlphabet{\mathpzc}{OT1}{pzc}{m}{it}

\def\R{{\mathbb R}}
\def\C{{\mathbb C}}
\def\N{{\mathbb N}}

\def\Z{{\mathbb Z}}

\def\KK{{\mathbb K}}

\def\MM{{\mathbb M}}

\def\OO{{\mathbb O}}

\def\rank{\mathrm{rank}}

\def\rd{\mathrm{d}}
\def\Om{\Omega}

\def\f{\frac}
\def\p{\partial}

\def\na{\nabla}
\def\la{\langle}
\def\ra{\rangle}

\def\a{{\boldsymbol a}}

\def\bb{{\boldsymbol b}}

\def\bc{{\boldsymbol c}}

\def\bsd{{\boldsymbol d}}
\def\e{{\boldsymbol e}}

\def\bsf{{\boldsymbol f}}

\def\bi{{\mathbf i}}

\def\bk{{\boldsymbol k}}
\def\bsl{{\boldsymbol l}}

\def\bm{\boldsymbol{m}}

\def\bu{\boldsymbol{u}}
\def\bsv{\boldsymbol{v}}

\def\w{{\boldsymbol w}}

\def\x{\boldsymbol{x}}

\def\bsA{{\boldsymbol A}}

\def\bB{{\boldsymbol B}}
\def\bC{{\boldsymbol C}}

\def\bsH{{\boldsymbol H}}

\def\bsI{{\boldsymbol I}}

\def\bsS{{\boldsymbol S}}
\def\bU{{\boldsymbol U}}
\def\bV{{\boldsymbol V}}

\def\bX{{\boldsymbol X}}
\def\bY{{\boldsymbol Y}}

\def\bsS{{\boldsymbol S}}
\def\bZ{{\boldsymbol Z}}
\def\mA{{\mathcal A}}

\def\mC{{\mathcal C}}
\def\mD{{\mathcal D}}

\def\mH{{\mathcal H}}
\def\mI{{\mathcal I}}

\def\mP{{\mathcal P}}

\def\mR{{\mathcal R}}
\def\mS{{\mathcal S}}
\def\mT{{\mathcal T}}

\def\mV{{\mathcal V}}
\def\mW{{\mathcal W}}

\def\bmA{{\boldsymbol \mA}}

\def\bmD{{\boldsymbol \mD}}

\def\bmH{{\boldsymbol \mH}}
\def\bmI{{\boldsymbol \mI}}

\def\bmP{{\boldsymbol \mP}}

\def\bmR{{\boldsymbol \mR}}
\def\bmS{{\boldsymbol \mS}}
\def\bmT{{\boldsymbol \mT}}

\def\bmW{{\boldsymbol \mW}}

\def\msF{{\mathscr F}}

\def\msH{{\mathscr H}}

\def\msN{{\mathscr N}}

\def\msR{{\mathscr R}}

\def\dde{\boldsymbol{\delta}}
\def\gga{\boldsymbol{\gamma}}

\def\bvphi{\boldsymbol{\varphi}}

\def\bi{\begin{itemize}} \def\ei{\end{itemize}}
\def\be{\begin{eqnarray*}}
\def\ee{\end{eqnarray*}}

\def\0{{\mathbf 0}}

\newcommand{\beq}{\begin{equation}}
\newcommand{\eeq}{\end{equation}}

\def\xxi{{\boldsymbol{\xi}}}

\def\zze{{\boldsymbol\zeta}}

\def\Ga{{\boldsymbol\Gamma}}
\def\Ph{\boldsymbol{\Phi}}
\def\Sig{\boldsymbol{\Sigma}}

\def\wt{\widetilde}
\def\wh{\widehat}

\def\Na{\boldsymbol \nabla}

\newcommand{\eps}{\varepsilon}
\def\la{\langle}
\def\ra{\rangle}

\def\XXint#1#2#3{{\setbox0=\hbox{$#1{#2#3}{\int}$ }
\vcenter{\hbox{$#2#3$ }}\kern-.55\wd0}}

\begin{document}
\maketitle

\centerline{\scshape Jian-Feng Cai}
\medskip
{\footnotesize
	\centerline{Department of Mathematics, Hong Kong University of Science and Technology}
	\centerline{Clearwater Bay Kowloon, Hong Kong}
}

\medskip

\centerline{\scshape Jae Kyu Choi$^*$}
\medskip
{\footnotesize
	\centerline{School of Mathematical Sciences, Tongji University}
	\centerline{Shanghai 200092, People's Republic of China}
}

\medskip

\centerline{\scshape Chenchuan Zhou}
\medskip
{\footnotesize
	\centerline{School of Mathematical Sciences, Tongji University}
	\centerline{Shanghai 200092, People's Republic of China}
}

\bigskip

 \centerline{(Communicated by the associate editor name)}

\begin{abstract}
The finite-rate-of-innovation (FRI) framework which corresponds a signal/image to a structured low-rank matrix is emerging as an alternative to the traditional sparse regularization. This is because such an off-the-grid approach is able to alleviate the basis mismatch between the true support in the continuous domain and the discrete grid. In this paper, we propose a two-stage off-the-grid regularization model for the image restoration. Given that the discontinuities/edges of the image lie in the zero level set of a band-limited periodic function, we can derive that the Fourier samples of the gradient of the image satisfy an annihilation relation, resulting in a low-rank two-fold Hankel matrix. In addition, since the singular value decomposition of a low-rank Hankel matrix corresponds to an adaptive tight frame system which can represent the image with sparse canonical coefficients, our approach consists of the following two stages. The first stage learns the tight wavelet frame system from a given measurement, and the second stage restores the image via the analysis approach based sparse regularization. The numerical results are presented to demonstrate that the proposed approach is compared favorably against several popular discrete regularization approaches and structured low-rank matrix approaches.
\end{abstract}

\section{Introduction}\label{Introduction}

Image restoration, including image denoising, deblurring, inpainting, etc., is one of the most fundamental processes in imaging sciences with a wide range of applications. Image restoration aims to obtain an image of high-quality from a given measurement which is degraded during the process of imaging, acquisition, communication, etc. Mathematically, it is in general modeled by the following linear inverse problem:
\begin{align}\label{Linear_IP}
\bsf=\bmA\bu+\zze
\end{align}
where $\bsf$ is the degraded measurement or the observed image, $\zze$ is a certain additive noise or a measurement error, and $\bmA$ is some linear operator which takes different forms depending on the specific image restoration problems.

Since the operator $\bmA$ is in general ill-conditioned or non-invertible, a proper regularization on the images to be recovered is needed to obtain a high quality recovery from the ill-posed linear inverse problem \cref{Linear_IP}. One of the most successful example is the sparse regularization based on the following $\ell_1$ norm minimization
\begin{align}\label{Model:L1}
\min_{\bu}~\lambda\left\|\Ph\bu\right\|_1+\f{1}{2}\left\|\bmA\bu-\bsf\right\|_2^2.
\end{align}
In \cref{Model:L1}, we design a linear transform $\Ph$ so that $\Ph\bu$ is close to zero in the smooth region. Since the inherent thresholding to solve \cref{Model:L1} is connected with the nonlinear evolution PDEs \cite{B.Dong2017}, we can say that \cref{Model:L1} aims to maintain the sharpness of image singularities while regularizing smooth regions at the same time.

Meanwhile, if we know the exact locations of image singularities as oracle, we can restore the image $\bu$ with sharp edges by solving the following least squares problem:
\begin{align}\label{Tikhonov_Oracle}
\min_{\bu}~\lambda\left\|\left(\Ph\bu\right)_{\Ga^c}\right\|_2^2+\left\|\bmA\bu-\bsf\right\|_2^2,
\end{align}
which can be understood as the variational formulation of
\begin{align*}
\bsf=\bmA\bu+\zze~~~\text{subject to}~~\big(\Ph\bu\big)_{\Ga^c}=\0.
\end{align*}
In \cref{Tikhonov_Oracle}, $\Ga$ is the set of indices corresponding to the image singularities; $\Ph\bu$ is close to zero on $\Ga^c$. Therefore, as long as $|\Ga^c|$ (the cardinality of $\Ga^c$) is sufficiently large, \cref{Tikhonov_Oracle} is well-posed, and the image singularities will be well kept in the restored image \cite{H.Ji2016}.

Even though both \cref{Model:L1,Tikhonov_Oracle} are able to restore an image with sharp edges, the major challenge lies in how to estimate $\Ga$ as accurately as possible. In many image restoration tasks, the degraded measurements are available only, which makes it difficult to obtain the exact information on the image singularities. Most importantly, in many practical applications, the true singularities can lie in the continuous domain or at least may not be exactly aligned with the discrete grid, and the discretization will lead to the \emph{basis mismatch} \cite{Y.Chi2010,G.Ongie2016,J.Ying2017} between the true singularities and the discrete image grid. Such a basis mismatch can distort information on the image singularities, and thus can yield the performance bottleneck \cite{J.Ying2017}.

In contrast to the conventional discrete sparse regularization, the continuous domain regularization exploits the sparsity prior in continuous domain, thereby alleviating the basis mismatch \cite{G.Ongie2018}. To the best of our knowledge, this ``off-the-grid'' approach stems from the Prony's method \cite{Prony1795} which corresponds a superposition of a few sinusoids to a structured low-rank matrix for the Dirac stream retrieval with unknown knots. Hence, we can adopt the so-called \emph{structured low-rank matrix (SLRM) approach} \cite{Y.Chen2014,G.Ongie2017,J.C.Ye2017} for the restoration of a superposition of a few sinusoids whose Fourier transform (or inverse transform) is a Dirac stream \cite{J.F.Cai2016,J.F.Cai2018a,J.F.Cai2019}, together with the subspace method for the support retrieval \cite{Schmidt1986}. In addition to the spectrally sparse signals, the SLRM can also be used to restore Fourier samples of a one dimensional piecewise constant signal \cite{T.Blu2008,M.Vetterli2002} as in this case the Fourier samples of a derivative becomes the superposition of a few sinusoids. Finally, the works in \cite{G.Ongie2018,G.Ongie2015a,G.Ongie2015,G.Ongie2016,H.Pan2014} have successfully extended the SLRM to the two dimensional images under a mild assumption on the image singularity curves to overcome the nonisolated image singularity curves in the two dimensional space. In addition, such an assumption on the curves has also enabled the estimation of the image singularity curves via the subspace method \cite{G.Ongie2016}.

In this paper, we propose a two-stage continuous domain regularization approach for the piecewise constant image restoration. Our two-stage approach follows the previous works \cite{J.F.Cai2020,G.Ongie2016} which correspond the Fourier samples of gradient to the structured low-rank matrices. Here we note that the SVD of a Hankel matrix induces tight frame filter banks due to its underlying convolutional structure. In other words, if we can associate a signal/image with a low-rank Hankel matrix, its right singular vectors form tight frame filter banks under which the canonical coefficients have a group sparsity according to the index of filters \cite[Theorem 3.2]{J.F.Cai2020}. Further motivated by the annihilating subspace method in \cite{G.Ongie2016} and the data-driven tight frame in \cite{J.F.Cai2014}, we first construct the tight frame filter banks from the low frequency Fourier samples of a degraded measurement. Once we have constructed the tight frame filter banks, we adopt a \emph{sparse regularization} via the wavelet frame analysis approach \cite{J.F.Cai2009/10} as an image restoration model via the continuous domain regularization. Even though the wavelet frame based analysis (e.g. \cite{J.F.Cai2012}) approach has been widely used for the sparse approximation of a discrete image, our analysis approach comes from the relaxation of the structured low-rank matrices generated by the Fourier samples for the continuous domain regularization.

\subsection{Overviews of related works}

Our two-stage approach for the piecewise constant image restoration is closely related to the recent SLRM framework for the two dimensional functions in \cite{G.Ongie2018,G.Ongie2015a,G.Ongie2015,G.Ongie2016,H.Pan2014}. Briefly speaking, if the singularity curves of a piecewise constant/holomorphic function, i.e. the supports of the first order (real/complex) derivatives of a target image, lie in the zero level set of a band-limited function called the \emph{annihilating polynomial}, the Fourier transform of the derivatives can be annihilated by the convolution with the annihilating filter. This annihilation relation in turn corresponds the Fourier samples of derivatives to the structured low-rank matrices. Based on this framework, the continuous domain regularization for the piecewise constant image restoration is proposed and studied in \cite{G.Ongie2015a,G.Ongie2015,G.Ongie2016,G.Ongie2017}, together with a restoration guarantee \cite{G.Ongie2018}.

Apart from the signal/image restoration, the SLRM framework can also be applied to the source/support estimation from a few (contagious) discrete samples. In the pioneer work \cite{Prony1795}, Prony's method retrieves the point sources of a Dirac stream from the roots of the annihilating polynomial. Since the Prony's method is numerically unstable in spite of the theoretical guarantee in the noiseless setting \cite{P.Guillaume1989}, several annihilating subspace methods have been proposed as stable alternatives in the literature. Most widely used methods include MUSIC \cite{Schmidt1986}, ESPIRIT \cite{R.Roy1989}, and the matrix pencil \cite{Y.Hua1990}. In the two dimensional piecewise constant image setting, the similar method has been introduced under the assumption that the singularities lie in the zero set of an annihilating polynomial. Using the tool of algebraic geometries, the authors in \cite{G.Ongie2016} present the necessary and the sufficient conditions on the low frequency samples for the singularity estimation. Based on these conditions, the subspace method for the singularity estimation is further introduced in \cite{G.Ongie2016}, which can be viewed as a two dimensional parallel of MUSIC algorithm.

Based on the SLRM framework and the corresponding subspace method, a similar two-stage approach for the piecewise constant image restoration is already presented in \cite{G.Ongie2016}. More precisely, the authors in \cite{G.Ongie2016} introduces a Cadzow denoising (e.g. \cite{Cadzow1988}) for the annihilating subspace estimation/the tight frame filter banks construction, and a least squares linear prediction (LSLP) model penalizing the $\ell_2$ norm of the canonical coefficients corresponding to annihilating subspaces for the image restoration. While the Cadzow denoising in \cite{G.Ongie2016} computes the SVD of a two-fold Hankel matrix at each iteration, our tight frame filter bank construction adopts the data-driven tight frame construction in \cite{J.F.Cai2014} which computes the SVD of a smaller matrix. In addition, compared to the LSLP model in \cite{G.Ongie2016}, our restoration model is less sensitive to the rank of the Hankel matrix, as we implicitly exploit the annihilating property via the sparsity promoting $\ell_1$ norm. Finally, since we consider the sparsity in the entire transform domain rather than the group sparsity according to the index of filters, our restoration model is expected to be more robust to the noise/the image modeling.

Finally, a similar wavelet frame relaxation can also be found in the previous work \cite{J.F.Cai2020}, where the data-driven tight frame is used for the piecewise constant image restoration via the continuous domain regularization. While the previous data-driven tight frame model in \cite{J.F.Cai2020} is a \emph{balanced approach} \cite{J.F.Cai2008,R.H.Chan2003} which assumes the distance between the canonical coefficients and the sparse coefficients, our proposed model is an \emph{analysis approach} \cite{J.F.Cai2009/10} which directly penalizes the sparse canonical coefficients. Notice that it is well known that the analysis approach reflects the structure of a target image better than the synthesis approach and the balanced approach (e.g. \cite{J.F.Cai2012}), and the similar arguments can be applied to our case as well; our proposed analysis approach is expected to reflect the SLRM structure than the previous balanced approach in \cite{J.F.Cai2020}.

\subsection{Organization and notation of paper}

The rest of this paper is organized as follows. The proposed two-stage approach based on the continuous domain regularization is presented in \cref{ProposedApproach}. To be more specific, we begin with the brief review on the structured low-rank matrix framework for the piecewise constant images and the corresponding edge estimation from the annihilating relations. Then we present the proposed image restoration model based on the wavelet frame, together with the alternating minimization algorithm, and we introduce the data-driven tight frame model for constructing the tight frame filter banks from degraded measurements. In \cref{Experiments}, we present some numerical results, and \cref{Conclusion} concludes this paper with a few future directions.

Throughout this paper, all two dimensional discrete images will be denoted by the bold faced lower case letters while the functions in a two dimensional continuous domain will be denoted by regular characters. Note that a two dimensional image can also be identified with a vector whenever convenient. All matrices will be denoted by the bold faced upper case letters, and the $j$th row and the $k$th column of a matrix $\bZ$ will be denoted by $\bZ^{(j,:)}$ and $\bZ^{(:,k)}$, respectively. We denote by
\begin{align}\label{ImageGrid}
\OO=\left\{-\lfloor N/2\rfloor,\cdots,\lfloor (N-1)/2\rfloor\right\}^2,
\end{align}
with $N\in\N$, the set of $N\times N$ grid. Throughout this paper, we only consider the square grid $\OO$ for simplicity. Note, however, that it is not difficult to generalize the setting into an arbitrary $N_1\times N_2$ grid
\begin{align*}
\OO=\left\{-\lfloor N_1/2\rfloor,\cdots,\lfloor (N_1-1)/2\rfloor\right\}\times\left\{-\lfloor N_2/2\rfloor,\cdots,\lfloor (N_2-1)/2\rfloor\right\}
\end{align*}
with $N_1,N_2\in\N$. The space of complex valued functions on $\OO$ is denoted by $\mV\simeq\C^{|\OO|}$. Given two symmetric rectangular grids $\KK$ and $\MM$, the contraction $\KK:\MM$ and the Minkowski sum $\KK+\MM$ are defined as
\begin{align*}
\KK:\MM&=\left\{\bk\in\KK:\bk+\MM\subseteq\KK\right\}=\left\{\bk\in\KK:\bk+\bsl\in\KK~\text{for all}~\bsl\in\MM\right\},\\
\KK+\MM&=\left\{\bk+\bsl:\bk\in\KK,~\text{and}~\bsl\in\MM\right\},
\end{align*}
respectively. Note that $2\KK=\KK+\KK$, $3\KK=2\KK+\KK$, etc. Finally, operators are denoted as the bold faced calligraphic letters. For instance, $\bmI$ denotes the identity operator on some Hilbert space, whereas $\bsI$ denotes the identity matrix. In particular, for a rectangular $K_1\times K_2$ grid $\KK$, $\bmH:\mV\to\C^{M_1\times M_2}$ ($M_1=|\OO:\KK|$ and $M_2=|\KK|$) is an operator which corresponds $\bsv\in\mV$ to the Hankel matrix $\bmH\bsv$ by concatenating $K_1\times K_2$ patches of $\bsv$ into row vectors. In the sense of two dimensional multi-indices, $\bmH\bsv$ is defined as
\begin{align*}
\left(\bmH\bsv\right)(\bk,\bsl)=\bsv(\bk+\bsl),~~~~~~\bk\in\OO:\KK~~~\text{and}~~~\bsl\in\KK.
\end{align*}
With a slight abuse of notation, we also use
\begin{align}\label{TwoFoldHankel}
\bmH\w=\left[\begin{array}{cc}
(\bmH\w_1)^T&(\bmH\w_2)^T
\end{array}\right]^T,
\end{align}
to denote the $2M_1\times M_2$ two-fold Hankel matrix constructed from $\w=\left(\w_1,\w_2\right)\in\mV\times\mV$.

\section{Proposed image restoration approach}\label{ProposedApproach}

\subsection{Structured low-rank matrix for piecewise constant images}\label{SLRMFramework}

In this section, we briefly review the continuous domain approach for two dimensional images. Interested readers can refer to \cite{G.Ongie2016} and references therein for more detailed surveys.

Throughout this paper, we consider the SLRM framework for the following piecewise constant function
\begin{align}\label{uModel}
u(\x)=\sum_{j=1}^J\alpha_j1_{\Om_j}(\x),~~~~~~~\x\in\R^2,
\end{align}
from its Fourier samples
\begin{align}\label{Forward}
\wh{u}(\bk)=\msF(u)(\bk)=\int_{\R^2}u(\x)e^{-2\pi i\bk\cdot\x}\rd\x,~~~~~~\bk\in\OO,
\end{align}
where the sample grid $\OO$ is defined as \cref{ImageGrid}. In \cref{uModel}, $\alpha_j\in\C$ and $1_{\Om}$ denotes the characteristic function on a set $\Om$: $1_{\Om}(\x)=1$ if $\x\in\Om$, and $0$ otherwise. Without loss of generality, we assume that each $\Om_j$ lies in $[-1/2,1/2)^2$ for simplicity. (However, it is not difficult to generalize the setting into an arbitrary rectangular region $[-L_1/2,L_1/2)\times[-L_2/2,L_2/2)$.) We further assume that \cref{uModel} is expressed with the smallest number of characteristic functions such that $\Om_j$'s are pairwise disjoint. Then the discontinuities of $u$ agrees with $\Gamma=\bigcup_{j=1}^J\p\Om_j$, called the \emph{singularity set} of $u$ \cite{G.Ongie2016}.

In general, it is difficult to directly establish the SLRM framework without any further information on $\Gamma$. Nevertheless, the authors in \cite{G.Ongie2016} assume that there exists a finite rectangular and symmetric index set $\KK$ such that
\begin{align}\label{MajorAssumption}
\Gamma\subseteq\left\{\x\in\R^2:\varphi(\x)=0\right\}~~~~~\text{with}~~~~~\varphi(\x)=\sum_{\bk\in\KK}\a(\bk)e^{-2\pi i\bk\cdot\x}.
\end{align}
Following the previous works in \cite{J.F.Cai2022,G.Ongie2015,G.Ongie2016}, any function $\varphi(\x)$ in the form of \cref{MajorAssumption} will be called the \emph{trigonometric polynomial}, and the zero level set $\left\{\x\in\R^2:\varphi(\x)=0\right\}$ the \emph{trigonometric curve}. For a trigonometric polynomial $\varphi$ in \cref{MajorAssumption}, the degree of $\varphi$ is defined as an ordered pair of degrees in each coordinate, and is denoted by $\deg(\varphi)$. In particular, a trigonometric polynomial with the smallest degree is called the \emph{minimal polynomial}. Under this setting, the authors in \cite{G.Ongie2015,G.Ongie2016} derived that the Fourier transform of $\na u=\left(\p_1u,\p_2u\right)$ is annihilated by the convolution with the Fourier coefficients $\a:=\left\{\a(\bk):\bk\in\KK\right\}$ of $\varphi$. More precisely, we have
\begin{align}\label{FouGradConvSemiDiscrete}
\left(\msF(\na u)\ast\wh{\varphi}\right)(\xxi)=\sum_{\bk\in\KK}\msF(\na u)(\xxi+\bk)\a(\bk)=0,~~~~~~\xxi\in\R^2,
\end{align}
where the convolution acts on $\msF(\p_1 u)$ and $\msF(\p_2 u)$ separately.

In this sense, we call the trigonometric polynomial $\varphi$ in \cref{MajorAssumption} the \emph{annihilating polynomial}, and its Fourier coefficient $\a=\left\{\a(\bk):\bk\in\KK\right\}$ the \emph{annihilating filter} with support $\KK$. Using the tool of algebraic geometries, the authors in \cite{G.Ongie2016} proved that there exists the unique minimal polynomial associated with a given trigonometric curve, and that the degree of a trigonometric polynomial bounds the number of connected components (i.e. $J$ in \cref{uModel}), from which we can estimate the degrees of freedom of an image to be restored \cite{G.Ongie2016}.

In our setting, the Fourier transform of $u$ is sampled on the grid $\OO$ in \cref{ImageGrid} with $N\in\N$ large enough to guarantee a high image resolution. Hence, \cref{FouGradConvSemiDiscrete} becomes the finite system of linear equations
\begin{align}\label{FouGradConvDiscrete}
\sum_{\bk\in\KK}\msF(\na u)(\bsl+\bk)\a(\bk)=0,~~~~~~~\bsl\in\OO:\KK.
\end{align}
In the matrix-vector multiplication form, we have
\begin{align}\label{RankDeficientHankel}
\bmH\left(\msF(\na u)\big|_{\OO}\right)\a=\0.
\end{align}
Hence, the two-fold Hankel matrix $\bmH\left(\msF(\na u)\big|_{\OO}\right)$ has a nontrivial nullspace. Indeed, as long as the filter support $\KK'$ defining $\bmH\left(\msF(\na u)\big|_{\OO}\right)$ is sufficiently large, $\bmH\left(\msF(\na u)\big|_{\OO}\right)$ will have a nontrivial nullspace. To see this, let $\varphi$ in \cref{MajorAssumption} be the minimal polynomial with the support $\KK$, and let $\KK'$ be the assumed filter support strictly containing $\KK$.  Then we have
\begin{align}\label{RankUB}
\rank\left(\bmH\left(\msF(\na u)\big|_{\OO}\right)\right)\leq|\KK'|-|\KK':\KK|,
\end{align}
which shows that the two-fold Hankel matrix constructed by $\msF(\na u)\big|_{\OO}$ is \emph{rank deficient} \cite{G.Ongie2016}. Roughly speaking, \cref{RankUB} means that if $\varphi$ in \cref{MajorAssumption} is the minimal polynomial for the image singularities $\Gamma$ with Fourier coefficients $\a$ supported on $\KK$, $e^{-2\pi i\bm\cdot\x}\varphi(\x)$ is also an annihilating polynomial for $\bm\in\KK':\KK$, or equivalently, $\a(\cdot-\bm)$ is also an annihilating filter \cite{J.F.Cai2020}.

Apart from the low-rank property \cref{RankUB} of the two-fold Hankel matrix $\bmH\left(\msF(\na u)\big|_{\OO}\right)$, we can also derive the necessary and the sufficient conditions on the number of Fourier samples (i.e. the size of $\OO$) to guarantee the unique recovery of the singularities from the annihilating relation \cref{FouGradConvDiscrete} (or \cref{RankDeficientHankel}) \cite{G.Ongie2016}. To be more precise, assume that $u(\x)$ be as in \cref{uModel} with the singularity set $\Gamma$ satisfying \cref{MajorAssumption} with the minimal polynomial $\varphi$ of the degree $(K_1,K_2)$. If we can restore $\Gamma$ using the samples $\msF(u)\big|_{\OO}$ on the $N\times N$ grid $\OO$, then it must be
\begin{align}\label{Necessary}
2\left(N-K_1\right)\left(N-K_2\right)\geq\left(K_1+1\right)\left(K_2+1\right)-1,
\end{align}
which gives the necessary condition on $\OO$ for the estimation of $\Gamma$ \cite[Proposition 3.1]{G.Ongie2016}. Conversely, let $\varphi$ in \cref{MajorAssumption} be the minimal polynomial with the Fourier coefficients $\a$ supported on $\KK$, and let $\KK'$ be a rectangular grid strictly containing $\KK$. If a filter $\bb=\left\{\bb(\bk):\bk\in\KK'\right\}$ on $\KK'$ satisfies
\begin{align}\label{Sufficient}
\sum_{\bk\in\KK'}\msF(\na u)(\bsl+\bk)\bb(\bk)=0,~~~~~~~\bsl\in\KK'+\KK,
\end{align}
there exists a trigonometric polynomial $\eta$ with Fourier coefficients supported on $\KK':\KK$ such that
\begin{align*}
\left(\eta\varphi\right)(\x)=\sum_{\bk\in\KK'}\bb(\bk)e^{-2\pi i\bk\cdot\x}.
\end{align*}
Hence, the trigonometric polynomial corresponding to the filter $\bb$ is also an annihilating polynomial, which gives a sufficient condition on the sample grid $\OO$ \cite[Theorem 3.4]{G.Ongie2016}.

In general, it is challenging to identify the minimal annihilating filter from the nullspace of $\bmH\left(\msF(\na u)\big|_{\OO}\right)$ \cite{G.Ongie2016}. Nevertheless, it is still possible to apply the subspace algorithm such as MUSIC \cite{Schmidt1986} to our setting. For the sake of clarity, let $\bB=\left[\begin{array}{ccc}
\bb_1&\cdots&\bb_R
\end{array}\right]$ be an orthonormal basis for $\msN\left(\bmH\left(\msF(\na u)\big|_{\OO}\right)\right)$, and let $\varphi_m$ be the trigonometric polynomial corresponding to $\bb_m$. Then it follows that $\varphi_m=\eta_m\varphi$ for some trigonometric polynomial $\eta_m$. We define the following \emph{sum-of-squares average}
\begin{align}\label{SumofSquares}
\overline{\varphi}(\x)=\left(\sum_{m=1}^R\left|\varphi_m(\x)\right|^2\right)^{1/2}=\left(\sum_{m=1}^R\left|\eta_m(\x)\varphi(\x)\right|^2\right)^{1/2}.
\end{align}
Obviously, we have $\overline{\varphi}(\x)=0$ for $\x\in\Gamma$. In addition, it is also proved in \cite{G.Ongie2016} that $\overline{\varphi}(\x)\neq0$ almost everywhere in $\Gamma^c$, which makes it possible to estimate $\Gamma$ from $\msN\left(\bmH\left(\msF(\na u)\big|_{\OO}\right)\right)$. Finally, $\overline{\varphi}(\x)$ can be expressed as the reciprocal of pseudospectra:
\begin{align}\label{Pseudospectra}
\overline{\varphi}(\x)=\left\|\bB\bB^*\e_{\x}\right\|_2
\end{align}
where $\e_{\x}(\bk)=e^{-2\pi i\bk\cdot\x}$ with $\bk\in\KK'$. Note that \cref{Pseudospectra} can be understood as a two dimensional parallel to MUSIC. See \cite{G.Ongie2015,G.Ongie2016} for details.

\subsection{Image restoration model}\label{ImageRestorationModel}

Let $\bsf$ be a degraded measurement modeled as
\begin{align}\label{Linear_IP2}
\bsf=\bmA\bsv+\zze,
\end{align}
where $\bsv=\msF(u)\big|_{\OO}\in\mV$ with $u$ defined as in \cref{uModel}, and $\zze$ is some additive noise/measurement error. In what follows, with a slight abuse of notation, we assume the linear operator $\bmA$ acts on the Fourier samples.

Since the two-fold Hankel matrix $\bmH\left(\msF(\na u)\big|_{\OO}\right)\in\C^{2M_1\times M_2}$ defined as \cref{TwoFoldHankel} is low-rank, we can consider
\begin{align}\label{RankMinimization}
\min_{\bsv}~\rank\left(\bmH\left(\bmD\bsv\right)\right)~~~~~\text{subject to}~~~~\bsf=\bmA\bsv+\zze
\end{align}
as a continuous domain regularization for the piecewise constant image restoration. In \cref{RankMinimization}, $\bmD:\mV\to\mV\times\mV$ is defined as
\begin{align}\label{Ddef}
\bmD\bsv(\bk)=\left(2\pi ik_1\bsv(\bk),2\pi ik_2\bsv(\bk)\right),~~~~~\bk=(k_1,k_2)\in\OO,
\end{align}
which reflects the samples of $\msF(\na u)(\xxi)=2\pi i\xxi\wh{u}(\xxi)=\left(2\pi i\xi_1\wh{u}(\xxi),2\pi i\xi_2\wh{u}(\xxi)\right)$.

In the literature, there are numerous tractable relaxations of \cref{RankMinimization}, which is an NP-hard problem in general. Some of them include the convex nuclear norm relaxation (e.g. \cite{J.F.Cai2016,M.Fazel2013}), the iterative reweighted least squares (IRLS) for the Schatten $p$-norm minimization \cite{M.Fornasier2011,Y.Hu2019,K.Mohan2012,G.Ongie2017}, and the nonconvex alternating projections \cite{J.F.Cai2018a,J.F.Cai2019} etc. In this paper, considering the SVD of a low-rank Hankel matrix and inspired by the previous work in \cite{J.F.Cai2022} for the piecewise smooth image restoration, we propose another relaxation of \cref{RankMinimization}, by presenting another interpretation of \cite[Theorem 3.2]{J.F.Cai2020}.

To begin with, we let $\a_1,\cdots,\a_{M_2}$ be $K_1\times K_2$ filters supported on $\KK$, and we introduce
\begin{align}
\bmW&=\big[\bmS_{\a_1(-\cdot)}^T,\bmS_{\a_2(-\cdot)}^T,\cdots,\bmS_{\a_{M_2}(-\cdot)}^T\big]^T,\label{OurAnalysis}\\
\bmW^*&=\big[\bmS_{\overline{\a}_1},\bmS_{\overline{\a}_2},\cdots,\bmS_{\overline{\a}_{M_2}}\big],\label{OurSynthesis}
\end{align}
where $\bmS_{\a}$ is a discrete convolution under the periodic boundary condition:
\begin{align*}
\left(\bmS_{\a}\bsv\right)(\bk)=\left(\a\ast\bsv\right)(\bk)=\sum_{\bm\in\Z^2}\a(\bk-\bm)\bsv(\bm).
\end{align*}
In other words, both $\bmW$ and $\bmW^*$ are concatenations of discrete convolutions.

Under this setting, we assume that $\rank\left(\bmH\left(\bmD\bsv\right)\right)=r\ll2M_1\wedge M_2$. Considering its full SVD $\bmH\left(\bmD\bsv\right)=\bX\Sig\bY^*$, we let
\begin{align}\label{OurFilters}
\a_m=M_2^{-1/2}\bY^{(:,m)}
\end{align}
by reformulating each $\bY^{(:,m)}\in\C^{M_2}$ into a $K_1\times K_2$ filter supported on $\KK$. Then \cite[Theorem 3.2]{J.F.Cai2020} tells us that $\bmW$ in \cref{OurAnalysis} defined by using the filters $\a_1,\cdots,\a_{M_2}$ in \cref{OurFilters} satisfies
\begin{align}\label{LowRankHankelTightFrame}
\bmW^*\bmW\left(\bmD\bsv\right)=\sum_{m=1}^{M_2}\bmS_{\overline{\a}_m}\left(\bmS_{\a_m(-\cdot)}\left(\bmD\bsv\right)\right)=\bmD\bsv,
\end{align}
and for $m=r+1,\cdots,M_2$, we have
\begin{align}\label{TightFrameSparse}
\bmS_{\a_m(-\cdot)}\left(\bmD\bsv\right)(\bk)=\0,~~~~~\bk\in\OO:\KK,
\end{align}
where the discrete convolution $\bmS_{\a}$ acts on each component of $\bmD\bsv$.

In other words, if we know the SVD of a two-fold Hankel matrix $\bmH\left(\bmD\bsv\right)\in\C^{2M_1\times M_2}$ as an oracle, we can explicitly construct tight frame filter banks under which $\bmD\bsv$ admits a sparse canonical representation, and the sparsity of canonical coefficients can be grouped according to the index of filters \cite{J.F.Cai2020}. Here, motivated by \cite{J.F.Cai2022}, we assume that a priori estimation $\wt{\bsv}\in\mV$ is available with the SVD
\begin{align*}
\bmH\left(\bmD\wt{\bsv}\right)=\wt{\bX}\wt{\Sig}\wt{\bY}^*,
\end{align*}
whose detailed method will be discussed in the following subsection. We define the tight frame transform $\bmW$ in \cref{OurAnalysis} via
\begin{align*}
\a_m=M_2^{-1/2}\wt{\bY}^{(:,m)},~~~~~m=1,\cdots,M_2.
\end{align*}
Since then \cref{TightFrameSparse} holds approximately under this $\bmW$, we remove the group sparsity in the canonical coefficients for the better sparse approximation instead. This leads us to solve
\begin{align}\label{ProposedModel}
\min_{\bsv}~\f{1}{2}\left\|\bmA\bsv-\bsf\right\|_2^2+\left\|\gga\cdot\bmW\left(\bmD\bsv\right)\right\|_1
\end{align}
where the weighted $\ell_1$ norm is defined as
\begin{align*}
\left\|\gga\cdot\bmW\left(\bmD\bsv\right)\right\|_1=\sum_{m=1}^{M_2}\gamma_m\left\|\bmS_{\a_m(-\cdot)}\left(\bmD\bsv\right)\right\|_1
\end{align*}
to reflect the different weights related to the singular values. That is,
\begin{align*}
\gamma_m=\f{\nu}{\wt{\Sig}^{(m,m)}+\eps}
\end{align*}
with some $\nu>0$ and a small $\eps>0$ to avoid the division by zero.

Note that the similar idea can be found in the previous related works. More precisely, since the filters $\a_{r+1},\cdots,\a_{M_2}$ form an orthogonal basis for $\msN\left(\bmH\left(\bmD\bsv\right)\right)$, we can also consider the following $\ell_2$ minimization model
\begin{align}\label{LSLP}
\min_{\bsv}~\left\|\bmA\bsv-\bsf\right\|_2^2+\gamma\sum_{m=r+1}^{M_2}\left\|\bmS_{\a_m(-\cdot)}\left(\bmD\bsv\right)\right\|_2^2,
\end{align}
which is the least squares linear prediction (LSLP) model in \cite{G.Ongie2016} to solve the ideal low-pass filter deconvolution. Indeed, if we know the exact right singular vectors $\bmH\left(\bmD\bsv\right)=\bX\Sig\bY^*$ as an oracle, the above LSLR model will be more appropriate from the viewpoint of \cref{LowRankHankelTightFrame,TightFrameSparse}, as our model implicitly reflects \cref{TightFrameSparse} by the $\ell_1$ norm. However, since we only have the degraded measurement in many practical cases, the challenge lies in how to find the wavelet frame transform $\bmW$ (and the rank $r$) satisfying \cref{TightFrameSparse} as exactly as possible. Nevertheless, since we relax the sparsity of the tight frame coefficients over the entire range of $\bmW$, (i.e. not necessarily grouped as in \cref{TightFrameSparse}), we can expect that our model \cref{ProposedModel} is able to achieve more flexibility, leading to the performance gain in the image restoration.

We mention that \cref{LowRankHankelTightFrame,TightFrameSparse} are originally to propose the following data-driven tight frame (DDTF) relaxation
\begin{align}\label{LRHDDTF}
\begin{split}
&~~~\min_{\bsv,\bc,\bmW}~\f{1}{2}\left\|\bmA\bsv-\bsf\right\|_2^2+\f{\mu}{2}\left\|\bmW\left(\bmD\bsv\right)-\bc\right\|_2^2+\left\|\gga\cdot\bc\right\|_0\\
&\text{subject to}~~\bmW^*\bmW=\bmI,
\end{split}
\end{align}
for the compressed sensing MRI restoration via the continuous domain regularization. In \cref{LRHDDTF}, $\ell_0$ norm $\left\|\gga\cdot\bc\right\|_0$ encodes the number of nonzero elements in $\bc$. At this moment, we see that our model \cref{ProposedModel} has the following advantages over the data-driven tight frame model \cref{LRHDDTF}. First of all, given that $\bmW$ is estimated appropriately, we may not require to further learn it at each iteration with additional computational costs. Second, since \cref{ProposedModel} is convex while \cref{LRHDDTF} is nonconvex, we can easily expect better behavior and theoretical support for the numerical algorithms. Most importantly, from the viewpoint of wavelet frame based image restoration, the model \cref{ProposedModel} is an \emph{analysis approach} \cite{J.F.Cai2009/10} while the data-driven tight frame model \cref{LRHDDTF} is (a relaxation of) the \emph{balanced approach} \cite{J.F.Cai2008,R.H.Chan2003} (See \cref{PreliminariesTightFrame} for the brief description on the wavelet frame based image restoration models). Note that the analysis approach reflects the structure of a target image better than other approaches under some mild assumptions on the tight frame (e.g. \cite{J.F.Cai2012}). Hence, the same arguments can be applied to our case; even though both \cref{ProposedModel} and \cref{LRHDDTF} are derived from \cref{LowRankHankelTightFrame,TightFrameSparse}, the analysis approach \cref{ProposedModel} reflects the SLRM framework better than the balanced approach \cref{LRHDDTF}.

\begin{remark} We note that our image restoration model \cref{ProposedModel} mainly focuses on the restoration of images which can be modeled/approximated by piecewise constant functions. This means that, neither piecewise smooth images nor textures may not be suitable for \cref{ProposedModel}. Even though a structured low-rank matrix framework for the piecewise smooth functions are proposed in \cite{J.F.Cai2022}, further researches are required as the estimation of annihilating subspaces is involved with the blind deconvolutions. In addition, noting that the structured low-rank matrices in the image domain can reflect the textures \cite{Y.R.Fan2017,H.Schaeffer2013}, we may be able to jointly adopt the structured low-rank matrix framework in different domains for the cartoon-texture decomposition. Nevertheless, we forgo further discussions on these topics as both of them are beyond the scope of this work. We only focus on the restoration of piecewise constant image.
\end{remark}

To solve the convex model \cref{ProposedModel}, there are numerous numerical algorithms available in the literature. In this paper, we use the ADMM \cite{J.Eckstein1992} or the split Bregman algorithm \cite{J.F.Cai2009/10} as we can convert \cref{ProposedModel} into several subproblems with closed form solutions, together with the convergence guarantee \cite{J.F.Cai2009/10}. Since the derivation is quite standard, we forgo further details and directly present the overall algorithm in \cref{Alg1}.

\begin{algorithm}[t!]
\begin{algorithmic}
\STATE{\textbf{Initialization:} $\bsv_{(0)}$, $\bc_{(0)}$, $\bsd_{(0)}$}
\FOR{$n=0$, $1$, $2$, $\cdots$}
\STATE{\textbf{(1)} Update $\bsv$:
\begin{align}\label{vsubprob}
\bsv_{(n+1)}=\argmin_{\bsv}\f{1}{2}\left\|\bmA\bsv-\bsf\right\|_2^2+\f{\beta}{2}\left\|\bmW\left(\bmD\bsv\right)-\bc_{(n)}+\wt{\bc}_{(n)}\right\|_2^2
\end{align}
\textbf{(2)} Update $\bc$:
\begin{align}\label{csubprob}
\bc_{(n+1)}=\argmin_{\bc}\left\|\gga\cdot\bc\right\|_1+\f{\beta}{2}\left\|\bc-\bmW\left(\bmD\bsv_{(n+1)}\right)-\bsd_{(n)}\right\|_2^2
\end{align}
\textbf{(3)} Update $\bsd$:
\begin{align}
\bsd_{(n+1)}=\bsd_{(n)}+\bmW\left(\bmD\bsv_{(n+1)}\right)-\bc_{(n+1)}
\end{align}}
\ENDFOR
\end{algorithmic}
\caption{Split Bregman algorithm for \cref{ProposedModel}}\label{Alg1}
\end{algorithm}

For \cref{vsubprob}, since $\bmW^*\bmW=\bmI$, we solve the following linear equation:
\begin{align}\label{vsubexplicit}
\bsv_{(n+1)}=\left(\bmA^*\bmA+\beta\bmD^*\bmD\right)^{-1}\left[\bmA^*\bsf+\beta\bmD^*\bmW^*\left(\bc_{(n)}-\bsd_{(n)}\right)\right].
\end{align}
Notice that if $\bmA$ is a pointwise multiplication in the frequency domain such as image denoising, image deblurring, and CS restoration, so is $\bmA^*\bmA+\beta\bmD^*\bmD$, and thus, no matrix inversion is needed. For other cases, we can use the iterative solver such as the conjugate gradient method.

The closed form solution for \cref{csubprob} is expressed in terms of the soft thresholding:
\begin{align}\label{csubexplicit}
\bc_{(n+1)}=\bmT_{\gga/\beta}\left(\bmW\left(\bmD\bsv_{(n+1)}\right)+\bsd_{(n)}\right).
\end{align}
In \cref{csubexplicit}, the soft thresholding $\bmT_{\gga}\left(\bc\right)$ for the frame coefficients $\bc=(\bc_1,\bc_2)$ with $\bc_l=\left[\begin{array}{ccc}
\bc_{l1}^T,\cdots,\bc_{lM_2}^T
\end{array}\right]^T$ ($l=1,2$), and the thresholding parameter $\gga=\left[\begin{array}{ccc}
\gamma_1&\cdots&\gamma_{M_2}
\end{array}\right]^T$ is defined as
\begin{align*}
\bmT_{\gga}\left(\bc\right)_{lm}(\bk)=\max\left\{\left|\bc_{lm}(\bk)\right|-\gamma_j,0\right\}\f{\bc_{lm}(\bk)}{\left|\bc_{lm}(\bk)\right|}
\end{align*}
for $\bk\in\OO$, $l=1,2$ and $m=1,\cdots,M_2$. By convention, we set $0/0=0$.

\subsection{Construction of tight frame filter banks}

In our proposed image restoration model \cref{ProposedModel}, we need a pre-estimation of $\bsv$ to construct appropriate tight frame filter banks. However, since what we need is that the filter banks satisfy \cref{TightFrameSparse} in \cref{Th2}, a proper estimation of the null space of the two-fold Hankel matrix will be sufficient. Indeed, it is proved in \cite{G.Ongie2016} that we can estimate the null space given that a sufficient number of low frequency samples are available. More precisely, let $\wt{\OO}\subseteq\OO$ be an $\wt{N}\times\wt{N}$ grid defined as \cref{ImageGrid} with $\wt{N}<N$. We further assume that $\wt{N}\in\N$ satisfies the sufficient condition on the singularity estimation described in \cite{G.Ongie2016}. Let $\wt{\bsv}=\bsv\big|_{\wt{\OO}}\in\wt{\mV}$ where $\wt{\mV}$ denotes the space of complex valued functions defined on $\wt{\OO}$. Then it suffices to set
\begin{align}\label{OurTightFrame}
\a_m=M_2^{-1/2}\wt{\bY}^{(:,m)},~~~~~\text{where}~~~~\bmH\left(\wt{\bmD}\wt{\bsv}\right)=\wt{\bX}\wt{\Sig}\wt{\bY}^*~~\text{is the SVD}
\end{align}
with $\wt{\bmD}=\bmD\big|_{\wt{\mV}}$.

Hence, for the construction of filter banks, it suffices to estimate low frequency samples $\wt{\bsv}\in\wt{\mV}$ from a given degraded measurement $\bsf$. For this purpose, we observe that, in the ideal low-pass filter deconvolution, i.e. $\bmA=\bmR_{\wt{\OO}}$ being the restriction on $\wt{\OO}$, we can restore the Fourier samples $\wt{\bsv}$ on $\wt{\OO}$ via
\begin{align}\label{CadzowDen}
\min_{\wt{\bsv}}\left\|\wt{\bsv}-\bsf\right\|_2^2~~~~~\text{subject to}~~~\rank\left(\bmH\left(\wt{\bmD}\wt{\bsv}\right)\right)\leq r,
\end{align}
where $r>0$ denotes the rank bound estimated from \cref{RankUB}. For the pre-restoration of low frequency Fourier samples, we can extend this idea to the generic image restoration problem \cref{Linear_IP2}. Letting $\wt{\bsf}=\bsf\big|_{\wt{\OO}}$ and $\wt{\bmA}=\bmA\big|_{\wt{\mV}}$, it is natural to solve
\begin{align}\label{CadzowGen}
\min_{\wt{\bsv}}\left\|\wt{\bmA}\wt{\bsv}-\wt{\bsf}\right\|_2^2~~~~~\text{subject to}~~~\rank\left(\bmH\left(\wt{\bmD}\wt{\bsv}\right)\right)\leq r,
\end{align}
which is a maximum likelihood with a low-rank constraint \cite{G.Ongie2016}.

In general, solving \cref{CadzowGen} by the Cadzow denoising on a lifted matrix is slow to converge \cite{G.Ongie2016,H.Pan2014}, which may not be suitable for our purpose. This can be resolved by the following variable splitting method \cite{Haldar2014,G.Ongie2016}
\begin{align}\label{CadzowVS}
\begin{split}
&~~~\min_{\wt{\bsv},\bZ}\left\|\wt{\bmA}\wt{\bsv}-\wt{\bsf}\right\|_2^2+\beta\left\|\bmH\left(\wt{\bmD}\wt{\bsv}\right)-\bZ\right\|_F^2\\
&\text{subject to}~~\rank\left(\bZ\right)\leq r.
\end{split}
\end{align}
Notice that \cref{CadzowVS} finds $\wt{\bsv}$ such that $\bmH\left(\wt{\bmD}\wt{\bsv}\right)$ can be well approximated by a low-rank matrix.

To further reduce the computational costs which arises in computing the SVD of a two-fold Hankel matrix during the alternating minimization, we apply \cref{LowRankHankelTightFrame,TightFrameSparse} to relax the low-rank approximation in \cref{CadzowVS} into the sparse approximation via the data-driven tight frame \cite{J.F.Cai2014} (See \cref{DataDrivenTightFrame} for the brief introduction on the data-driven tight frames). This leads us to solve
\begin{align}\label{DDTFApproach}
\begin{split}
&~~~\min_{\wt{\bsv},\wt{\bc},\wt{\bmW}}\left\|\wt{\bmA}\wt{\bsv}-\wt{\bsf}\right\|_2^2+\beta\left\|\wt{\bmW}\left(\wt{\bmD}\wt{\bsv}\right)-\wt{\bc}\right\|_2^2\\
&\text{subject to}~~\wt{\bc}\in\mC,~~\text{and}~~\wt{\bmW}^*\wt{\bmW}=\bmI.
\end{split}
\end{align}
In \cref{DDTFApproach}, the constraint $\mC$ is defined as
\begin{align}\label{Constraint}
\mC=\left\{\wt{\bc}=\left(\wt{\bc}_1,\wt{\bc}_2\right):\wt{\bc}_m(\bk)=\left(\wt{\bc}_{1m}(\bk),\wt{\bc}_{2m}(\bk)\right)=\0~~~\text{for}~~m=r+1,\cdots,M_2,~\&~\bk\in\OO:\KK\right\},
\end{align}
which reflects \cref{TightFrameSparse}.

To solve \cref{DDTFApproach}, we use the proximal alternating minimization algorithm in \cite{H.Attouch2010}. More precisely, we initialize $\wt{\bsv}_{(0)}=\wt{\bmA}^*\wt{\bsf}$. From the SVD $\bmH\left(\wt{\bmD}\wt{\bsv}_{(0)}\right)=\wt{\bX}_{(0)}\wt{\Sig}_{(0)}\wt{\bY}_{(0)}^*$, we set $\wt{\bmW}_{(0)}$ via
\begin{align*}
\wt{\a}_{(0),m}=M_2^{-1/2}\wt{\bY}_{(0)}^{(:,m)}~~~~~m=1,\cdots,M_2.
\end{align*}
Using the initial filter banks, we initialize $\wt{\bc}_{(0)}$ by projecting $\bmW_{(0)}\left(\wt{\bmD}\wt{\bsv}_{(0)}\right)$ onto the constraint set $\mC$:
\begin{align*}
\wt{\bc}_{(0)}=\bmP_{\mC}\left(\wt{\bmW}_{(0)}\left(\wt{\bmD}\wt{\bsv}_{(0)}\right)\right).
\end{align*}
After the initializations, we optimize $(\wt{\bsv},\wt{\bc},\wt{\bmW})$ alternatively, as summarized in \cref{Alg2}.

\begin{algorithm}[t!]
\begin{algorithmic}
\STATE{\textbf{Initialization:} $\wt{\bsv}_{(0)}$, $\wt{\bc}_{(0)}$, $\wt{\bmW}_{(0)}$}
\FOR{$n=0$, $1$, $2$, $\cdots$}
\STATE{\textbf{(1)} Update $\wt{\bsv}$:
\begin{align}\label{vtsubprob}
\wt{\bsv}_{(n+1)}=\argmin_{\wt{\bsv}}\left\|\wt{\bmA}\wt{\bsv}-\wt{\bsf}\right\|_2^2+\beta\left\|\wt{\bmW}_{(n)}\left(\wt{\bmD}\wt{\bsv}\right)-\wt{\bc}_{(n)}\right\|_2^2+\beta_1\left\|\wt{\bsv}-\wt{\bsv}_{(n)}\right\|_2^2
\end{align}
\textbf{(2)} Update $\wt{\bc}$:
\begin{align}\label{ctsubprob}
\bc_{(n+1)}=\argmin_{\wt{\bc}\in\mC}\left\|\wt{\bc}-\wt{\bmW}_{(n)}\left(\wt{\bmD}\wt{\bsv}_{(n+1)}\right)\right\|_2^2+\beta_2\left\|\wt{\bc}-\wt{\bc}_{(n)}\right\|_2^2
\end{align}
\textbf{(3)} Update $\wt{\bmW}$:
\begin{align}\label{Wsub}
\bmW_{(n+1)}=\argmin_{\wt{\bmW}^*\wt{\bmW}=\bmI}\left\|\wt{\bmW}\left(\wt{\bmD}\wt{\bsv}_{(n+1)}\right)-\wt{\bc}_{(n+1)}\right\|_2^2+\beta_3\left\|\wt{\bmW}-\wt{\bmW}_{(n)}\right\|_F^2
\end{align}}
\ENDFOR
\end{algorithmic}
\caption{Proximal alternating minimization algorithm for \cref{DDTFApproach}}\label{Alg2}
\end{algorithm}

The solution to \cref{vtsubprob} is similar to \cref{vsubexplicit}; since $\wt{\bmW}_{(n)}^*\wt{\bmW}_{(n)}=\bmI$ for $n=0,1,2,\cdots$, we have
\begin{align}\label{vtsubexplicit}
\wt{\bsv}_{(n+1)}=\left(\wt{\bmA}^*\wt{\bmA}+\beta\wt{\bmD}^*\wt{\bmD}+\beta_1\bmI\right)^{-1}\left(\wt{\bmA}^*\wt{\bsf}+\beta\wt{\bmD}^*\wt{\bmW}_{(n)}^*\bc_{(n)}+\beta_1\wt{\bsv}_{(n)}\right).
\end{align}
To solve \cref{ctsubprob,Wsub}, we introduce $\wt{\bsH}_{(n+1)}=\bmH\left(\wt{\bmD}\wt{\bsv}_{(n+1)}\right)\in\C^{2\wt{M}_1\times M_2}$ with $\wt{M}_1=|\wt{\OO}:\KK|$, and let $\wt{\bsA}\in\C^{M_2\times M_2}$ be a matrix whose column vectors are concatenations of the filters $\wt{\a}_1,\cdots,\wt{\a}_{M_2}$. Then we can reformulate \cref{ctsubprob,Wsub} into
\begin{align}
\wt{\bC}_{(n+1)}&=\argmin_{\wt{\bC}\in\mC}\beta\left\|\wt{\bC}-\wt{\bsH}_{(n+1)}\wt{\bsA}_{(n)}\right\|_F^2+\beta_2\left\|\wt{\bC}-\wt{\bC}_{(n)}\right\|_F^2\label{Ctsub}\\
\wt{\bsA}_{(n+1)}&=\argmin_{\wt{\bsA}\wt{\bsA}^*=M_2^{-1}\bsI}\beta\left\|\wt{\bsH}_{(n+1)}\wt{\bsA}-\wt{\bC}_{(n+1)}\right\|_F^2+\beta_3\left\|\wt{\bsA}-\wt{\bsA}_{(n)}\right\|_F^2\label{Dupdate}.
\end{align}
In \cref{Ctsub,Dupdate}, $\left\|\cdot\right\|_F$ denotes the Frobenius norm of a matrix, and $\wt{\bC}=\left[\begin{array}{cc}
\wt{\bC}_1^T&\wt{\bC}_2^T
\end{array}\right]^T\in\C^{2\wt{M}_1\times M_2}$, and $\wt{\bC}_j\in\C^{\wt{M}_1\times M_2}$ ($j=1,2$) denotes the corresponding matrix formulations of frame coefficients $\wt{\bc}_j$. With a slight abuse of notation, we use the same notation $\mC$ for the $2\wt{M}_1\times M_2$ matrix reformulation of \cref{Constraint}:
\begin{align}\label{Constraint2}
\mC=\left\{\wt{\bC}\in\C^{2M_1\times M_2}:\wt{\bC}^{(:,m)}=\0~~\text{for}~~m=r+1,\cdots,M_2\right\}.
\end{align}
Under this reformulation, the solutions to \cref{Ctsub,Dupdate} have the closed form solutions; namely, we have
\begin{align}
\wt{\bC}_{(n+1)}^{(:,m)}&=\left\{\begin{array}{cl}
\left(\beta\wt{\bsH}_{(n+1)}\wt{\bsA}_{(n)}^{(:,m)}+\beta_2\wt{\bC}_{(n)}^{(:,m)}\right)/\left(\beta+\beta_2\right)~&\text{if}~m=1,\cdots,r\vspace{0.25em}\\
\0~&\text{otherwise,}
\end{array}\right.\label{Ctsubexplicit}\\
\wt{\bsA}_{(n+1)}&=M_2^{-1/2}\wt{\bU}\wt{\bV}^*~~\text{where}~~\wt{\bsH}_{(n+1)}^*\wt{\bC}_{(n+1)}+\beta_3\beta^{-1}\wt{\bsA}_{(n)}=\wt{\bU}\wt{\bsS}\wt{\bV}^*~~\text{is the SVD}.\label{Wsubexplicit}
\end{align}
Therefore, we only compute the SVD of $\wt{\bsH}_{(n+1)}^*\wt{\bC}_{(n+1)}+\beta_3\beta^{-1}\wt{\bsA}_{(n)}\in\C^{M_2\times M_2}$, leading to the computational efficiency over the SVD of $\bmH\left(\wt{\bmD}\wt{\bsv}_{(n+1)}\right)\in\C^{2\wt{M}_1\times M_2}$ required to solving \cref{CadzowVS}.

Notice that, based on the framework of \cite{H.Attouch2010,C.Bao2016,J.Bolte2014,Y.Xu2013}, we can verify that the sequence generated by \cref{Alg2} globally converges to a stationary point provided that a linear operator $\wt{\bmA}$ satisfies the following:
\begin{align}\label{Assumption}
\wt{\bmA}\dde_0\neq\0,
\end{align}
where $\dde_0(\bk)=1$ if $\bk=\0$; $\dde_0(\bk)=0$ if $\bk\neq\0$. In words, since $\wh{1}=\delta$, \cref{Assumption} is equivalent to the assumption that every constant does not belong to the nullspace of an operator in the image domain, which is standard in the theoretical analysis of variational problems (e.g. \cite{G.Aubert2006}). In addition, \cref{Assumption} holds in many practical cases including the image deblurring and the compressed sensing MRI, as, in particular, the compressed sensing MRI is required to sample $\0$ frequency to avoid the DC-off set (e.g. \cite{M.Lustig2007}).

For the computational cost, we first mention that we do not explicitly formulate $\wt{\bsH}_{(n+1)}=\bmH\left(\wt{\bmD}\wt{\bsv}_{(n+1)}\right)$ at each iteration step. To simplify the discussion, we consider $K_1=K_2=K$, i.e. the square patch, without loss of generality. Then both $\wt{\bsH}_{(n+1)}\wt{\bsA}_{(n)}$ and $\wt{\bsH}_{(n+1)}^*\wt{\bC}_{(n+1)}$ are computed by using $4r$ two dimensional convolutions/fast Hankel matrix-vector multiplications \cite{L.Lu2015} directly from $\wt{\bmD}\wt{\bsv}_{(n+1)}$, requiring $O(r\wt{N}^2\log\wt{N})$ operations. In addition, we can update $\wt{\bsH}_{(n+1)}^*\wt{\bC}_{(n+1)}^{(:,m)}+\beta_3\beta^{-1}\wt{\bsA}_{(n)}^{(:,m)}$ directly after computing $\wt{\bmH}_{(n+1)}\wt{\bsA}_{(n)}^{(:,m)}$ in a single loop. Hence, together with the SVD of $\wt{\bsH}_{(n+1)}^*\wt{\bC}_{(n+1)}+\beta_3\beta^{-1}\wt{\bsA}_{(n)}\in\C^{K^2\times K^2}$ requiring $O(K^6)$ operations, we need $O(r\wt{N}^2\log\wt{N}+K^6)$ operations for \cref{Ctsubexplicit,Wsubexplicit}, which requires less computational costs than the direct computation of SVD in \cref{CadzowVS}.

For the memory storage, we do not need the full storage of $\wt{\bC}\in\C^{2\wt{M}_1\times M_2}$ and $\wt{\bsH}_{(n+1)}^*\wt{\bC}_{(n+1)}\in\C^{M_2\times M_2}$ either. From \cref{Ctsubexplicit}, both $\wt{\bC}_{(n+1)}$ and $\wt{\bsH}_{(n+1)}^*\wt{\bC}_{(n+1)}$ eventually contain at least $M_2-r$ zero column vectors. Hence, it suffices to store at most $r$ nonzero column vectors for these two matrices. Since $\wt{\bC}_{(n+1)}$ is the largest matrix to be stored, \cref{Alg2} requires $O(r\wt{N}^2)$ memory storages.

Finally, we mention that, we obtain low frequency samples $\wt{\bsv}$ by solving \cref{DDTFApproach} and set $\bmW$ via \cref{OurTightFrame} rather than directly setting $\wt{\bmW}$ in \cref{DDTFApproach} as an output. For the explanation of the reason, we consider the SVD
\begin{align}\label{HankelSVD}
\bmH\left(\wt{\bmD}\wt{\bsv}\right)=\wt{\bX}\wt{\Sig}\wt{\bY}^*=\left[\begin{array}{cc}
\wt{\bX}_{\|}&\wt{\bX}_{\perp}
\end{array}\right]\left[\begin{array}{cc}
\wt{\Sig}_{\|}&\0\\
\0&\wt{\Sig}_{\perp}
\end{array}\right]\left[\begin{array}{c}
\wt{\bY}_{\|}^*\\
\wt{\bY}_{\perp}^*
\end{array}\right].
\end{align}
In \cref{HankelSVD}, we have $\wt{\bX}_{\|}\in\C^{2\wt{M}_1\times r}$, $\wt{\bX}_{\perp}\in\C^{2\wt{M}_1\times(M_2-r)}$, $\wt{\Sig}_{\|}\in\R^{r\times r}$, $\wt{\Sig}_{\perp}\in\R^{(M_2-r)\times(M_2-r)}$, $\wt{\bY}_{\perp}\in\C^{M_2\times r}$, and $\wt{\bY}_{\perp}\in\C^{M_2\times(M_2-r)}$. If we set
\begin{align*}
\wt{\bC}=M_2^{-1/2}\left[\begin{array}{cc}
\wt{\bX}_{\|}\wt{\Sig}_{\|}&\0
\end{array}\right]:=M_2^{-1/2}\left[\begin{array}{cc}
\wt{\bX}_{\|}&\wt{\bX}_{\perp}
\end{array}\right]\left[\begin{array}{cc}
\wt{\Sig}_{\|}&\0\\
\0&\0
\end{array}\right],
\end{align*}
\cref{Wsub} in \cref{Alg2} (with $\beta_3=0$) gives the following SVD
\begin{align*}
\bmH\left(\wt{\bmD}\wt{\bsv}\right)^*\wt{\bC}=\left[\begin{array}{cc}
\wt{\bY}_{\|}&\wt{\bY}_{\perp}
\end{array}\right]\left[\begin{array}{cc}
M_2^{-1/2}\wt{\Sig}_{\|}^2&\0\\
\0&\0
\end{array}\right]=\wt{\bY}\left[\begin{array}{cc}
M_2^{-1/2}\wt{\Sig}_{\|}^2&\0\\
\0&\0
\end{array}\right]\bsI.
\end{align*}
This means that \cref{ctsubprob,Wsub} is closely related to the rank $r$ approximation of $\bmH\left(\wt{\bmD}\wt{\bsv}\right)$ with the right singular vectors and the singular values as outputs. However, the assumption that the underlying image is piecewise constant with edges described by a trigonometric curve is approximately true in practice \cite{G.Ongie2016}. Since the low-rank approximation will strongly force edges to be given as a trigonometric curve in our setting, the use of such truncated singular values will be too restrictive for the practical applications. Hence, we restore low frequency samples $\wt{\bsv}$ from \cref{DDTFApproach} for the sake of flexibility.

\section{Numerical results}\label{Experiments}

In this section, we conduct some numerical simulations on the following image restoration tasks:
\begin{enumerate}
\item image restoration from partial random Fourier samples (RS),
\item ideal low-pass filter deconvolution (ILF).
\end{enumerate}
Specifically, to compare the performance of the proposed model \cref{ProposedModel} over the existing approaches, we choose to compare the proposed model
\begin{align}\label{ProposedModel2}
\min_{\bsv}~\f{1}{2}\left\|\bmA\bsv-\bsf\right\|_2^2+\left\|\gga\cdot\bmW\left(\bmD\bsv\right)\right\|_1,
\end{align}
with the LSLP model \cite{G.Ongie2016}
\begin{align}\label{LSLP2}
\min_{\bsv}~\left\|\bmA\bsv-\bsf\right\|_2^2+\gamma\sum_{m=r+1}^{M_2}\left\|\bmS_{\a_m(-\cdot)}\left(\bmD\bsv\right)\right\|_2^2,
\end{align}
the LRHDDTF model \cite{J.F.Cai2020}
\begin{align}\label{LRHDDTF2}
\begin{split}
&~~~\min_{\bsv,\bc,\bmW}~\f{1}{2}\left\|\bmA\bsv-\bsf\right\|_2^2+\f{\mu}{2}\left\|\bmW\left(\bmD\bsv\right)-\bc\right\|_2^2+\left\|\gga\cdot\bc\right\|_0\\
&\text{subject to}~~\bmW^*\bmW=\bmI,
\end{split}
\end{align}
and the following Schatten $0$-norm relaxation \cite{G.Ongie2017} of a direct rank minimization
\begin{align}\label{Schatten0}
\min_{\bsv}~\f{1}{2}\left\|\bmA\bsv-\bsf\right\|_2^2+\gamma\left\|\bmH\left(\bmD\bsv\right)\right\|_0.
\end{align}
In \cref{Schatten0},  $\left\|\bZ\right\|_0$ is the Schatten $0$-norm of a matrix $\bZ$ defined as
\begin{align*}
\left\|\bZ\right\|_0=\ln\det\left(\left(\bZ^*\bZ\right)^{1/2}+\eps\bsI\right)
\end{align*}
with a small $\eps>0$. Finally, to further study the improvements over the conventional sparse regularizations, we compare with the total variation (TV) model \cite{L.I.Rudin1992}
\begin{align}\label{TV}
\min_{\bu}~\f{1}{2}\left\|\bmA\bu-\bsf\right\|_2^2+\left\|\gga\cdot\Na\bu\right\|_1,
\end{align}
and the wavelet frame model (e.g. \cite{J.F.Cai2009/10})
\begin{align}\label{Frame}
\min_{\bu}\f{1}{2}\left\|\bmA\bu-\bsf\right\|_2^2+\left\|\gga\cdot\bmW\bu\right\|_1.
\end{align}
Recall that, with a slight abuse of notation, the forward operator $\bmA$ in \cref{ProposedModel2,LSLP2,LRHDDTF2,Schatten0} acts on the frequency domain while that in \cref{TV,Frame} acts on the image domain, depending on the context. All experiments are implemented on MATLAB $\mathrm{R}2014\mathrm{a}$ running on a laptop with $64\mathrm{GB}$ RAM and Intel(R) Core(TM) CPU $\mathrm{i}7$-$8750\mathrm{H}$ at $2.20\mathrm{GHz}$ with $6$ cores.

\begin{figure}[t]
\centering
\subfloat[E]{\label{Ellipse}\begin{minipage}{2.250cm}
\includegraphics[width=2.250cm]{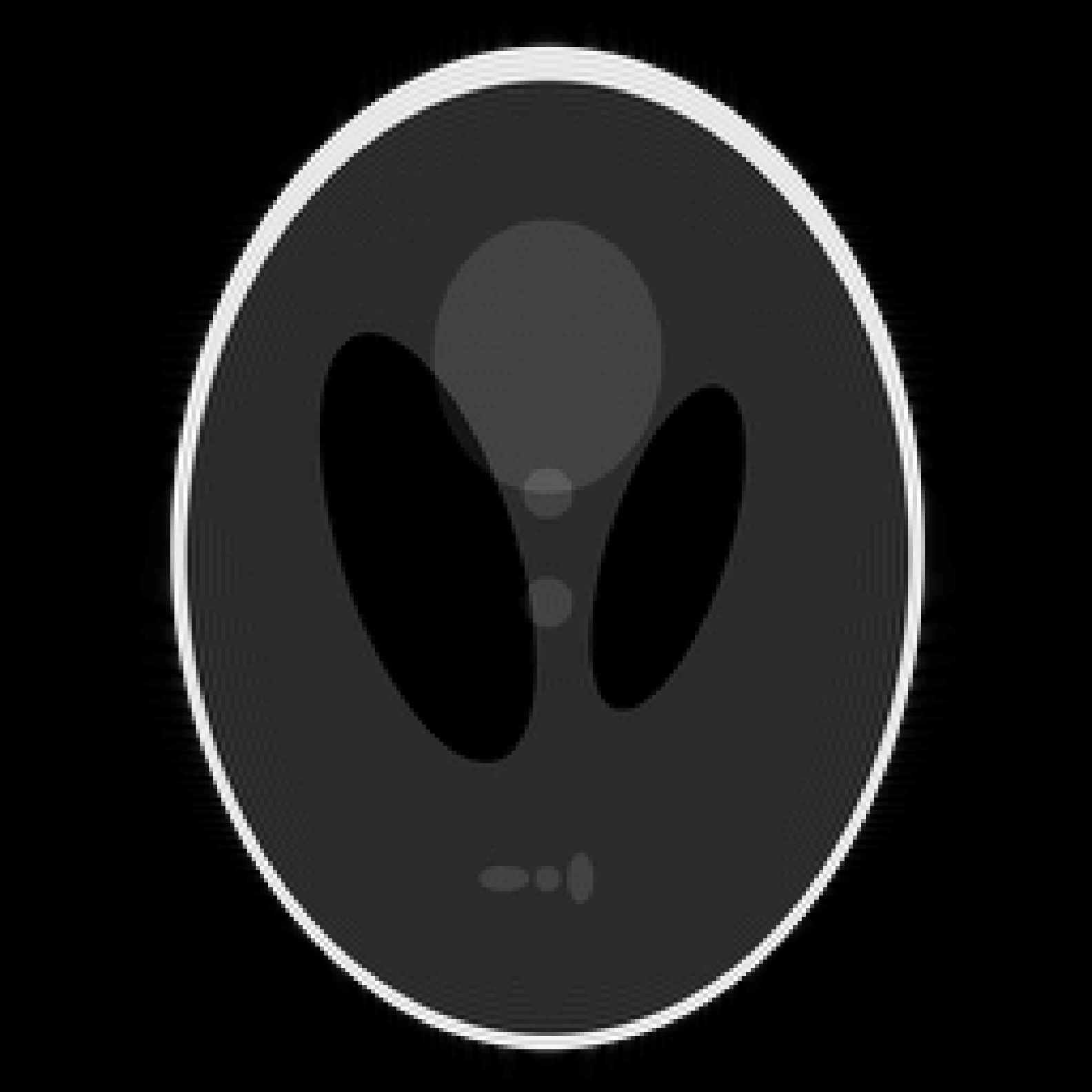}\vspace{0.025cm}\\
\includegraphics[width=2.250cm]{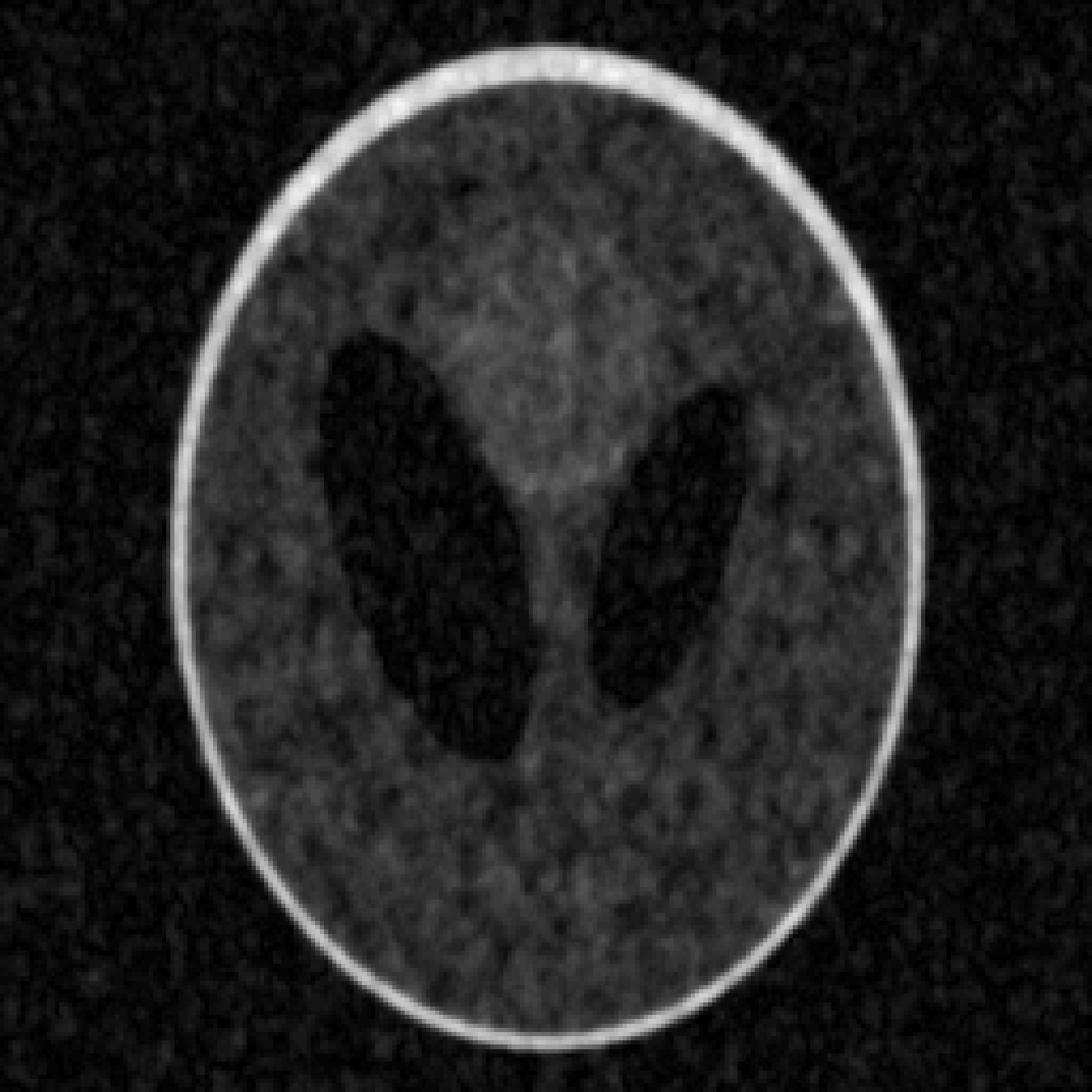}\vspace{0.025cm}\\
\includegraphics[width=2.250cm]{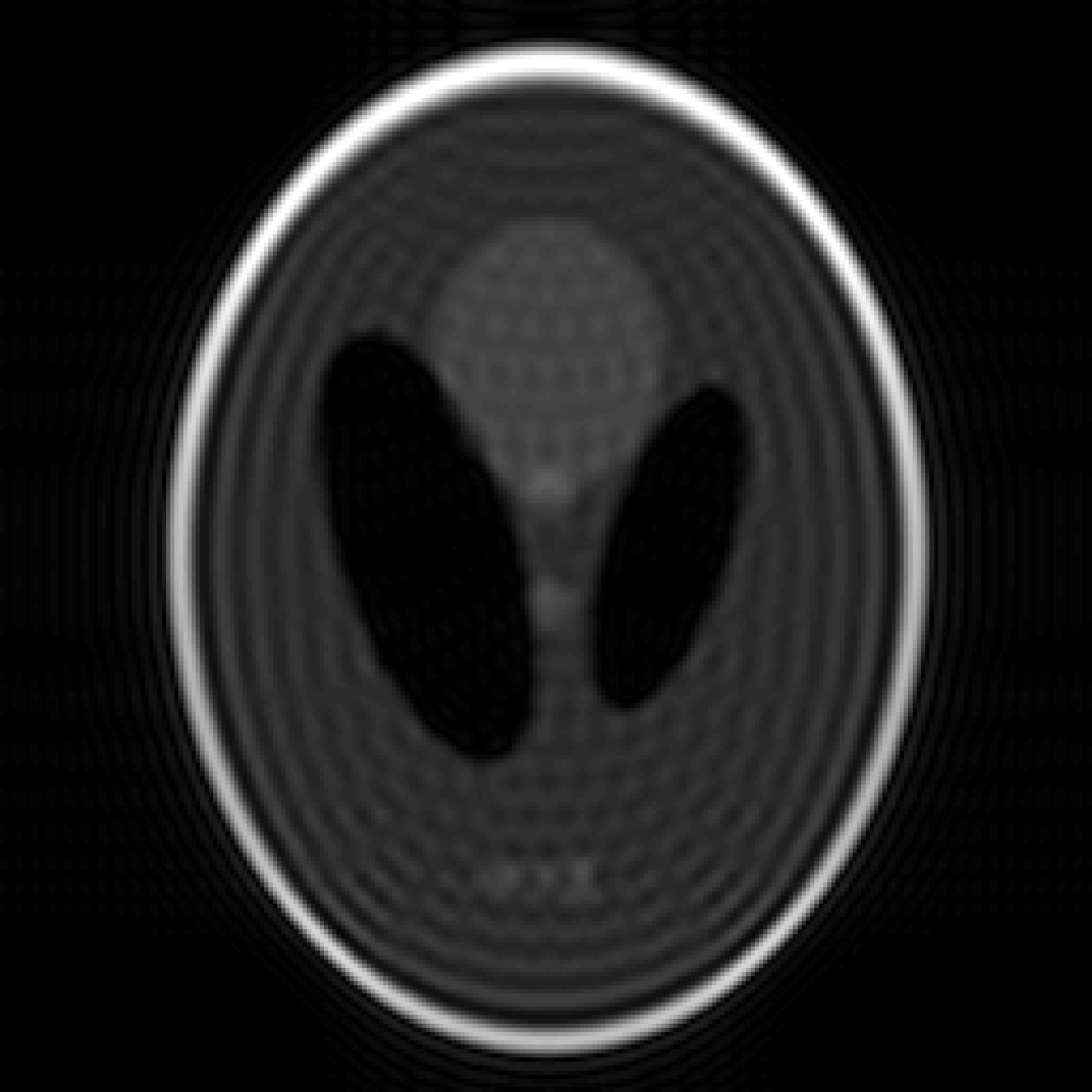}
\end{minipage}}\hspace{-0.05cm}
\subfloat[R]{\label{Rectangle}\begin{minipage}{2.250cm}
\includegraphics[width=2.250cm]{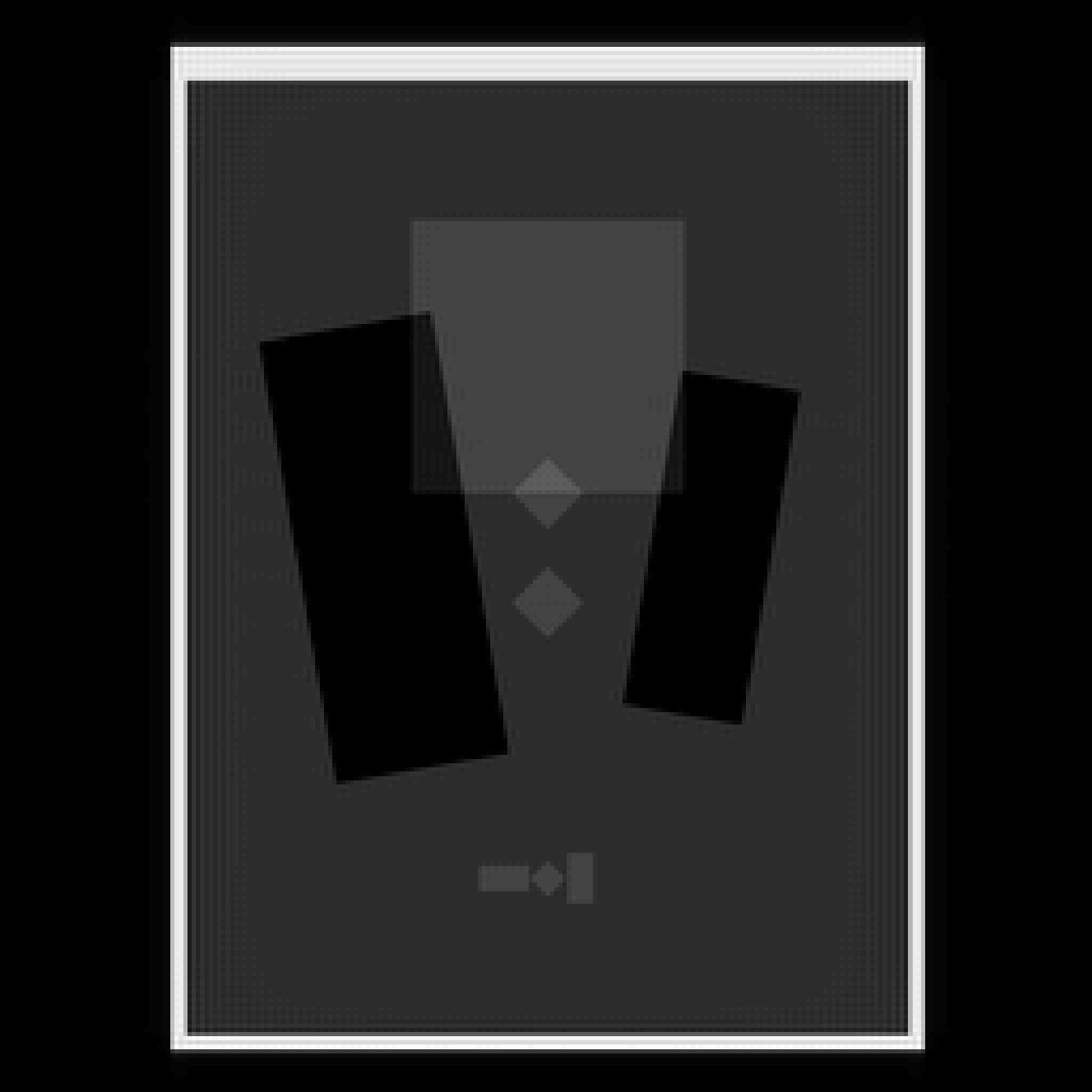}\vspace{0.025cm}\\
\includegraphics[width=2.250cm]{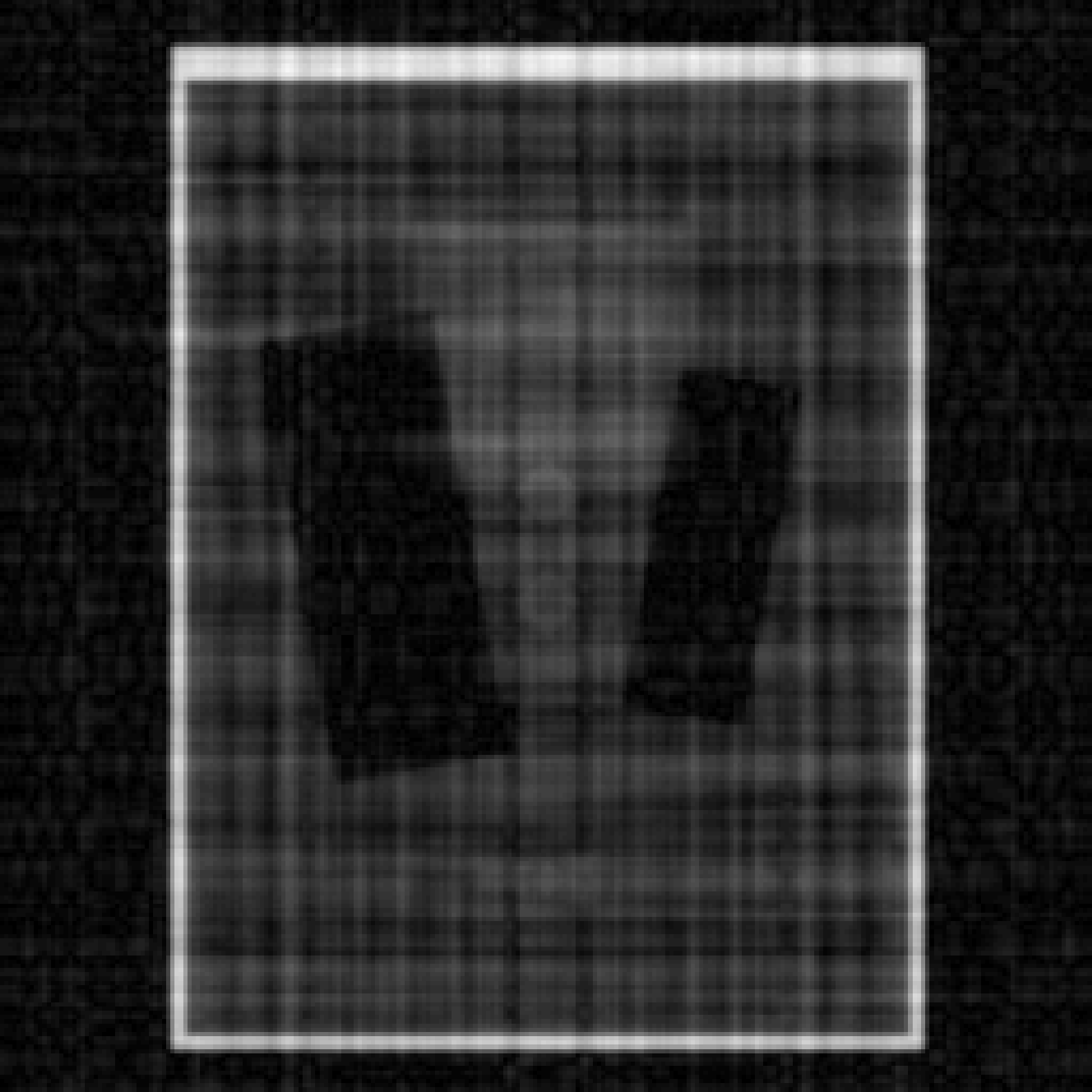}\vspace{0.025cm}\\
\includegraphics[width=2.250cm]{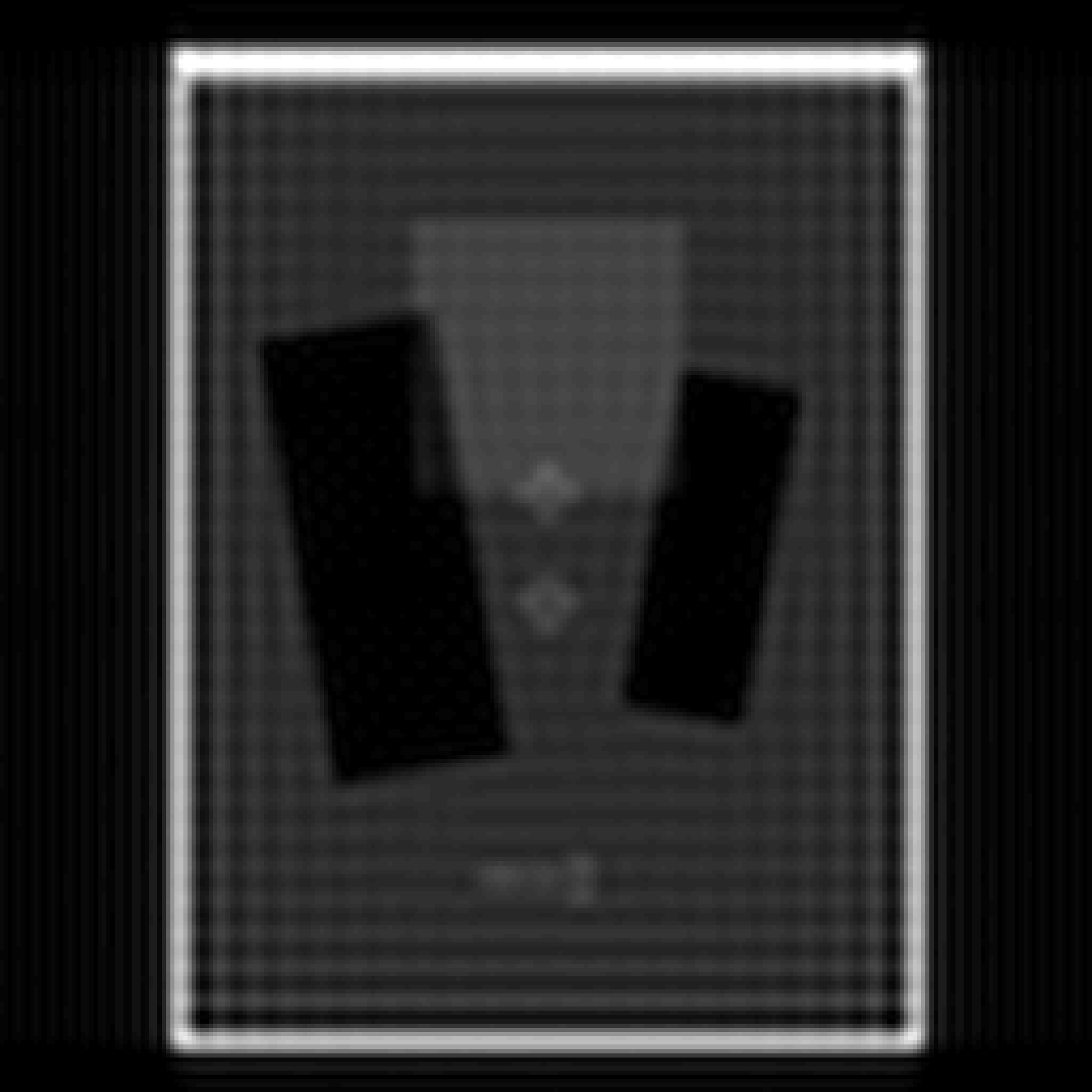}
\end{minipage}}\hspace{-0.05cm}
\subfloat[AB]{\label{AngryBirds}\begin{minipage}{2.250cm}
\includegraphics[width=2.250cm]{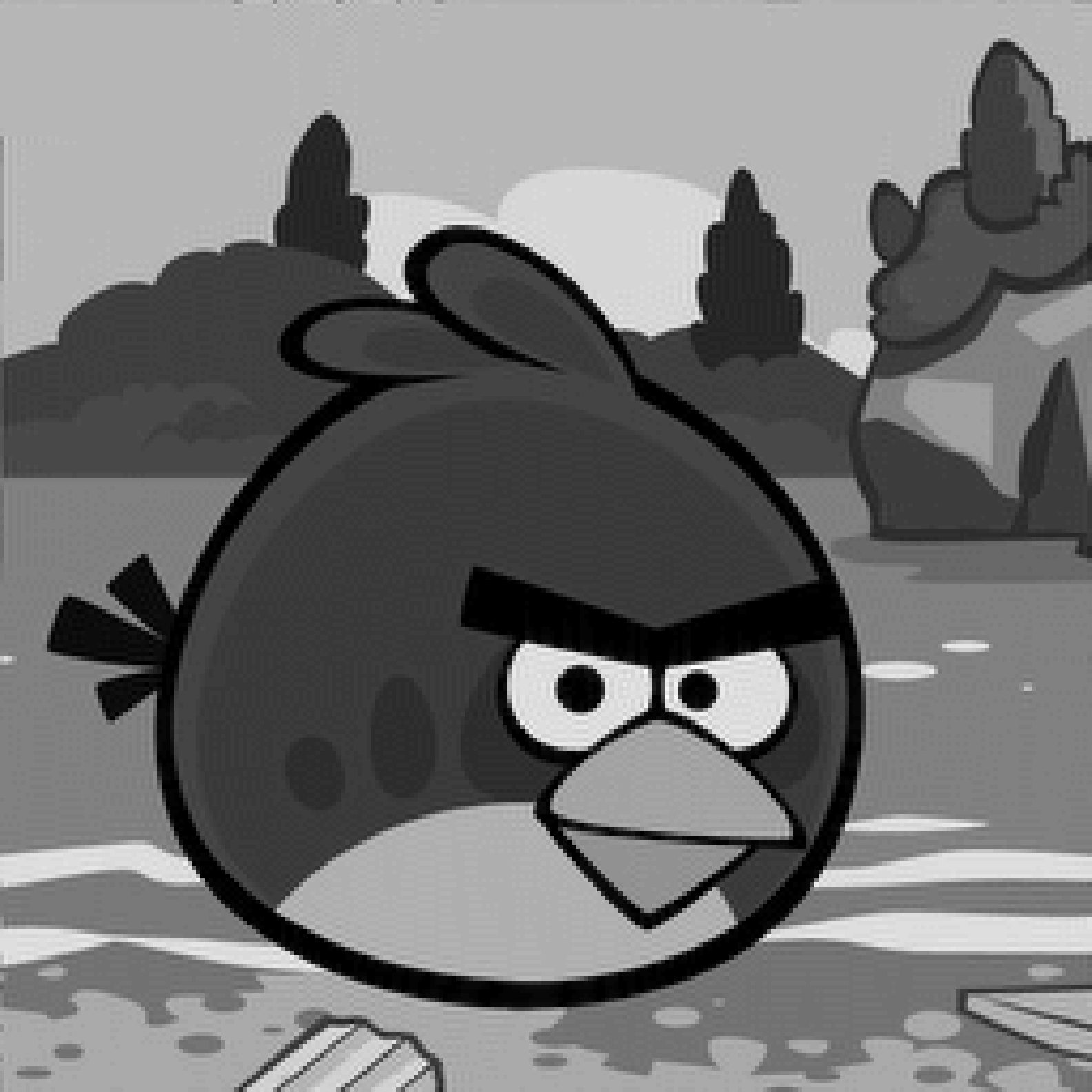}\vspace{0.025cm}\\
\includegraphics[width=2.250cm]{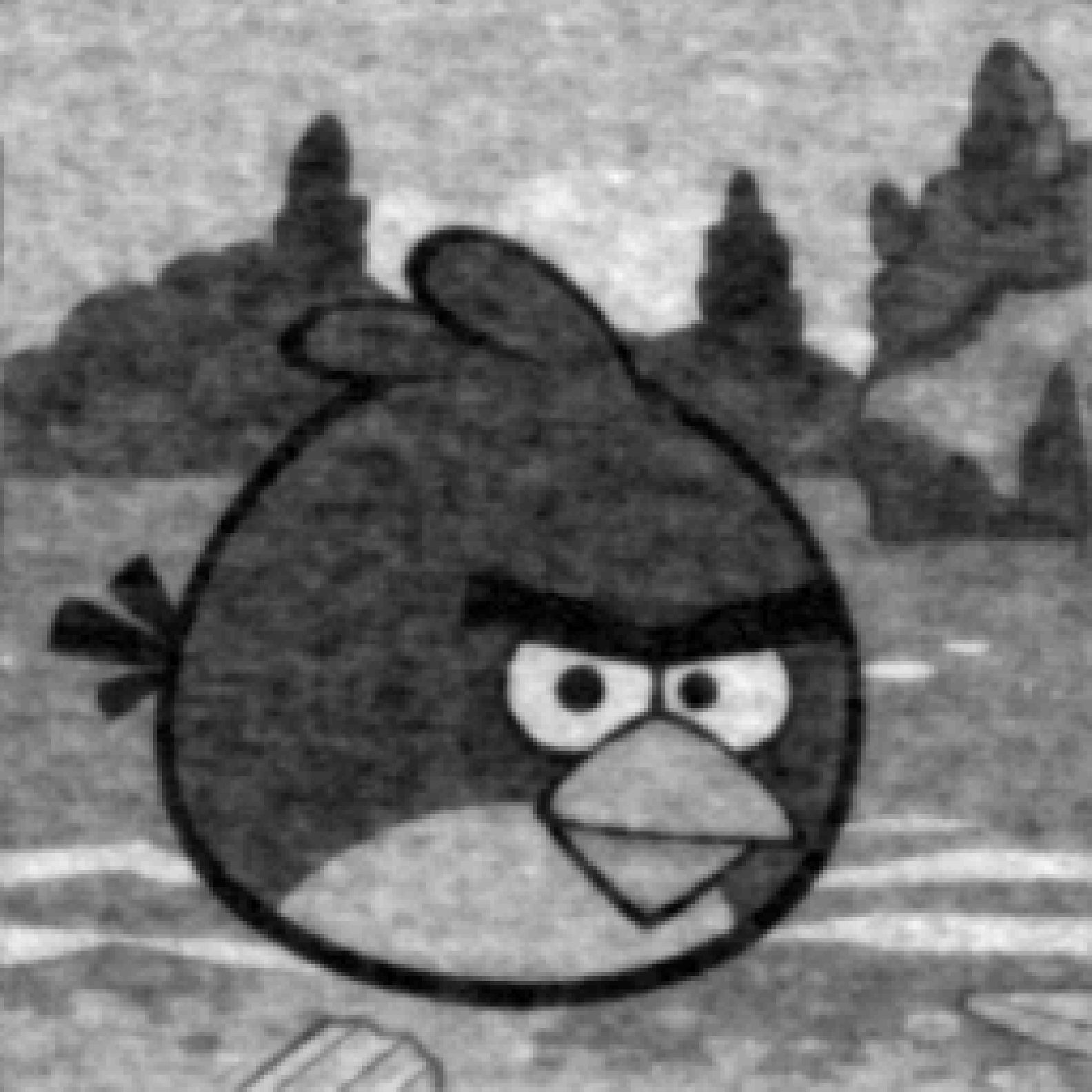}\vspace{0.025cm}\\
\includegraphics[width=2.250cm]{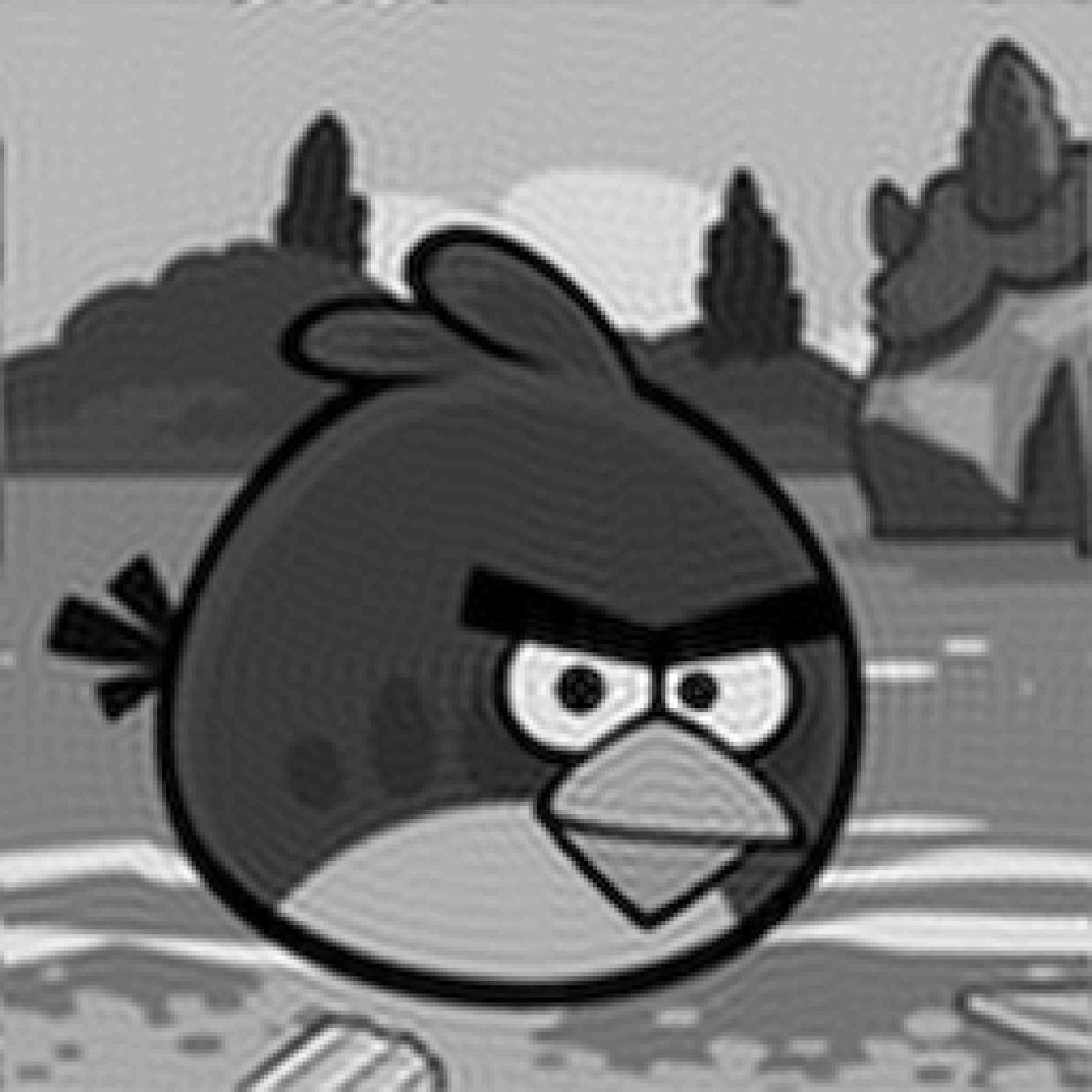}
\end{minipage}}\hspace{-0.05cm}
\subfloat[SM]{\label{Mario}\begin{minipage}{2.250cm}
\includegraphics[width=2.250cm]{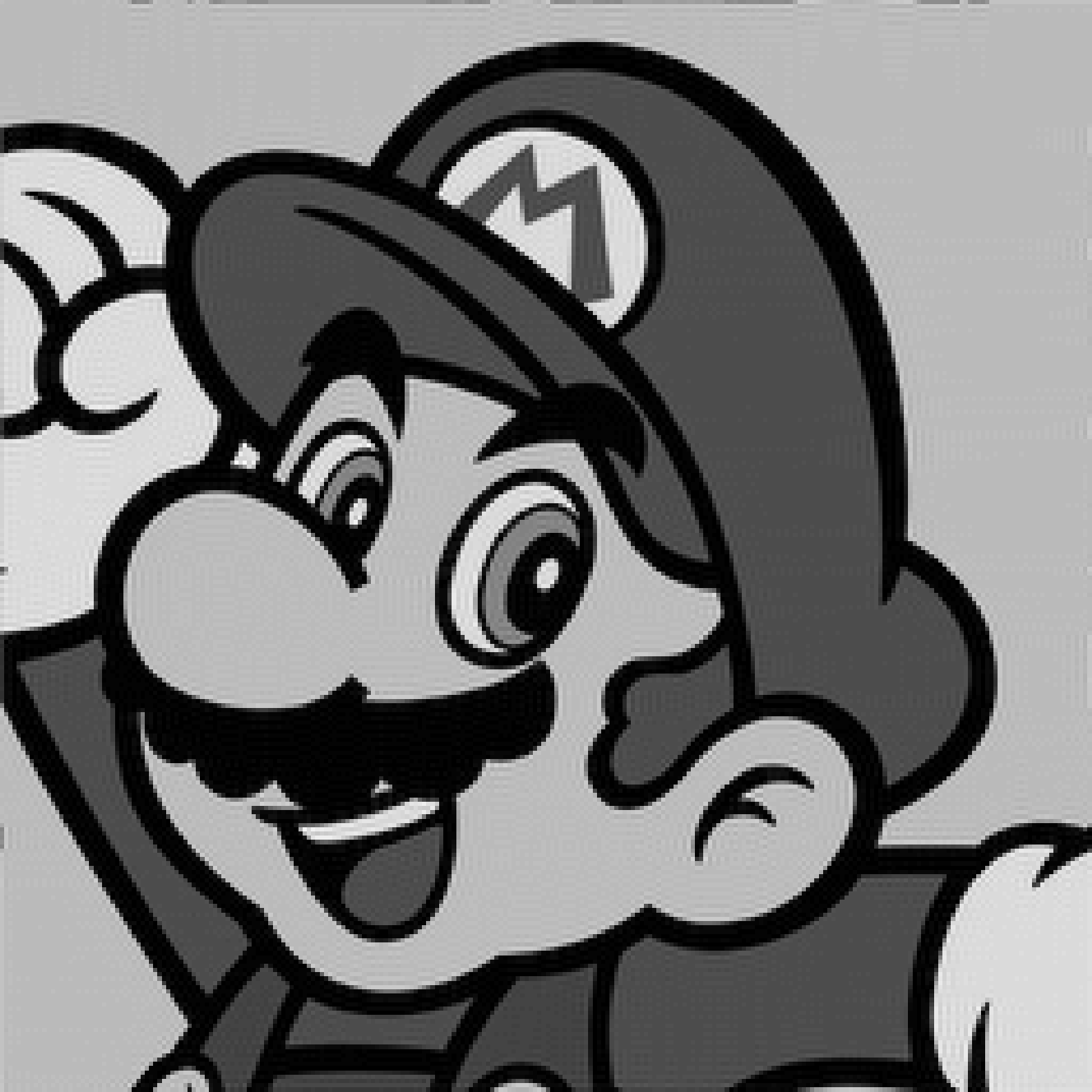}\vspace{0.025cm}\\
\includegraphics[width=2.250cm]{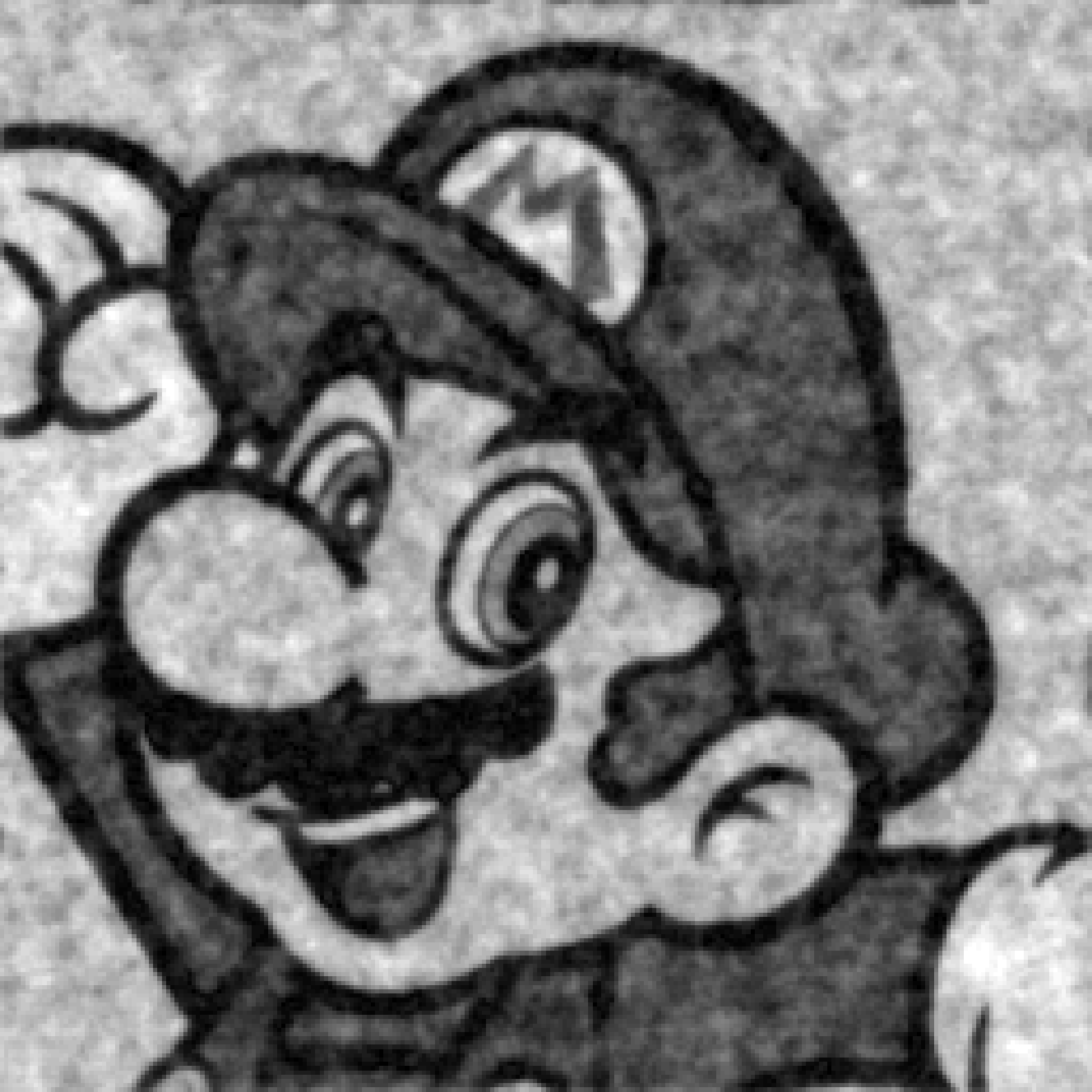}\vspace{0.025cm}\\
\includegraphics[width=2.250cm]{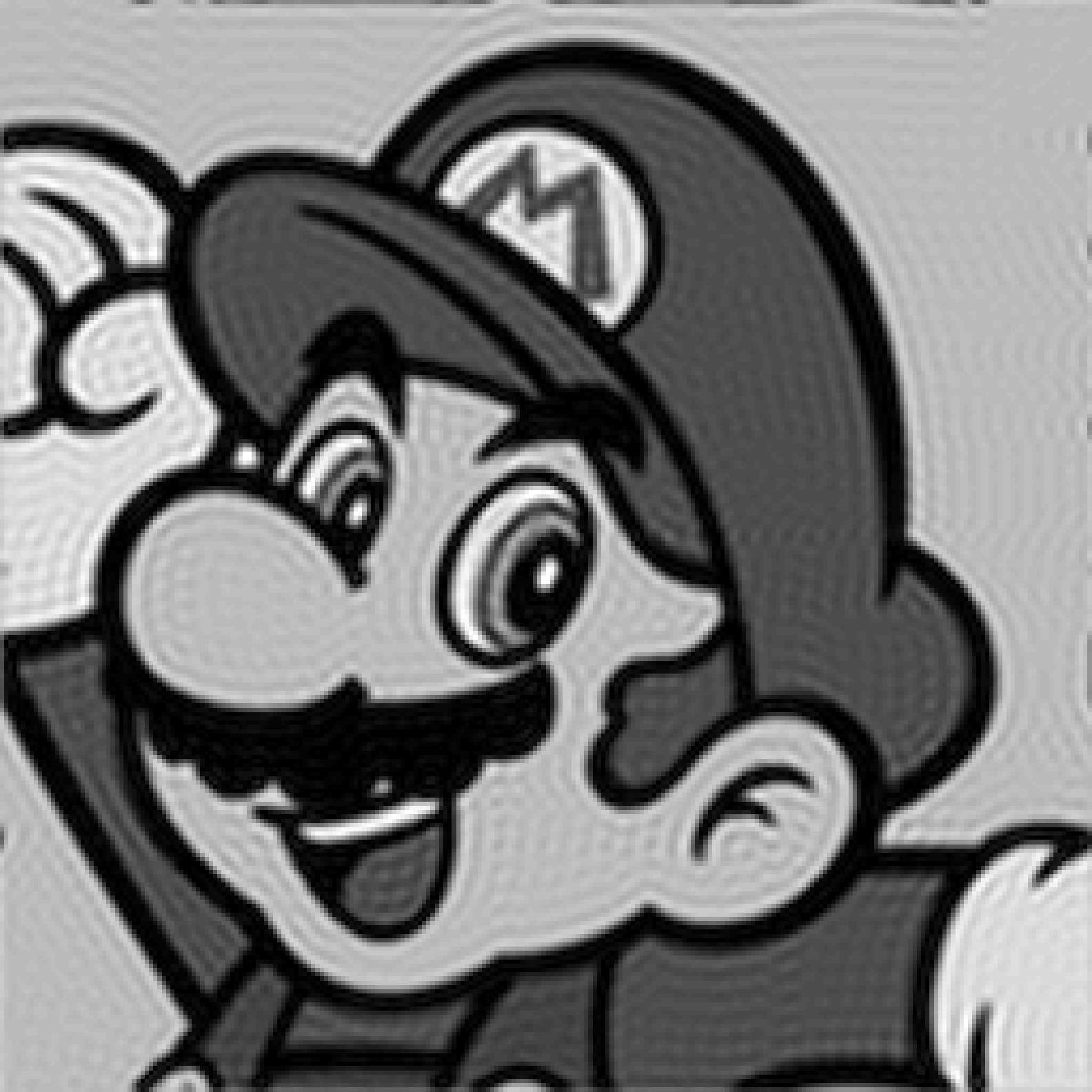}
\end{minipage}}\hspace{-0.05cm}
\subfloat[MC]{\label{MC}\begin{minipage}{2.250cm}
\includegraphics[width=2.250cm]{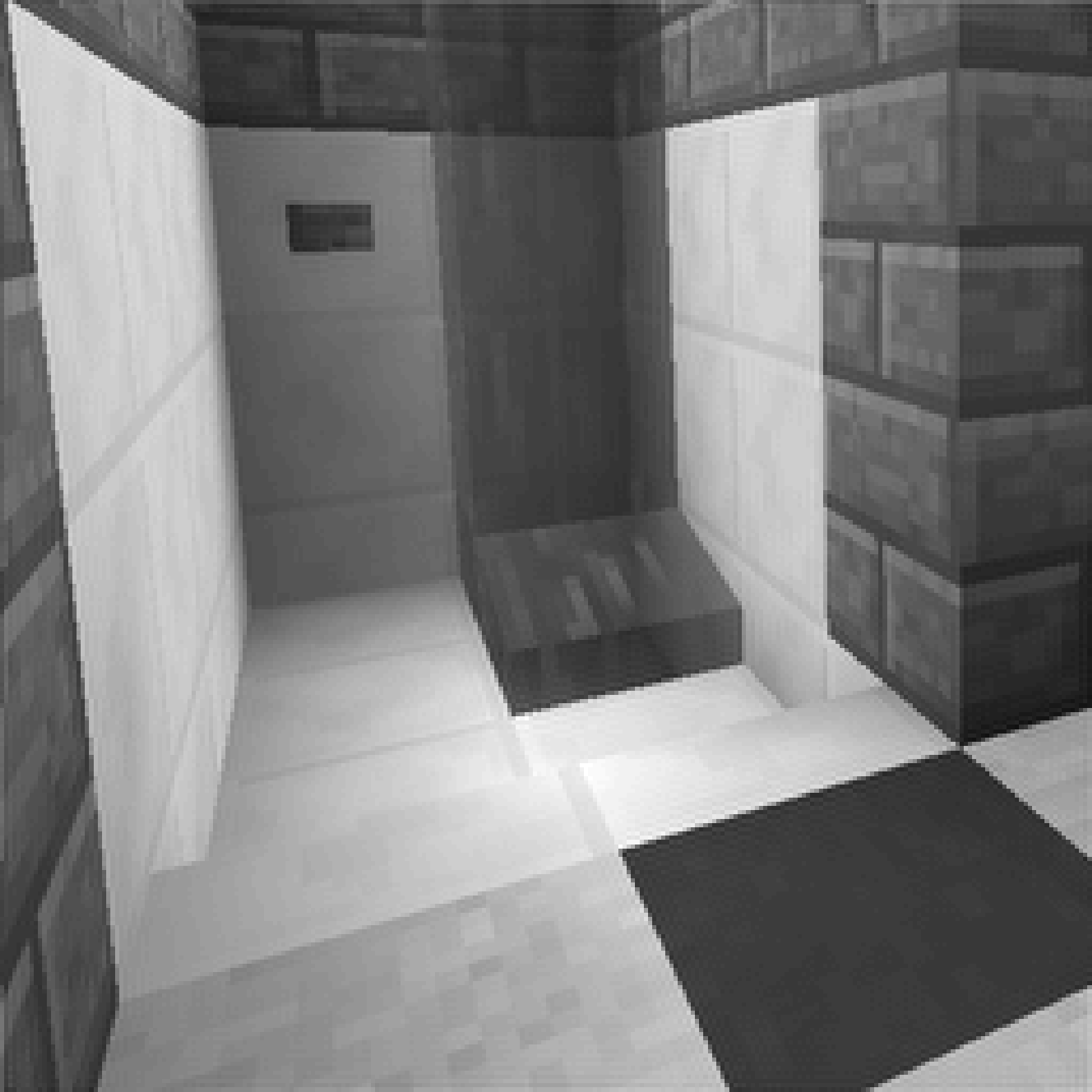}\vspace{0.025cm}\\
\includegraphics[width=2.250cm]{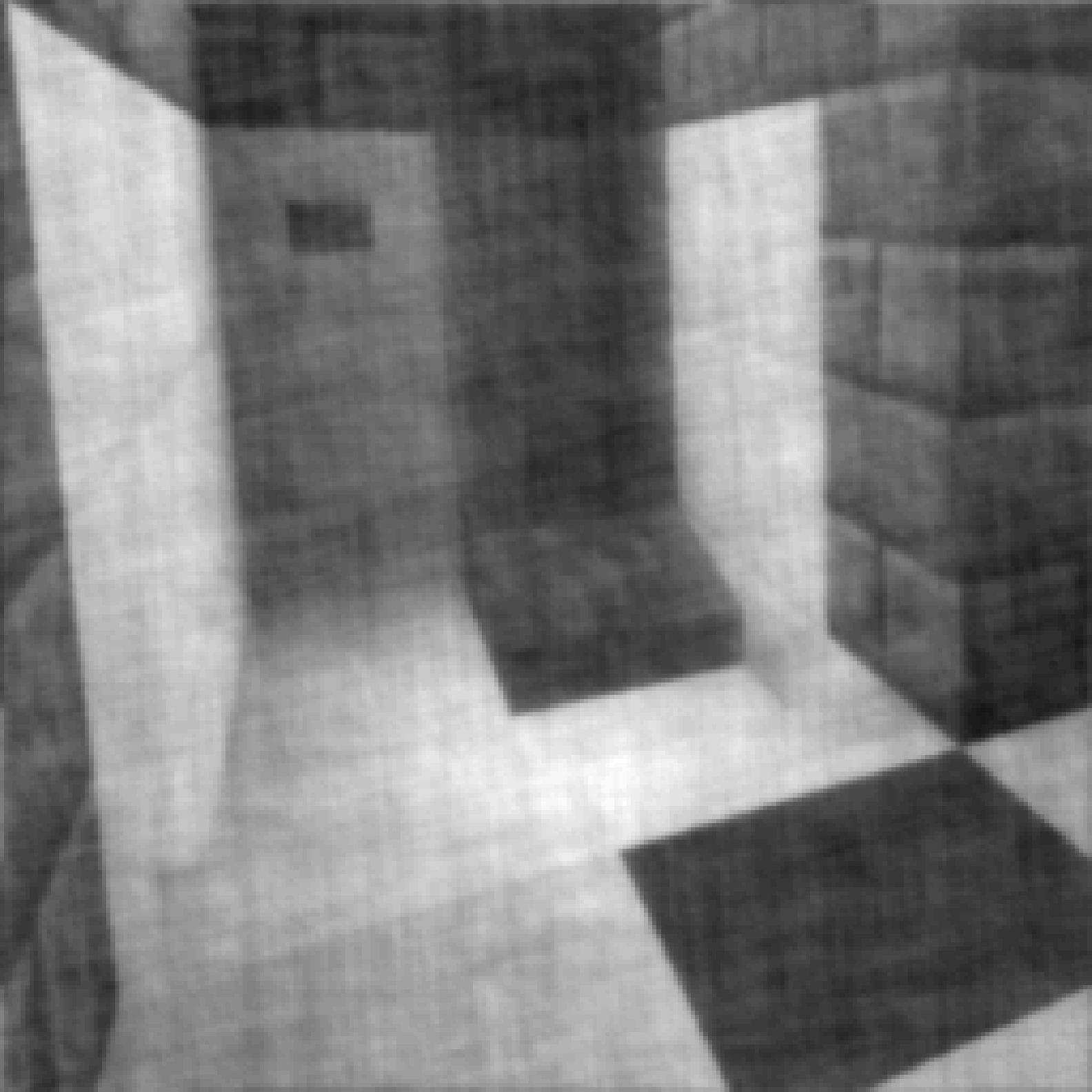}\vspace{0.025cm}\\
\includegraphics[width=2.250cm]{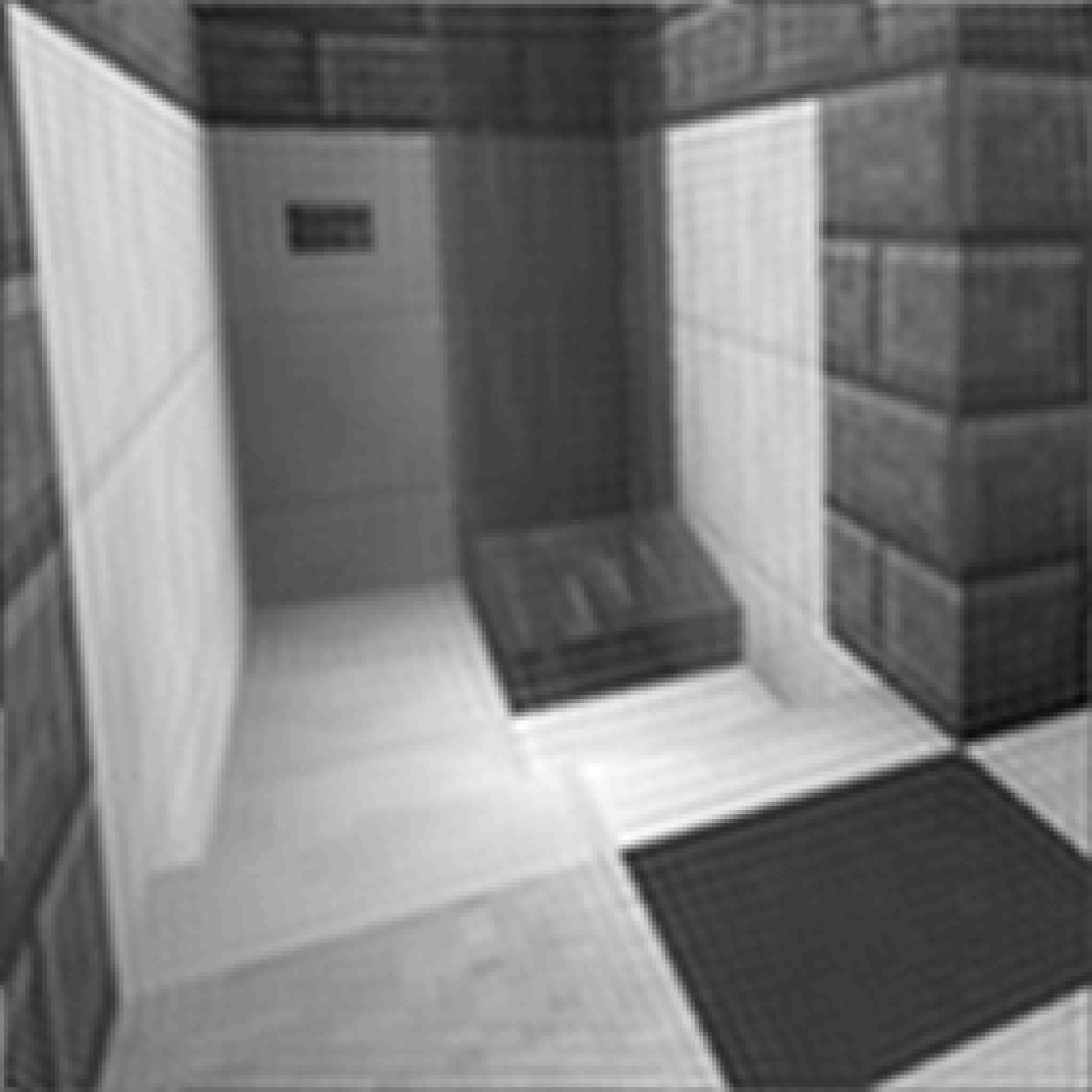}
\end{minipage}}
\caption{Visualization of test images (first row) and observed images (second to fourth rows). Second row corresponds to the random sampling, third row to the ideal low-pass filtering. All images are displayed in the window level $[0,1]$ for the fair comparisons.}\label{OriginalImages}
\end{figure}

Throughout the experiments, we test two synthetic images (``Ellipse (E)'' and ``Rectangle (R)'') and three non-synthetic images (``Angry Birds (AB)'', ``Super Mario (SM)'', and ``Minecraft (MC)'') with the size $256\times256$, taking the values in $[0,1]$, as shown in \cref{OriginalImages}. For the random sampling, the data $\bsf$ is synthesized by randomly sampling $20\%$ of $256\times256$ Fourier samples via the variable density sampling method described in \cite{M.Lustig2007}. For the ideal low-pass filtering, we take $65\times65$ low frequency Fourier samples for the synthetic images and $101\times101$ for the non-synthetic images. In any image restoration task, a mild white Gaussian noise is also added to generate noisy measurements. For \cref{ProposedModel2,LSLP2,LRHDDTF2,Schatten0}, we use the $K\times K$ square patch for simplicity. Specifically, we choose $K=33$ and $r=330$ for ``Ellipse'', $K=25$ and $r=285$ for the ``Rectangles'', and $K=51$ and $r$ to be around $60\sim70\%$ of the size of filters for the non-synthetic images (``Angry Birds'', ``Super Mario'', and ``Minecraft''), depending on the geometry of the target images. The tight frame transform in \cref{ProposedModel2} and the filters in \cref{LSLP2}, including the initial tight frame transform in \cref{LRHDDTF2} are obtained by solving \cref{Alg2}. Specifically, we use \cref{Alg2} to restore $65\times65$ low frequency samples for ``Ellipse'', $49\times49$ for ``Rectangle'', and $101\times101$ for non-synthetic images. Then using the pre-restored low frequency samples $\wt{\bsv}$, we compute the SVD of $\bmH\left(\wt{\bmD}\wt{\bsv}\right)$, which will be used for the filter banks in \cref{ProposedModel2,LSLP2} and the initialization of \cref{LRHDDTF2}. For \cref{TV}, we use the forward difference with the periodic boundary condition for the difference operator $\Na$. In \cref{Frame}, two types of $\bmW$ are used; one is the undecimated tensor product Haar framelet transform with $1$ level of decomposition \cite{B.Dong2013}, and the other is the data-driven tight frame (DDTF) in \cite{W.Zhou2013}. In particular, we use the method in \cite{J.F.Cai2014} by using $4\times4$ local discrete cosine transform filters and a reference image restored by the Haar framelet model. In all of the experiments, we have manually tuned the regularization parameters to achieve the optimal restoration results in each scenario. In particular, since it has been demonstrated in \cite{G.Ongie2016} that the regularization term weighted by the edge information shows advantages over the one without the weights, we compute \cref{Pseudospectra} using the filters obtained from \cref{Alg2} for the weights in \cref{TV,Frame}.In all of the experiments, we have manually tuned the regularization parameters to achieve the optimal restoration results in each scenario. For the quantitative comparison, we compute the signal-to-noise ratio (SNR), the high frequency error norm (HFEN) \cite{S.Ravishankar2011}, and the structure similarity index map (SSIM) \cite{Z.Wang2004}. Recall that for \cref{ProposedModel2,LSLP2,LRHDDTF2,Schatten0}, the restored image is computed via the inverse DFT of the restored Fourier samples.

\subsection{Edge estimation}\label{EdgeEstimation}

We investigate the edge estimation of \cref{CadzowVS,DDTFApproach} for the synthetic images (``Ellipse'' and ``Rectangle''). Specifically, for each image restoration task, we solve \cref{CadzowVS} by using the alternating minimization method in \cite{G.Ongie2016} and \cref{DDTFApproach} by \cref{Alg2} to obtain low frequency samples $\wt{\bsv}$. Specifically, we record the number of iterations and time to reach the stopping criterion
\begin{align}\label{SCEdgeEstimation}
\f{\left\|\wt{\bsv}_{(n+1)}-\wt{\bsv}_{(n)}\right\|_2}{\left\|\wt{\bsv}_{(n+1)}\right\|_2}\leq5\times10^{-4}.
\end{align}
Using the singular value decompositions $\bmH\left(\wt{\bmD}\wt{\bsv}\right)=\wt{\bX}\wt{\Sig}\wt{\bY}^*$, we compute \cref{Pseudospectra} to compare the edge sets.

\begin{table}[t]
\centering
\begin{tabular}{|c||c|c|c|}
\hline
\multicolumn{4}{|c|}{Ellipse}\\ \hline
Method&Index&Random Sampling&Ideal Low Pass Filter\\ \hline
\multirow{2}{*}{Cadzow \cref{CadzowVS}}&NIters&$59$&$4$\\ \cline{2-4}
&Time&$547.4\mathrm{sec}$&$46.9\mathrm{sec}$\\ \hline
\multirow{2}{*}{DDTF \cref{DDTFApproach}}&NIters&$61$&$4$\\ \cline{2-4}
&Time&$91.7\mathrm{sec}$&$14.6\mathrm{sec}$\\ \hline\hline
\multicolumn{4}{|c|}{Rectangle}\\ \hline
Method&Index&Random Sampling&Ideal Low Pass Filter\\ \hline
\multirow{2}{*}{Cadzow \cref{CadzowVS}}&NIters&$59$&$3$\\ \cline{2-4}
&Time&$215.6\mathrm{sec}$&$15.2\mathrm{sec}$\\ \hline
\multirow{2}{*}{DDTF \cref{DDTFApproach}}&NIters&$62$&$3$\\ \cline{2-4}
&Time&$27.7\mathrm{sec}$&$5.0\mathrm{sec}$\\ \hline
\end{tabular}
\caption{The number of iterations and the time to reach the stopping criterion for the edge estimation}\label{TableEdge}
\end{table}

\begin{figure}[t]
\centering
\subfloat[Ref.]{\label{EllipseEdgeOriginal}\begin{minipage}{3cm}
\includegraphics[width=3cm]{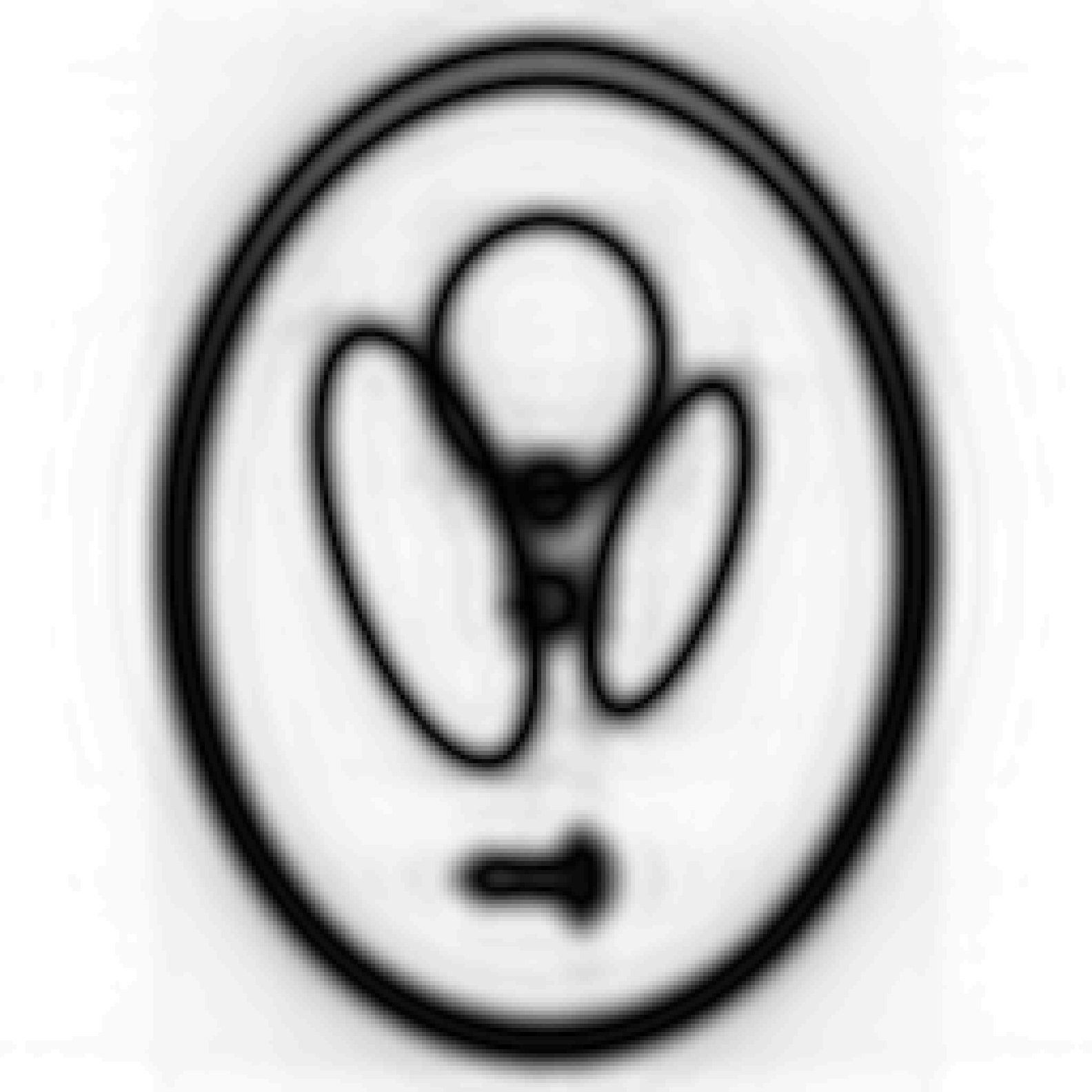}\vspace{0.025cm}\\
\includegraphics[width=3cm]{EllipsesEdgeOriginal.pdf}
\end{minipage}}\hspace{-0.05cm}
\subfloat[RS]{\label{EllipseRandomEdge}\begin{minipage}{3cm}
\includegraphics[width=3cm]{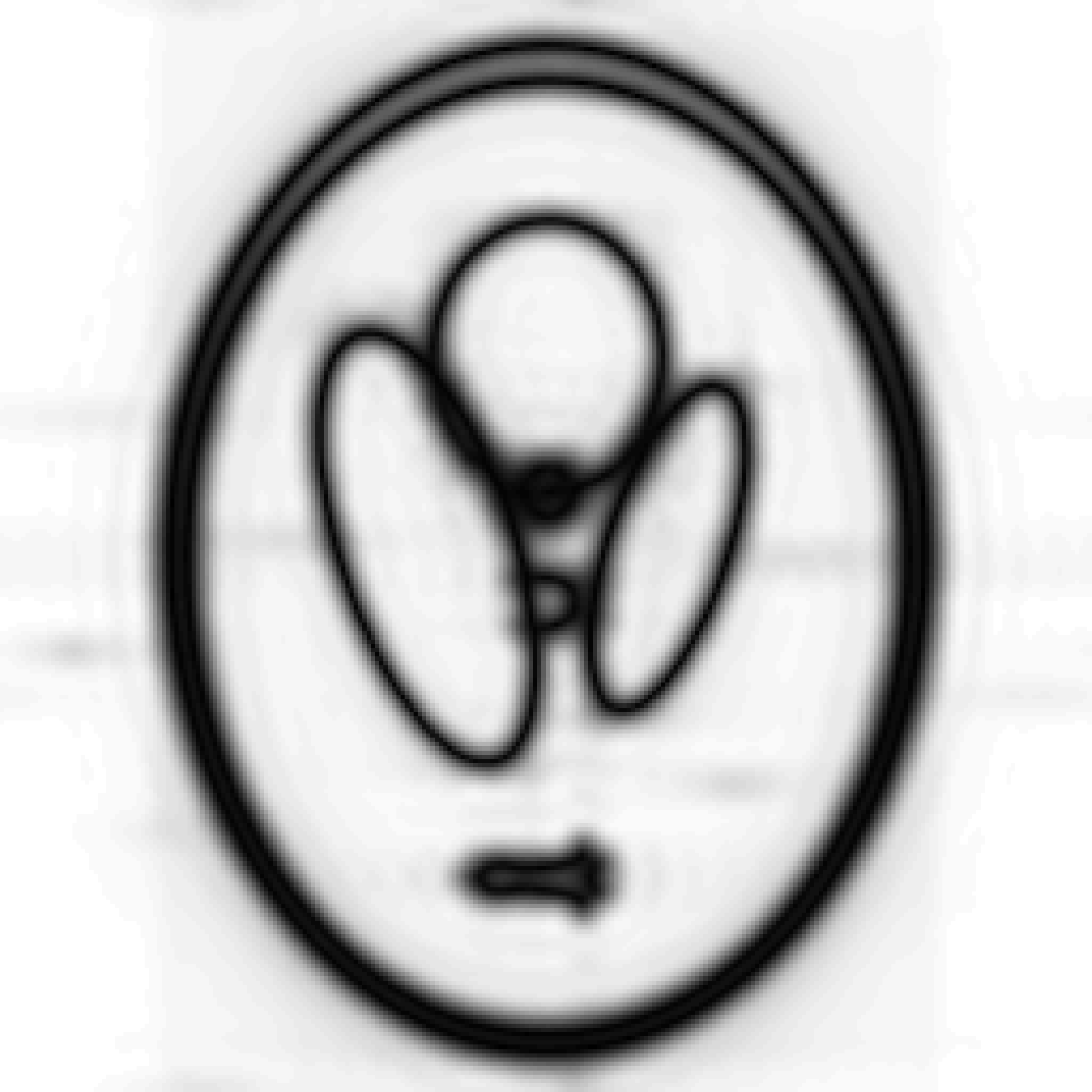}\vspace{0.025cm}\\
\includegraphics[width=3cm]{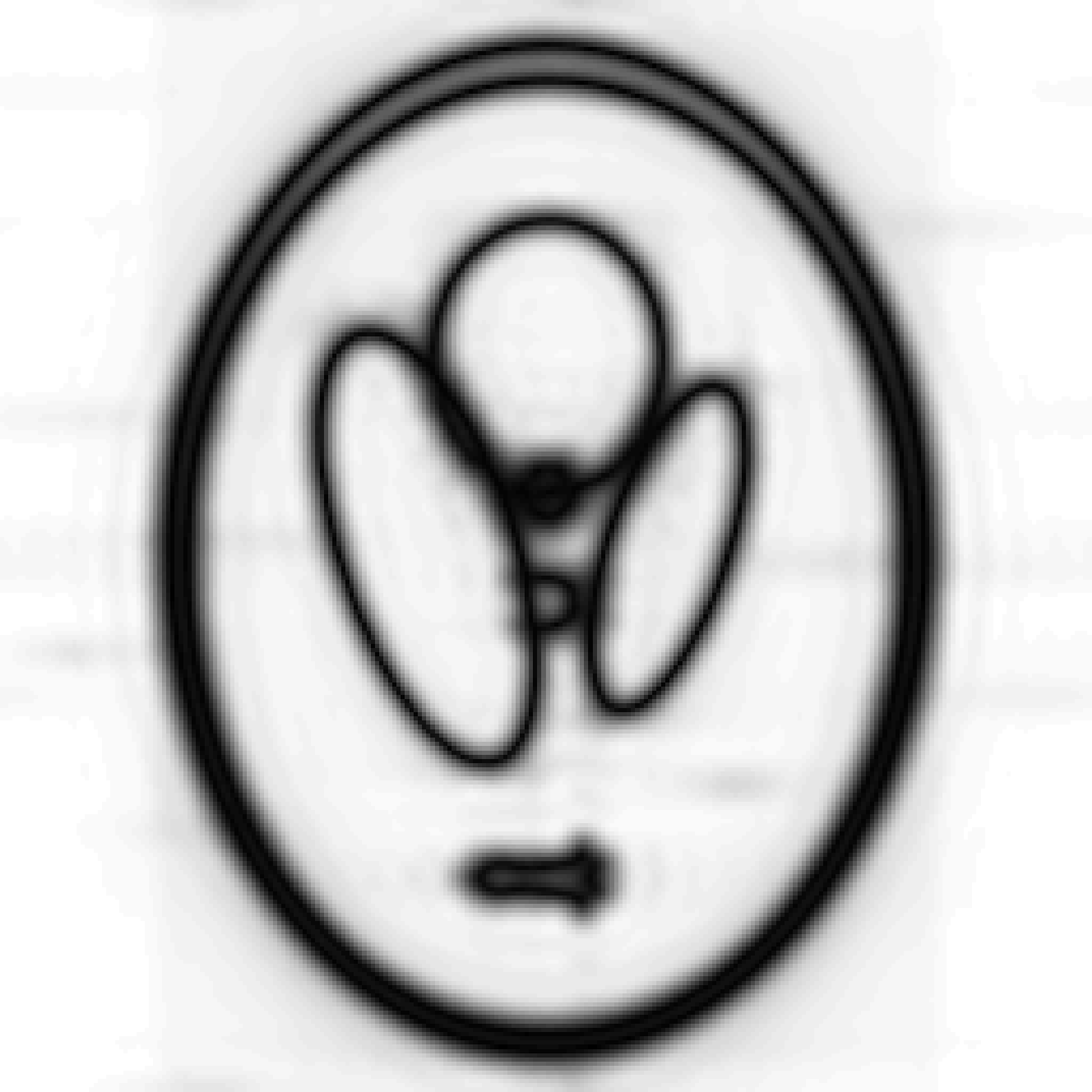}
\end{minipage}}\hspace{-0.05cm}
\subfloat[ILF]{\label{EllipseILFEdge}\begin{minipage}{3cm}
\includegraphics[width=3cm]{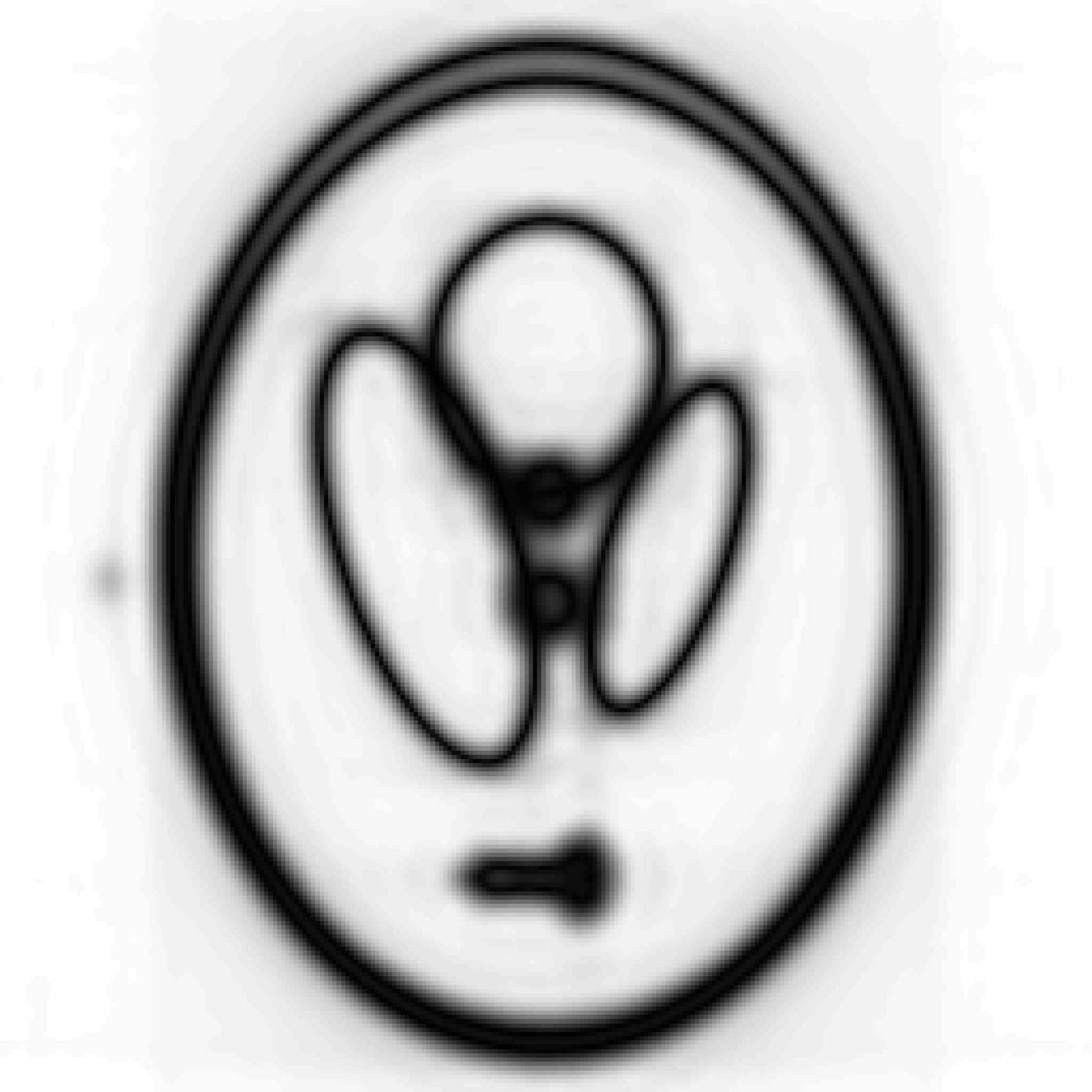}\vspace{0.025cm}\\
\includegraphics[width=3cm]{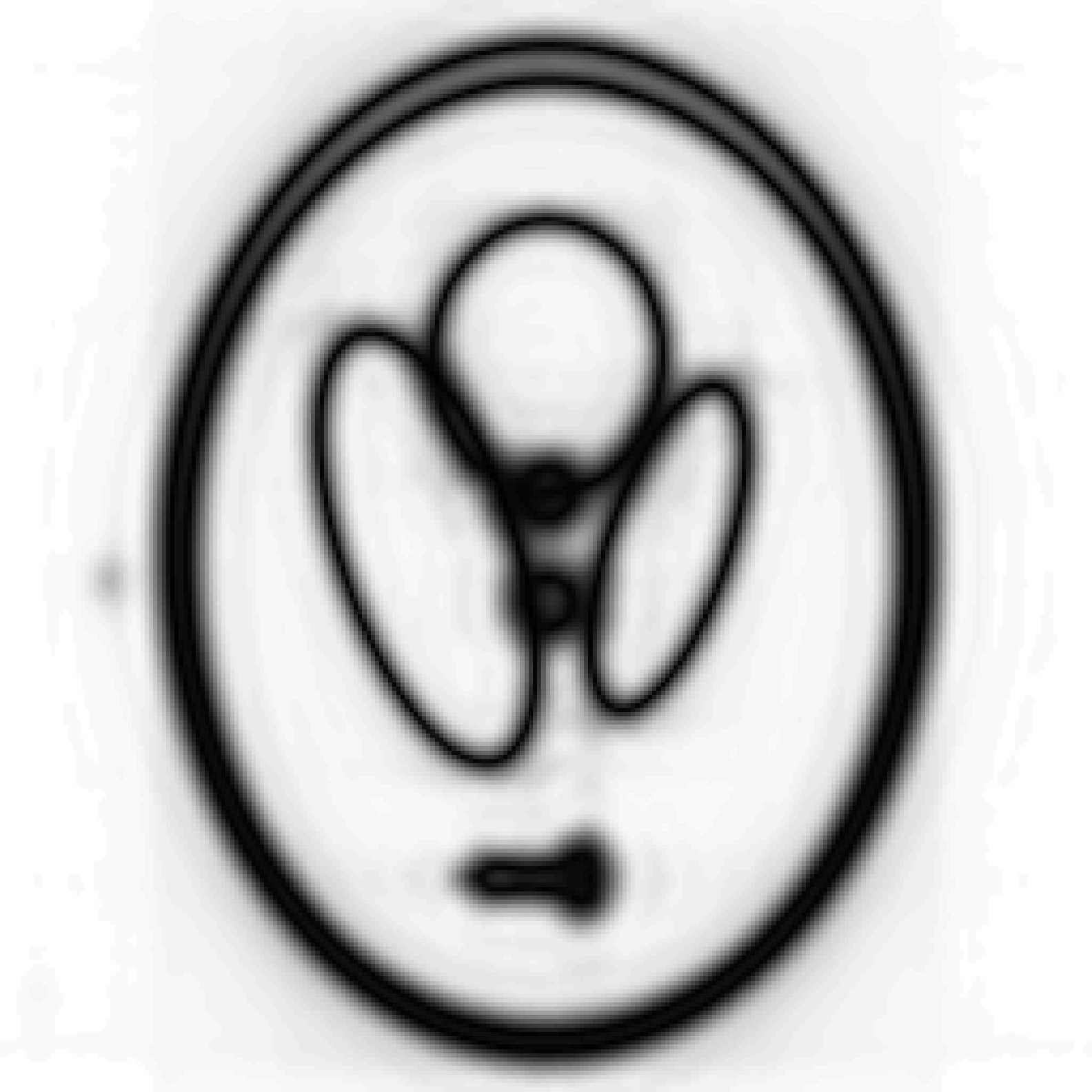}
\end{minipage}}
\caption{Comparison of estimated edges for ``Ellipse''. In \cref{EllipseRandomEdge,EllipseILFEdge}, the first row corresponds to \cref{CadzowVS}, and the second row to \cref{DDTFApproach}.}\label{EllipseEdge}
\end{figure}

\begin{figure}[t]
\centering
\subfloat[Ref.]{\label{RectangleEdgeOriginal}\begin{minipage}{3cm}
\includegraphics[width=3cm]{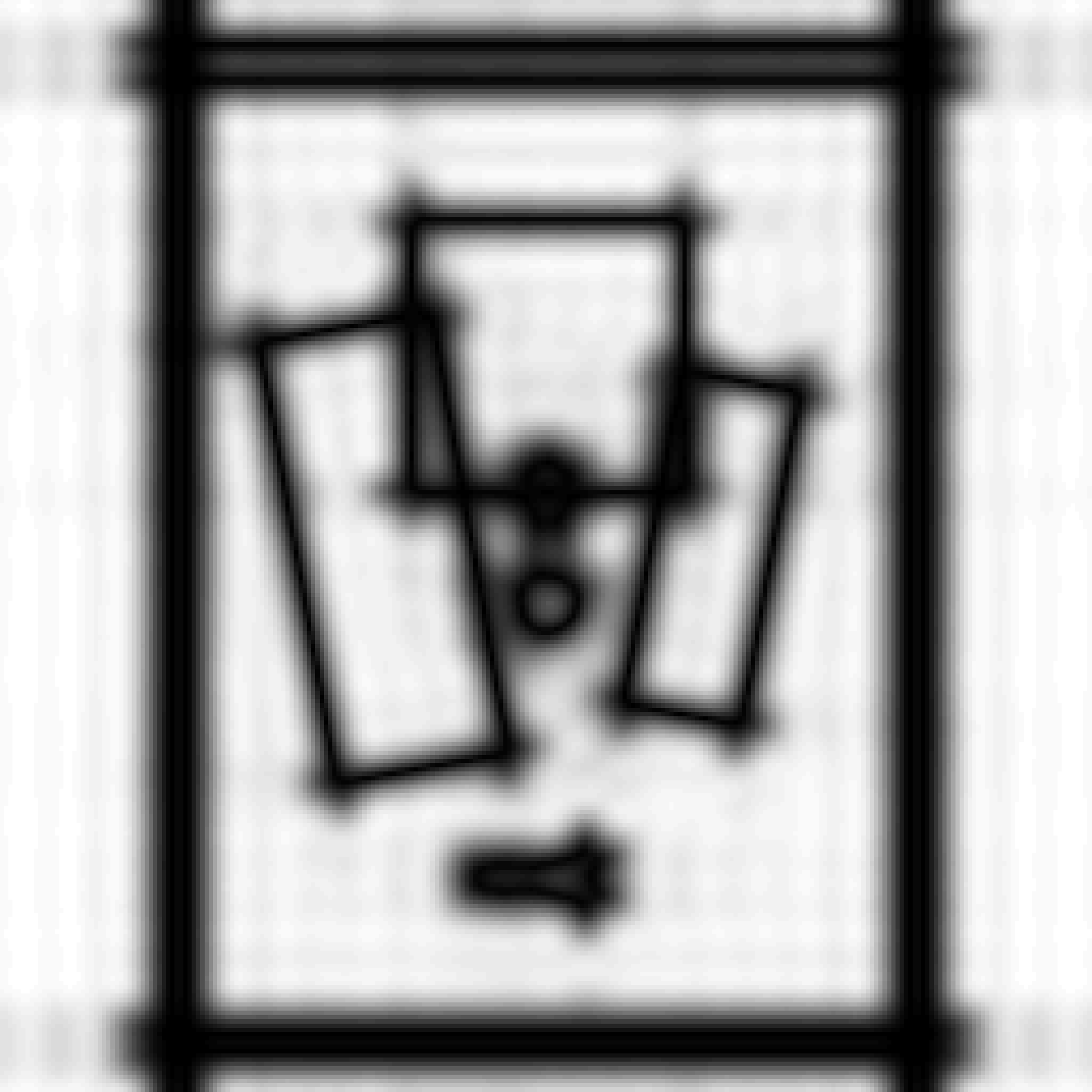}\vspace{0.025cm}\\
\includegraphics[width=3cm]{RectanglesEdgeOriginal.pdf}
\end{minipage}}\hspace{-0.05cm}
\subfloat[RS]{\label{RectangleRandomEdge}\begin{minipage}{3cm}
\includegraphics[width=3cm]{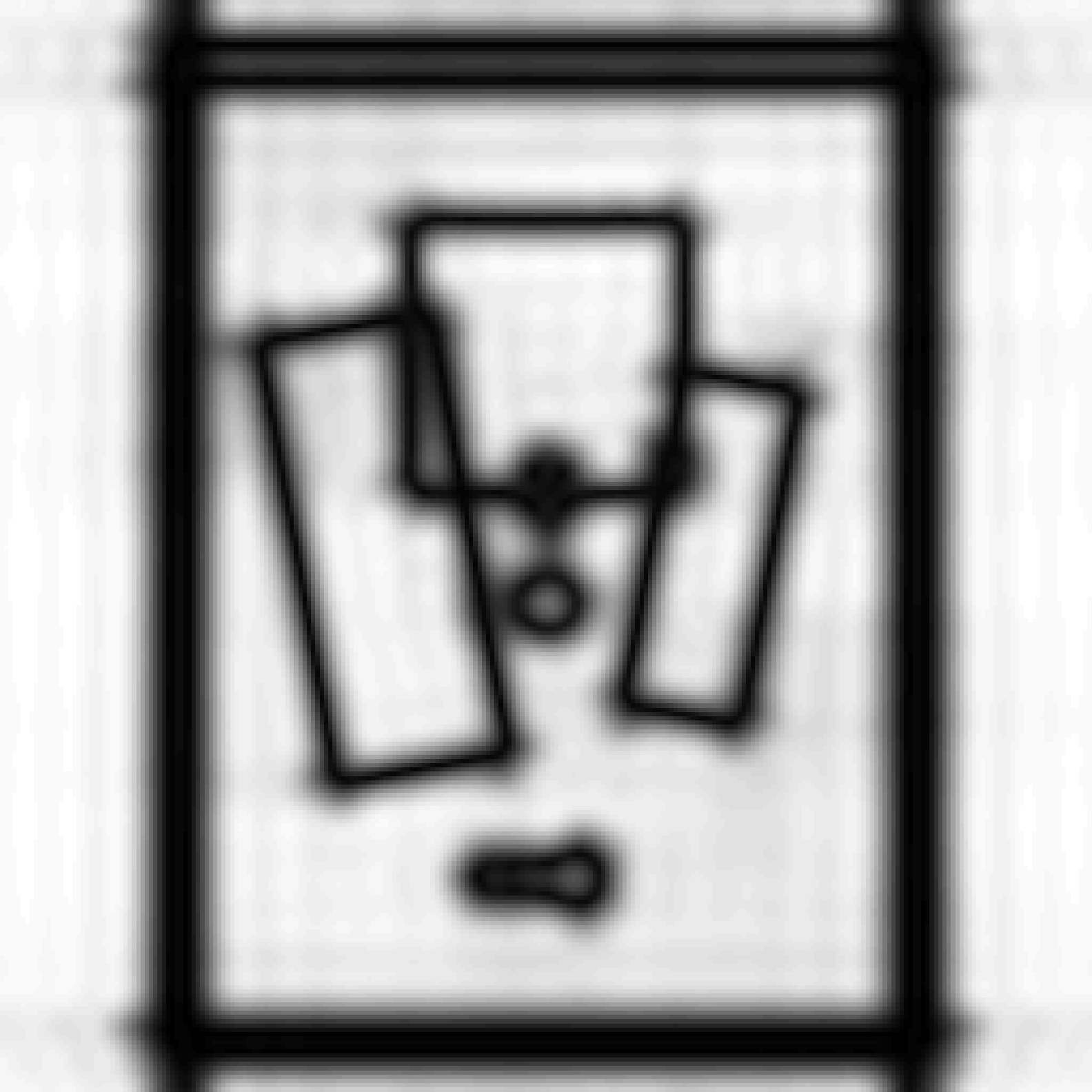}\vspace{0.025cm}\\
\includegraphics[width=3cm]{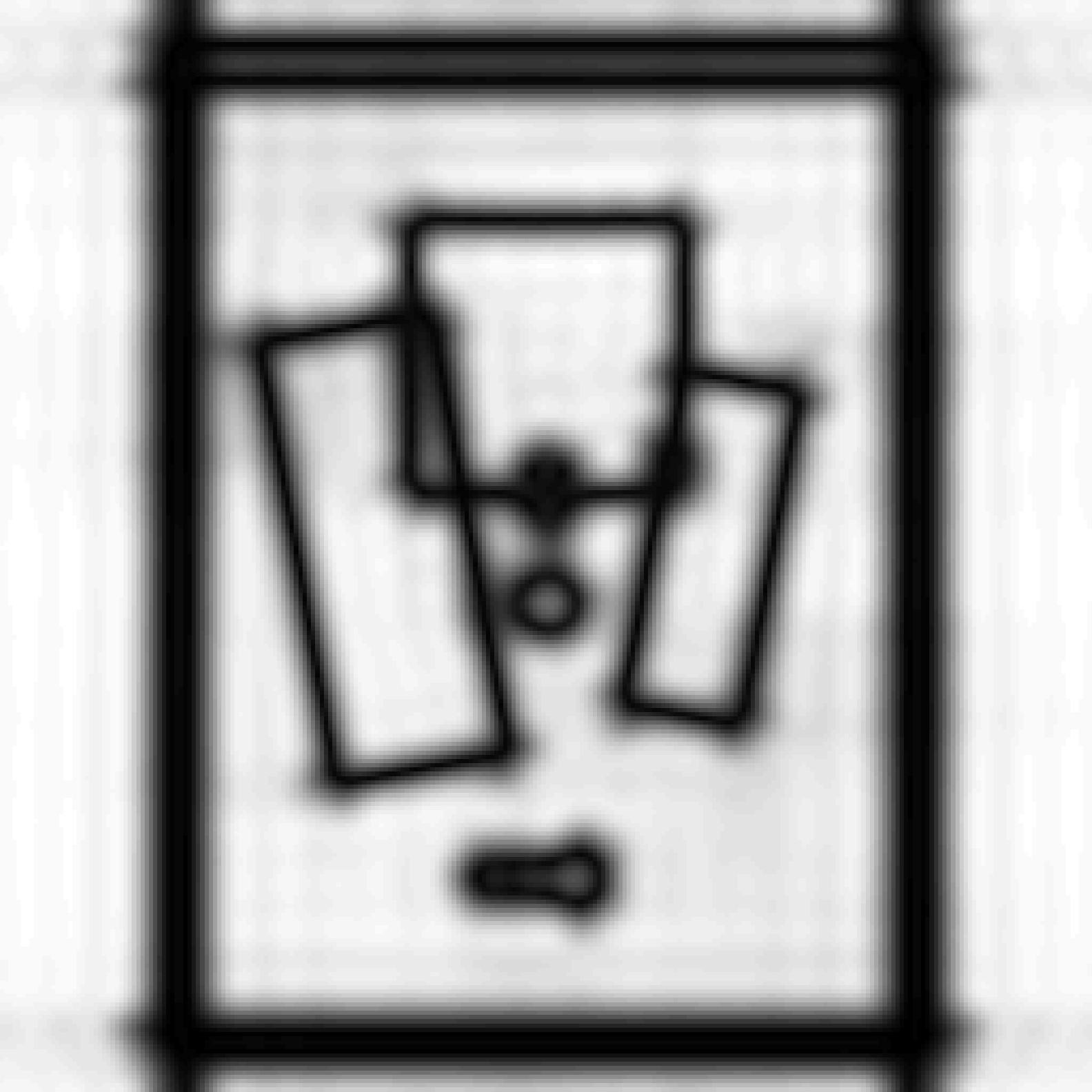}
\end{minipage}}\hspace{-0.05cm}
\subfloat[ILF]{\label{RectangleILFEdge}\begin{minipage}{3cm}
\includegraphics[width=3cm]{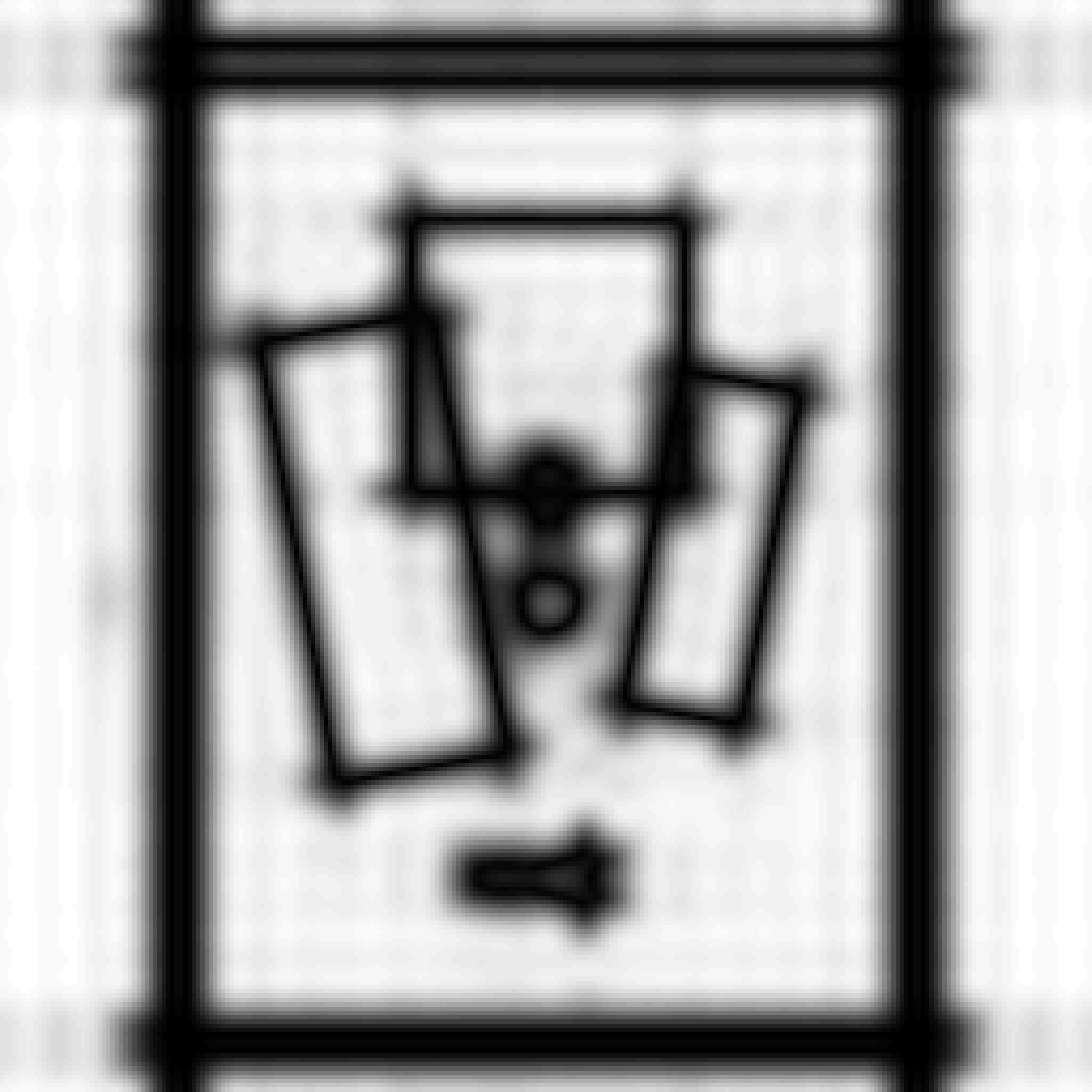}\vspace{0.025cm}\\
\includegraphics[width=3cm]{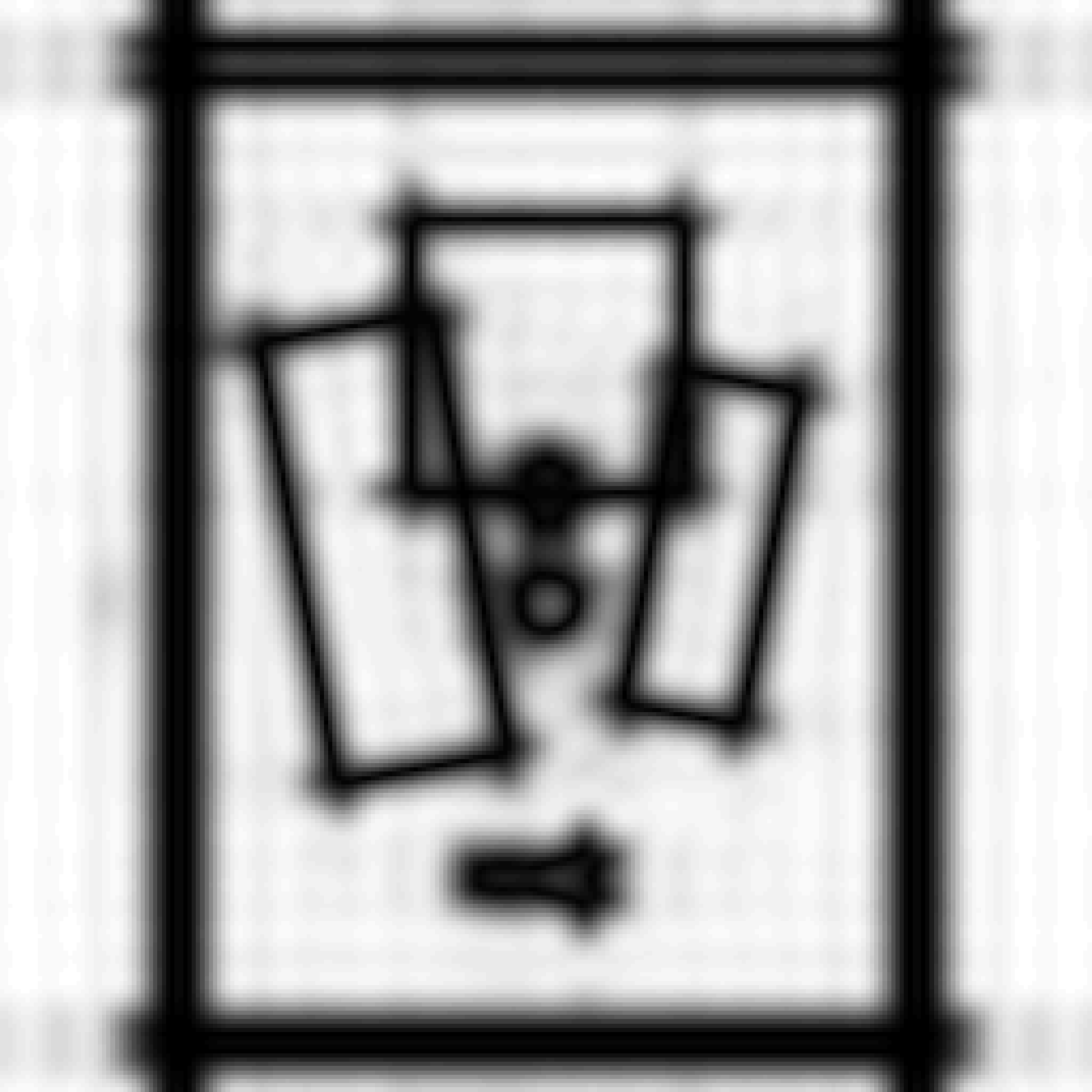}
\end{minipage}}
\caption{Comparison of estimated edges for ``Rectangle''. In \cref{RectangleRandomEdge,RectangleILFEdge}, the first row corresponds to \cref{CadzowVS}, and the second row to \cref{DDTFApproach}.}\label{RectangleEdge}
\end{figure}

\cref{TableEdge} summarizes the number of iterations and the time to reach \cref{SCEdgeEstimation}, and \cref{EllipseEdge,RectangleEdge} display the visual comparisons of edge set obtained by \cref{CadzowVS,DDTFApproach} for each scenario. At first glance, we can see that edge sets obtained from both \cref{CadzowVS,DDTFApproach} are visually similar to the reference edge set, as shown in \cref{EllipseEdge,RectangleEdge}. In addition, both methods show similar iteration numbers to meet \cref{SCEdgeEstimation} as shown in \cref{TableEdge}. However, it is evident that \cref{DDTFApproach} shows much shorter computational time than \cref{CadzowVS}. Notice that in \cref{Wsubexplicit} of \cref{Alg2}, we only compute the SVD of a $K^2\times K^2$ matrix $\wt{\bsH}_{(n+1)}^*\wt{\bC}_{(n+1)}+\beta_3\beta^{-1}\wt{\bsA}_{(n)}$, which is much smaller than the two-fold Hankel matrix $\bmH\left(\wt{\bmD}\wt{\bsv}_{(n+1)}\right)$ arising in the alternating minimization to solve \cref{CadzowVS}. Hence, we can conclude that, even though it is numerically demonstrated that there exists an equivalence between \cref{CadzowVS,DDTFApproach}, our approach \cref{DDTFApproach} shows an advantage over \cref{CadzowVS} in the fast computation.

\subsection{Image restoration results: synthetic images}\label{Synthetic}

\begin{table}[t]
\centering
\scriptsize
\begin{tabular}{|c|c|c|c|c|c|c|c|c|c|}
\hline
\multicolumn{10}{|c|}{Random sampling}\\ \hline
\multirow{2}{*}{Image}&\multirow{2}{*}{Index}&\multirow{2}{*}{IFFT}&\multirow{2}{*}{Proposed \cref{ProposedModel2}}&\multirow{2}{*}{LSLP \cref{LSLP2}}&\multirow{2}{*}{LRHDDTF \cref{LRHDDTF2}}&\multirow{2}{*}{Schatten $0$ \cref{Schatten0}}&\multirow{2}{*}{TV \cref{TV}}&\multicolumn{2}{|c|}{Frame \cref{Frame}}\\ \cline{9-10}
&&&&&&&&Haar&DDTF\\ \hline
\multirow{3}{*}{Ellipse}&SNR&$10.26$&$\textbf{36.05}$&$26.46$&$29.74$&$26.09$&$24.77$&$25.13$&$26.52$\\ \cline{2-10}
&HFEN&$0.4608$&$\textbf{0.0127}$&$0.0311$&$0.0266$&$0.0314$&$0.0482$&$0.0485$&$0.0300$\\ \cline{2-10}
&SSIM&$0.3521$&$\textbf{0.9965}$&$0.9757$&$0.9805$&$0.9749$&$0.9246$&$0.9390$&$0.9341$\\ \hline
\multirow{3}{*}{Rectangle}&SNR&$10.88$&$\textbf{33.65}$&$26.41$&$24.49$&$25.67$&$26.11$&$26.45$&$26.24$\\ \cline{2-10}
&HFEN&$0.4461$&$\textbf{0.0196}$&$0.0302$&$0.0446$&$0.0372$&$0.0402$&$0.0403$&$0.0322$\\ \cline{2-10}
&SSIM&$0.4126$&$\textbf{0.9846}$&$0.9455$&$0.9495$&$0.9526$&$0.8809$&$0.9244$&$0.8963$\\ \hline\hline
\multicolumn{10}{|c|}{Ideal low-pass filter}\\ \hline
\multirow{2}{*}{Image}&\multirow{2}{*}{Index}&\multirow{2}{*}{IFFT}&\multirow{2}{*}{Proposed \cref{ProposedModel2}}&\multirow{2}{*}{LSLP \cref{LSLP2}}&\multirow{2}{*}{LRHDDTF \cref{LRHDDTF2}}&\multirow{2}{*}{Schatten $0$ \cref{Schatten0}}&\multirow{2}{*}{TV \cref{TV}}&\multicolumn{2}{|c|}{Frame \cref{Frame}}\\ \cline{9-10}
&&&&&&&&Haar&DDTF\\ \hline
\multirow{3}{*}{Ellipse}&SNR&$11.46$&$\textbf{36.91}$&$22.53$&$14.52$&$18.53$&$20.65$&$20.77$&$14.27$\\ \cline{2-10}
&HFEN&$0.5102$&$\textbf{0.0183}$&$0.0622$&$0.2681$&$0.0965$&$0.0799$&$0.0827$&$0.2314$\\ \cline{2-10}
&SSIM&$0.6932$&$\textbf{0.9973}$&$0.9639$&$0.9287$&$0.9454$&$0.9196$&$0.9217$&$0.9231$\\ \hline
\multirow{3}{*}{Rectangle}&SNR&$10.74$&$\textbf{30.21}$&$21.85$&$13.84$&$18.07$&$22.73$&$21.07$&$14.79$\\ \cline{2-10}
&HFEN&$0.5960$&$\textbf{0.0409}$&$0.0897$&$0.3343$&$0.1221$&$0.731$&$0.1010$&$0.2137$\\ \cline{2-10}
&SSIM&$0.6401$&$\textbf{0.9944}$&$0.9045$&$0.8890$&$0.9073$&$0.8755$&$0.8780$&$0.8842$\\ \hline
\end{tabular}
\caption{Comparison of SNR, HFEN, and SSIM for synthetic images}\label{PhantomResults}
\end{table}

\begin{figure}[t]
\centering
\subfloat[Ref.]{\label{EllipseRandomOriginal}\begin{minipage}{3cm}
\includegraphics[width=3cm]{EllipsesOriginal.pdf}
\end{minipage}}\hspace{-0.05cm}
\subfloat[Proposed]{\label{EllipseRandomProposed}\begin{minipage}{3cm}
\includegraphics[width=3cm]{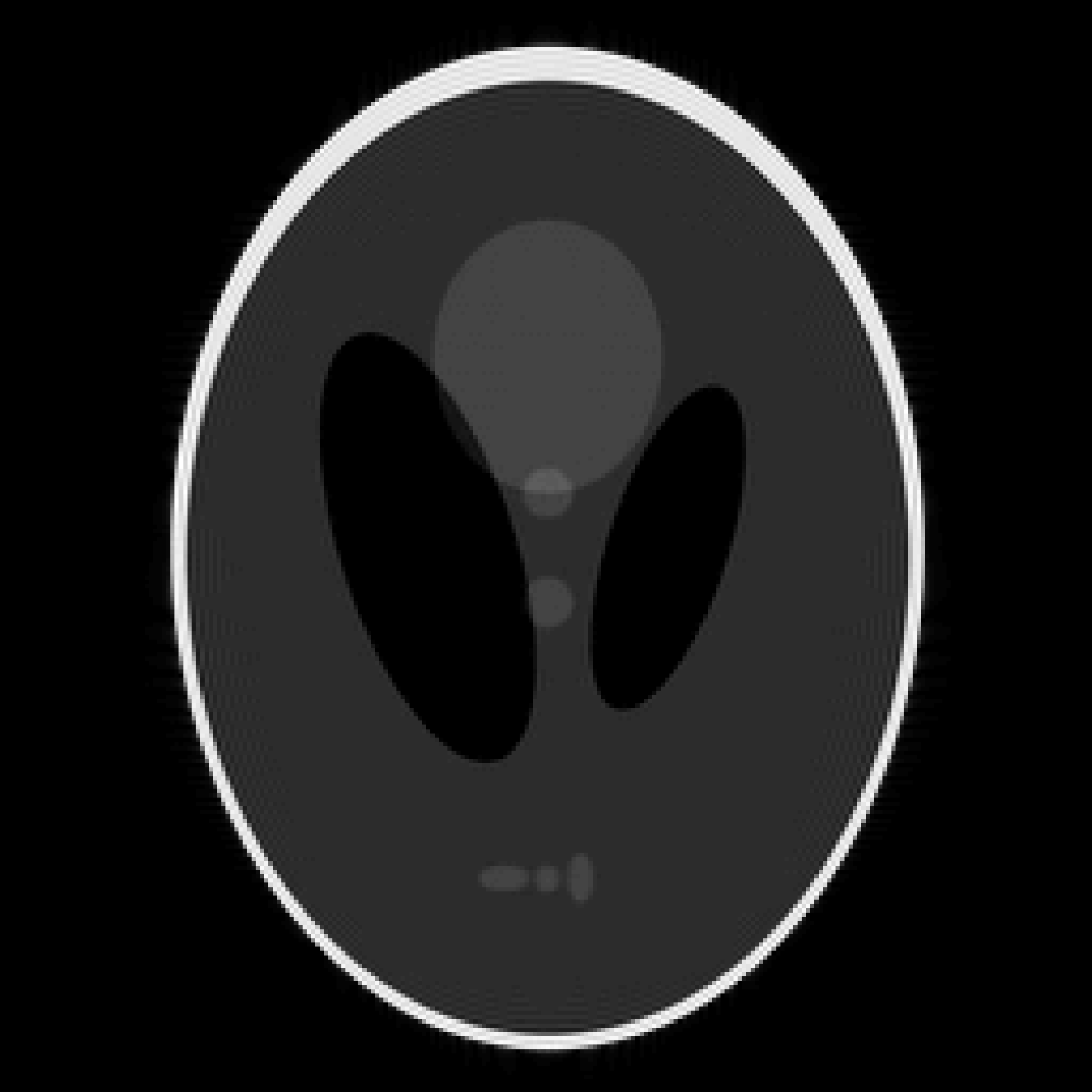}
\end{minipage}}\hspace{-0.05cm}
\subfloat[LSLP]{\label{EllipseRandomLSLP}\begin{minipage}{3cm}
\includegraphics[width=3cm]{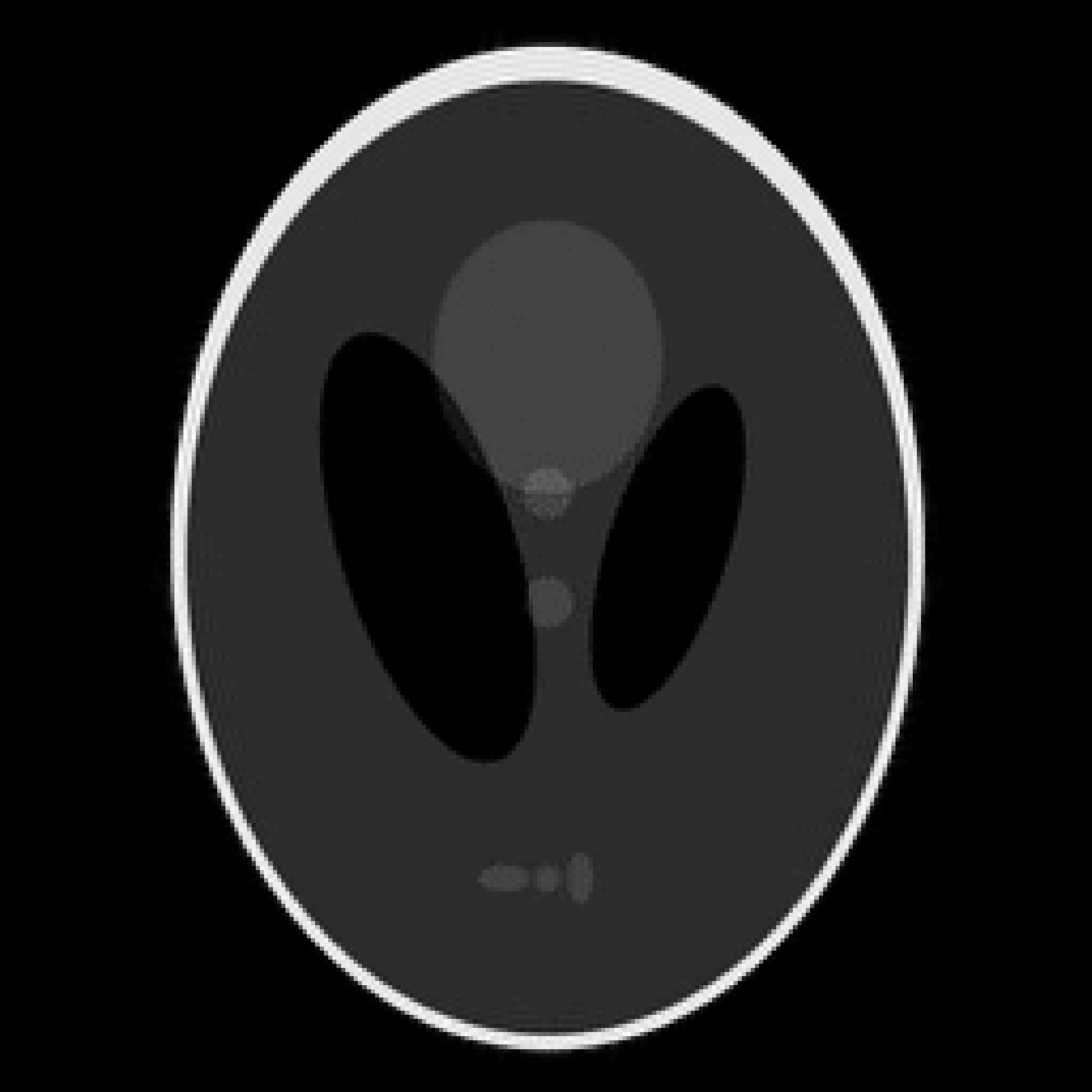}
\end{minipage}}\hspace{-0.05cm}
\subfloat[LRHDDTF]{\label{EllipseRandomLRHDDTF}\begin{minipage}{3cm}
\includegraphics[width=3cm]{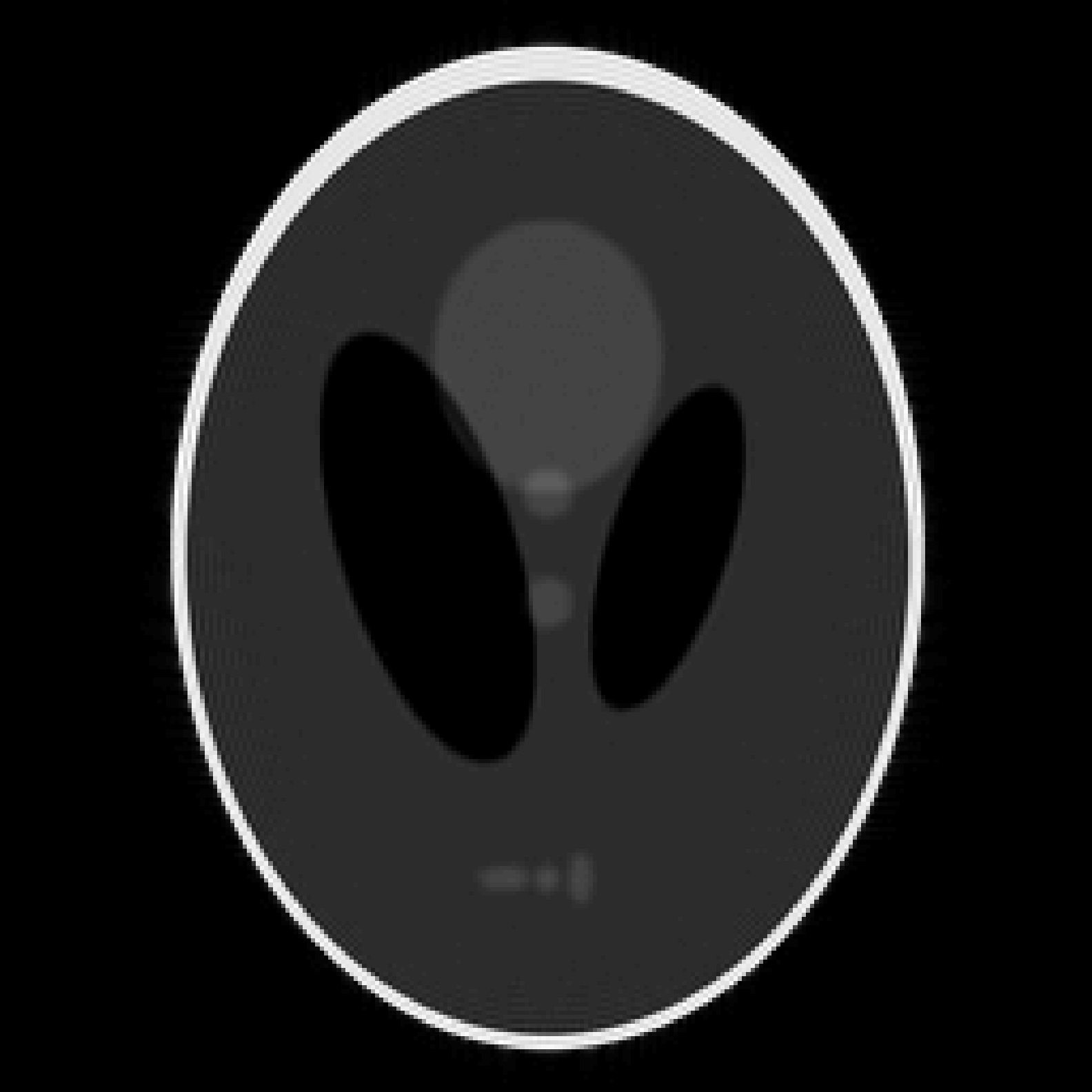}
\end{minipage}}\vspace{-0.25cm}\\
\subfloat[Schatten $0$]{\label{EllipseRandomGIRAF}\begin{minipage}{3cm}
\includegraphics[width=3cm]{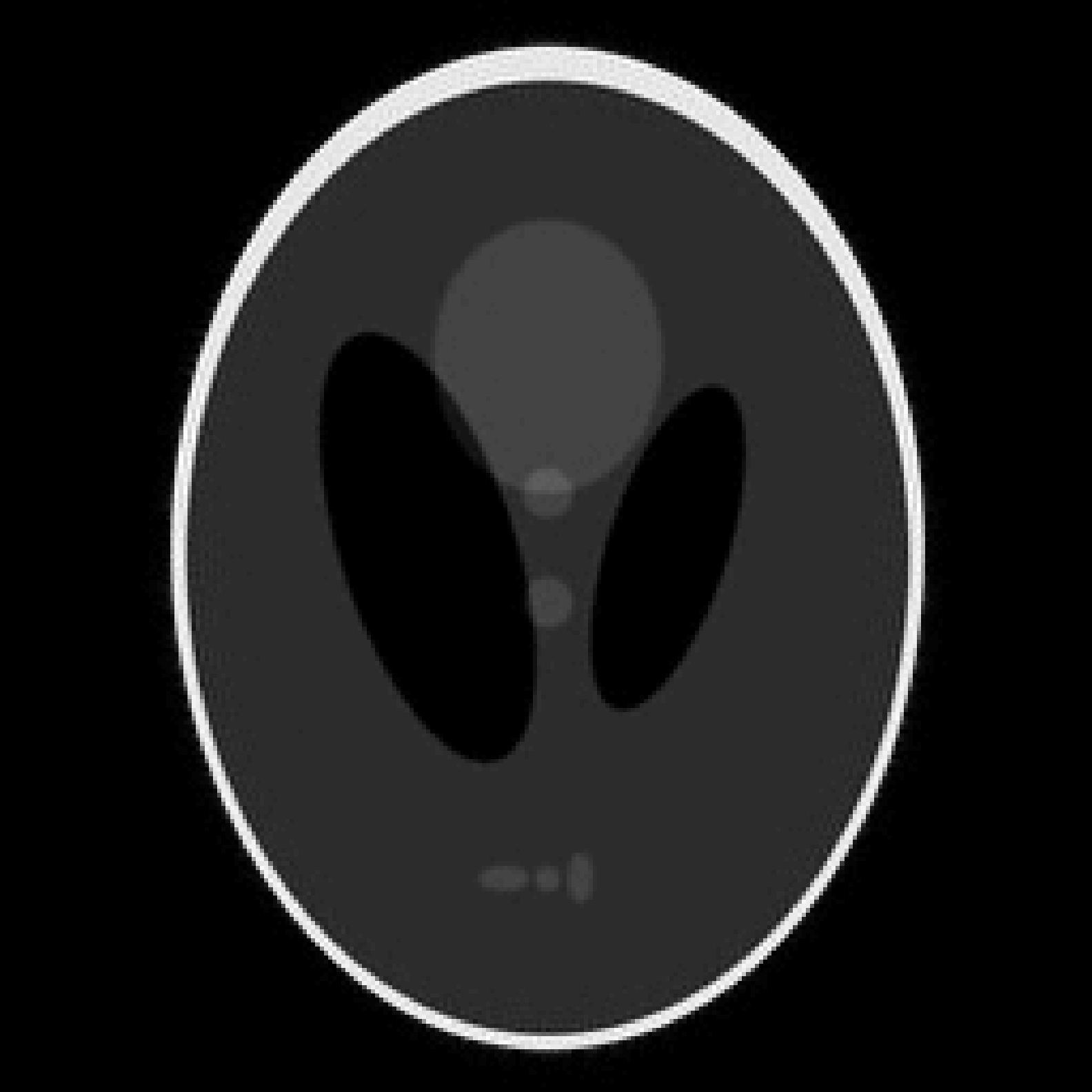}
\end{minipage}}\hspace{-0.05cm}
\subfloat[TV]{\label{EllipseRandomTV}\begin{minipage}{3cm}
\includegraphics[width=3cm]{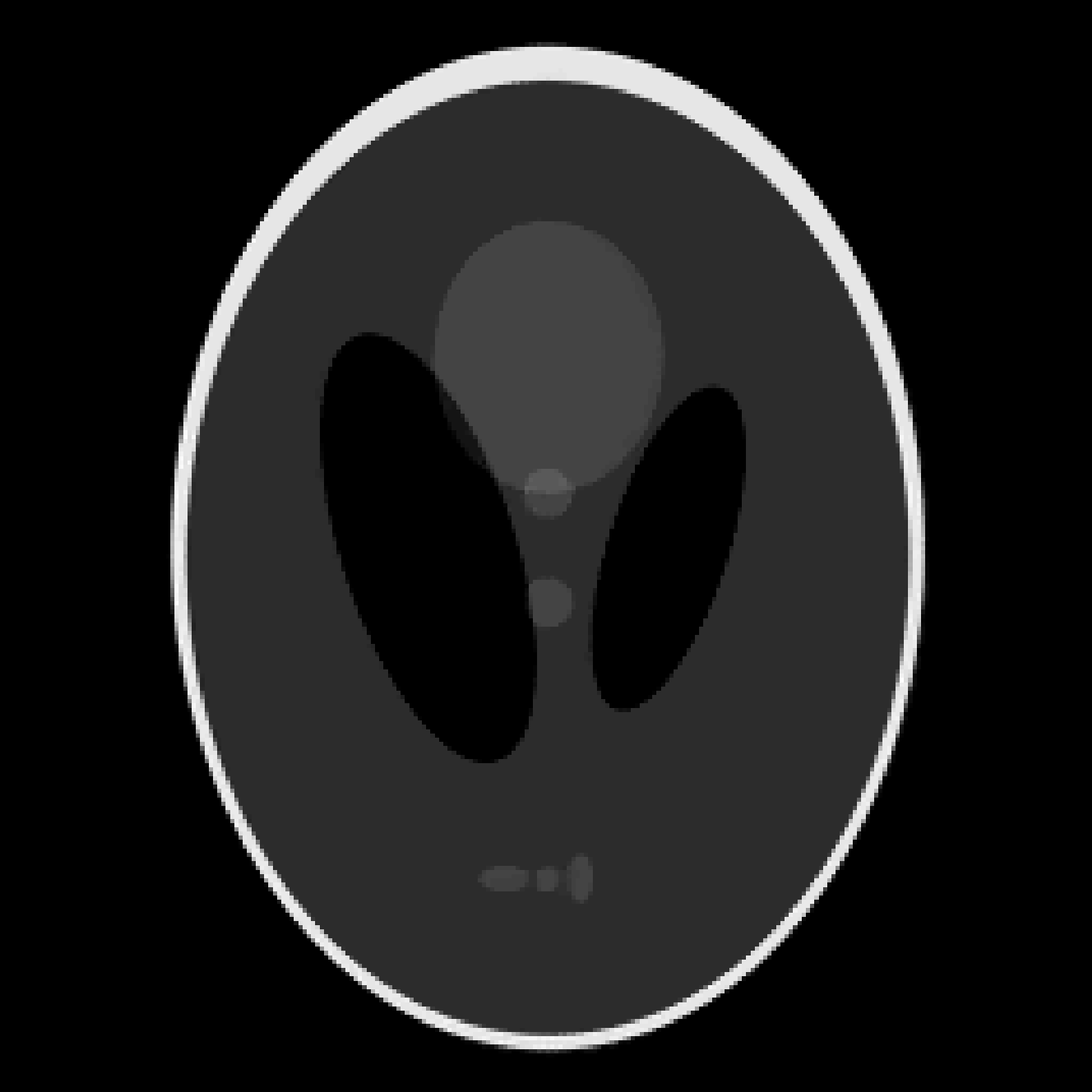}
\end{minipage}}\hspace{-0.05cm}
\subfloat[Haar]{\label{EllipseRandomHaar}\begin{minipage}{3cm}
\includegraphics[width=3cm]{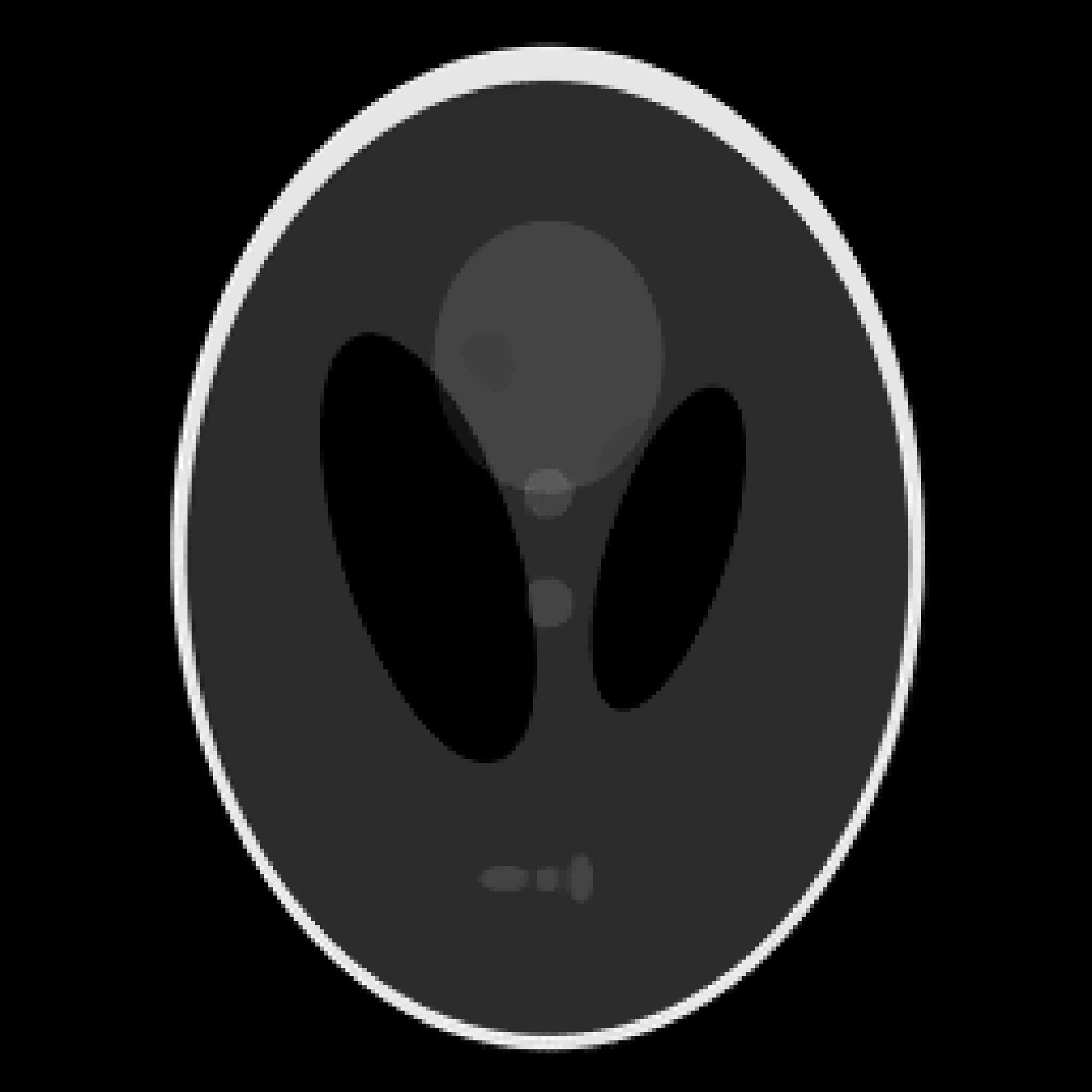}
\end{minipage}}\hspace{-0.05cm}
\subfloat[DDTF]{\label{EllipseRandomDDTF}\begin{minipage}{3cm}
\includegraphics[width=3cm]{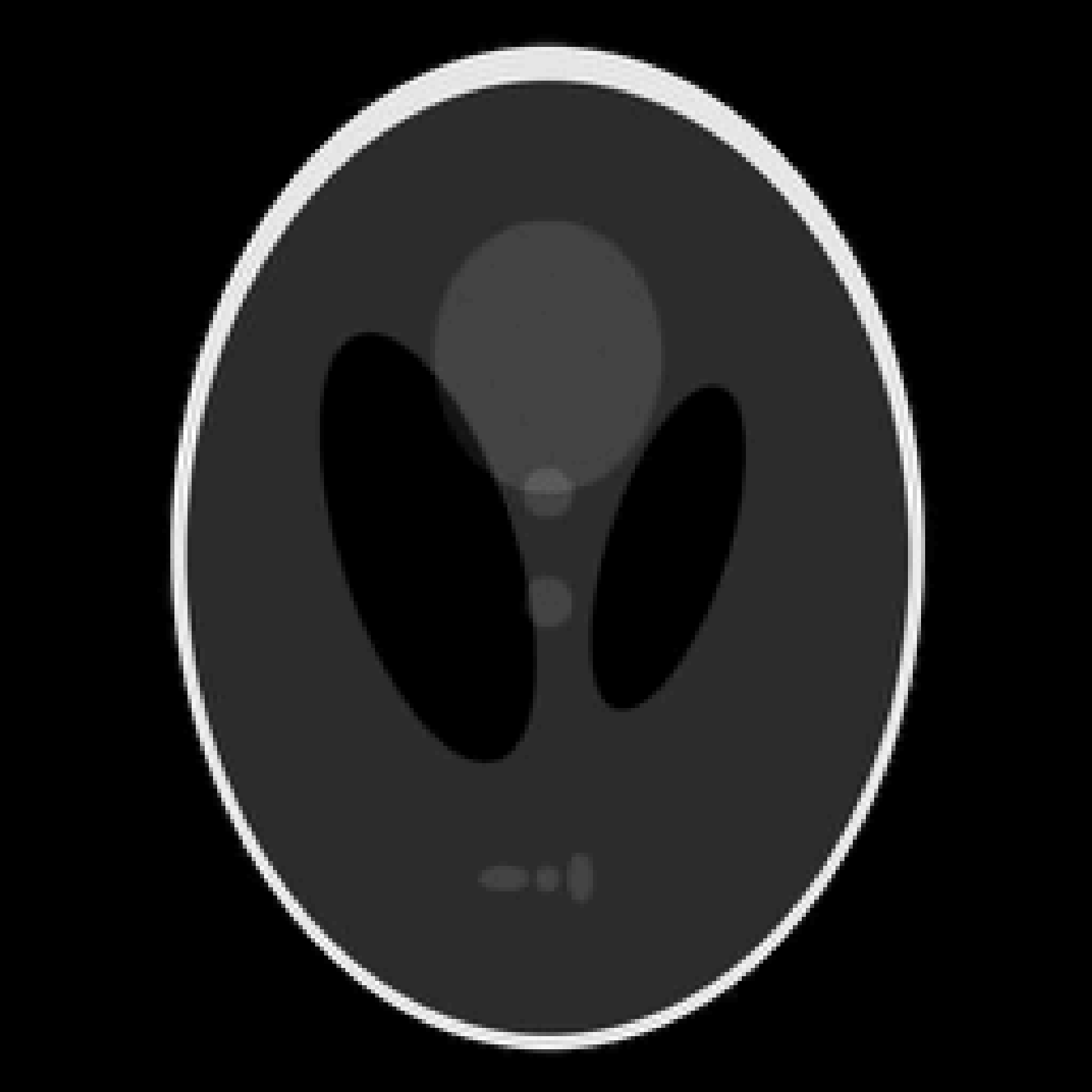}
\end{minipage}}
\caption{Visual comparison of ``Ellipse'' for random sampling. All restored images are displayed in the window level $[0,1]$ for the fair comparisons.}\label{EllipseRandomResults}
\end{figure}

\begin{figure}[ht]
\centering
\subfloat[Samples]{\label{EllipseRandomOriginal}\begin{minipage}{3cm}
\includegraphics[width=3cm]{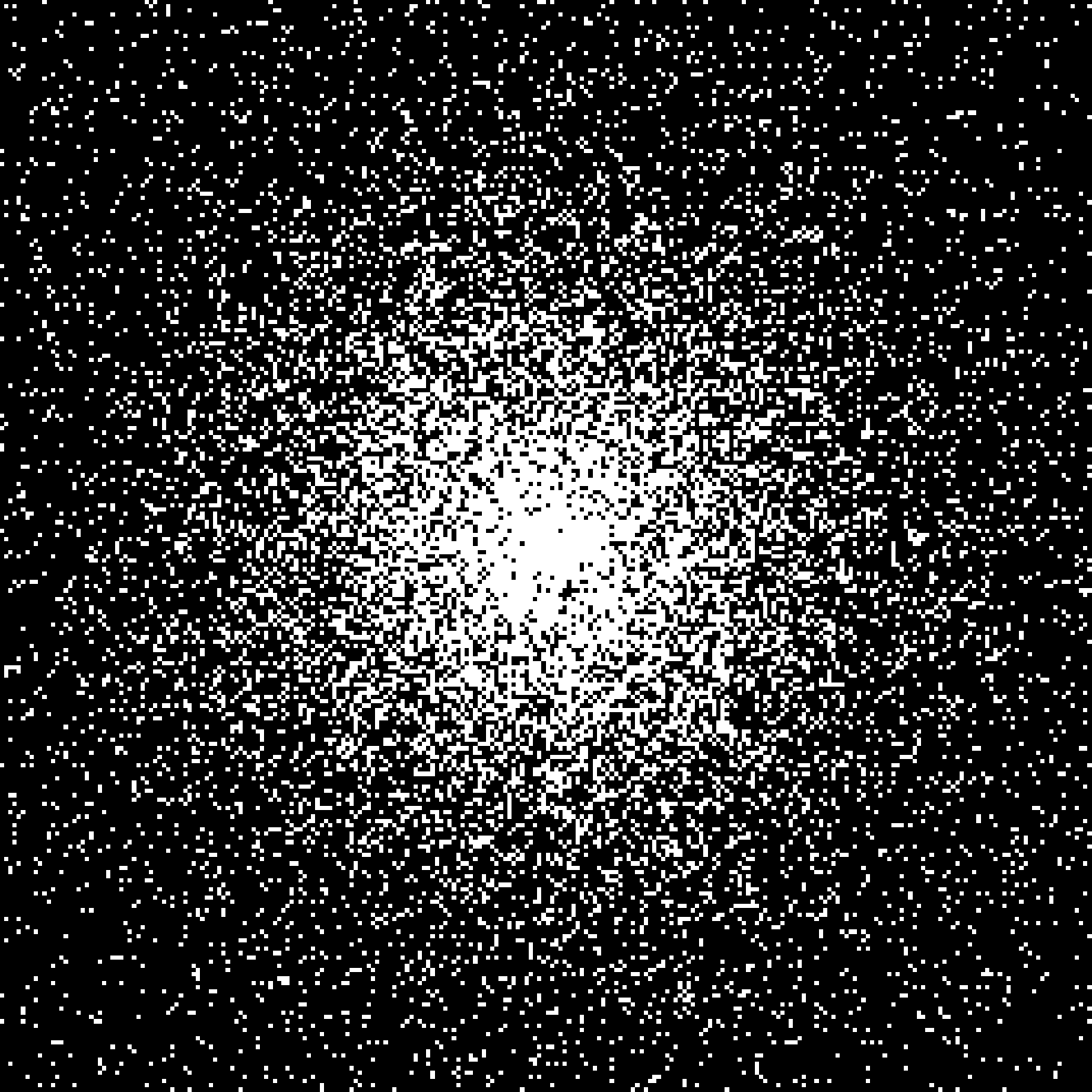}
\end{minipage}}\hspace{-0.05cm}
\subfloat[Proposed]{\label{EllipseRandomProposedError}\begin{minipage}{3cm}
\includegraphics[width=3cm]{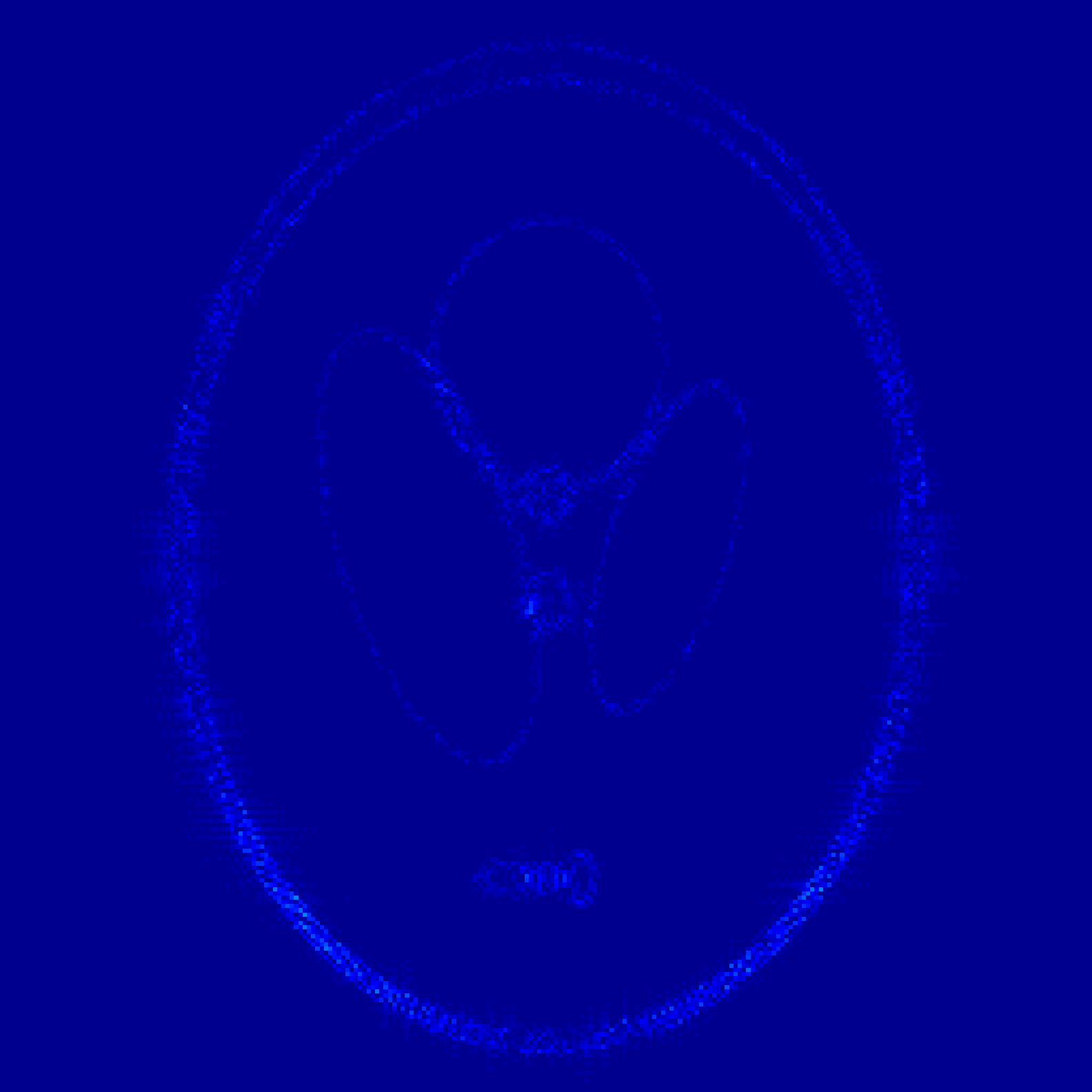}
\end{minipage}}\hspace{-0.05cm}
\subfloat[LSLP]{\label{EllipseRandomLSLPError}\begin{minipage}{3cm}
\includegraphics[width=3cm]{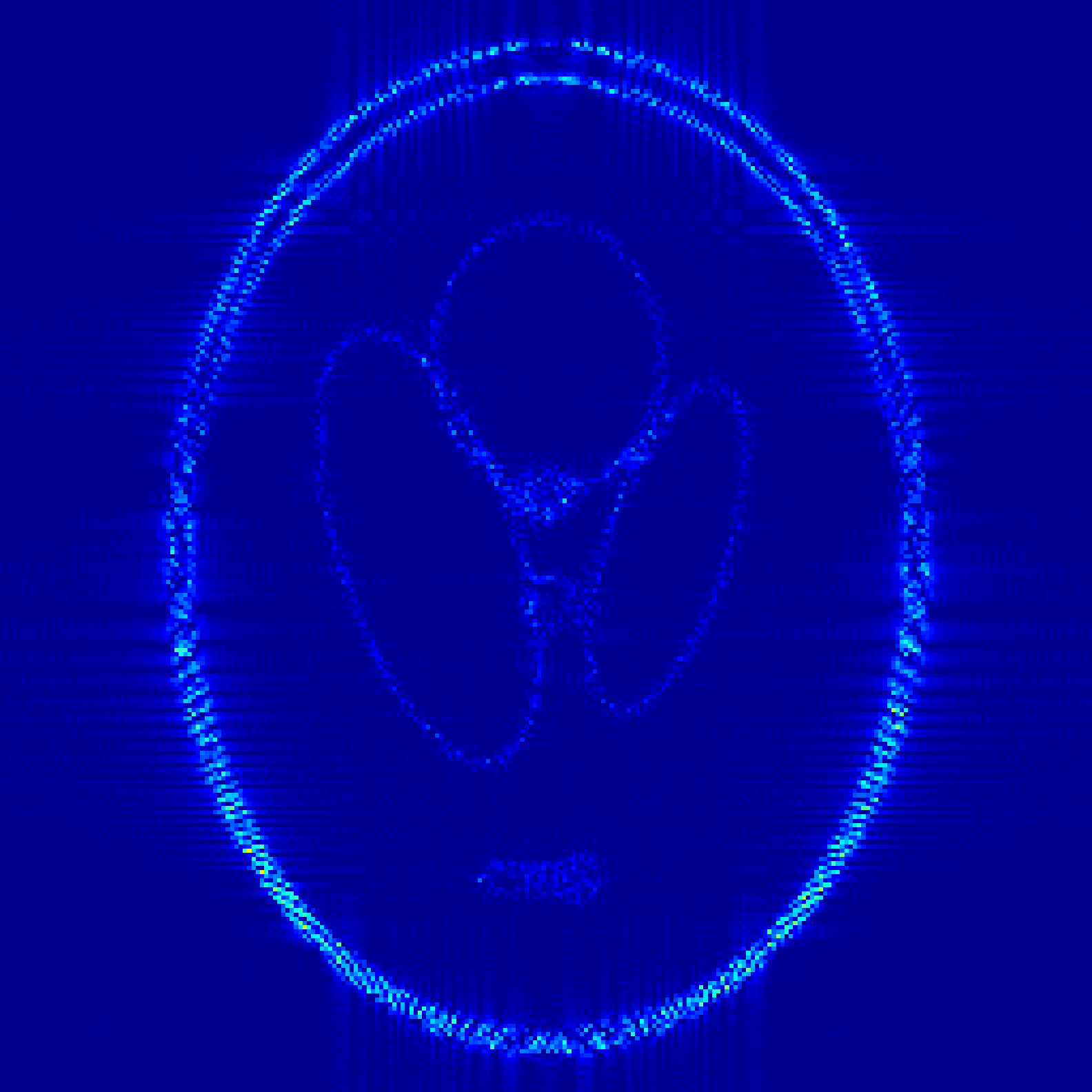}
\end{minipage}}\hspace{-0.05cm}
\subfloat[LRHDDTF]{\label{EllipseRandomLRHDDTFError}\begin{minipage}{3cm}
\includegraphics[width=3cm]{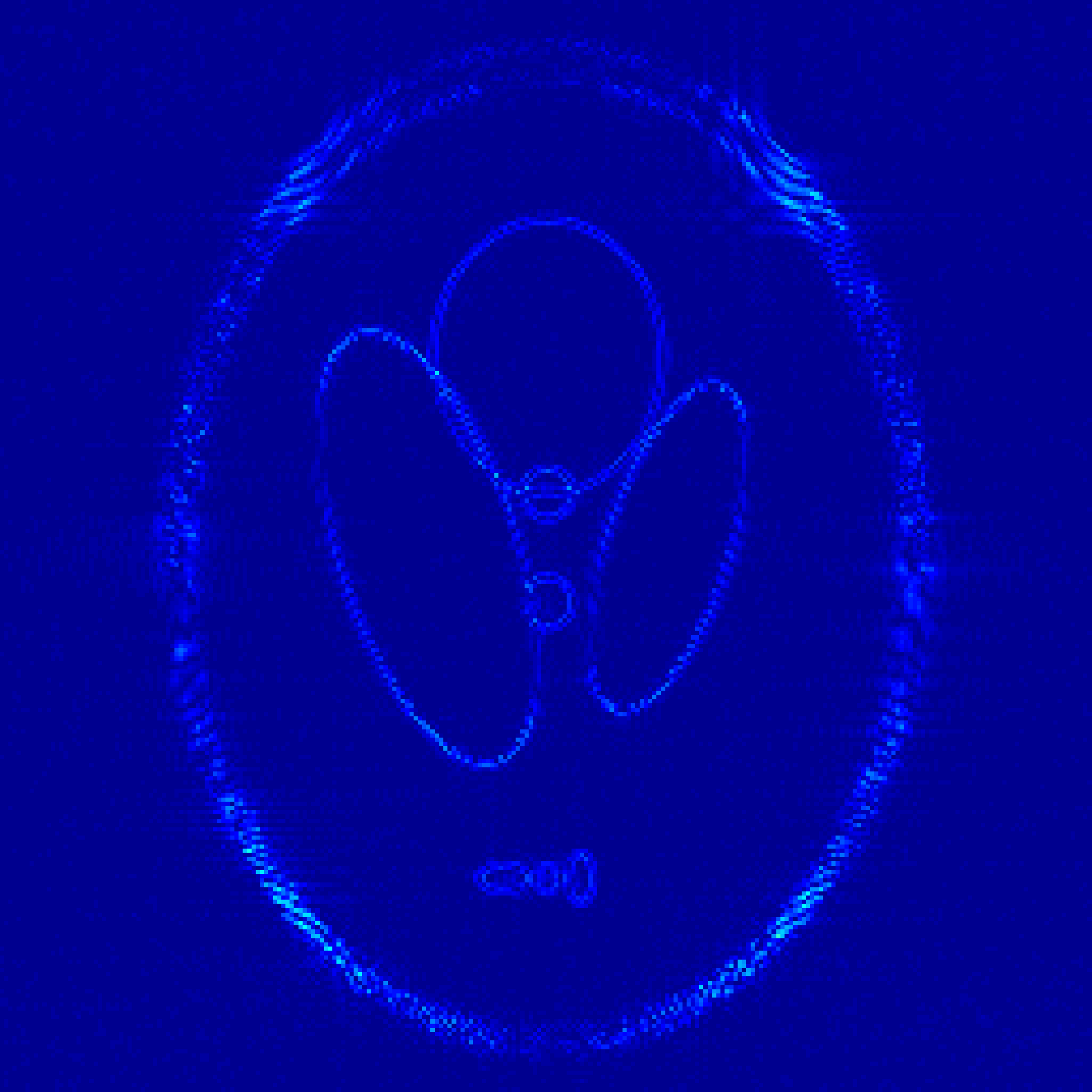}
\end{minipage}}\vspace{-0.25cm}\\
\subfloat[Schatten $0$]{\label{EllipseRandomGIRAFError}\begin{minipage}{3cm}
\includegraphics[width=3cm]{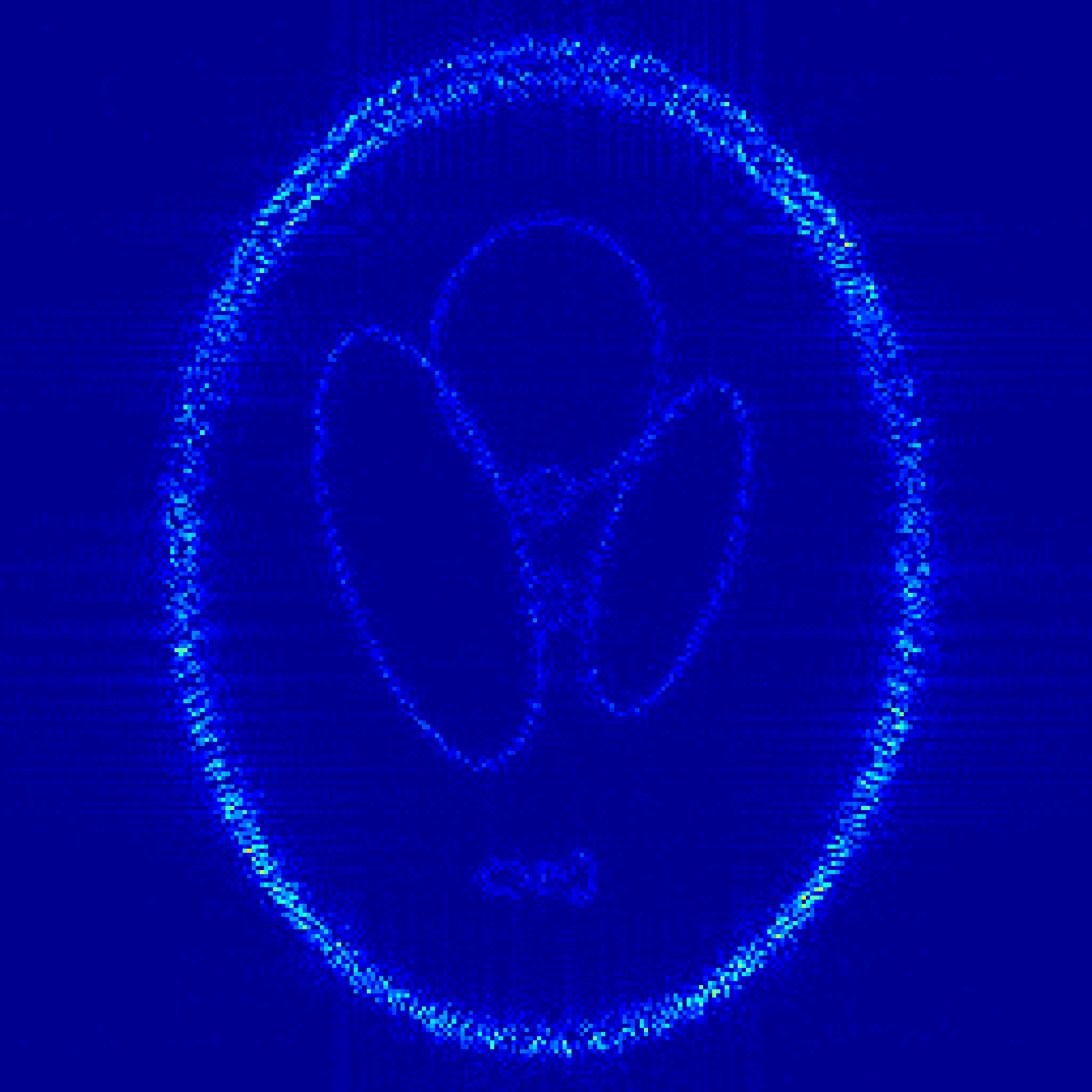}
\end{minipage}}\hspace{-0.05cm}
\subfloat[TV]{\label{EllipseRandomTVError}\begin{minipage}{3cm}
\includegraphics[width=3cm]{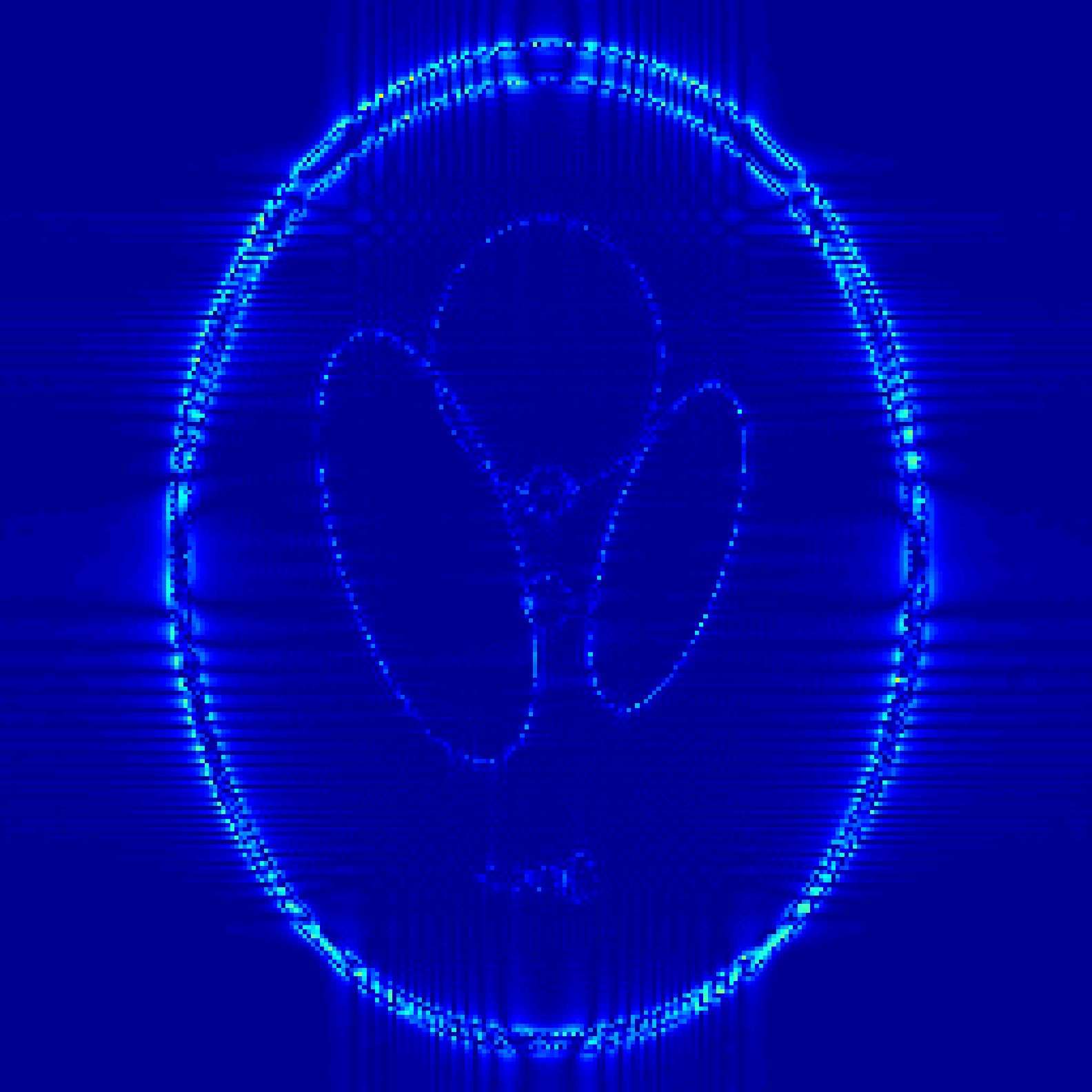}
\end{minipage}}\hspace{-0.05cm}
\subfloat[Haar]{\label{EllipseRandomHaarError}\begin{minipage}{3cm}
\includegraphics[width=3cm]{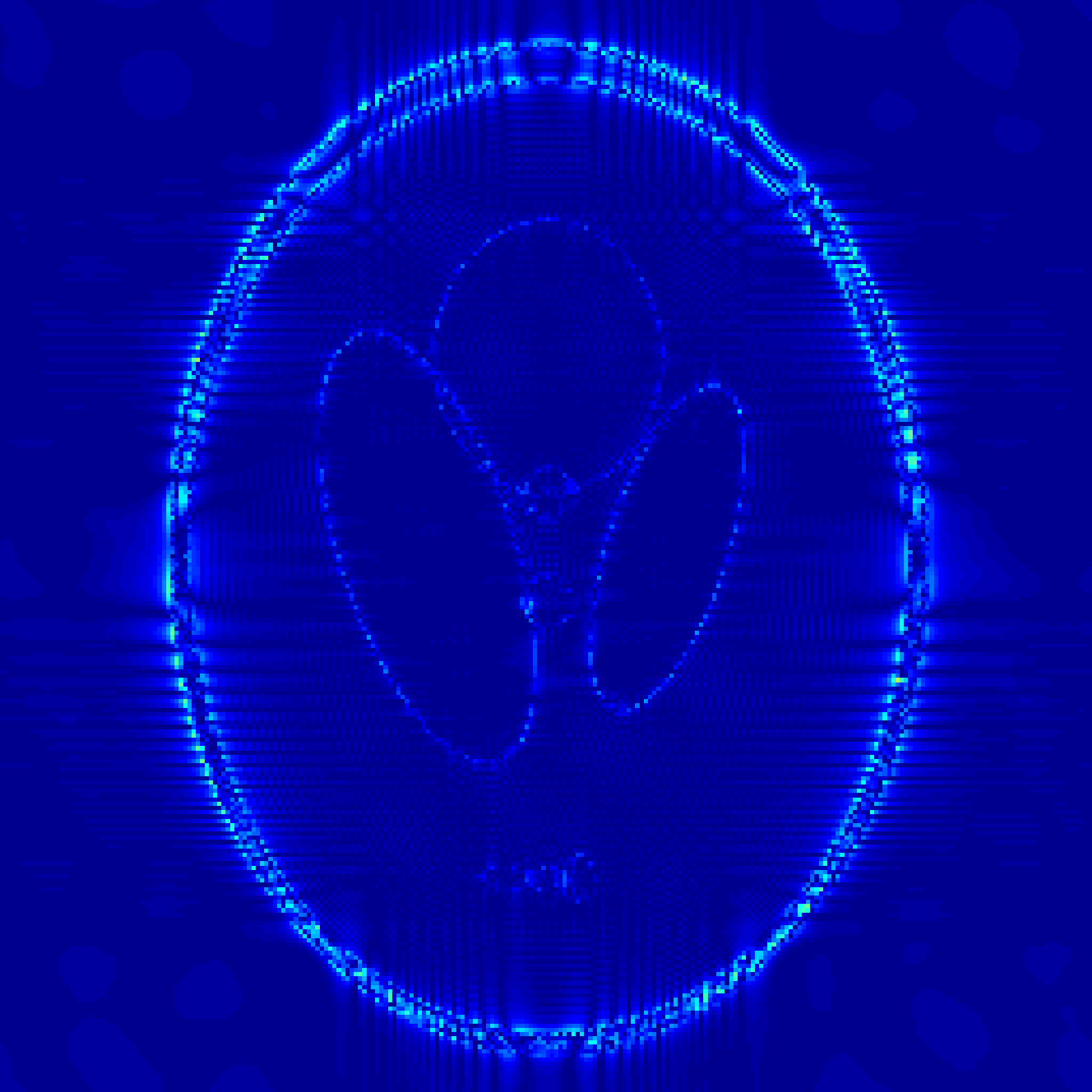}
\end{minipage}}\hspace{-0.05cm}
\subfloat[DDTF]{\label{EllipseRandomDDTFError}\begin{minipage}{3cm}
\includegraphics[width=3cm]{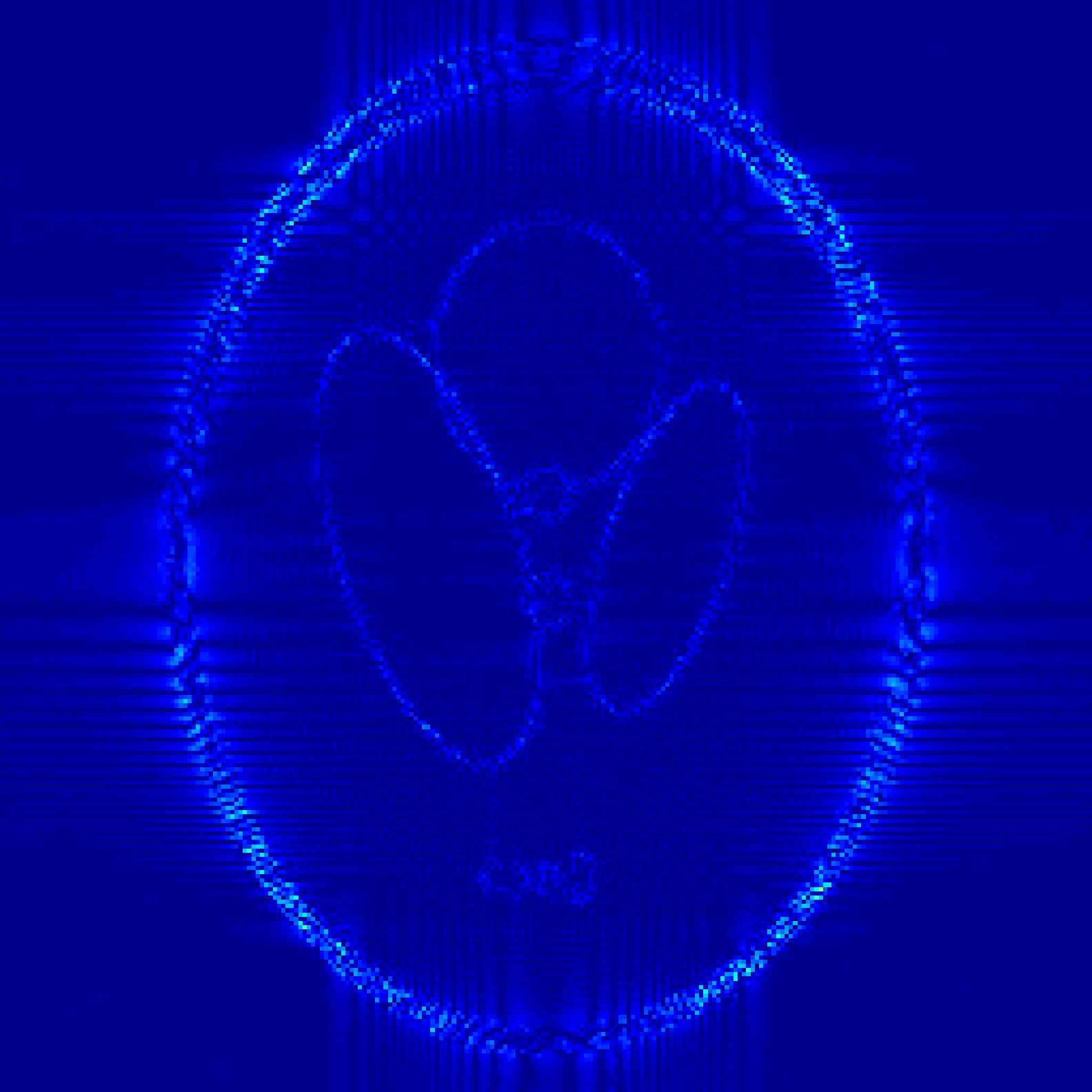}
\end{minipage}}
\caption{Error maps of \cref{EllipseRandomResults}. All error maps are displayed in the window level $[0,0.02]$ for the fair comparisons.}\label{EllipseRandomError}
\end{figure}

\begin{figure}[t]
\centering
\subfloat[Ref.]{\label{EllipseILFOriginal}\begin{minipage}{3cm}
\includegraphics[width=3cm]{EllipsesOriginal.pdf}
\end{minipage}}\hspace{-0.05cm}
\subfloat[Proposed]{\label{EllipseILFProposed}\begin{minipage}{3cm}
\includegraphics[width=3cm]{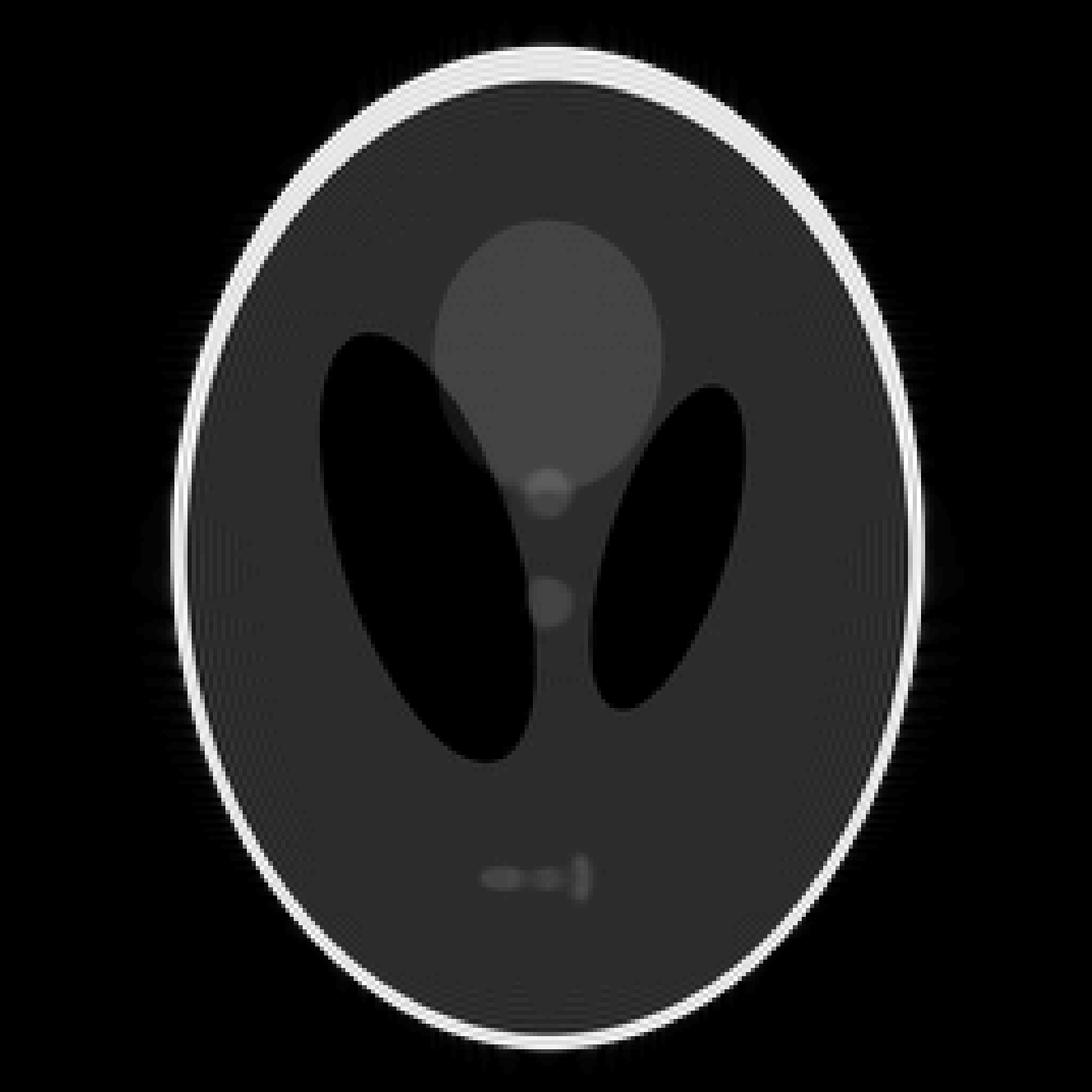}
\end{minipage}}\hspace{-0.05cm}
\subfloat[LSLP]{\label{EllipseILFLSLP}\begin{minipage}{3cm}
\includegraphics[width=3cm]{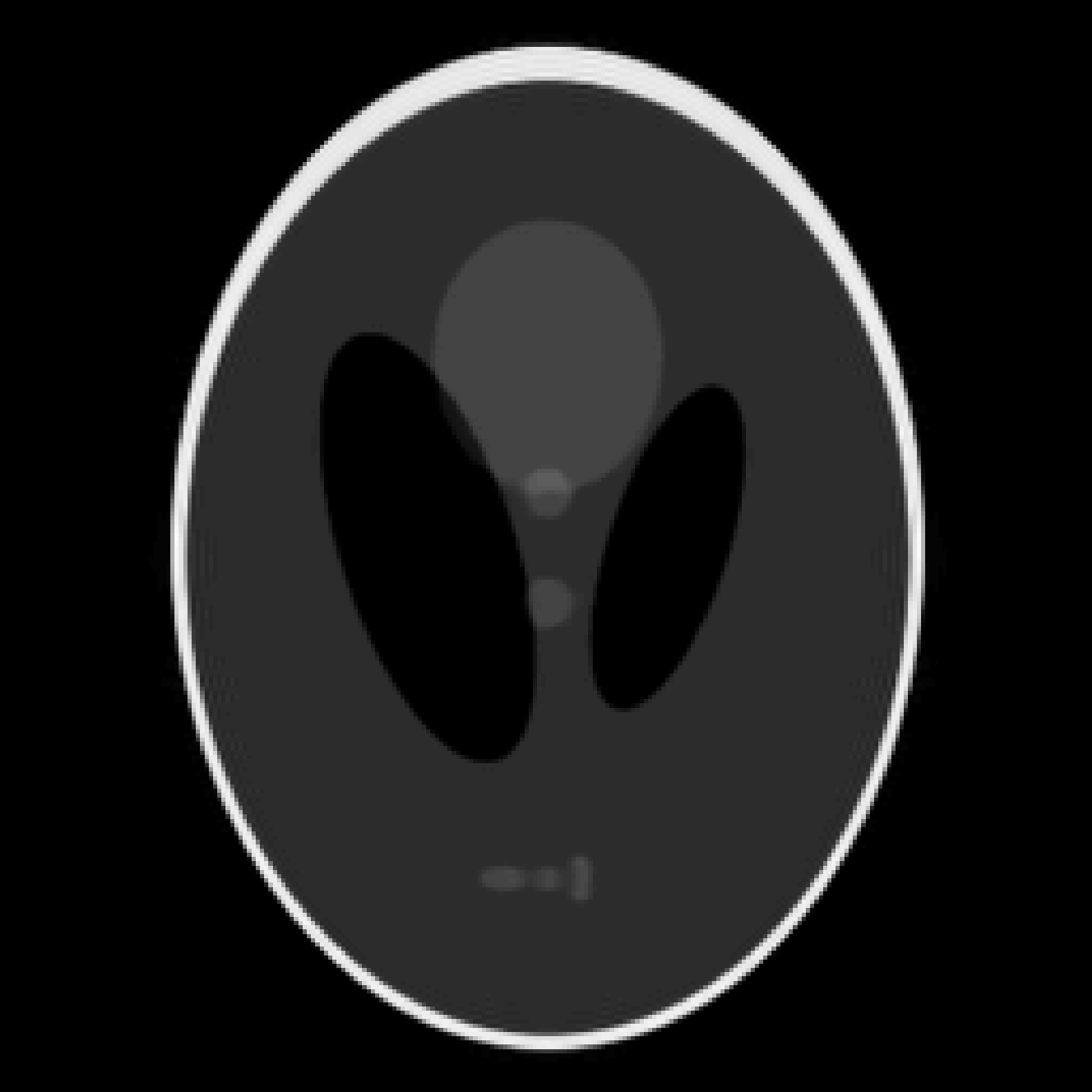}
\end{minipage}}\hspace{-0.05cm}
\subfloat[LRHDDTF]{\label{EllipseILFLRHDDTF}\begin{minipage}{3cm}
\includegraphics[width=3cm]{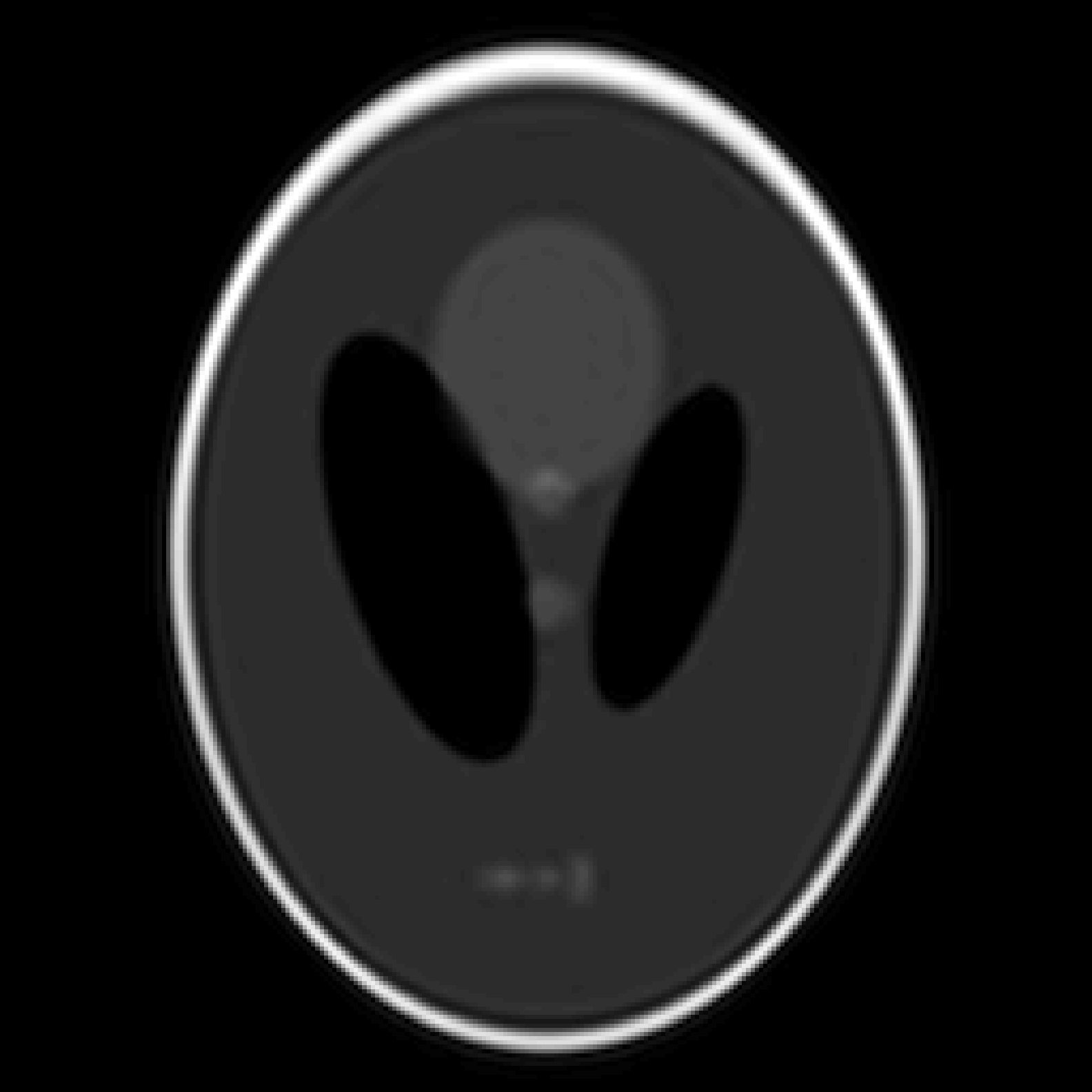}
\end{minipage}}\vspace{-0.25cm}\\
\subfloat[Schatten $0$]{\label{EllipseILFGIRAF}\begin{minipage}{3cm}
\includegraphics[width=3cm]{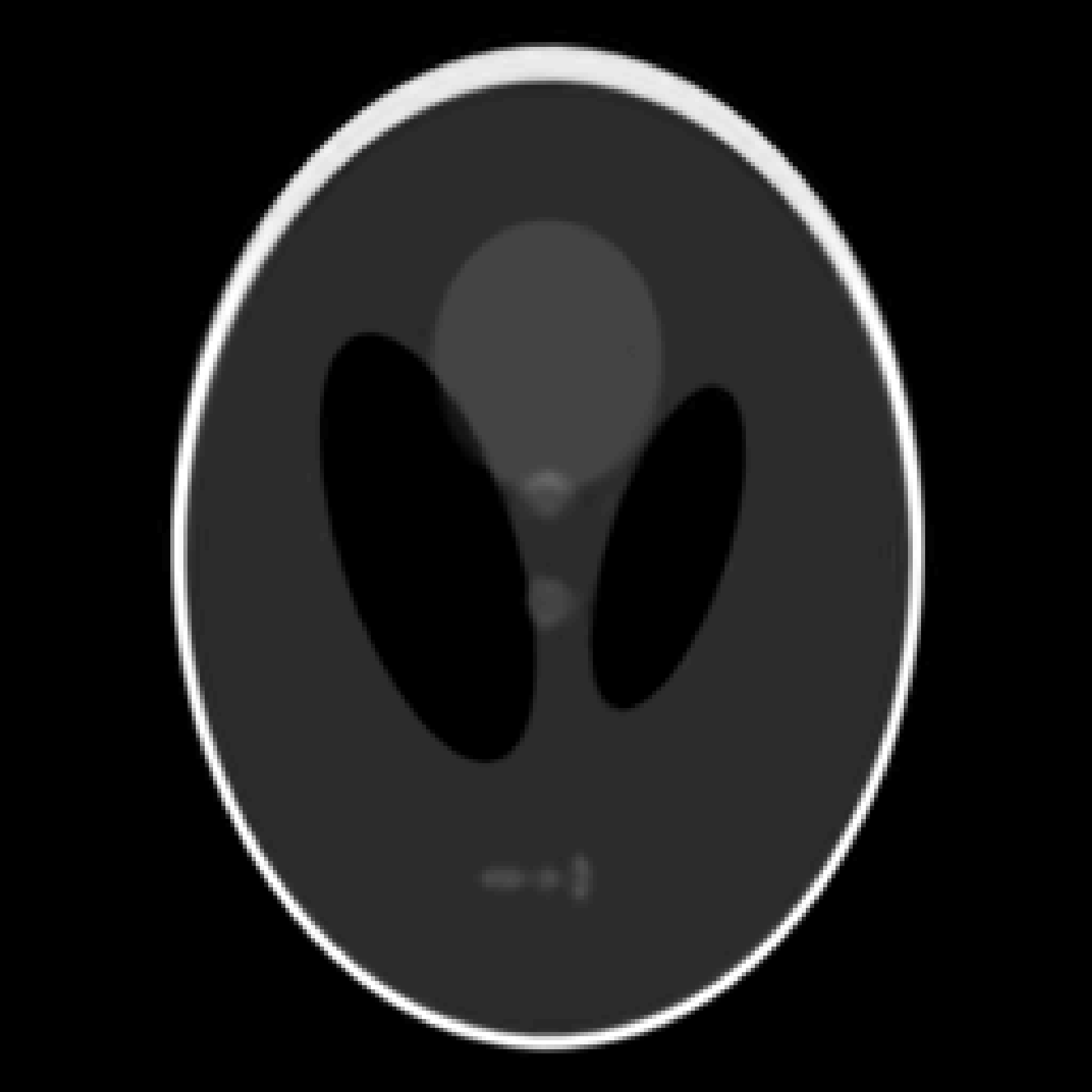}
\end{minipage}}\hspace{-0.05cm}
\subfloat[TV]{\label{EllipseILFTV}\begin{minipage}{3cm}
\includegraphics[width=3cm]{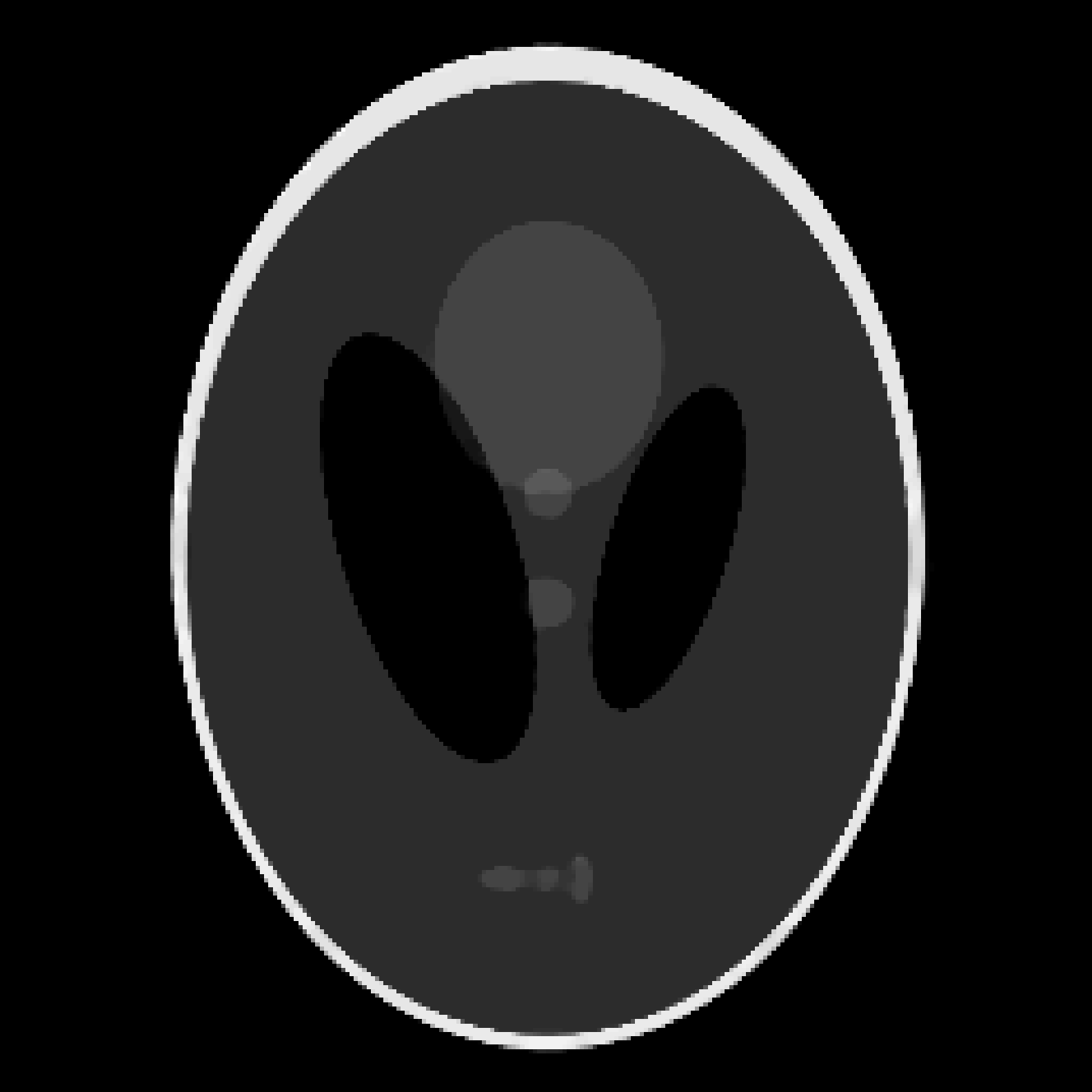}
\end{minipage}}\hspace{-0.05cm}
\subfloat[Haar]{\label{EllipseILFHaar}\begin{minipage}{3cm}
\includegraphics[width=3cm]{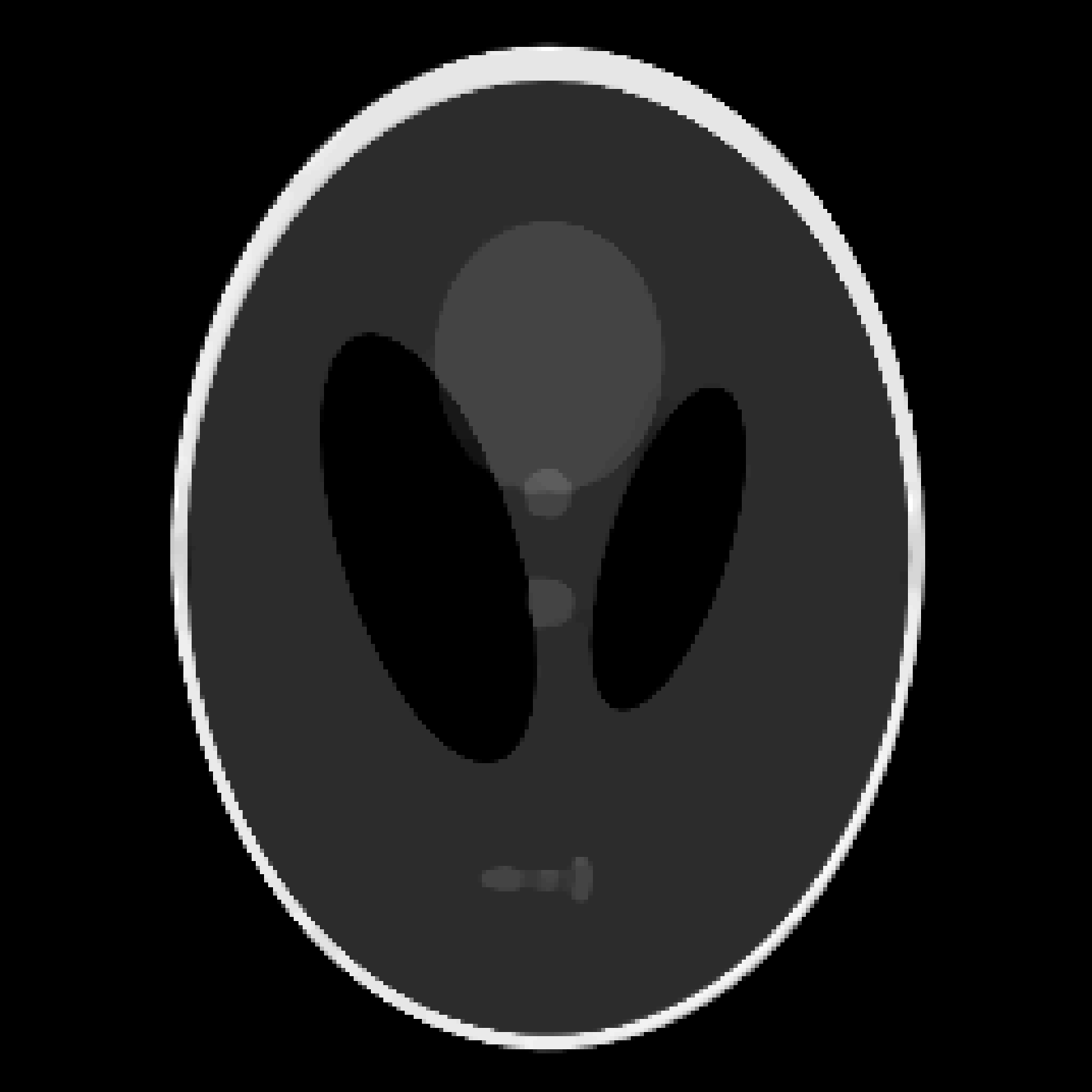}
\end{minipage}}\hspace{-0.05cm}
\subfloat[DDTF]{\label{EllipseILFDDTF}\begin{minipage}{3cm}
\includegraphics[width=3cm]{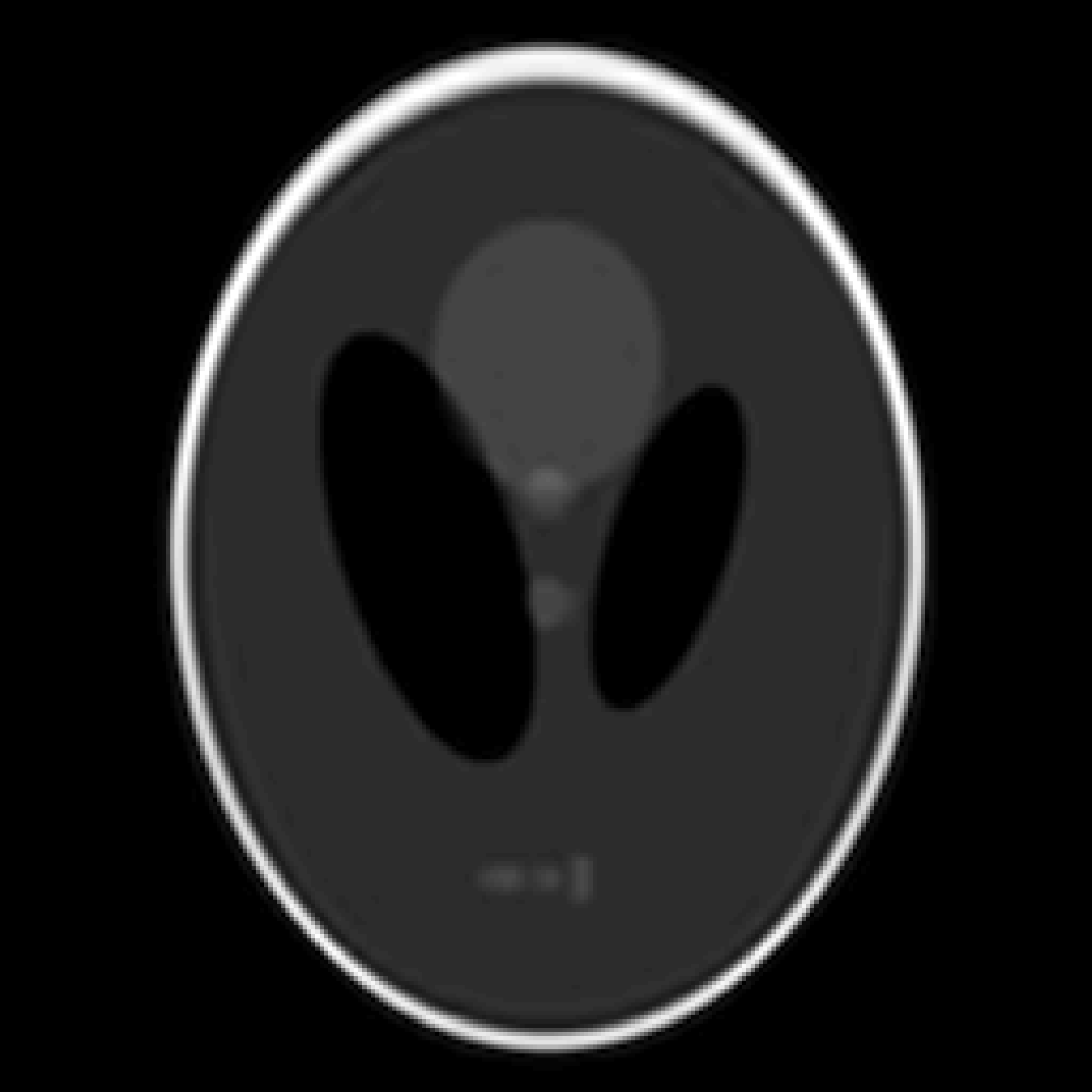}
\end{minipage}}
\caption{Visual comparison of ``Ellipse'' for ideal low-pass filter deconvolution.}\label{EllipseILFResults}
\end{figure}

\begin{figure}[t]
\centering
\subfloat[TF]{\label{EllipseILFTF}\begin{minipage}{3cm}
\includegraphics[width=3cm]{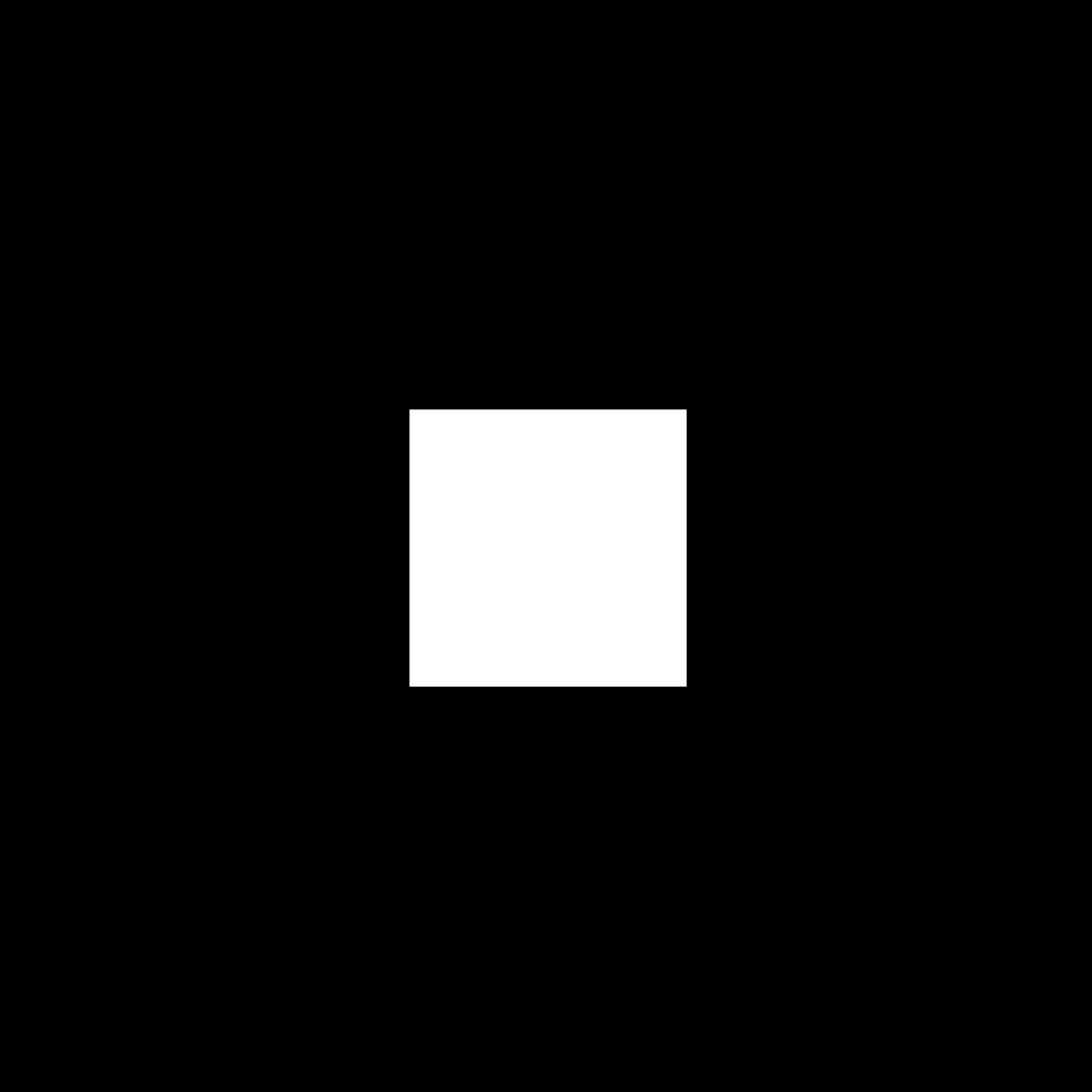}
\end{minipage}}\hspace{-0.05cm}
\subfloat[Proposed]{\label{EllipseILFProposedError}\begin{minipage}{3cm}
\includegraphics[width=3cm]{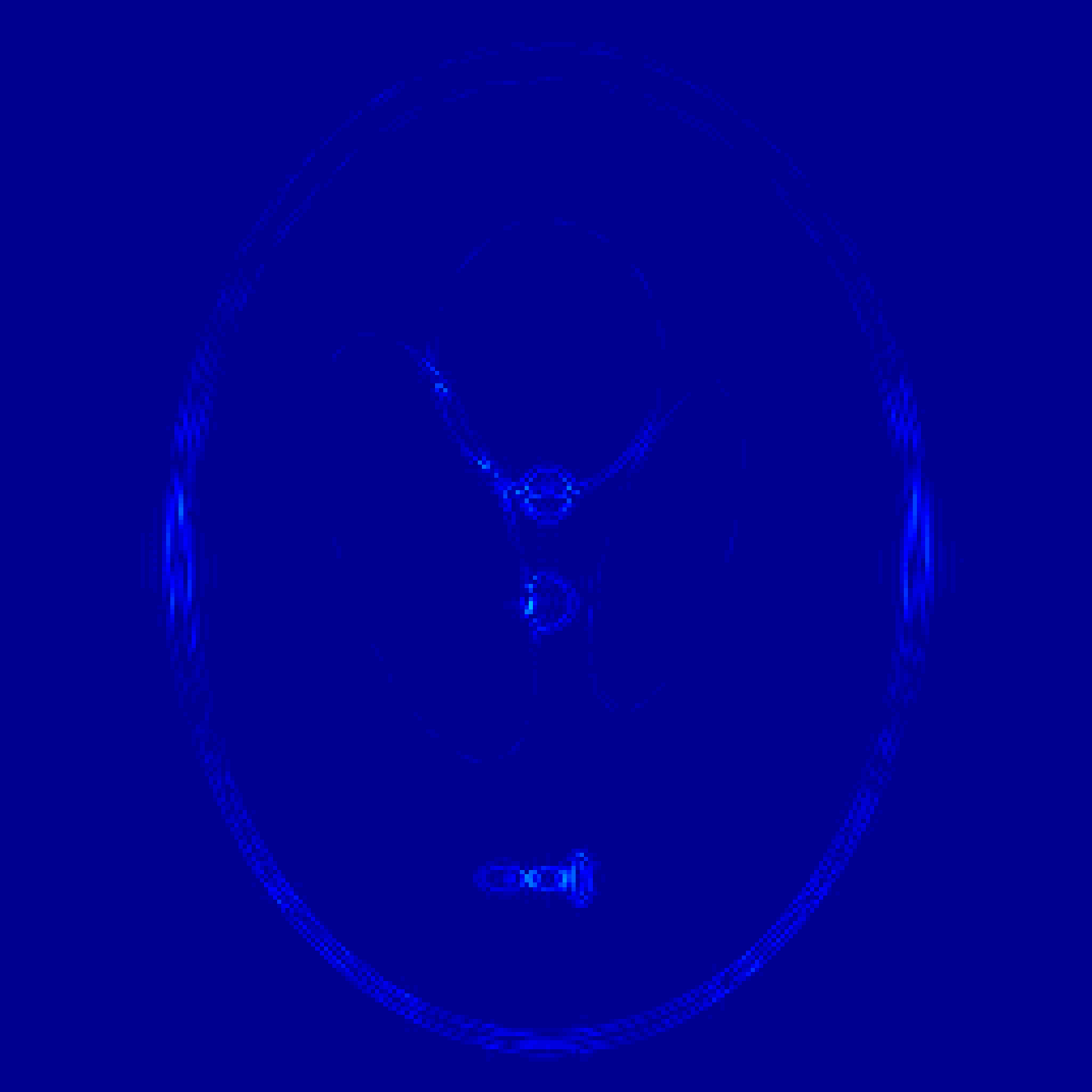}
\end{minipage}}\hspace{-0.05cm}
\subfloat[LSLP]{\label{EllipseILFLSLPError}\begin{minipage}{3cm}
\includegraphics[width=3cm]{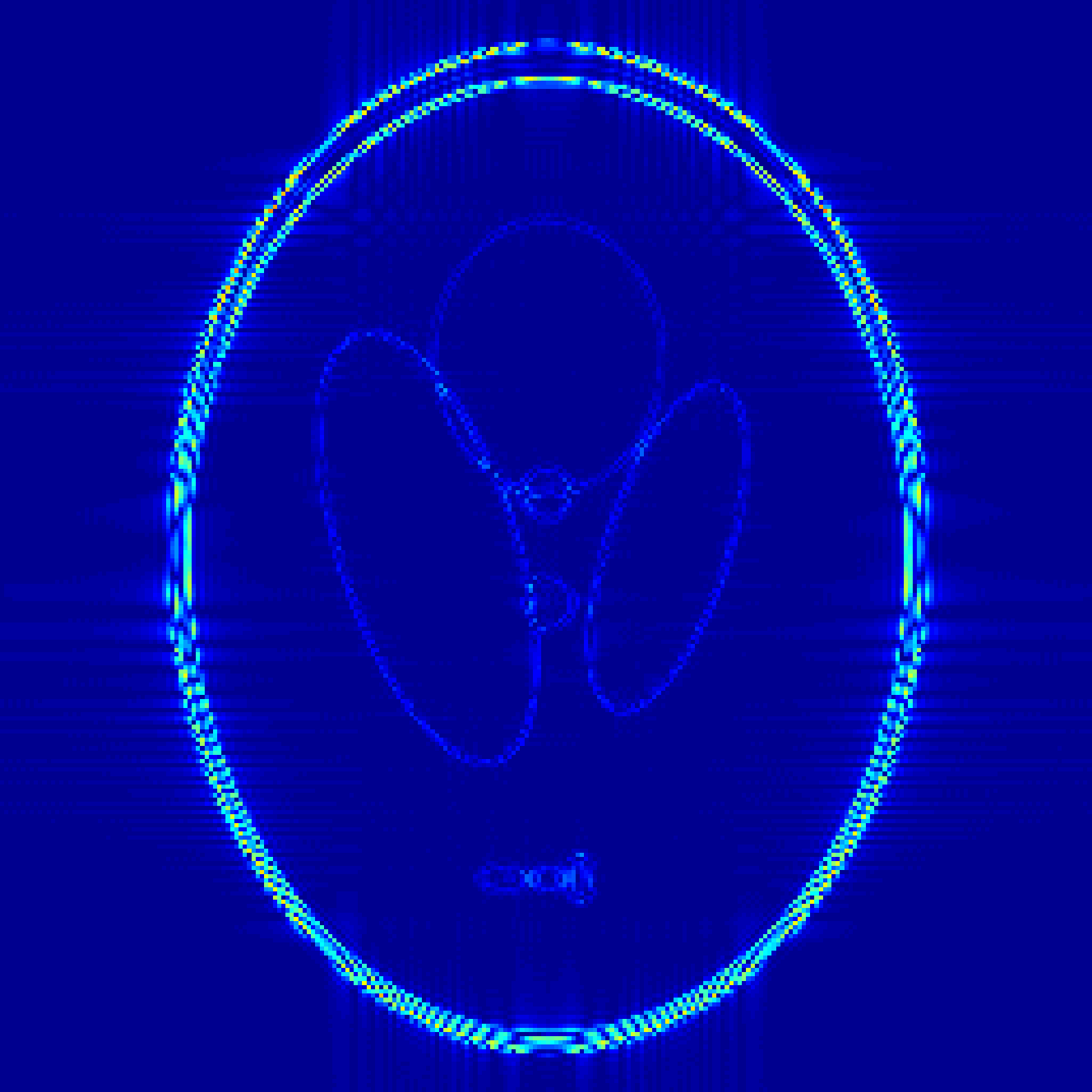}
\end{minipage}}\hspace{-0.05cm}
\subfloat[LRHDDTF]{\label{EllipseILFLRHDDTFError}\begin{minipage}{3cm}
\includegraphics[width=3cm]{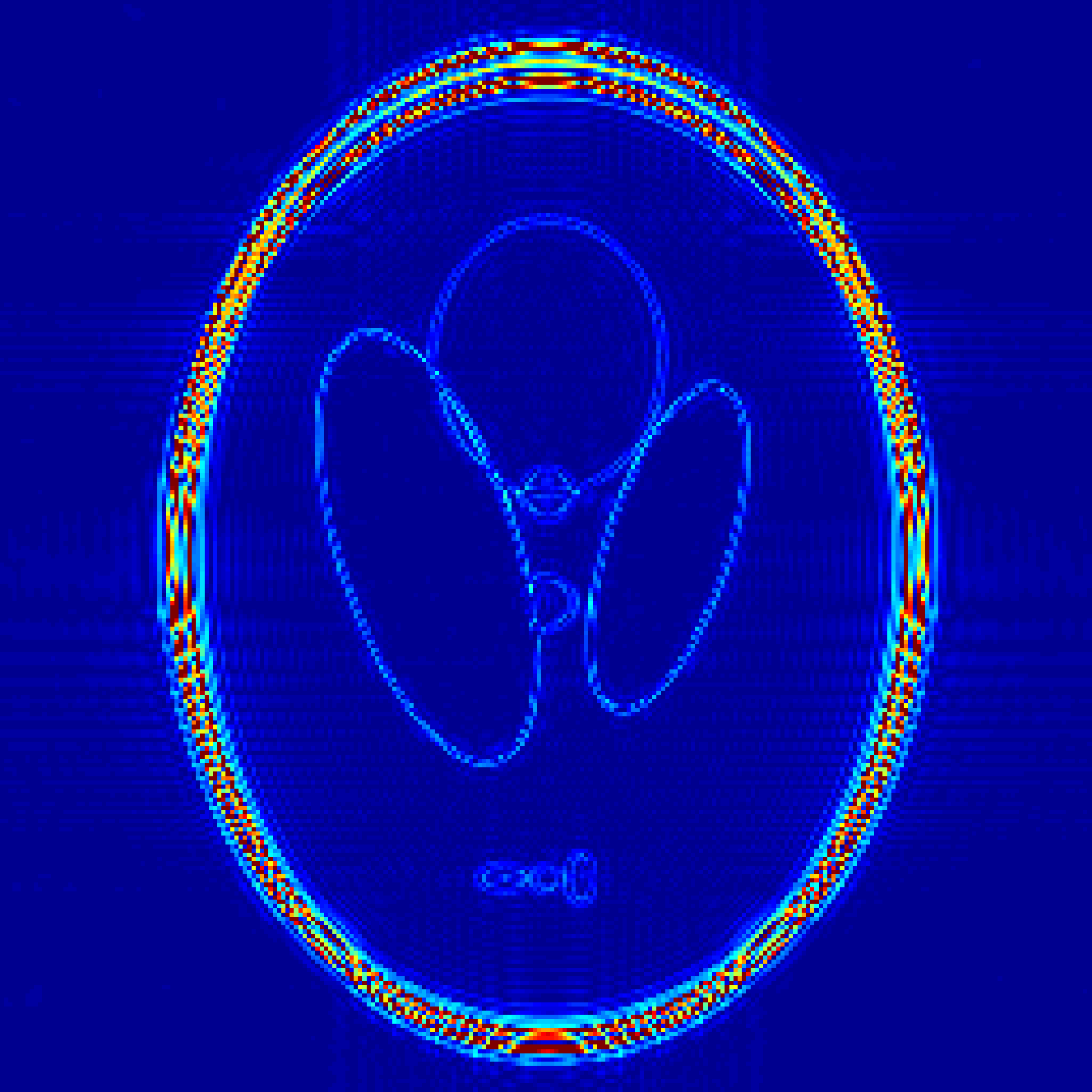}
\end{minipage}}\vspace{-0.25cm}\\
\subfloat[Schatten $0$]{\label{EllipseILFGIRAFError}\begin{minipage}{3cm}
\includegraphics[width=3cm]{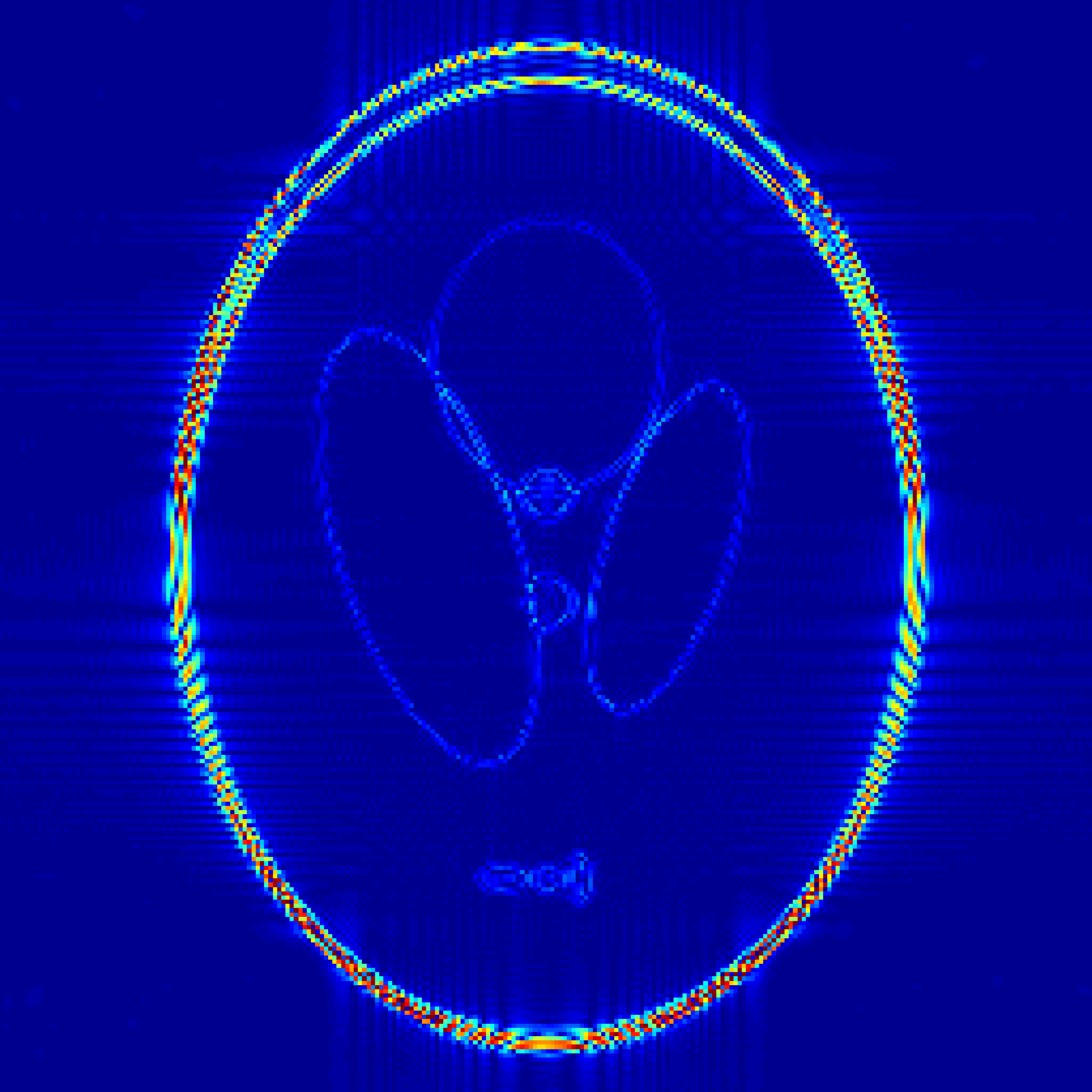}
\end{minipage}}\hspace{-0.05cm}
\subfloat[TV]{\label{EllipseILFTVError}\begin{minipage}{3cm}
\includegraphics[width=3cm]{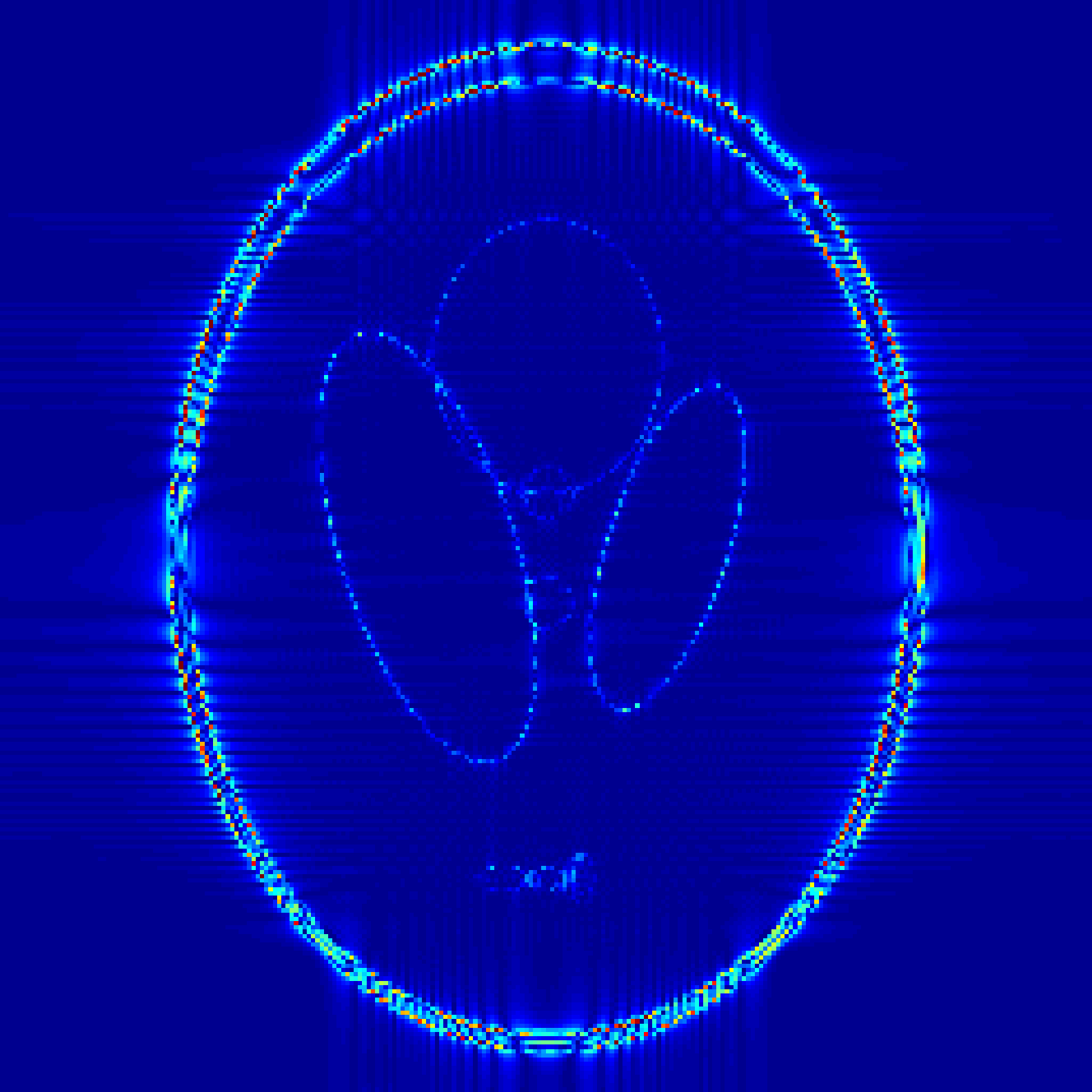}
\end{minipage}}\hspace{-0.05cm}
\subfloat[Haar]{\label{EllipseILFHaarError}\begin{minipage}{3cm}
\includegraphics[width=3cm]{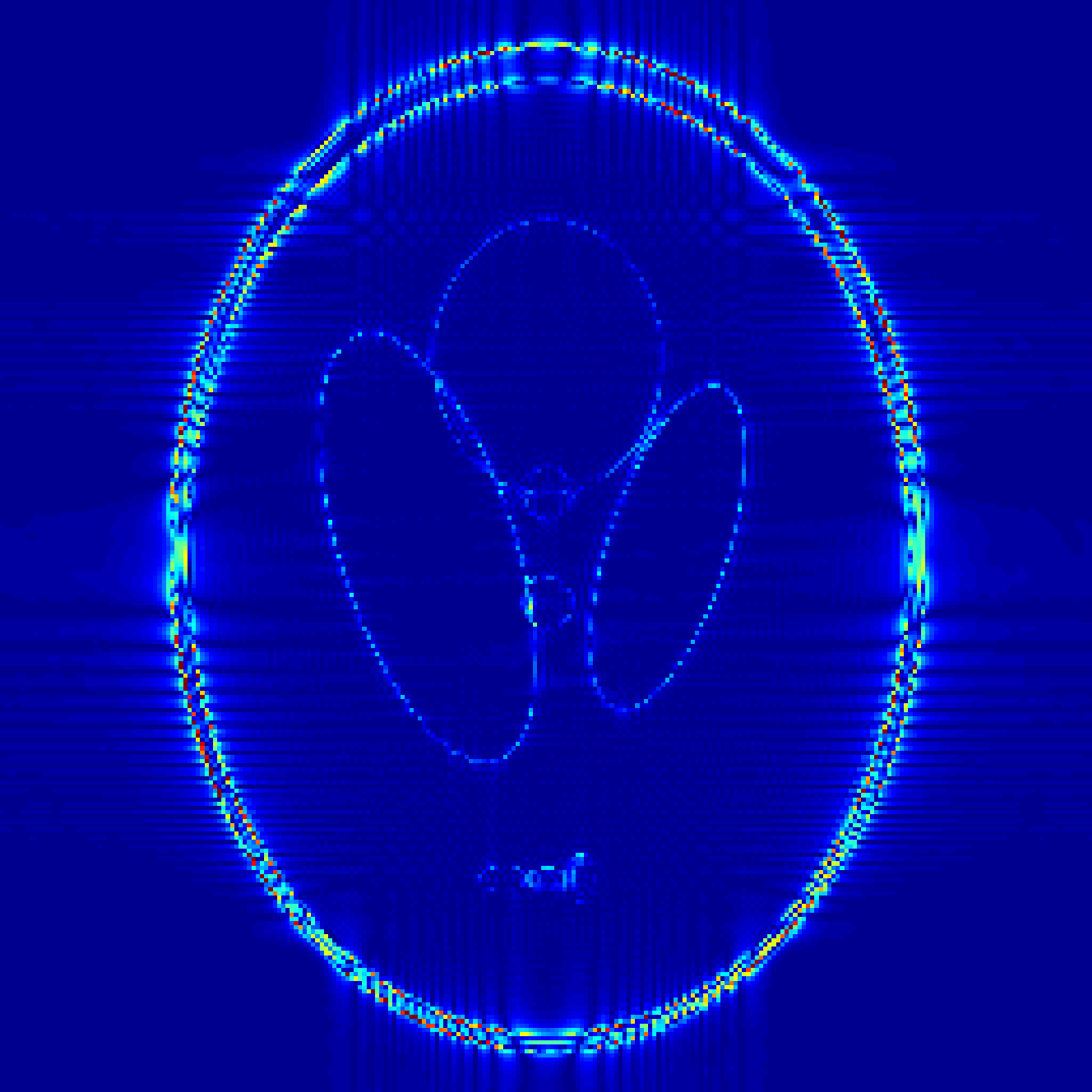}
\end{minipage}}\hspace{-0.05cm}
\subfloat[DDTF]{\label{EllipseILFDDTFError}\begin{minipage}{3cm}
\includegraphics[width=3cm]{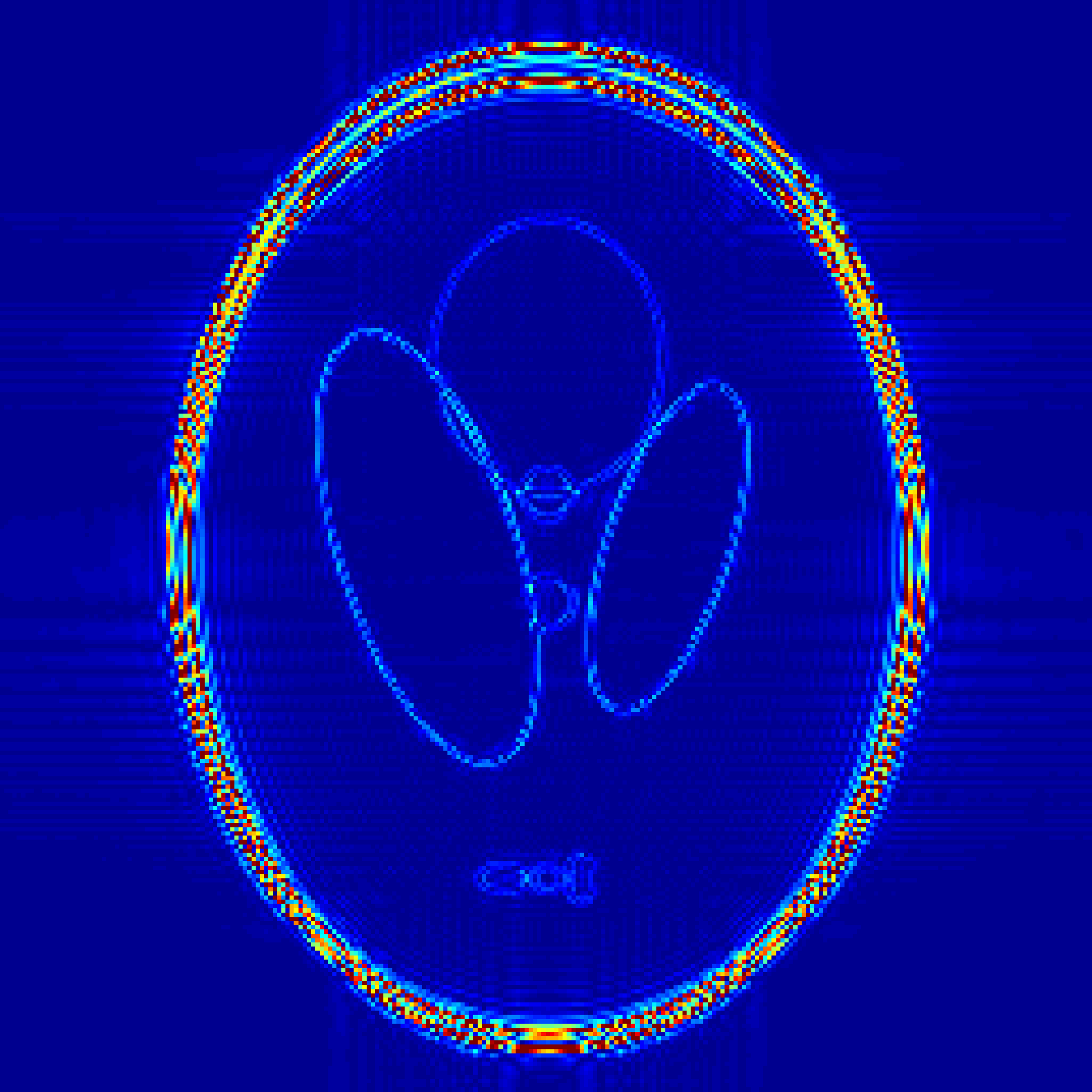}
\end{minipage}}
\caption{Error maps of \cref{EllipseILFResults}. In \cref{EllipseILFTF}, TF stands for the transfer function.}\label{EllipseILFError}
\end{figure}

\begin{figure}[t]
\centering
\subfloat[Ref.]{\label{RectangleRandomOriginal}\begin{minipage}{3cm}
\includegraphics[width=3cm]{RectanglesOriginal.pdf}
\end{minipage}}\hspace{-0.05cm}
\subfloat[Proposed]{\label{RectangleRandomProposed}\begin{minipage}{3cm}
\includegraphics[width=3cm]{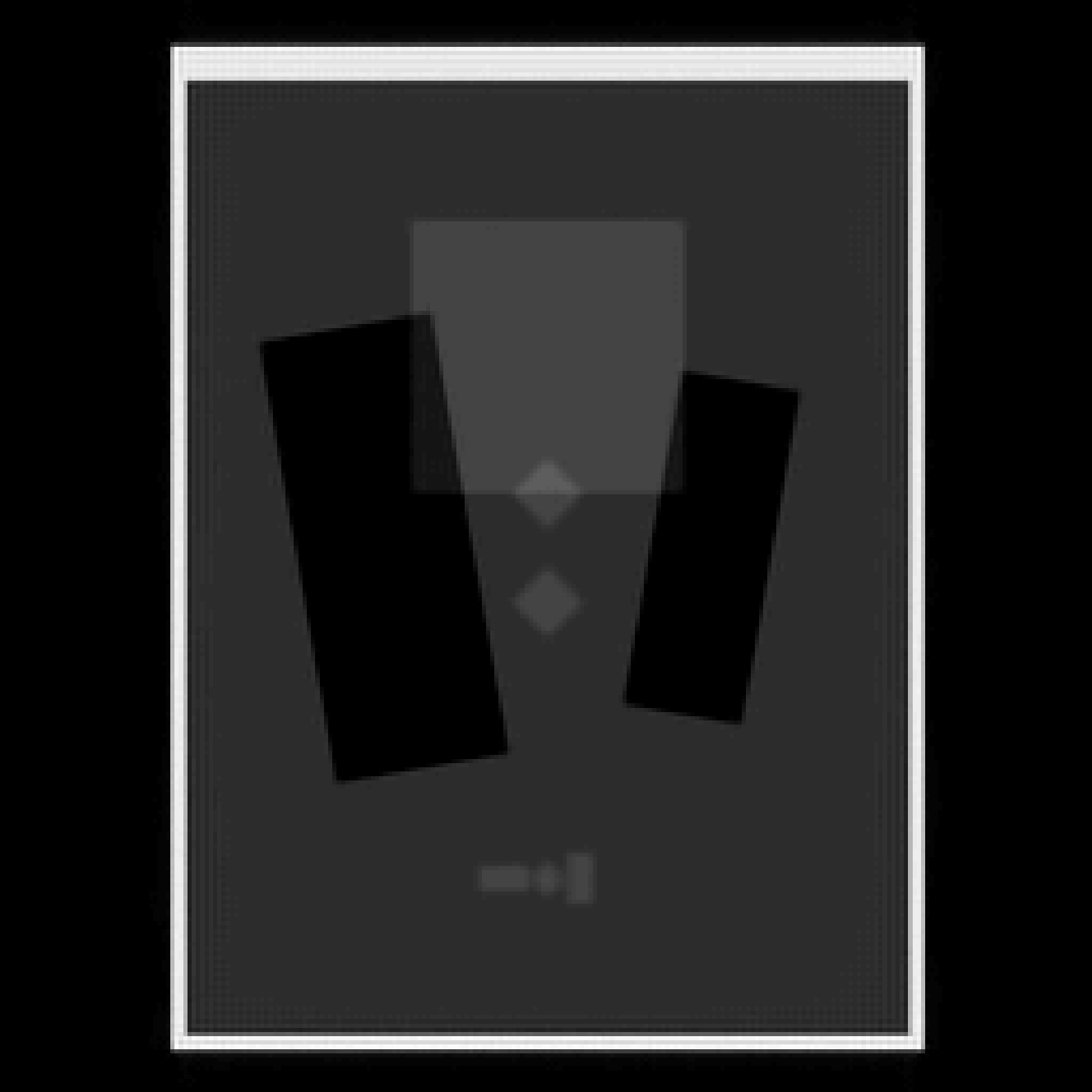}
\end{minipage}}\hspace{-0.05cm}
\subfloat[LSLP]{\label{RectangleRandomLSLP}\begin{minipage}{3cm}
\includegraphics[width=3cm]{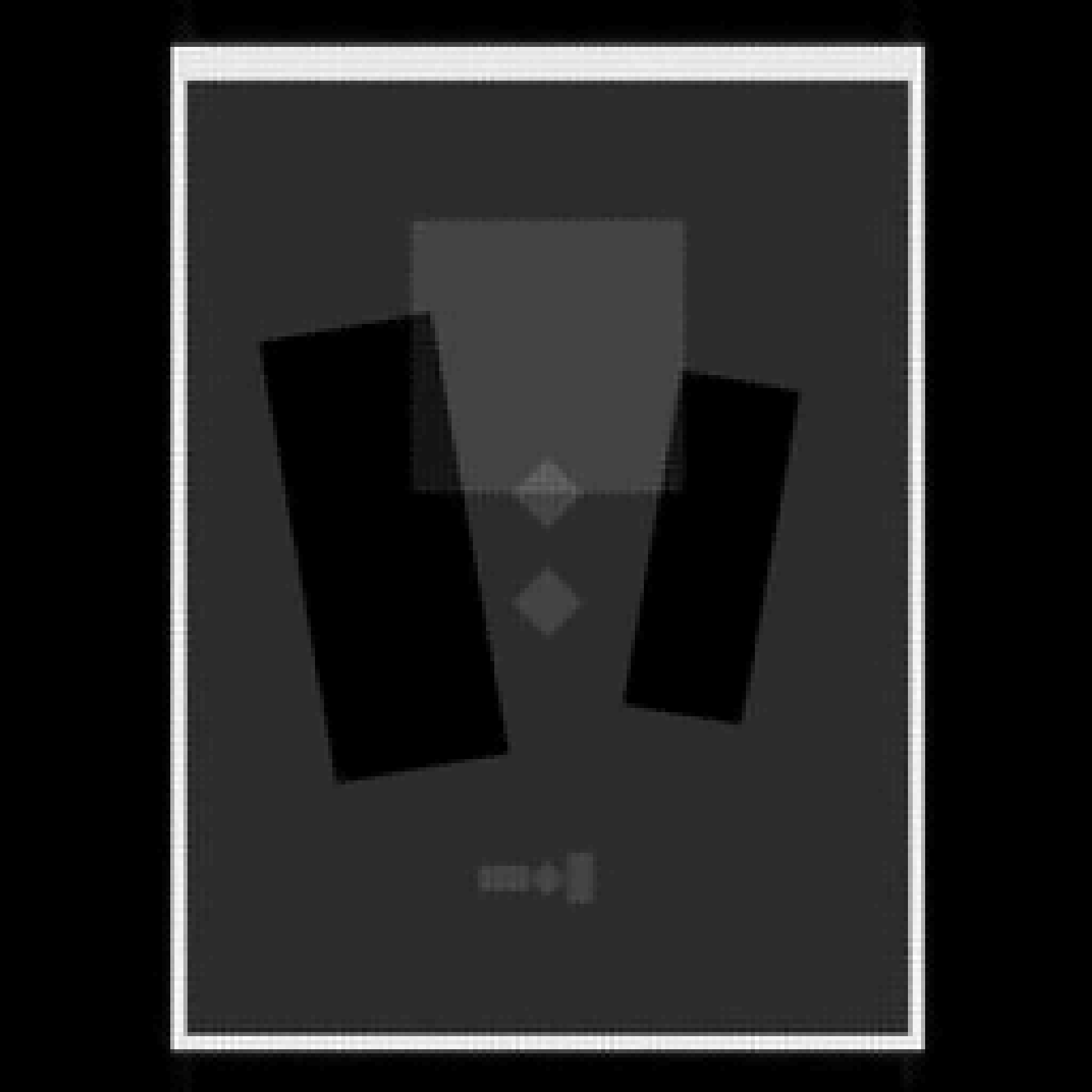}
\end{minipage}}\hspace{-0.05cm}
\subfloat[LRHDDTF]{\label{RectangleRandomLRHDDTF}\begin{minipage}{3cm}
\includegraphics[width=3cm]{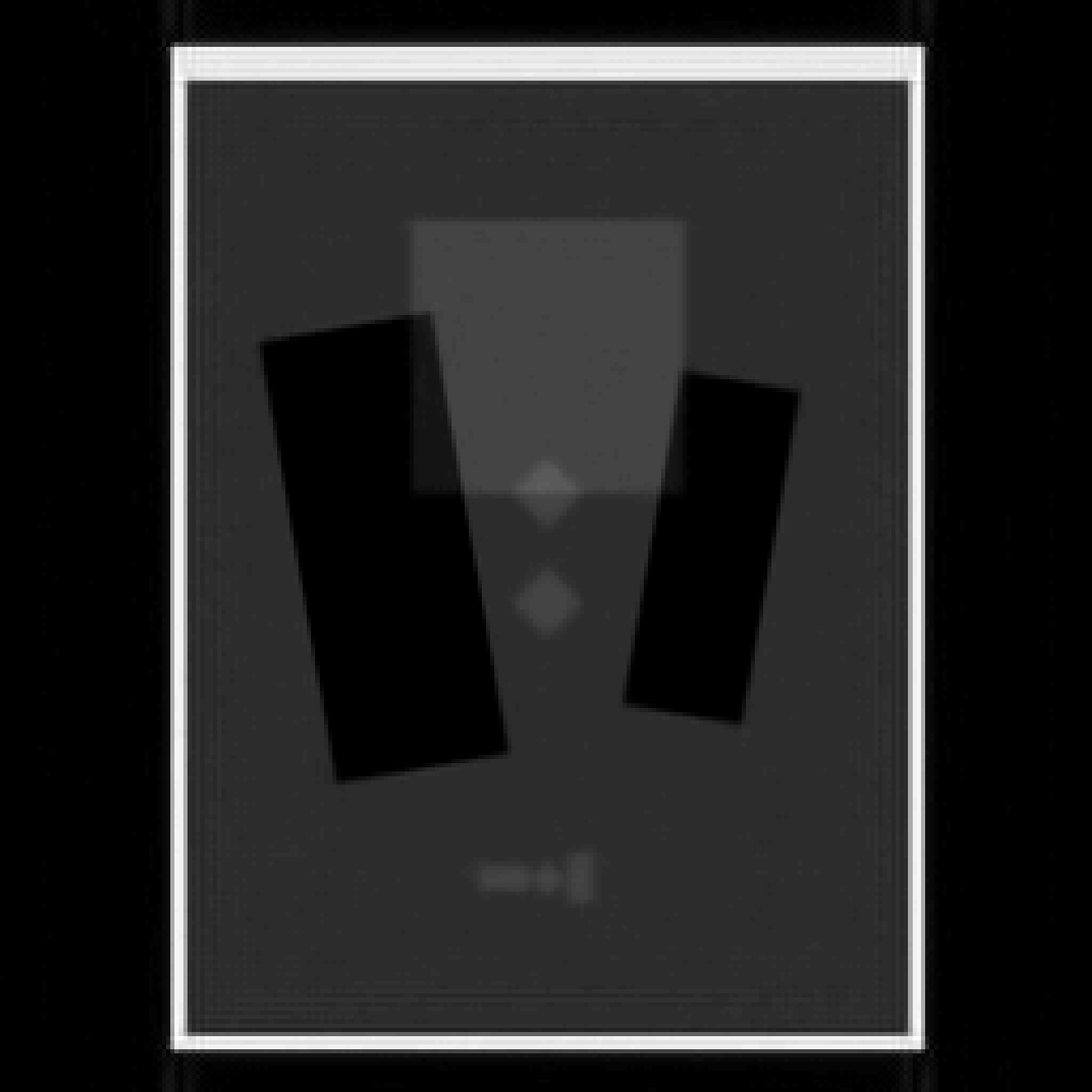}
\end{minipage}}\vspace{-0.25cm}\\
\subfloat[Schatten $0$]{\label{RectangleRandomGIRAF}\begin{minipage}{3cm}
\includegraphics[width=3cm]{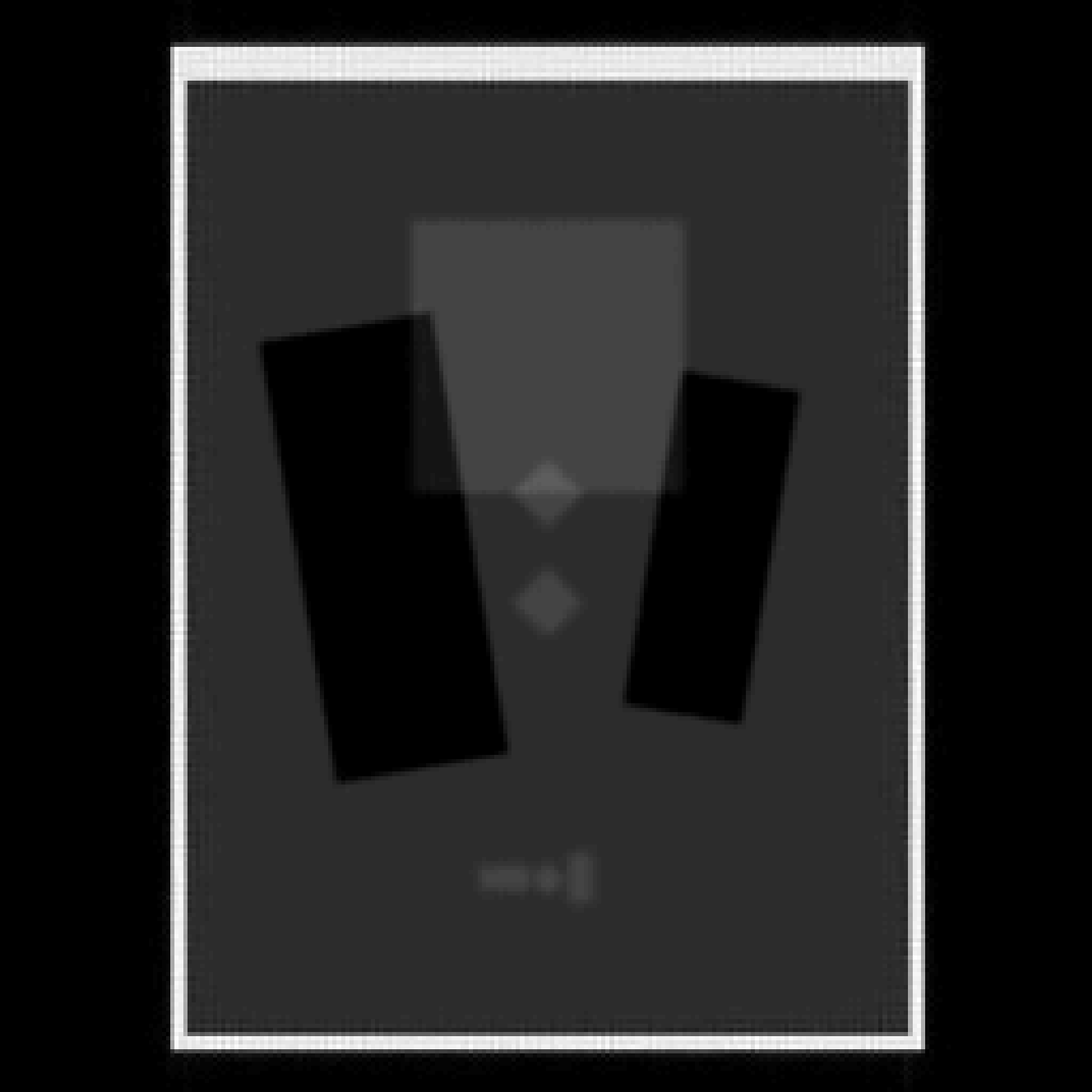}
\end{minipage}}\hspace{-0.05cm}
\subfloat[TV]{\label{RectangleRandomTV}\begin{minipage}{3cm}
\includegraphics[width=3cm]{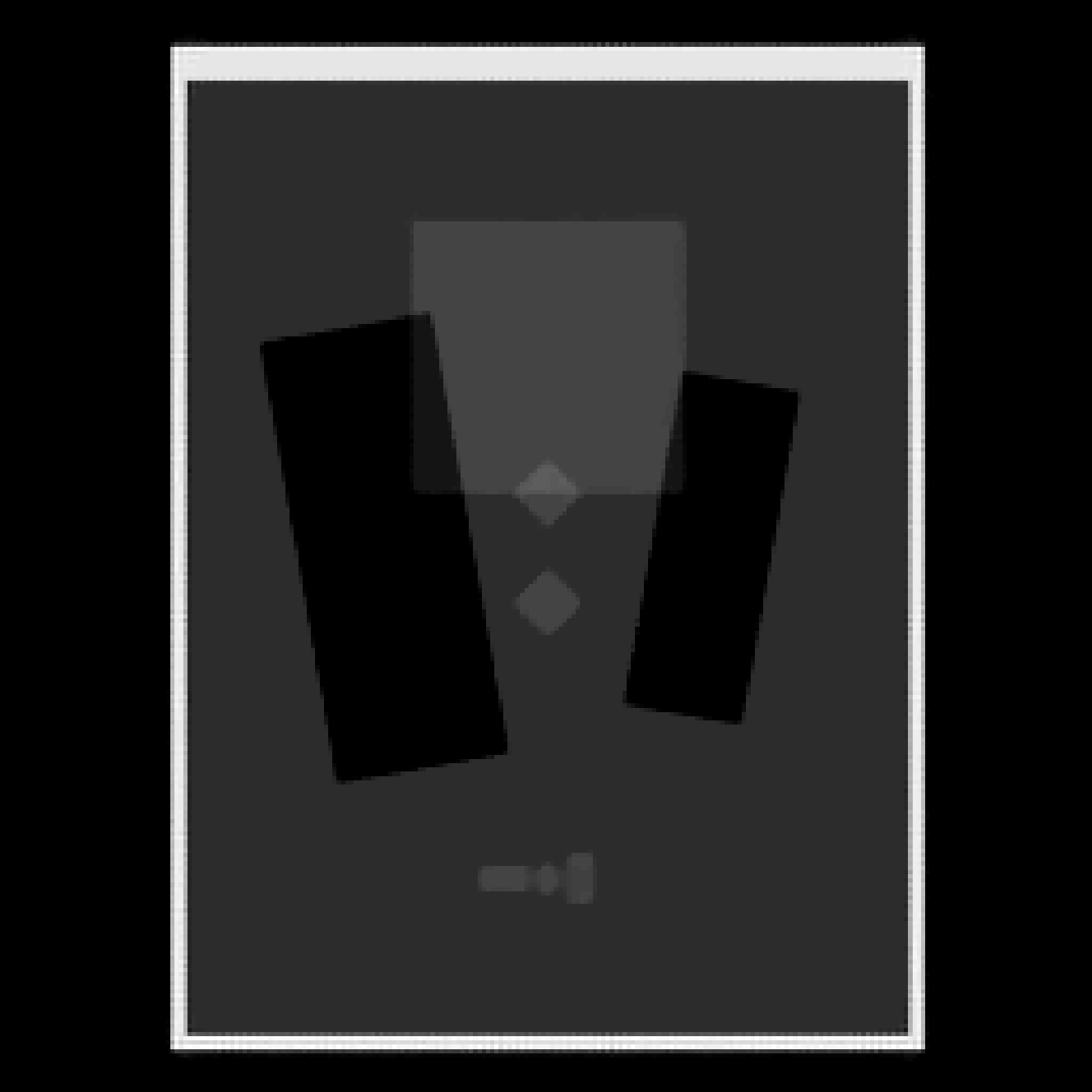}
\end{minipage}}\hspace{-0.05cm}
\subfloat[Haar]{\label{RectangleRandomHaar}\begin{minipage}{3cm}
\includegraphics[width=3cm]{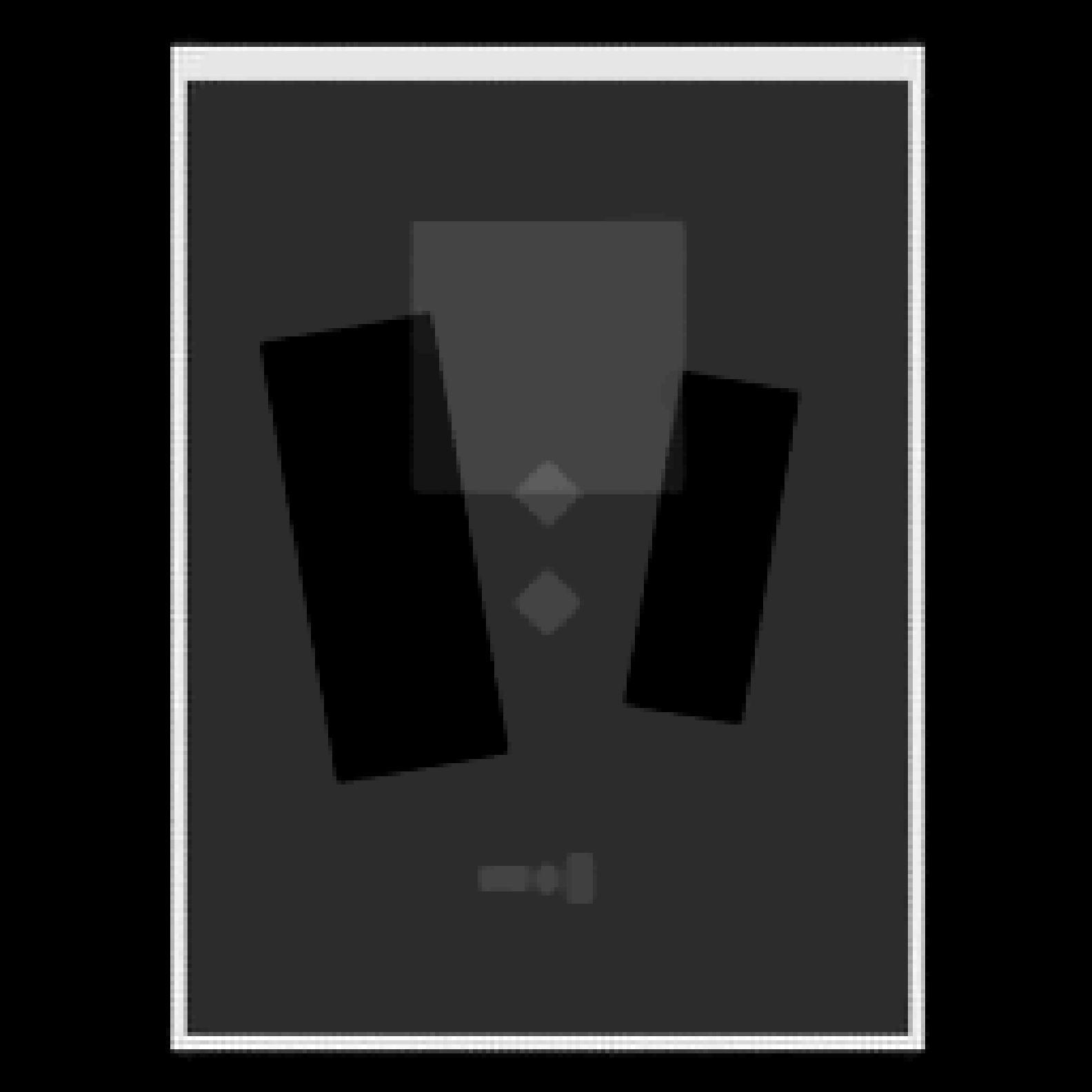}
\end{minipage}}\hspace{-0.05cm}
\subfloat[DDTF]{\label{RectangleRandomDDTF}\begin{minipage}{3cm}
\includegraphics[width=3cm]{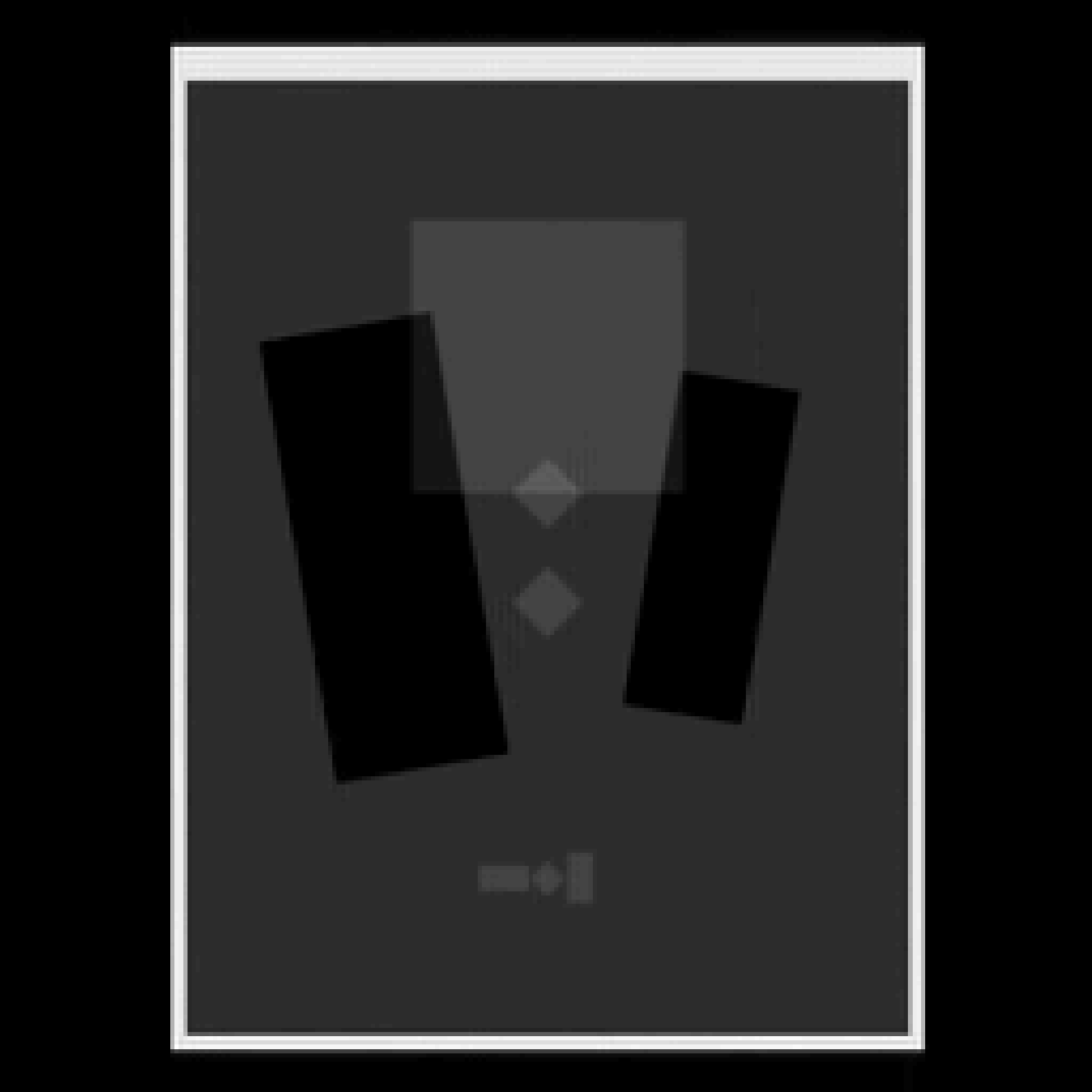}
\end{minipage}}
\caption{Visual comparison of ``Rectangle'' for random sampling.}\label{RectangleRandomResults}
\end{figure}

\begin{figure}[ht]
\centering
\subfloat[Samples]{\label{RectangleRandomOriginal}\begin{minipage}{3cm}
\includegraphics[width=3cm]{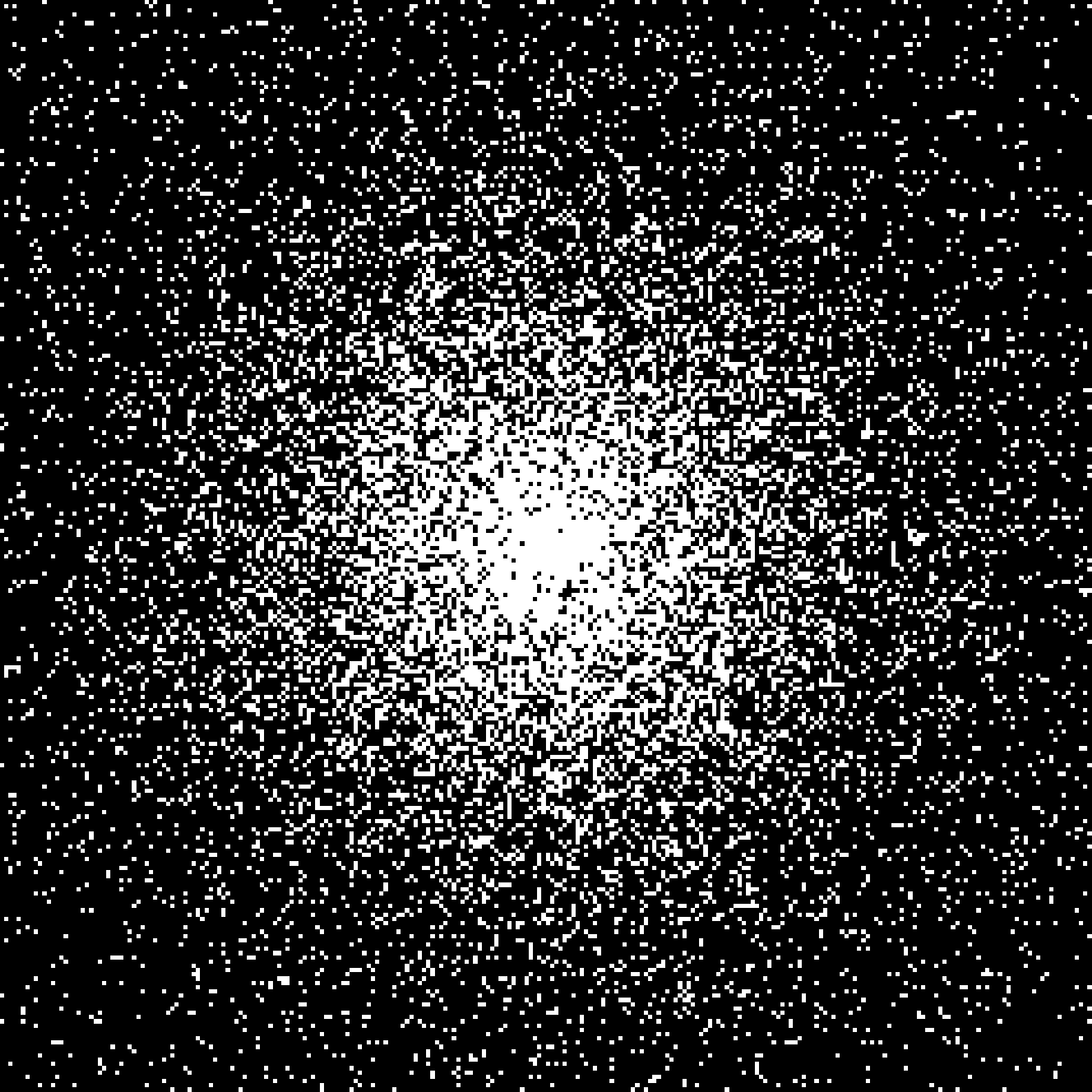}
\end{minipage}}\hspace{-0.05cm}
\subfloat[Proposed]{\label{RectangleRandomProposedError}\begin{minipage}{3cm}
\includegraphics[width=3cm]{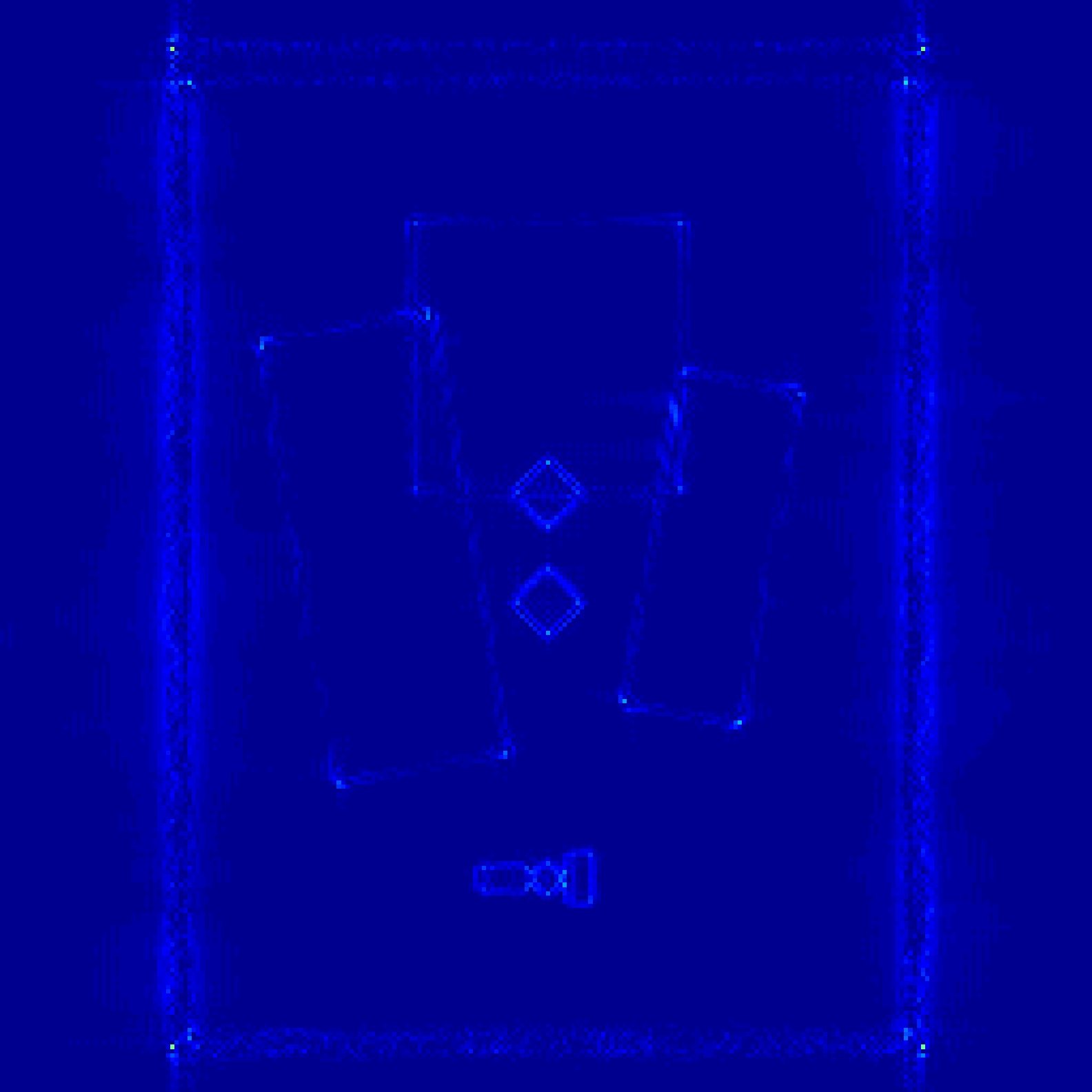}
\end{minipage}}\hspace{-0.05cm}
\subfloat[LSLP]{\label{RectangleRandomLSLPError}\begin{minipage}{3cm}
\includegraphics[width=3cm]{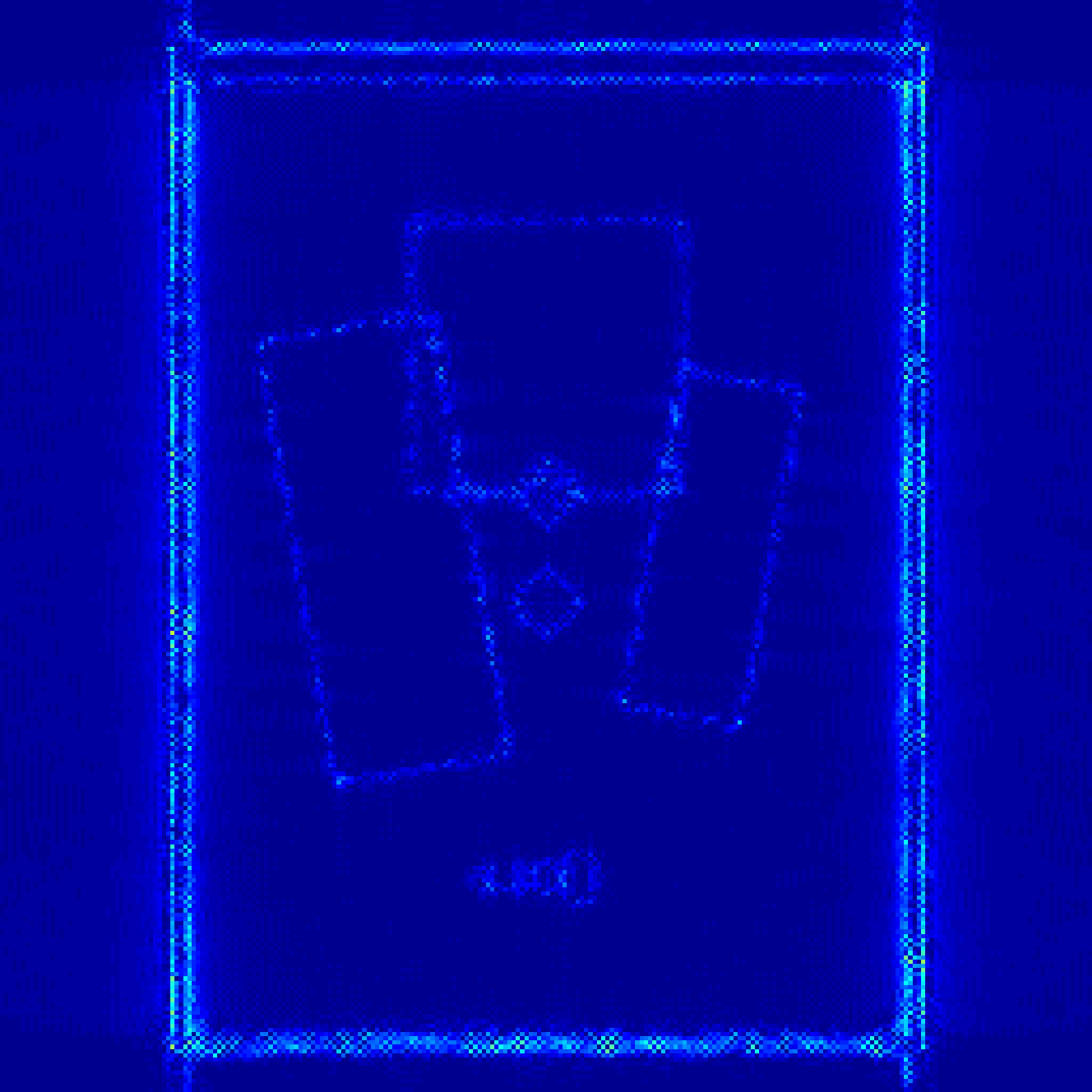}
\end{minipage}}\hspace{-0.05cm}
\subfloat[LRHDDTF]{\label{RectangleRandomLRHDDTFError}\begin{minipage}{3cm}
\includegraphics[width=3cm]{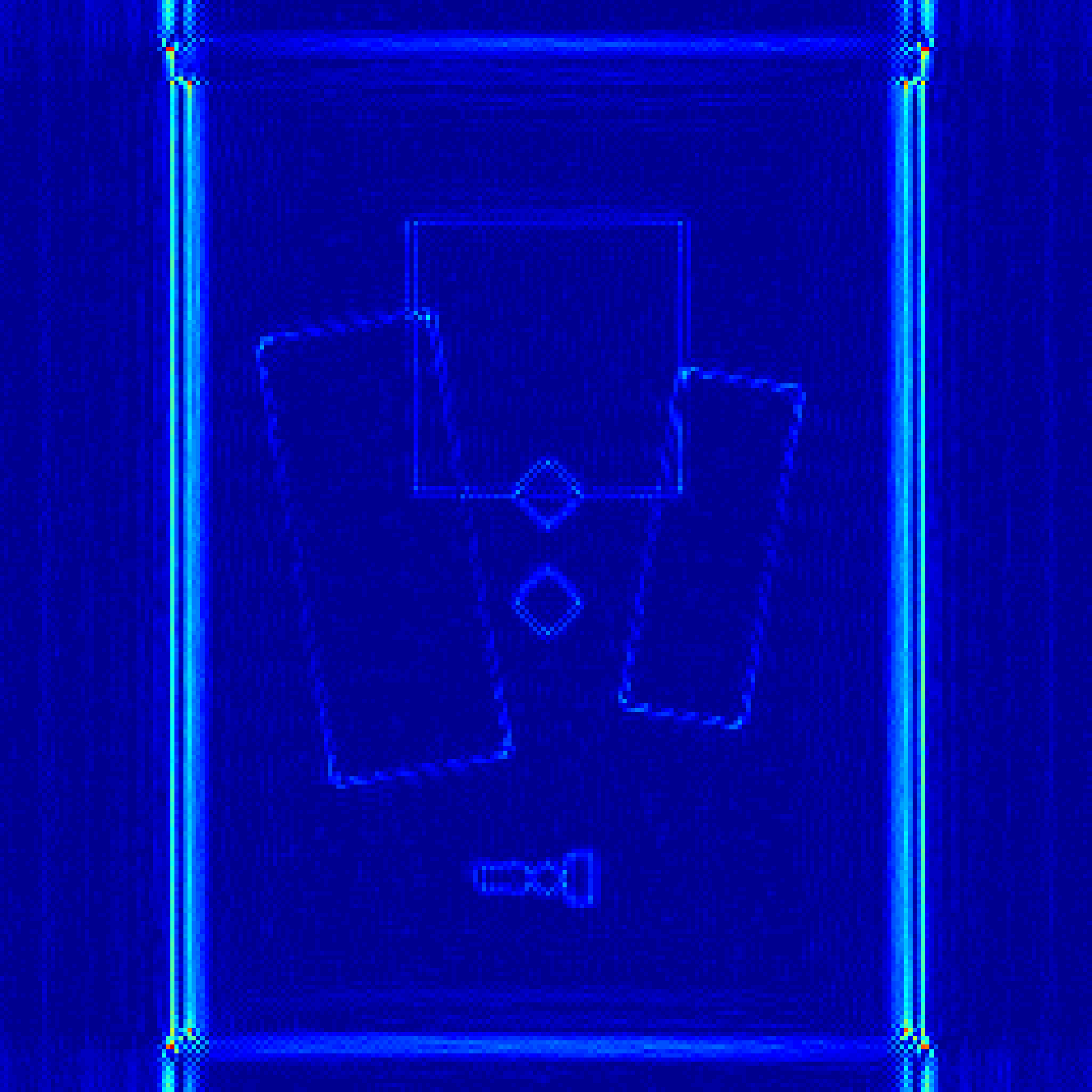}
\end{minipage}}\vspace{-0.25cm}\\
\subfloat[Schatten $0$]{\label{RectangleRandomGIRAFError}\begin{minipage}{3cm}
\includegraphics[width=3cm]{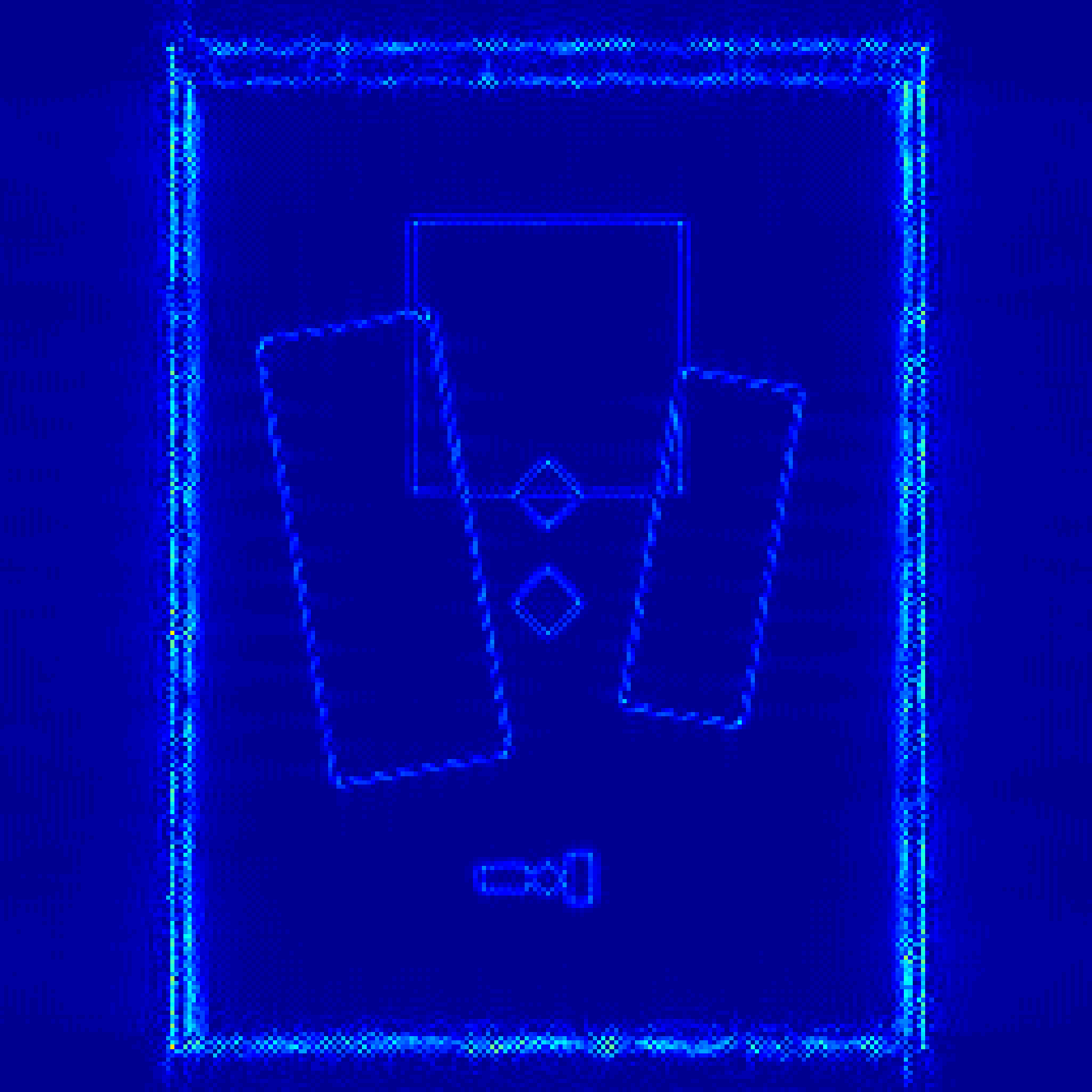}
\end{minipage}}\hspace{-0.05cm}
\subfloat[TV]{\label{RectangleRandomTVError}\begin{minipage}{3cm}
\includegraphics[width=3cm]{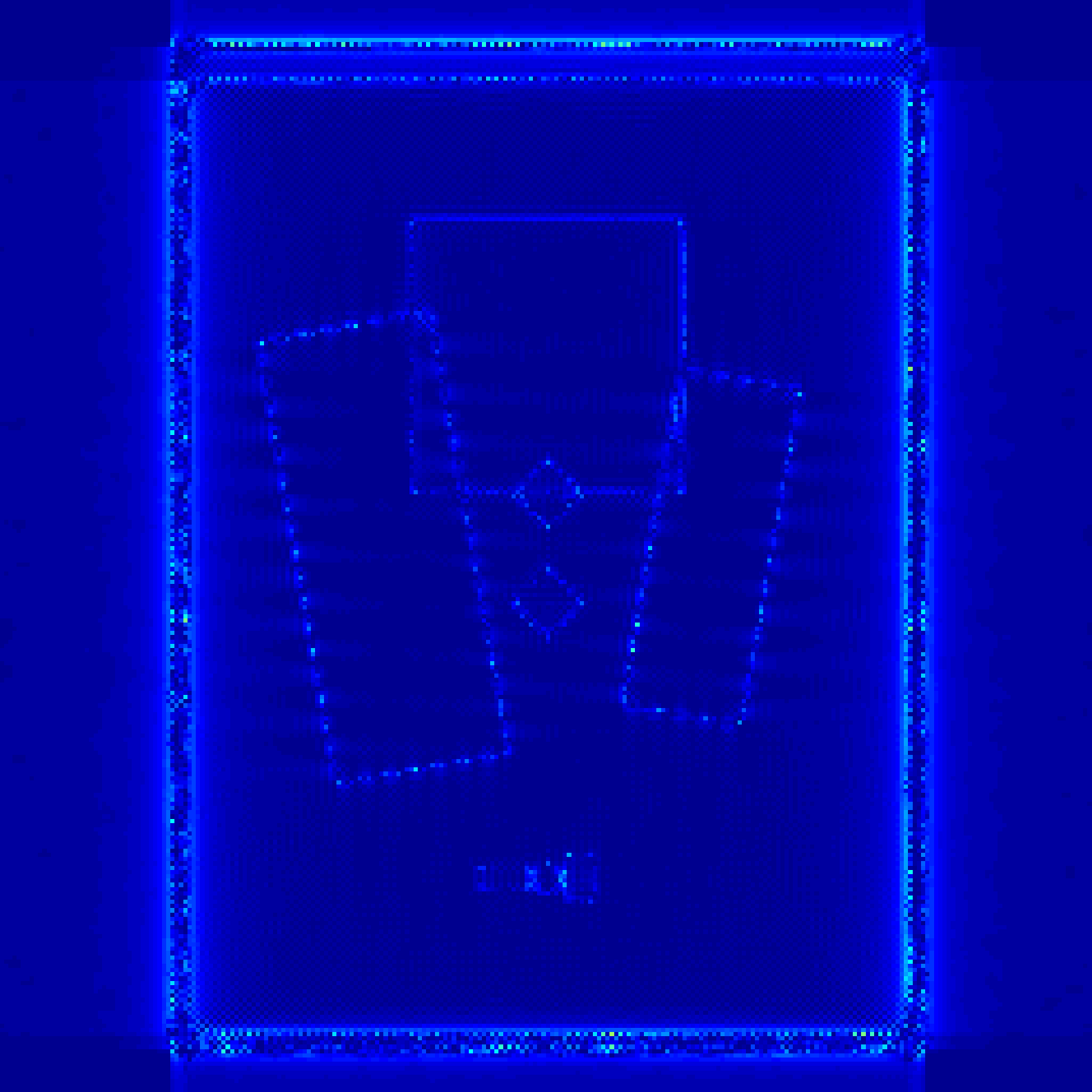}
\end{minipage}}\hspace{-0.05cm}
\subfloat[Haar]{\label{RectangleRandomHaarError}\begin{minipage}{3cm}
\includegraphics[width=3cm]{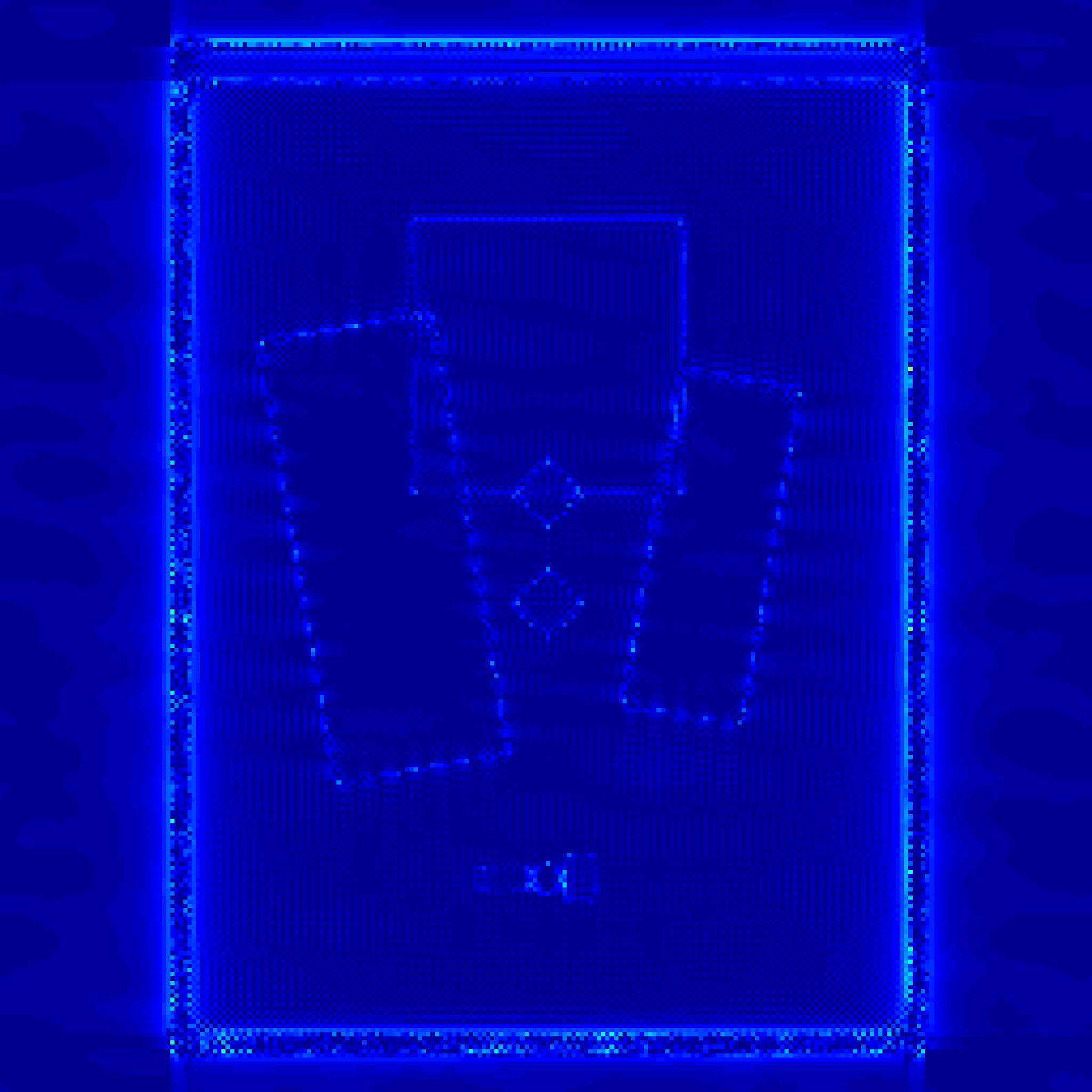}
\end{minipage}}\hspace{-0.05cm}
\subfloat[DDTF]{\label{RectangleRandomDDTFError}\begin{minipage}{3cm}
\includegraphics[width=3cm]{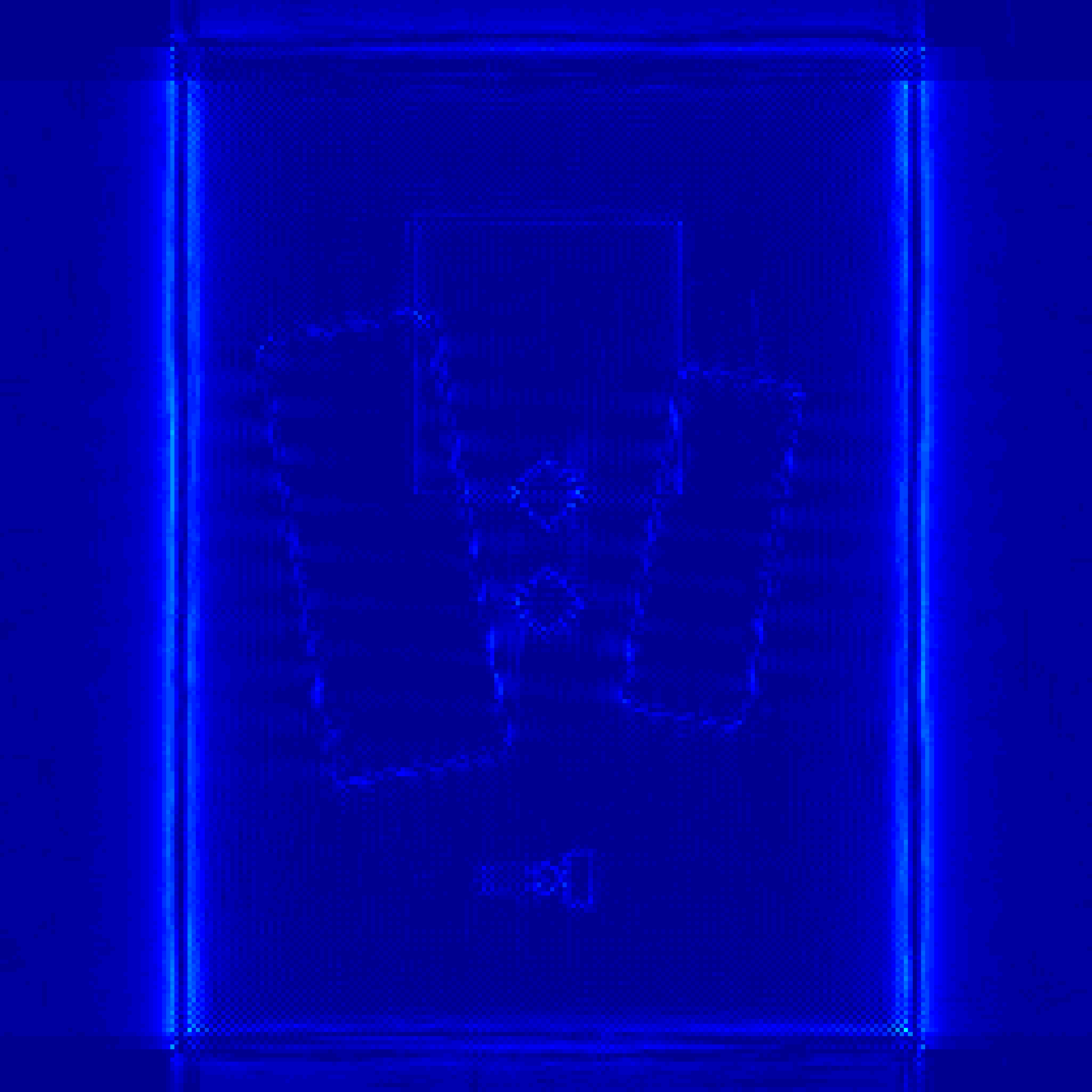}
\end{minipage}}
\caption{Error maps of \cref{RectangleRandomResults}. All error maps are displayed in the window level $[0,0.02]$ for the fair comparisons.}\label{RectangleRandomError}
\end{figure}

\begin{figure}[t]
\centering
\subfloat[Ref.]{\label{RectangleILFOriginal}\begin{minipage}{3cm}
\includegraphics[width=3cm]{RectanglesOriginal.pdf}
\end{minipage}}\hspace{-0.05cm}
\subfloat[Proposed]{\label{RectangleILFProposed}\begin{minipage}{3cm}
\includegraphics[width=3cm]{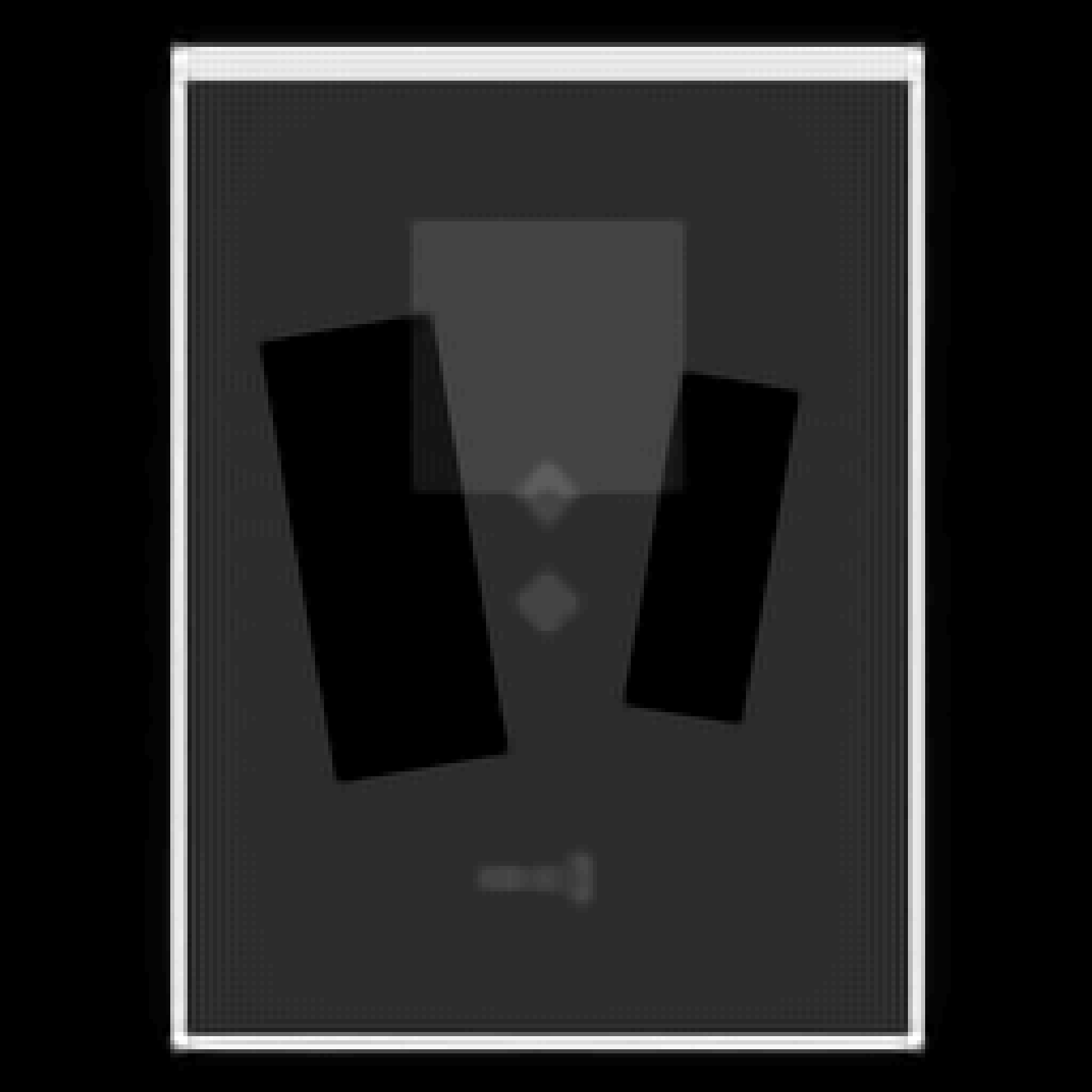}
\end{minipage}}\hspace{-0.05cm}
\subfloat[LSLP]{\label{RectangleILFLSLP}\begin{minipage}{3cm}
\includegraphics[width=3cm]{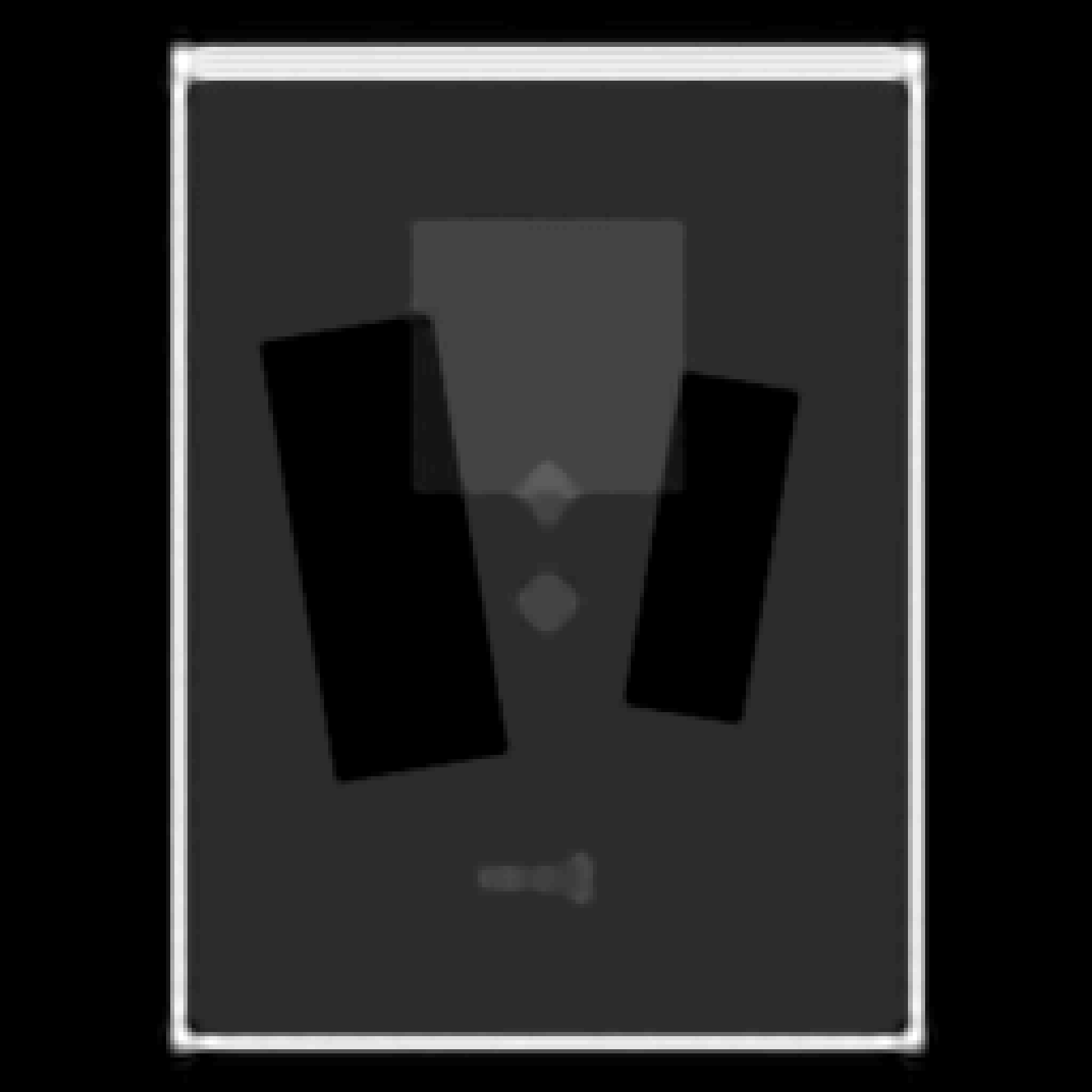}
\end{minipage}}\hspace{-0.05cm}
\subfloat[LRHDDTF]{\label{RectangleILFLRHDDTF}\begin{minipage}{3cm}
\includegraphics[width=3cm]{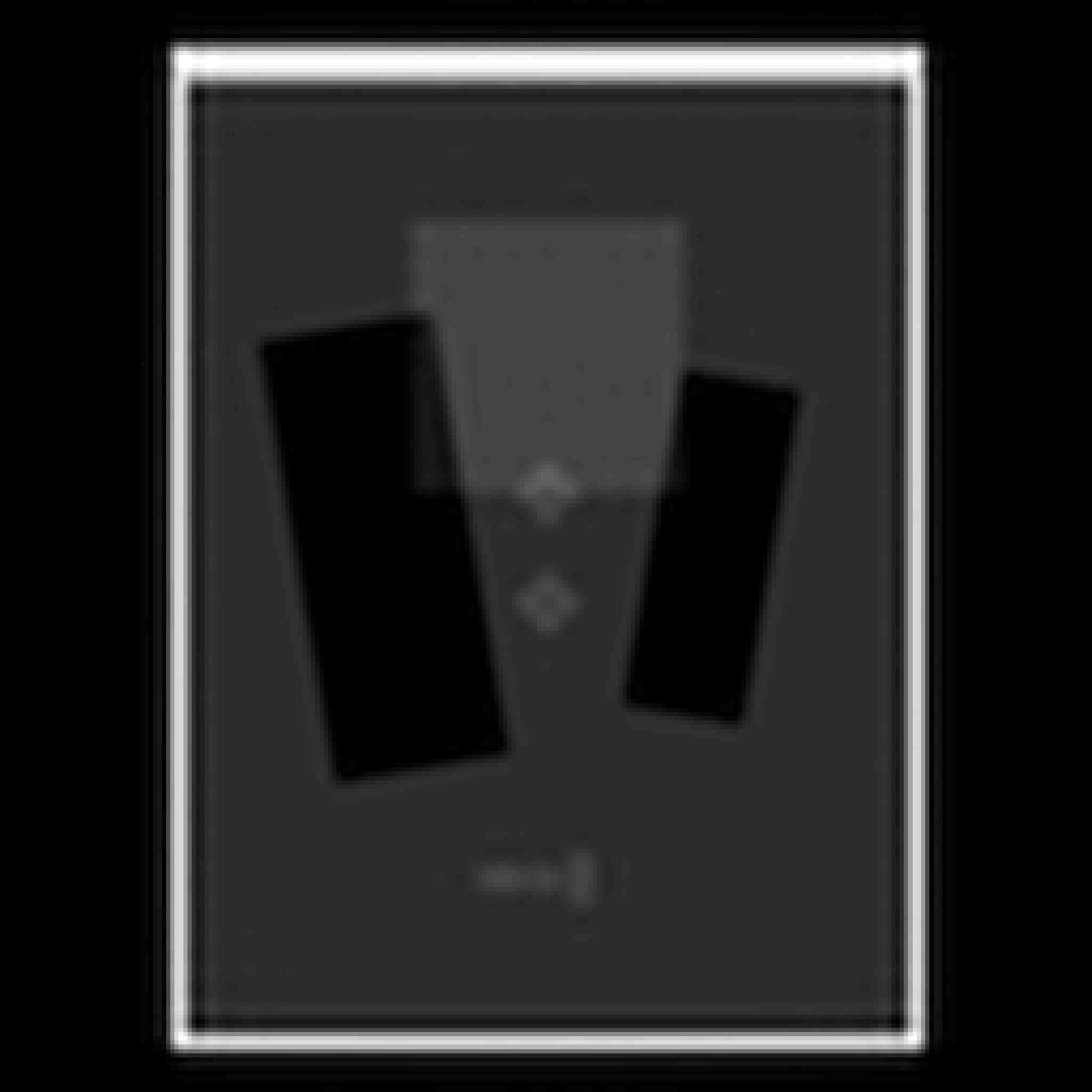}
\end{minipage}}\vspace{-0.25cm}\\
\subfloat[Schatten $0$]{\label{RectangleILFGIRAF}\begin{minipage}{3cm}
\includegraphics[width=3cm]{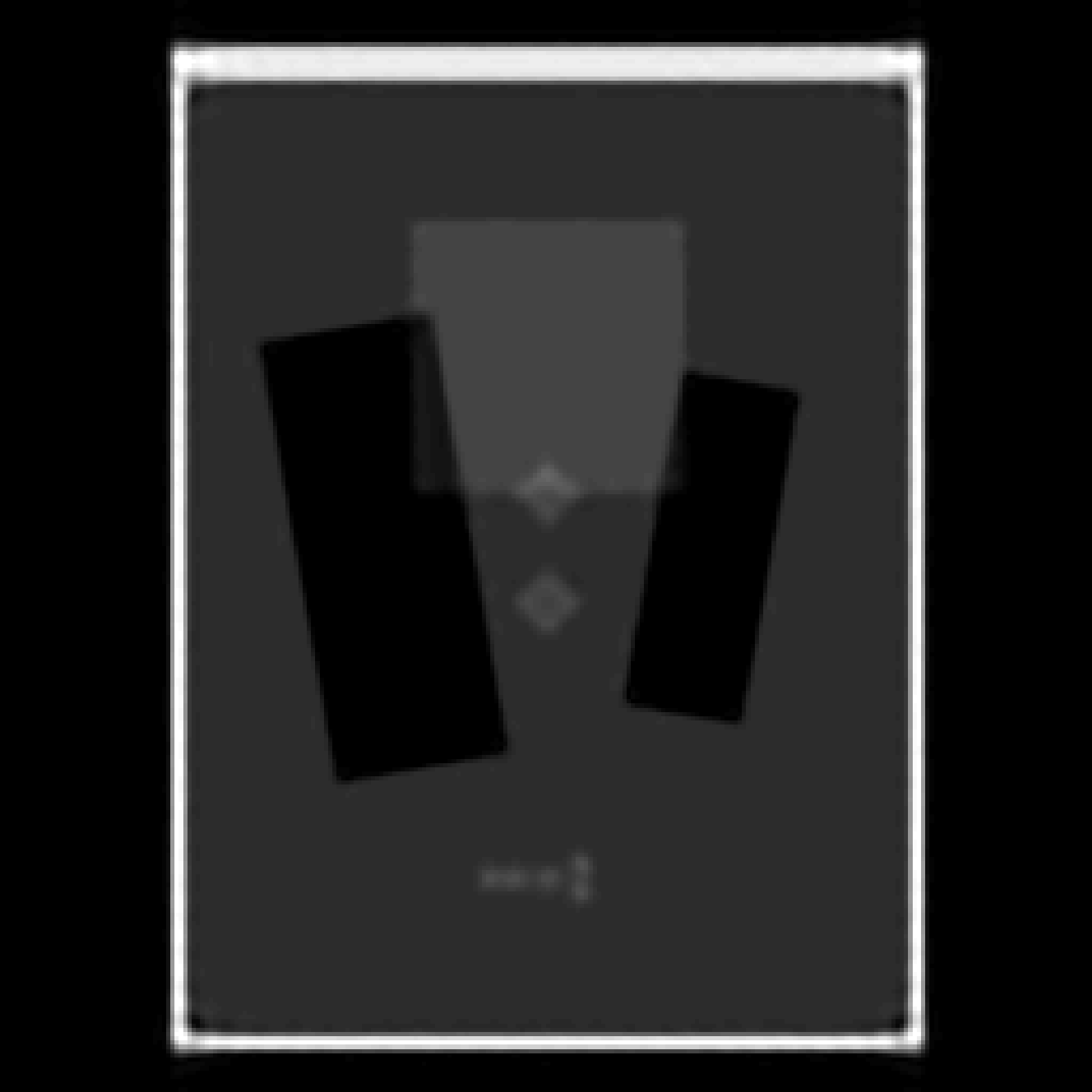}
\end{minipage}}\hspace{-0.05cm}
\subfloat[TV]{\label{RectangleILFTV}\begin{minipage}{3cm}
\includegraphics[width=3cm]{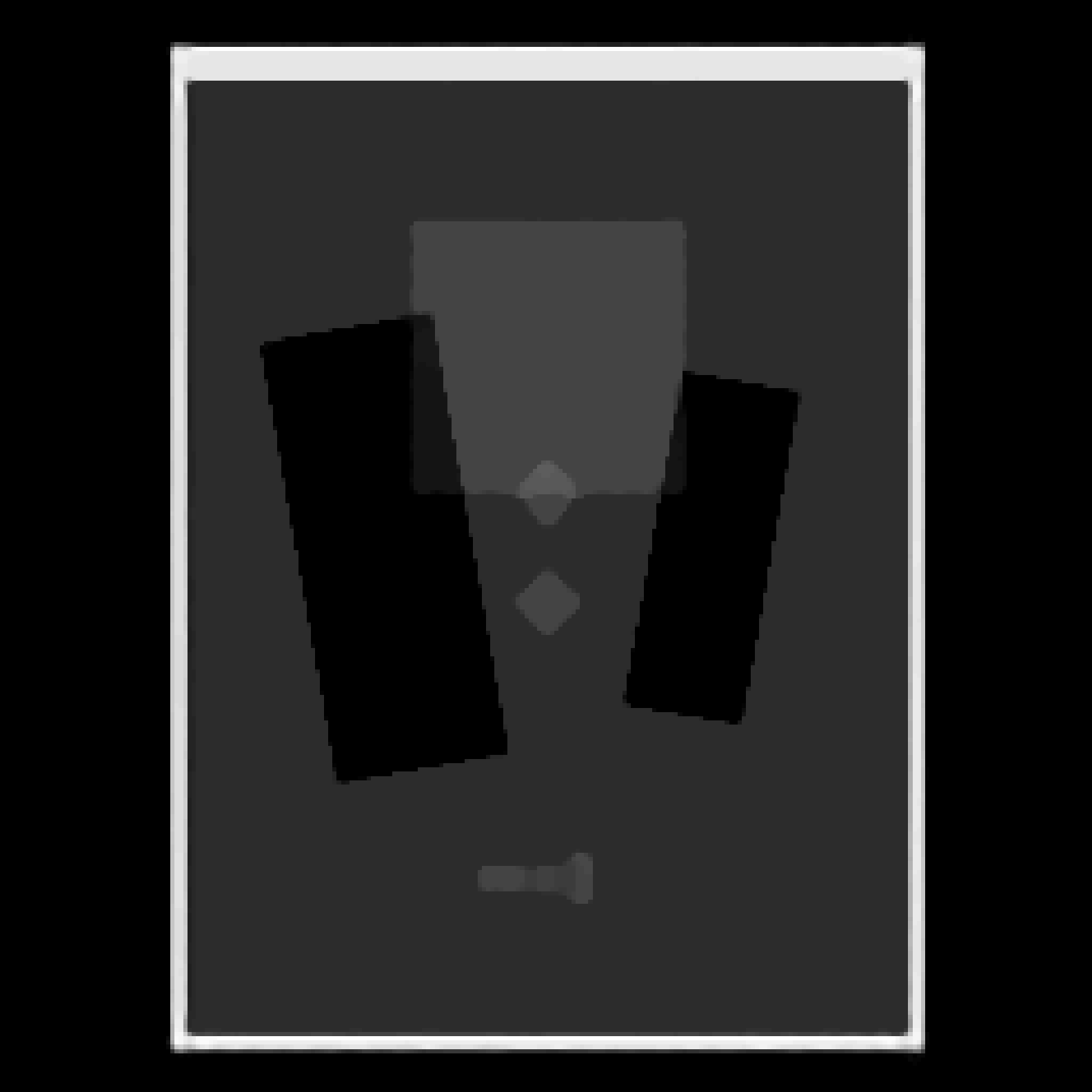}
\end{minipage}}\hspace{-0.05cm}
\subfloat[Haar]{\label{RectangleILFHaar}\begin{minipage}{3cm}
\includegraphics[width=3cm]{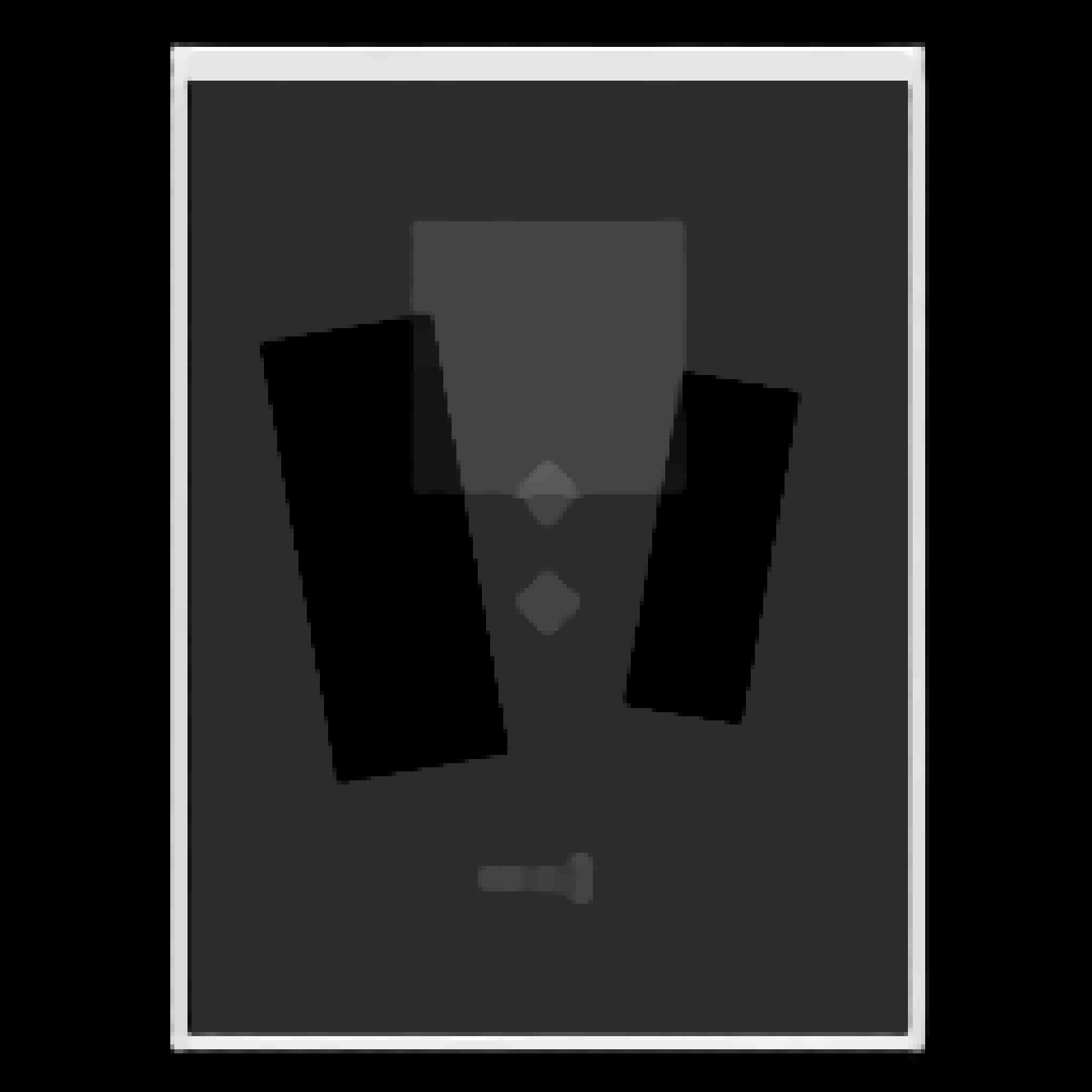}
\end{minipage}}\hspace{-0.05cm}
\subfloat[DDTF]{\label{RectangleILFDDTF}\begin{minipage}{3cm}
\includegraphics[width=3cm]{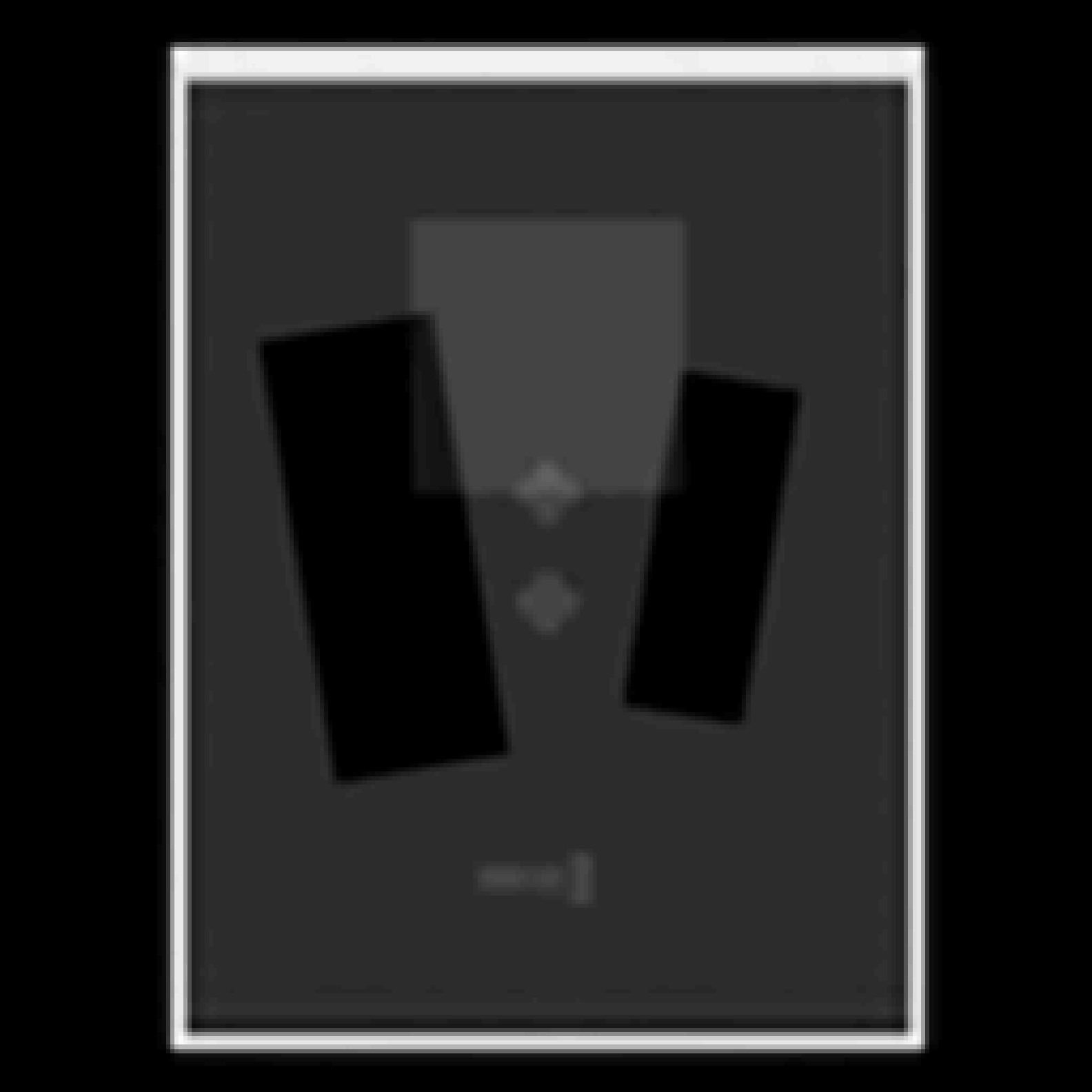}
\end{minipage}}
\caption{Visual comparison of ``Rectangle'' for ideal low-pass filter deconvolution.}\label{RectangleILFResults}
\end{figure}

\begin{figure}[t]
\centering
\subfloat[TF]{\label{RectangleILFTF}\begin{minipage}{3cm}
\includegraphics[width=3cm]{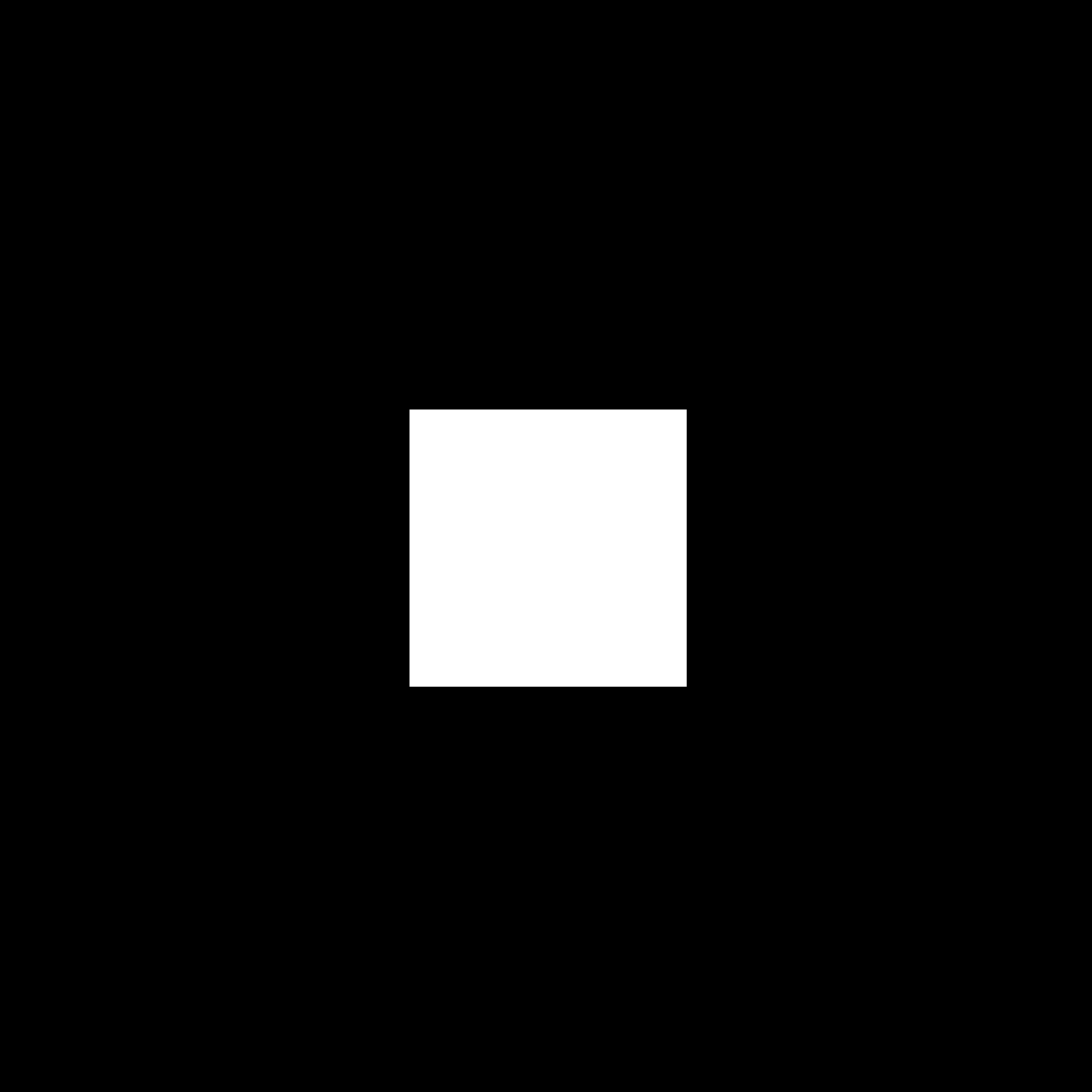}
\end{minipage}}\hspace{-0.05cm}
\subfloat[Proposed]{\label{RectangleILFProposedError}\begin{minipage}{3cm}
\includegraphics[width=3cm]{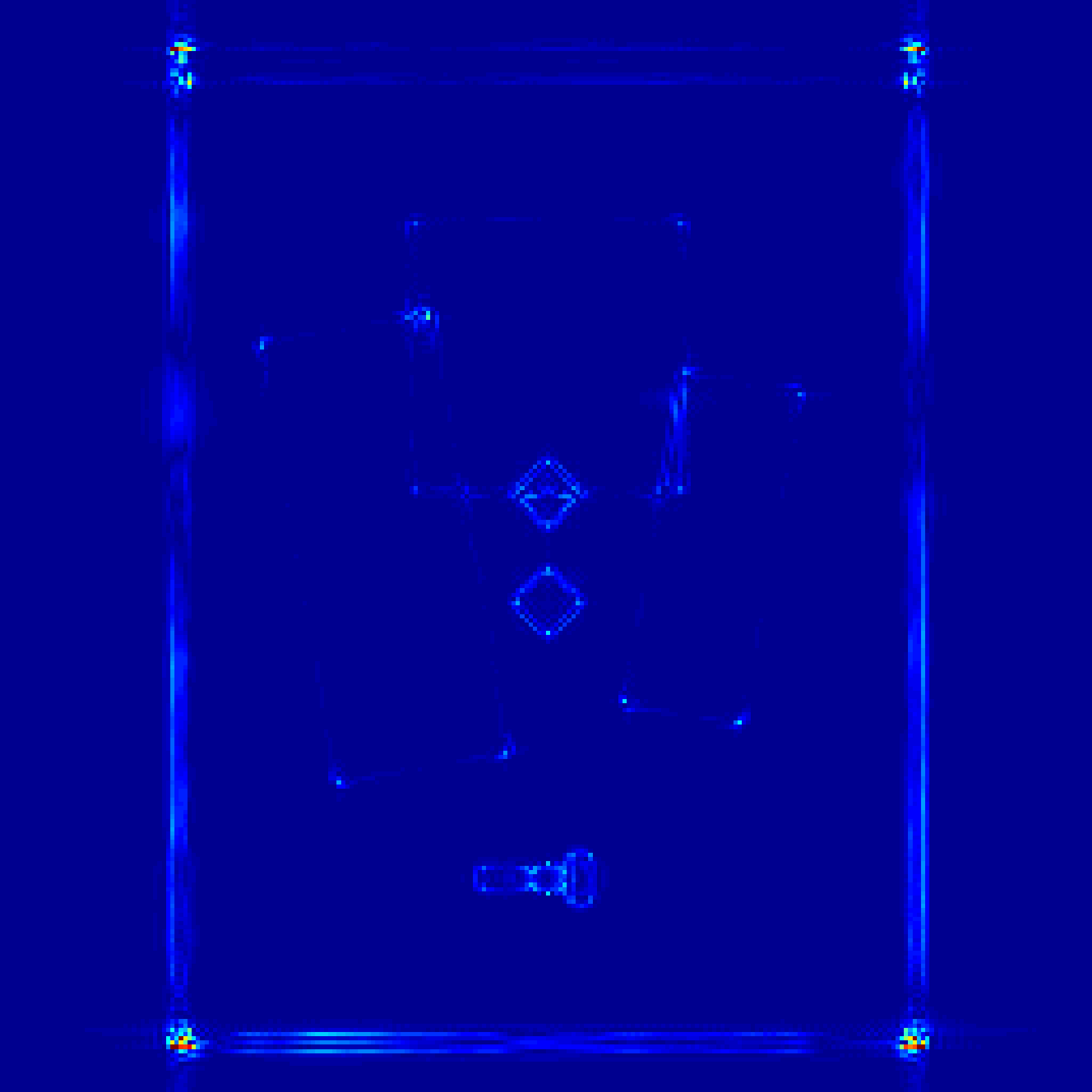}
\end{minipage}}\hspace{-0.05cm}
\subfloat[LSLP]{\label{RectangleILFLSLPError}\begin{minipage}{3cm}
\includegraphics[width=3cm]{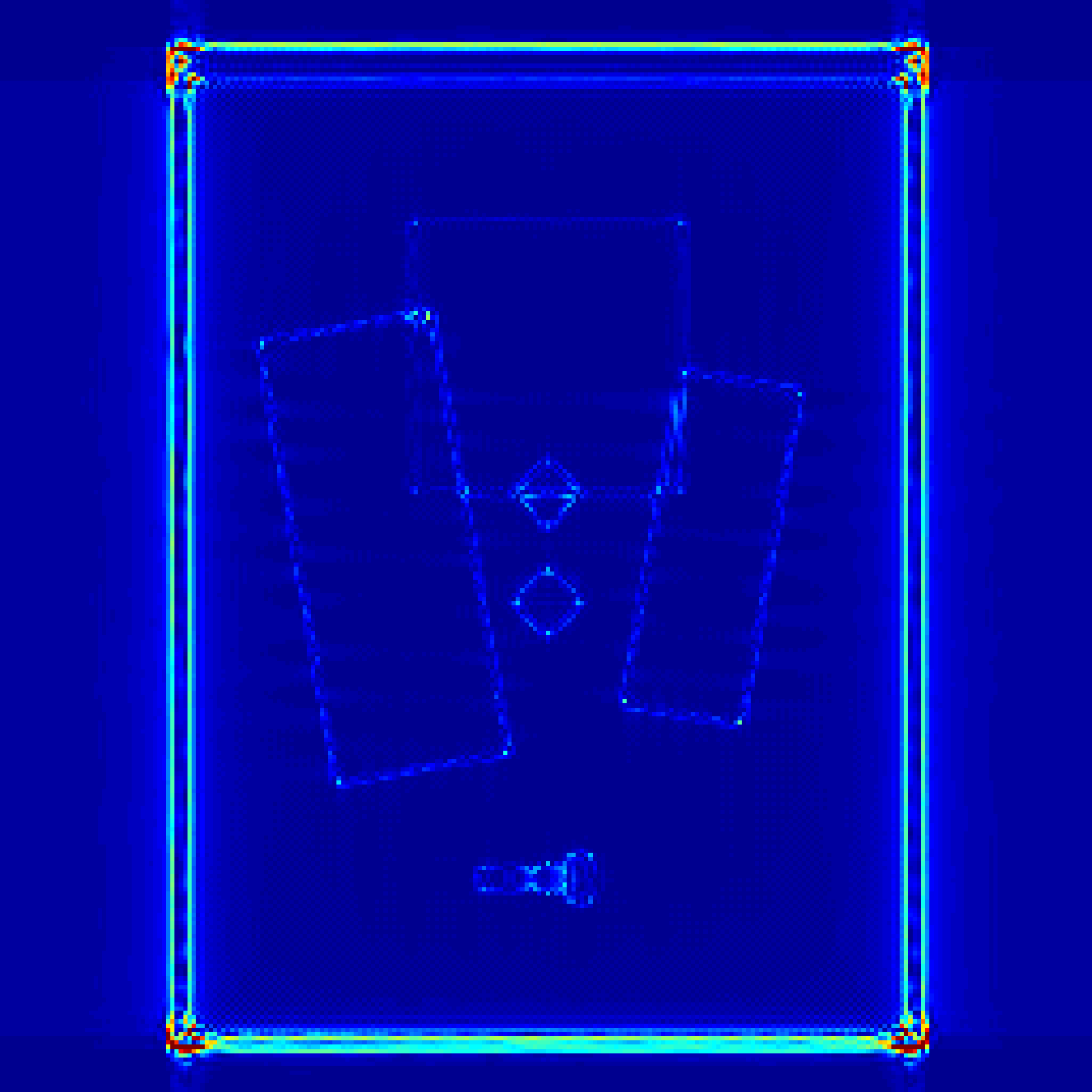}
\end{minipage}}\hspace{-0.05cm}
\subfloat[LRHDDTF]{\label{RectangleILFLRHDDTFError}\begin{minipage}{3cm}
\includegraphics[width=3cm]{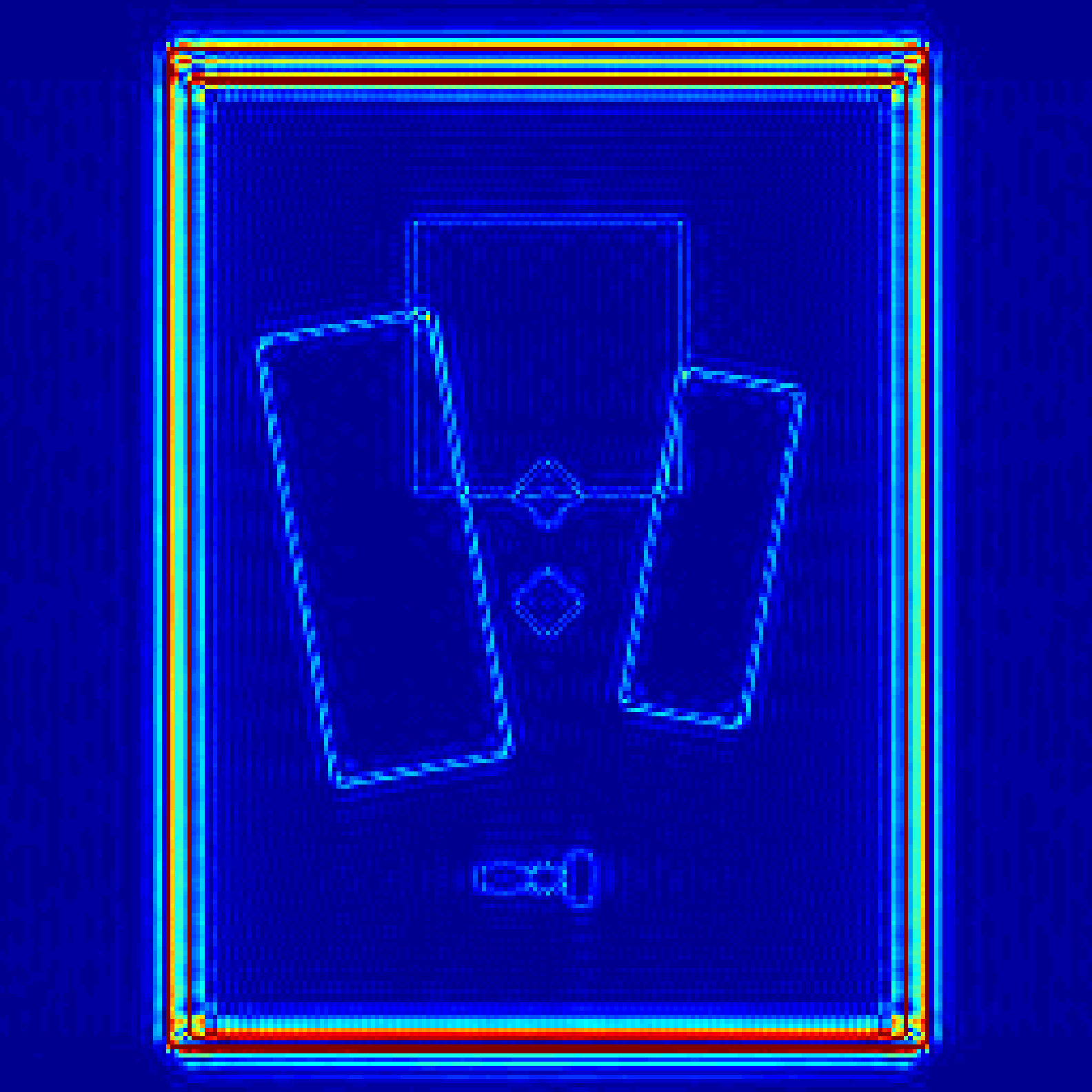}
\end{minipage}}\vspace{-0.25cm}\\
\subfloat[Schatten $0$]{\label{RectangleILFGIRAFError}\begin{minipage}{3cm}
\includegraphics[width=3cm]{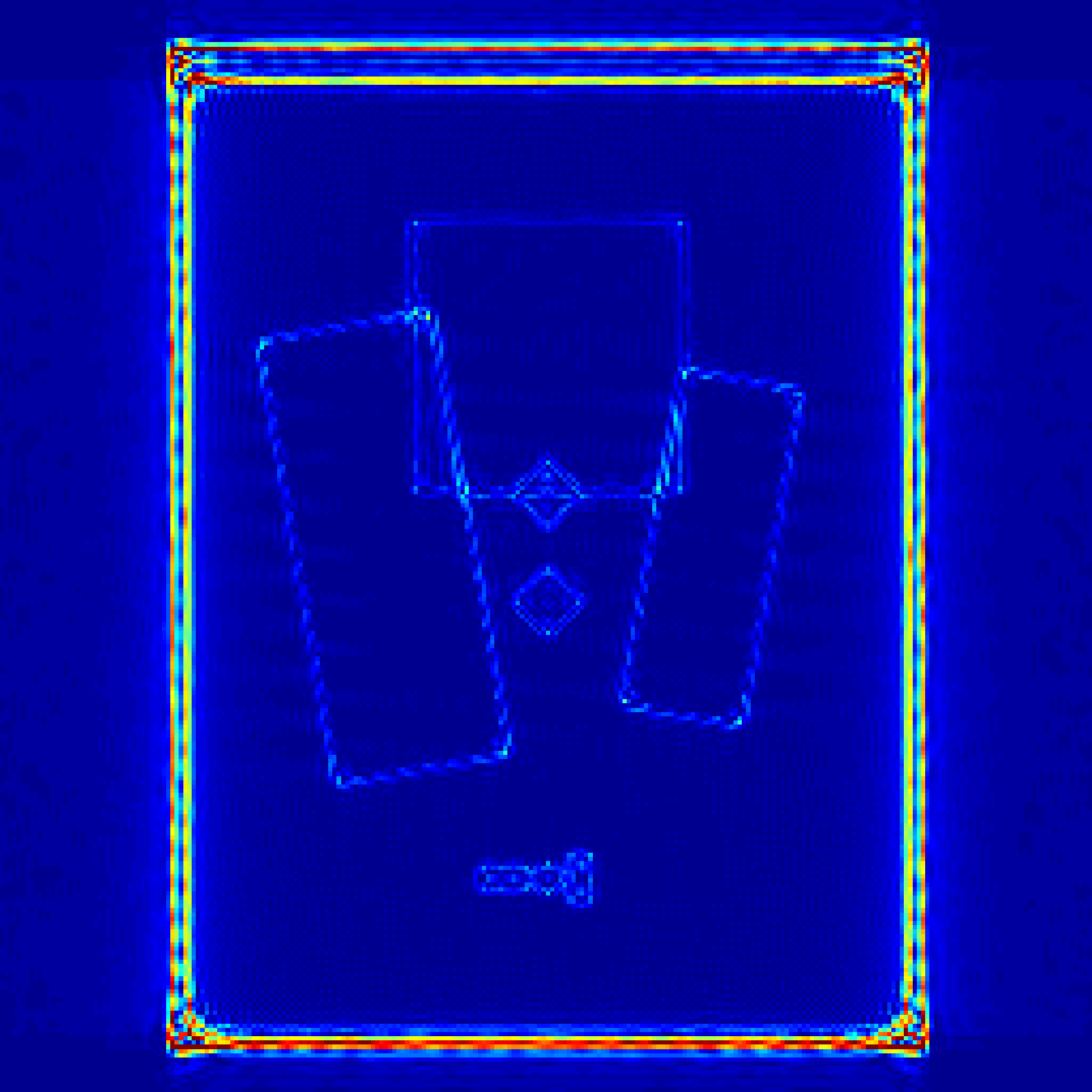}
\end{minipage}}\hspace{-0.05cm}
\subfloat[TV]{\label{RectangleILFTVError}\begin{minipage}{3cm}
\includegraphics[width=3cm]{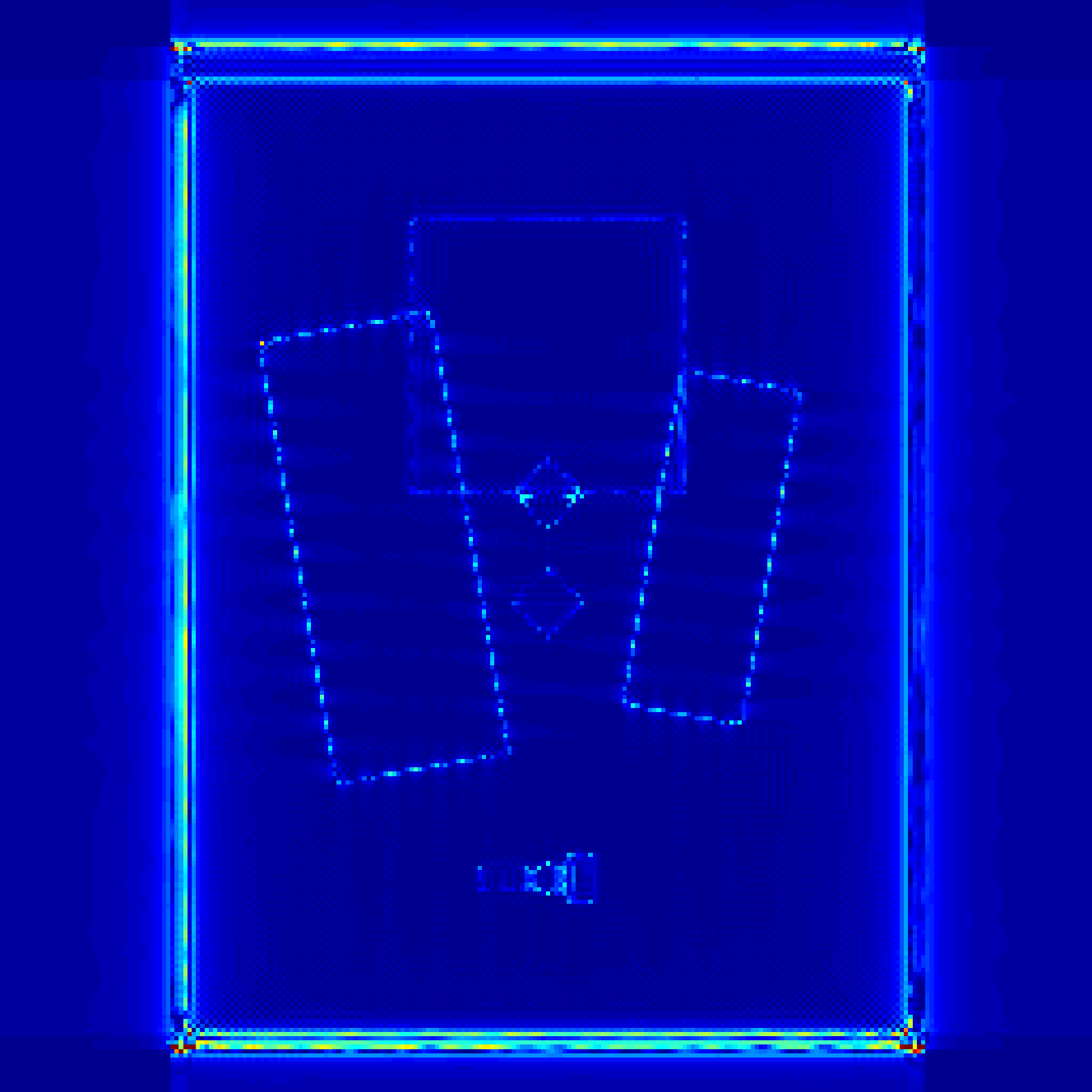}
\end{minipage}}\hspace{-0.05cm}
\subfloat[Haar]{\label{RectangleILFHaarError}\begin{minipage}{3cm}
\includegraphics[width=3cm]{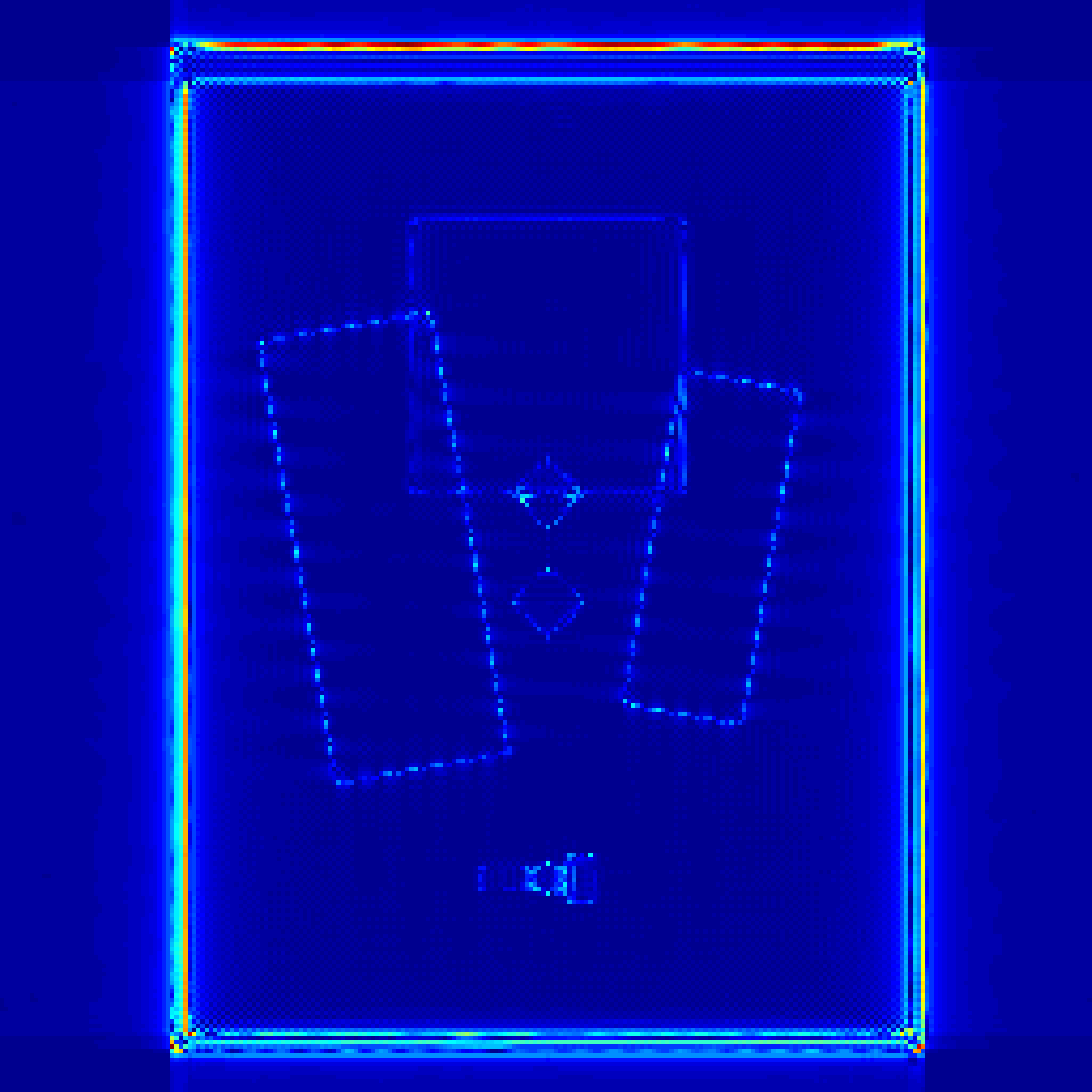}
\end{minipage}}\hspace{-0.05cm}
\subfloat[DDTF]{\label{RectangleILFDDTFError}\begin{minipage}{3cm}
\includegraphics[width=3cm]{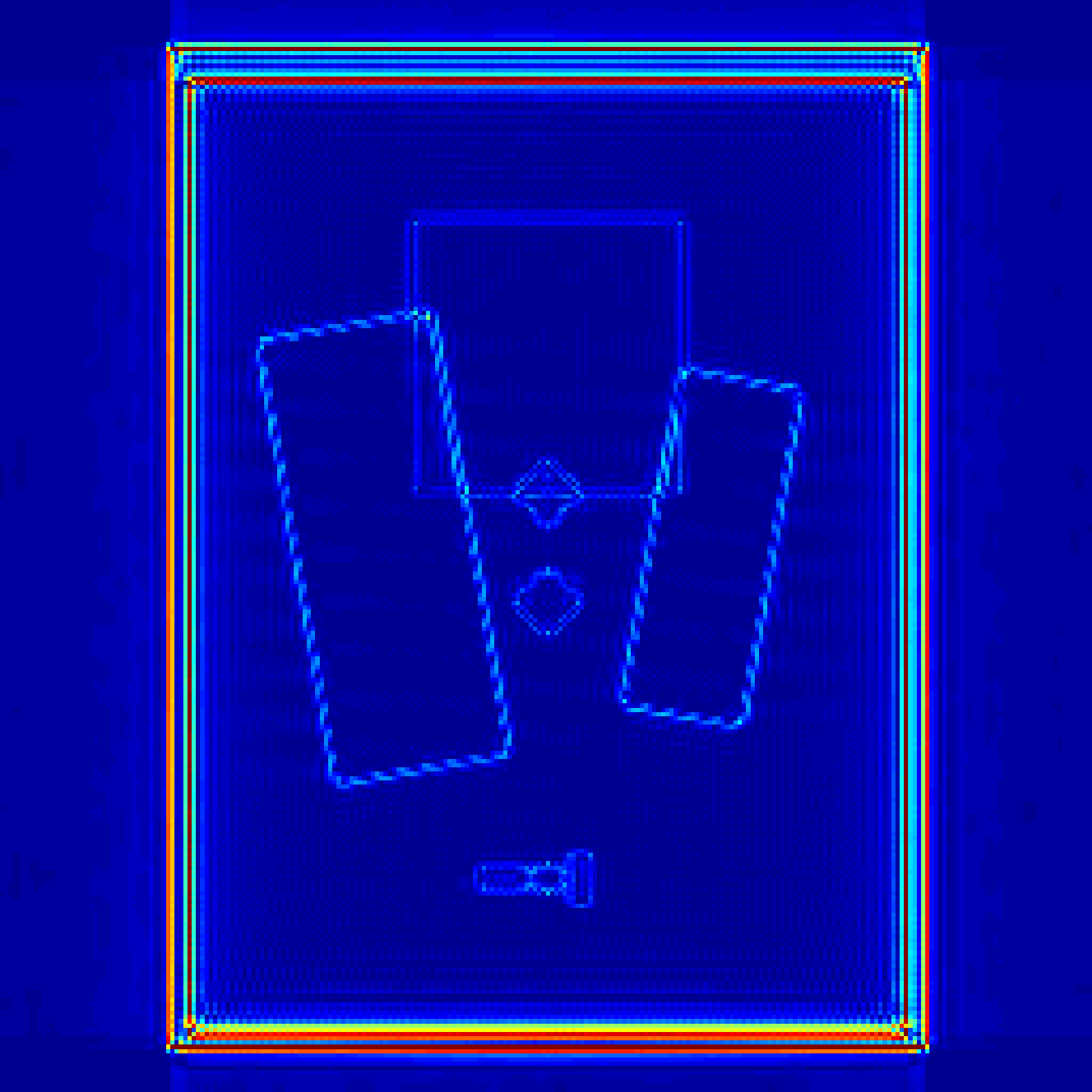}
\end{minipage}}
\caption{Error maps of \cref{RectangleILFResults}. In \cref{RectangleILFTF}, TF stands for the transfer function.}\label{RectangleILFError}
\end{figure}

The SNR, the HFEN, and the SSIM of the aforementioned approaches for the synthetic images are summarized in \cref{PhantomResults}, and \cref{EllipseRandomResults,EllipseILFResults,RectangleRandomResults,RectangleILFResults} display the visual comparisons with the corresponding error maps in \cref{EllipseRandomError,EllipseILFError,RectangleRandomError,RectangleILFError}. Throughout this section, all restored images, including the zoom-in views, are displayed in the window level $[0,1]$, and all error maps are displayed in the window level $[0,0.2]$ for fair comparisons. We can see that, our proposed model \cref{ProposedModel2} consistently outperforms other image restoration models in every index with visually observable improvements. These improvements demonstrate that the proposed approach is suitable for the piecewise constant image restoration.

Similar to the previous works \cite{J.F.Cai2022,J.F.Cai2020}, our proposed analysis approach \cref{ProposedModel2} is also based on the structured low-rank matrix framework for the continuous domain regularization. Hence, the similar arguments are still available in our case; the continuous domain regularization is able to reduce the basis mismatch between the true edge set in the continuum and the discrete grid, leading to the improvements in both indices and visual qualities. Indeed, even though the on-the-grid sparse regularization models \cref{TV,Frame} can manages to show improvements with the aid of the pre-estimated edge weights, it can be observed in the error maps \cref{EllipseRandomError,EllipseILFError,RectangleRandomError,RectangleILFError} that the edges are in general restricted to the discrete grids, leading to the performance bottleneck compared to our proposed model \cref{ProposedModel2}.

Most importantly, the results demonstrate that among the continuous domain regularizations \cref{ProposedModel2,LSLP2,LRHDDTF2,Schatten0}, our proposed model \cref{ProposedModel2} performs best. First of all, the LSLP model \cref{LSLP2} shows some artifacts near the edges. Notice that the operator $\bmD$ corresponding to the derivatives in the image domain amplifies the noise in the high frequencies. Since the LSLP model penalizes the canonical coefficients corresponding to the annihilating filters only, it is likely that such an amplified noise in the high frequencies may not be fully suppressed by \cref{LSLP2}. In contrast, since the proposed model promotes the sparsity on the entire range of $\bmW$, we can achieve a better denoising effect in spite of the amplified noise, leading to the better restoration results with less artifacts near the edges.

In addition, comparing \cref{ProposedModel2,LRHDDTF2}, we can see that the analysis approach \cref{ProposedModel2} consistently outperforms the balanced approach \cref{LRHDDTF2}. The results again verify the well-known literature; the analysis approach can reflect the structure of the target signal/image better than the synthesis/balanced approach (e.g. \cite{J.F.Cai2012}). Notice that \cref{TightFrameSparse} tells us that $\bmD\bsv$ is sparse under the right singular vectors of $\bmH\left(\bmD\bsv\right)$. This means that, the analysis approach \cref{ProposedModel2} is able to reflect the low-rank structure of $\bmH\left(\bmD\bsv\right)$ better than the balanced approach \cref{LRHDDTF2}, as the numerical results in \cref{EdgeEstimation} show that \cref{DDTFApproach} is able to obtain a reasonable approximation of right singular vectors.

Finally, the results also demonstrate that \cref{ProposedModel2} are less sensitive to the image restoration tasks compared to other continuous domain regularization approaches \cref{LSLP2,LRHDDTF2,Schatten0}. In particular, we can see that the LRHDDTF model \cref{LRHDDTF2} and the Schatten $0$ model \cref{Schatten0} perform relatively well in the case of random sampling where there are high frequency samples available in $\bsf$, as already illustrated in \cite{J.F.Cai2020} for the CS-MRI restoration. However, when the measurements $\bsf$ contain few/even no high frequency samples (ideal low-pass filter deconvolution), these two models \cref{LRHDDTF2,Schatten0} perform at most comparably to the discrete sparse regularization models \cref{TV,Frame}. In summary, our proposed analysis approach \cref{ProposedModel2} shows the overall better restoration quality in both the indices and the image quality, in addition to the less sensitiveness to the image restoration tasks.

\subsection{Image restoration results: non-synthetic images}

\begin{table}[t]
\centering
\scriptsize
\begin{tabular}{|c|c|c|c|c|c|c|c|c|c|}
\hline
\multicolumn{10}{|c|}{Random sampling}\\ \hline
\multirow{2}{*}{Image}&\multirow{2}{*}{Index}&\multirow{2}{*}{IFFT}&\multirow{2}{*}{Proposed \cref{ProposedModel2}}&\multirow{2}{*}{LSLP \cref{LSLP2}}&\multirow{2}{*}{LRHDDTF \cref{LRHDDTF2}}&\multirow{2}{*}{Schatten $0$ \cref{Schatten0}}&\multirow{2}{*}{TV \cref{TV}}&\multicolumn{2}{|c|}{Frame \cref{Frame}}\\ \cline{9-10}
&&&&&&&&Haar&DDTF\\ \hline
\multirow{3}{*}{AB}&SNR&$18.29$&$\textbf{31.95}$&$29.98$&$27.99$&$27.87$&$30.02$&$30.42$&$29.32$\\ \cline{2-10}
&HFEN&$0.4710$&$\textbf{0.0428}$&$0.0559$&$0.0877$&$0.0906$&$0.06172$&$0.0599$&$0.0571$\\ \cline{2-10}
&SSIM&$0.6265$&$\textbf{0.9831}$&$0.9758$&$0.9592$&$0.9576$&$0.9670$&$0.9689$&$0.9675$\\ \hline
\multirow{3}{*}{SM}&SNR&$15.34$&$\textbf{30.78}$&$27.36$&$27.00$&$24.57$&$28.74$&$29.01$&$26.20$\\ \cline{2-10}
&HFEN&$0.4461$&$\textbf{0.0196}$&$0.0302$&$0.0446$&$0.0372$&$0.0402$&$0.0403$&$0.0322$\\ \cline{2-10}
&SSIM&$0.4126$&$\textbf{0.9846}$&$0.9455$&$0.9495$&$0.9526$&$0.8809$&$0.9244$&$0.8963$\\ \hline
\multirow{3}{*}{MC}&SNR&$23.09$&$\textbf{33.56}$&$31.81$&$32.39$&$32.30$&$30.58$&$30.42$&$32.34$\\ \cline{2-10}
&HFEN&$0.4478$&$0.1235$&$0.1454$&$0.1345$&$0.1258$&$0.1699$&$0.1780$&$\textbf{0.1057}$\\ \cline{2-10}
&SSIM&$0.7871$&$\textbf{0.9557}$&$0.9431$&$0.9505$&$0.9488$&$0.9259$&$0.9226$&$0.9541$\\ \hline\hline
\multicolumn{10}{|c|}{Ideal low-pass filter}\\ \hline
\multirow{2}{*}{Image}&\multirow{2}{*}{Index}&\multirow{2}{*}{IFFT}&\multirow{2}{*}{Proposed \cref{ProposedModel2}}&\multirow{2}{*}{LSLP \cref{LSLP2}}&\multirow{2}{*}{LRHDDTF \cref{LRHDDTF2}}&\multirow{2}{*}{Schatten $0$ \cref{Schatten0}}&\multirow{2}{*}{TV \cref{TV}}&\multicolumn{2}{|c|}{Frame \cref{Frame}}\\ \cline{9-10}
&&&&&&&&Haar&DDTF\\ \hline
\multirow{3}{*}{AB}&SNR&$21.62$&$\textbf{28.88}$&$27.69$&$22.05$&$26.21$&$26.65$&$26.47$&$23.42$\\ \cline{2-10}
&HFEN&$0.2002$&$\textbf{0.0551}$&$0.0675$&$0.1769$&$0.0850$&$0.0677$&$0.0785$&$0.1109$\\ \cline{2-10}
&SSIM&$0.8472$&$\textbf{0.9788}$&$0.9726$&$0.8840$&$0.9596$&$0.9601$&$0.9587$&$0.9351$\\ \hline
\multirow{3}{*}{SM}&SNR&$19.60$&$\textbf{29.08}$&$26.13$&$20.06$&$24.17$&$25.52$&$25.86$&$20.77$\\ \cline{2-10}
&HFEN&$0.1910$&$\textbf{0.0453}$&$0.0683$&$0.1699$&$0.0892$&$0.0579$&$0.0610$&$0.1139$\\ \cline{2-10}
&SSIM&$0.8020$&$\textbf{0.9901}$&$0.9752$&$0.8538$&$0.9585$&$0.9587$&$0.9587$&$0.9207$\\ \hline
\multirow{3}{*}{MC}&SNR&$26.14$&$\textbf{32.35}$&$31.10$&$26.61$&$31.07$&$29.59$&$29.98$&$29.01$\\ \cline{2-10}
&HFEN&$0.2299$&$\textbf{0.0925}$&$0.1074$&$0.2041$&$0.0985$&$0.1364$&$0.1256$&$0.1047$\\ \cline{2-10}
&SSIM&$0.8820$&$\textbf{0.9576}$&$0.9483$&$0.8968$&$0.9482$&$0.9307$&$0.9367$&$0.9373$\\ \hline
\end{tabular}
\caption{Comparison of SNR, HFEN, and SSIM for non-synthetic images}\label{NaturalImageResults}
\end{table}

\begin{figure}[t]
\centering
\subfloat[Ref.]{\label{AngryBirdsRandomOriginal}\begin{minipage}{3cm}
\includegraphics[width=3cm]{AngryBirdsOriginal.pdf}
\end{minipage}}\hspace{-0.05cm}
\subfloat[Proposed]{\label{AngryBirdsRandomProposed}\begin{minipage}{3cm}
\includegraphics[width=3cm]{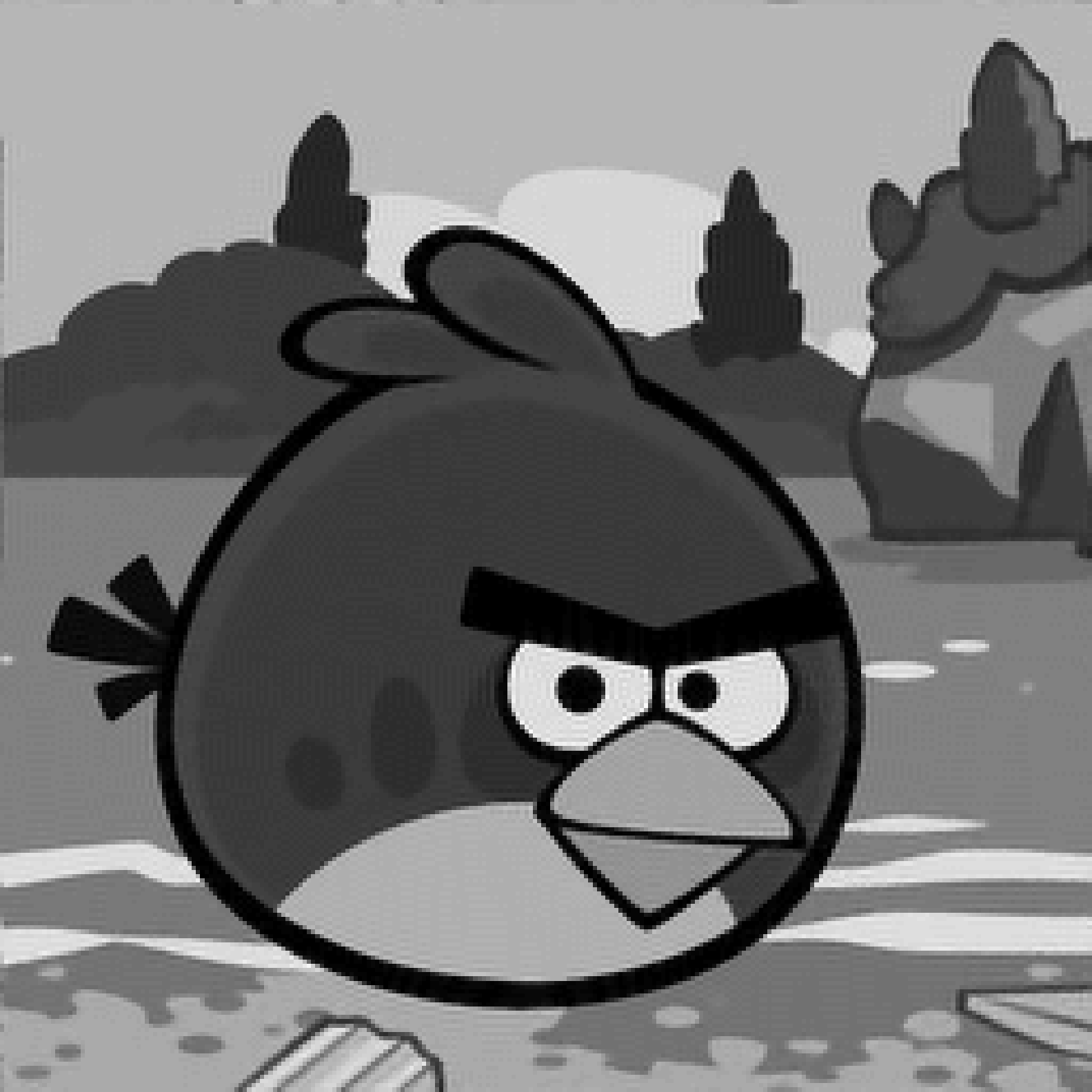}
\end{minipage}}\hspace{-0.05cm}
\subfloat[LSLP]{\label{AngryBirdsRandomLSLP}\begin{minipage}{3cm}
\includegraphics[width=3cm]{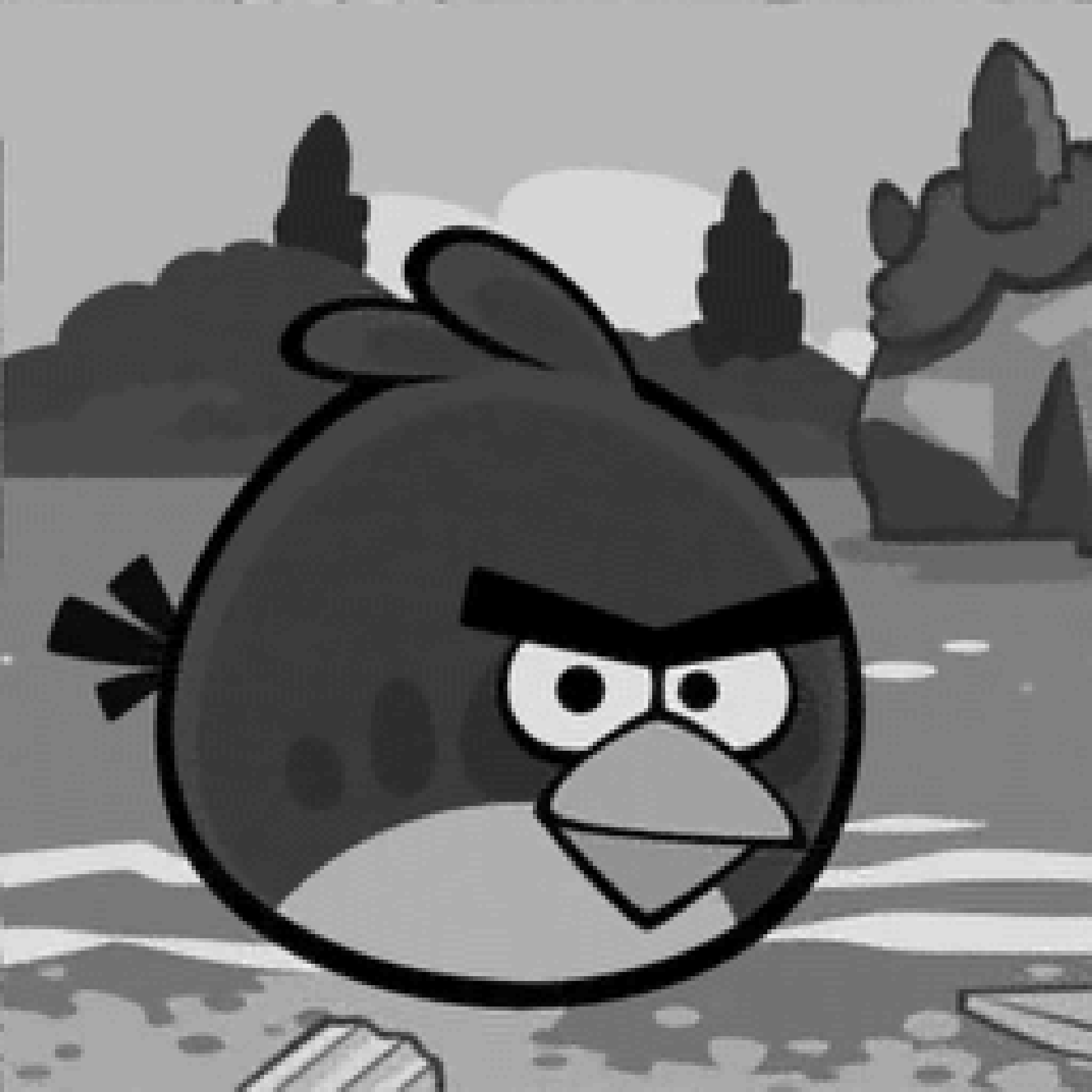}
\end{minipage}}\hspace{-0.05cm}
\subfloat[LRHDDTF]{\label{AngryBirdsRandomLRHDDTF}\begin{minipage}{3cm}
\includegraphics[width=3cm]{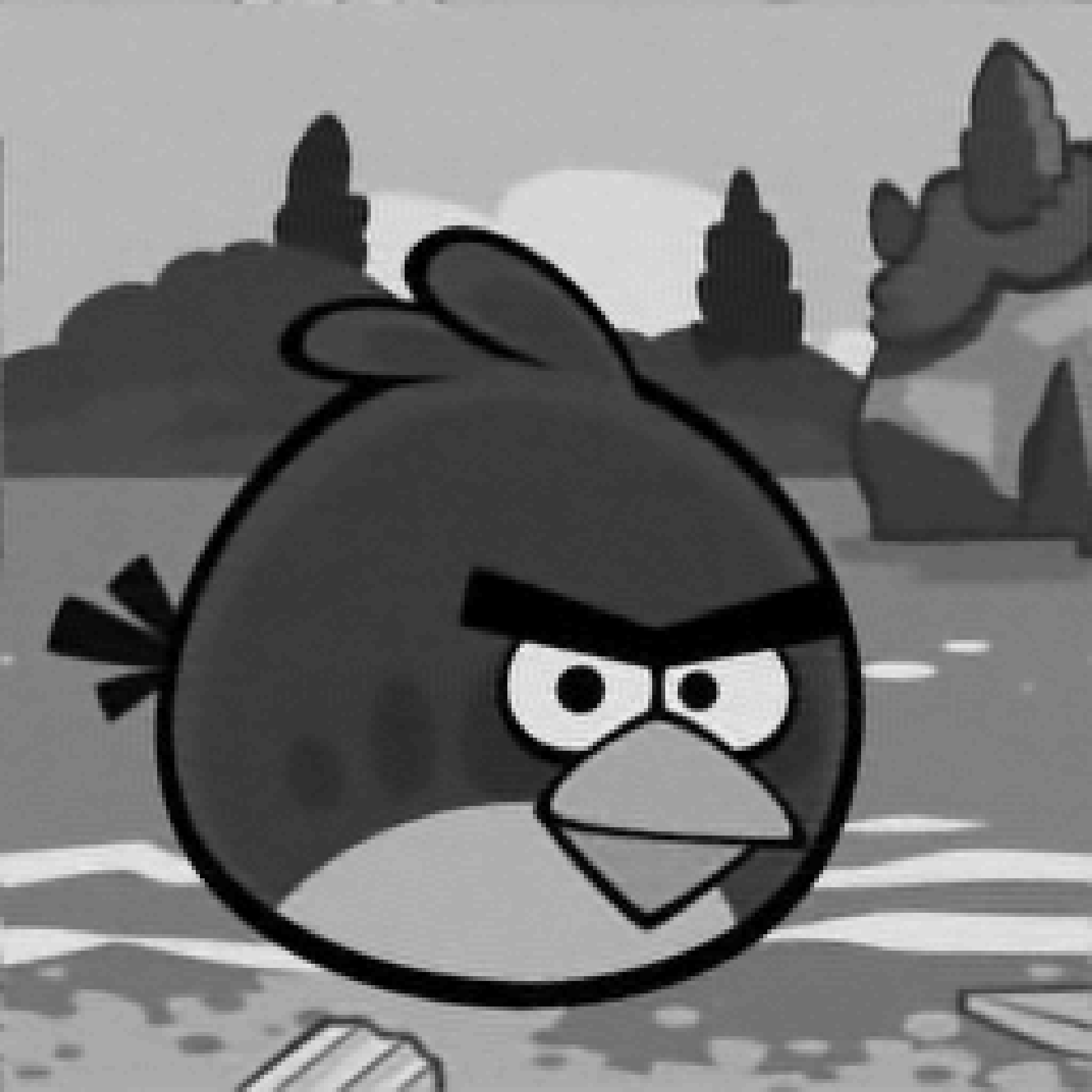}
\end{minipage}}\vspace{-0.25cm}
\subfloat[Schatten $0$]{\label{AngryBirdsRandomGIRAF}\begin{minipage}{3cm}
\includegraphics[width=3cm]{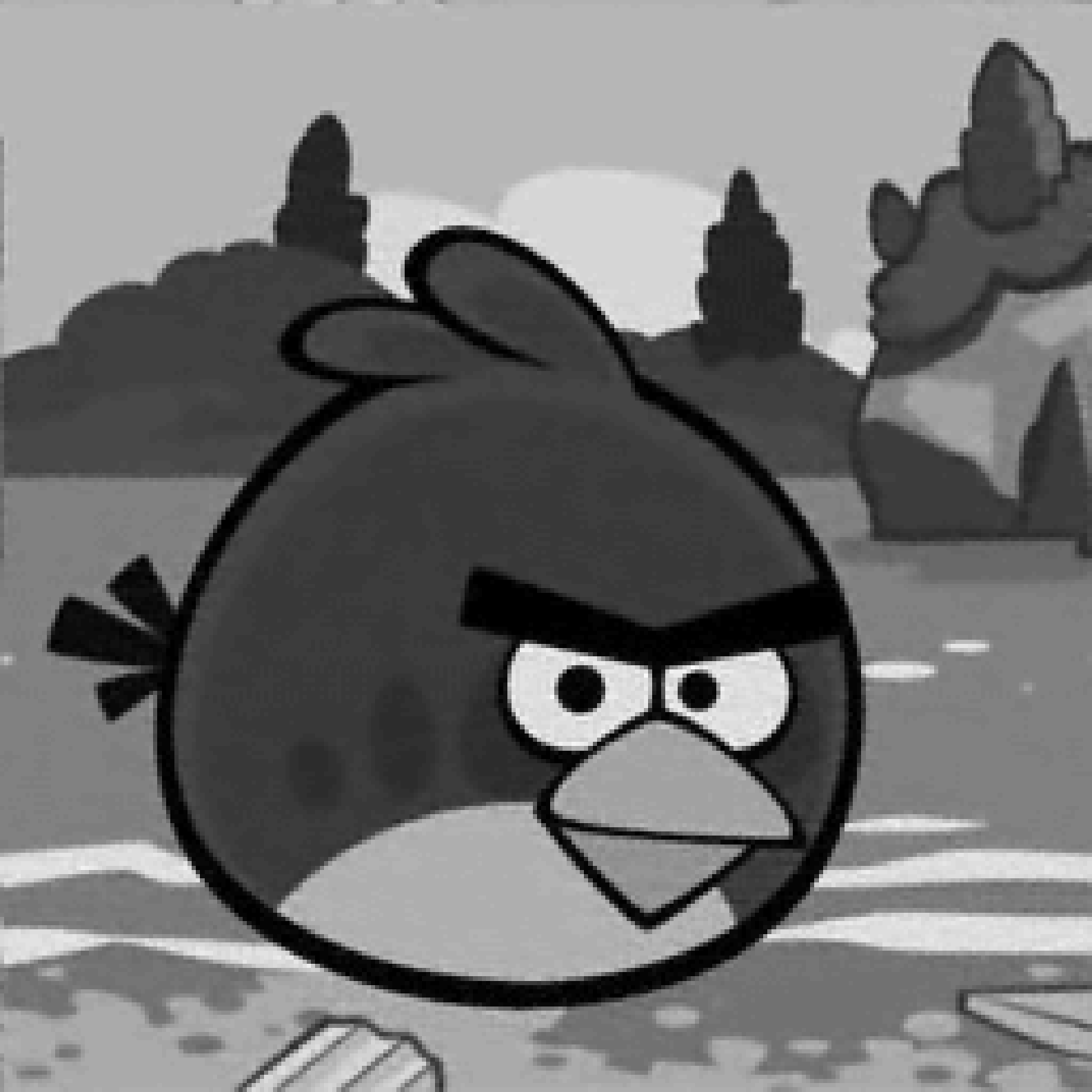}
\end{minipage}}\hspace{-0.05cm}
\subfloat[TV]{\label{AngryBirdsRandomTV}\begin{minipage}{3cm}
\includegraphics[width=3cm]{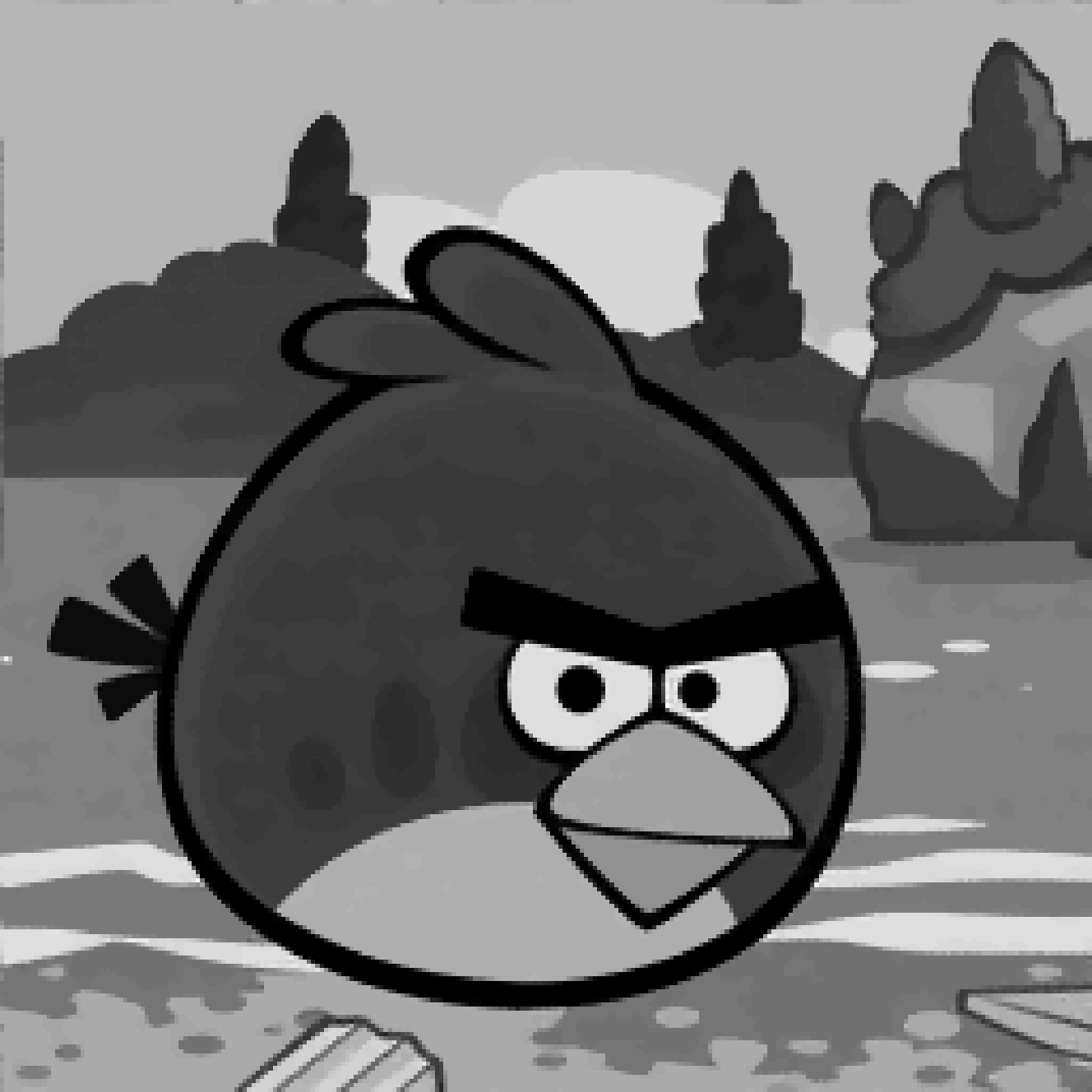}
\end{minipage}}\hspace{-0.05cm}
\subfloat[Haar]{\label{AngryBirdsRandomHaar}\begin{minipage}{3cm}
\includegraphics[width=3cm]{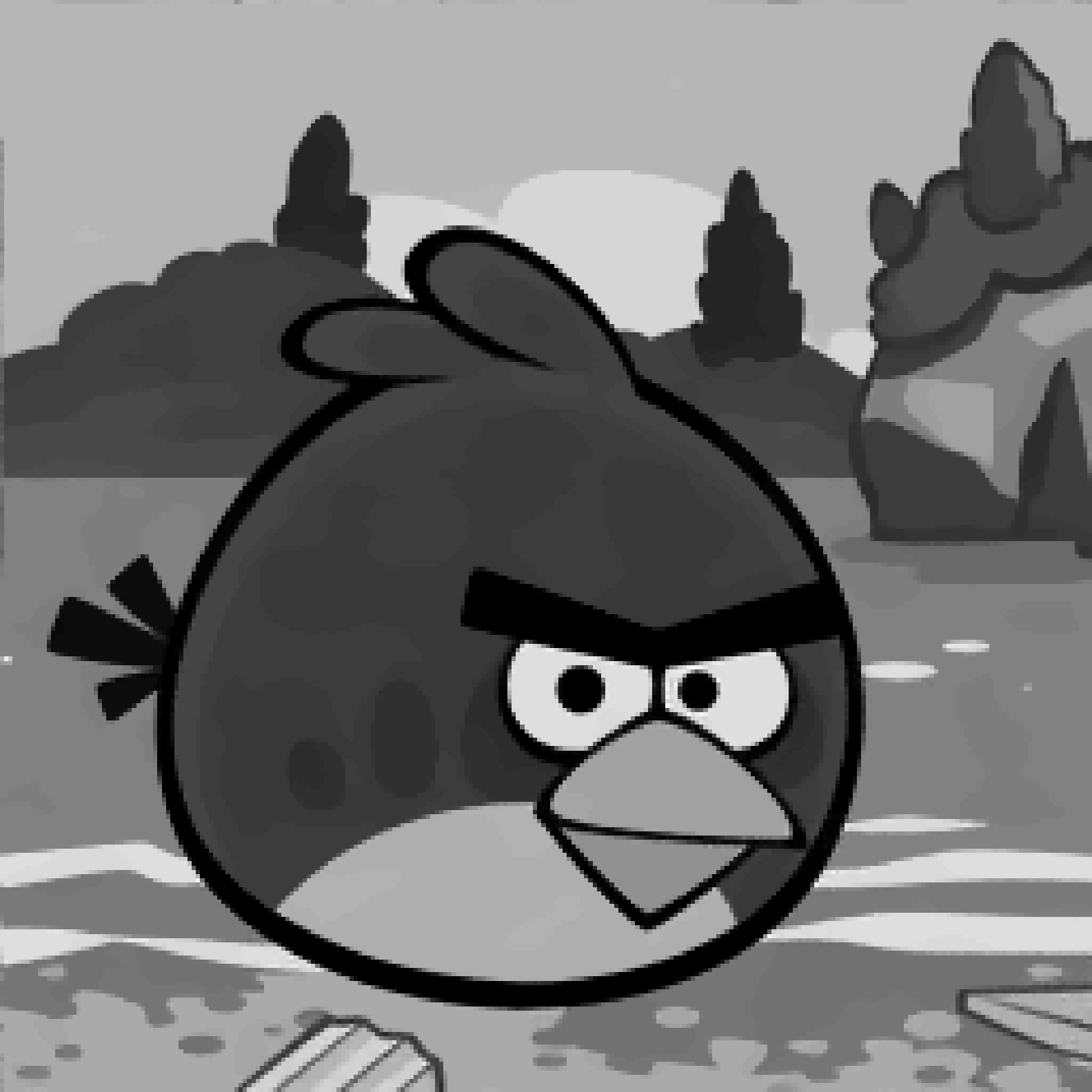}
\end{minipage}}\hspace{-0.05cm}
\subfloat[DDTF]{\label{AngryBirdsRandomDDTF}\begin{minipage}{3cm}
\includegraphics[width=3cm]{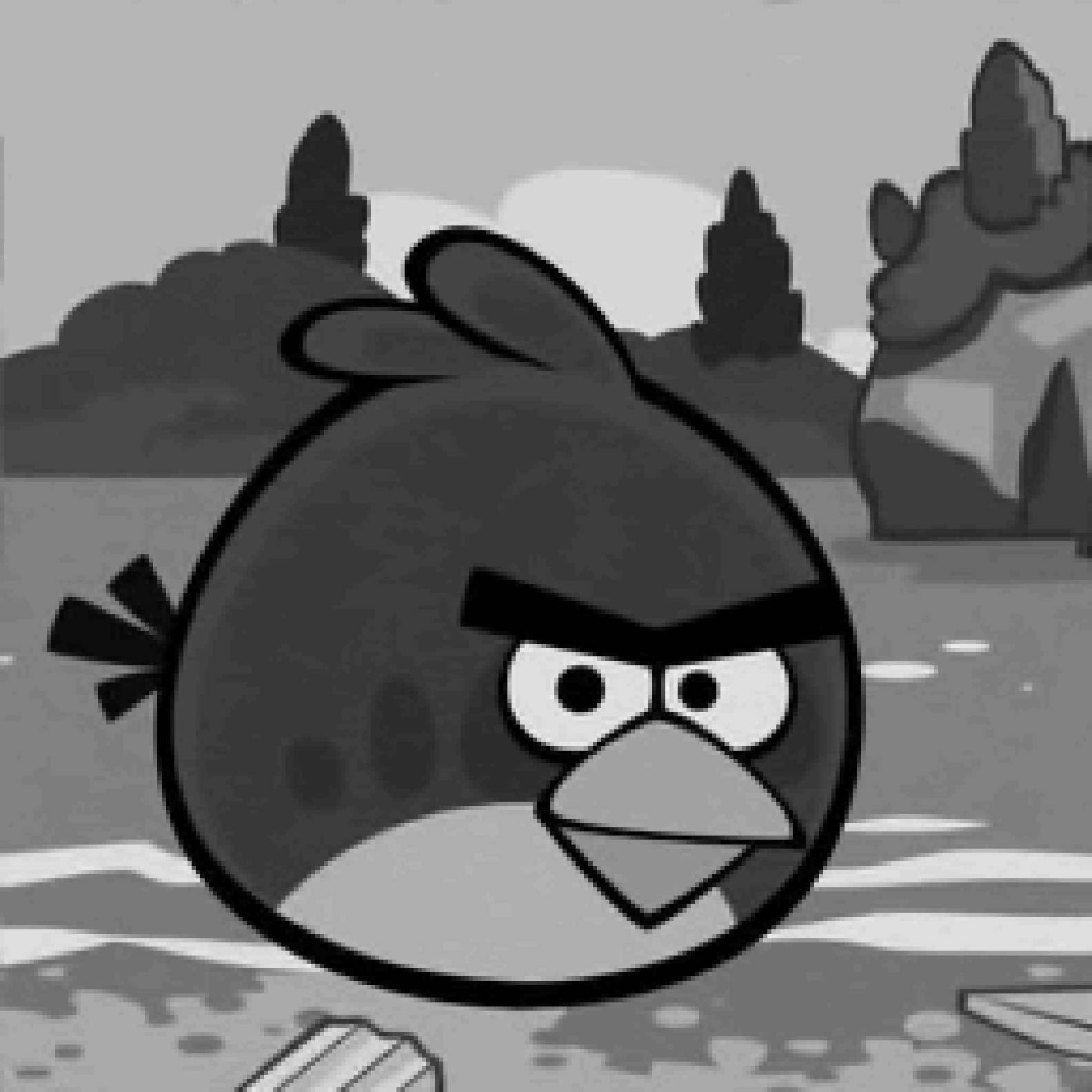}
\end{minipage}}
\caption{Visual comparison of ``Angry Birds'' for random sampling.}\label{AngryBirdsRandomResults}
\end{figure}

\begin{figure}[t]
\centering
\subfloat[Samples]{\label{AngryBirdsRandomSample}\begin{minipage}{3cm}
\includegraphics[width=3cm]{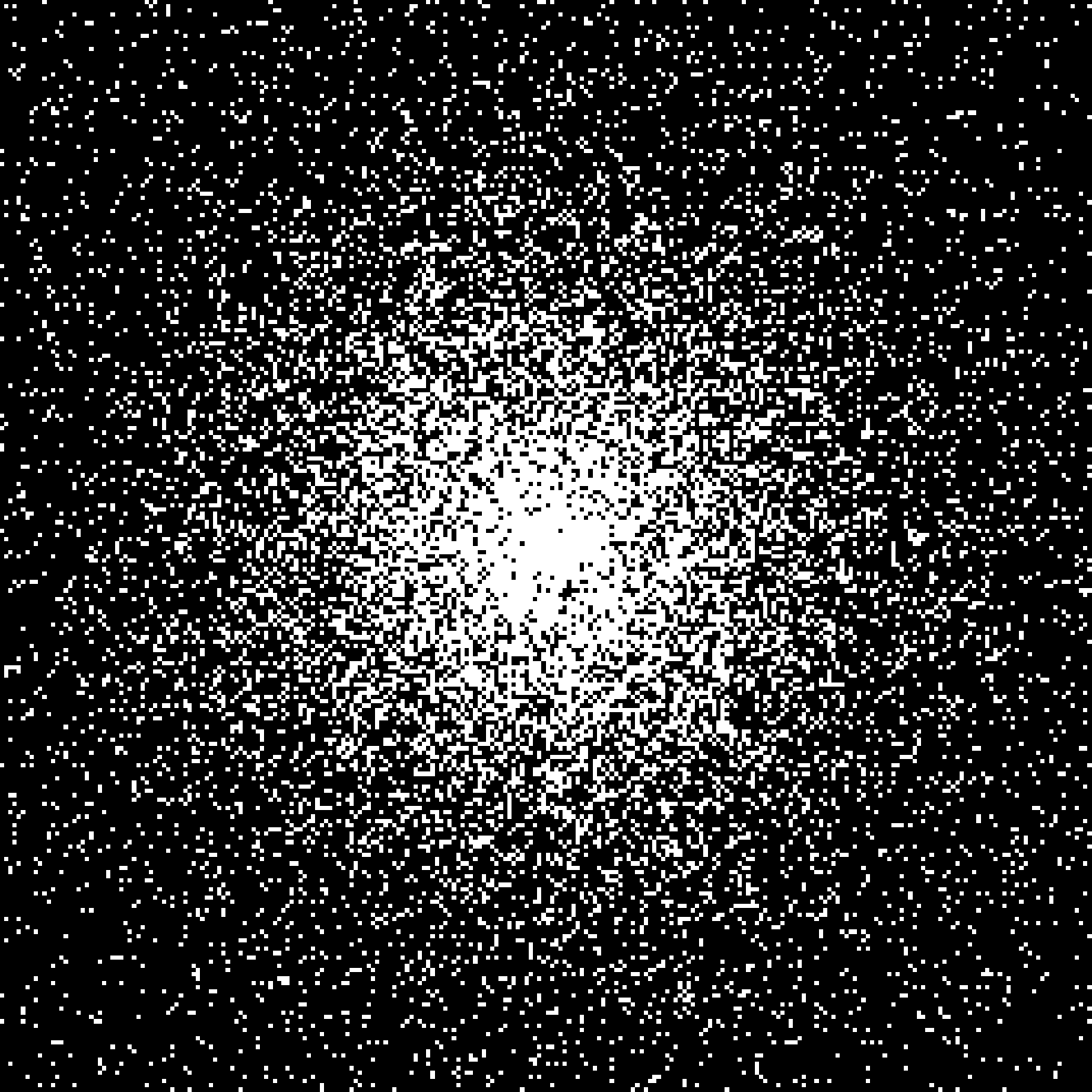}
\end{minipage}}\hspace{-0.05cm}
\subfloat[Proposed]{\label{AngryBirdsRandomProposedError}\begin{minipage}{3cm}
\includegraphics[width=3cm]{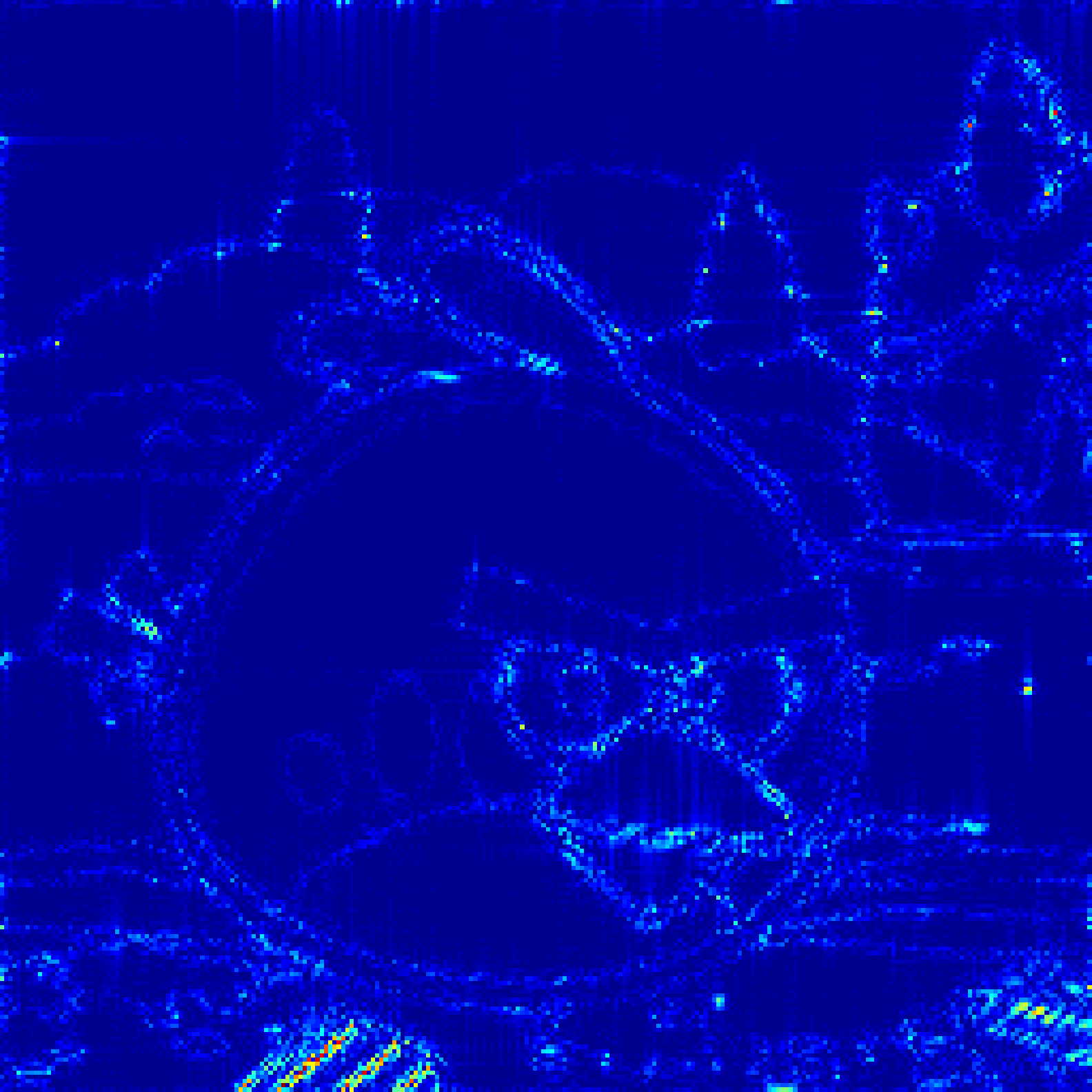}
\end{minipage}}\hspace{-0.05cm}
\subfloat[LSLP]{\label{AngryBirdsRandomLSLPError}\begin{minipage}{3cm}
\includegraphics[width=3cm]{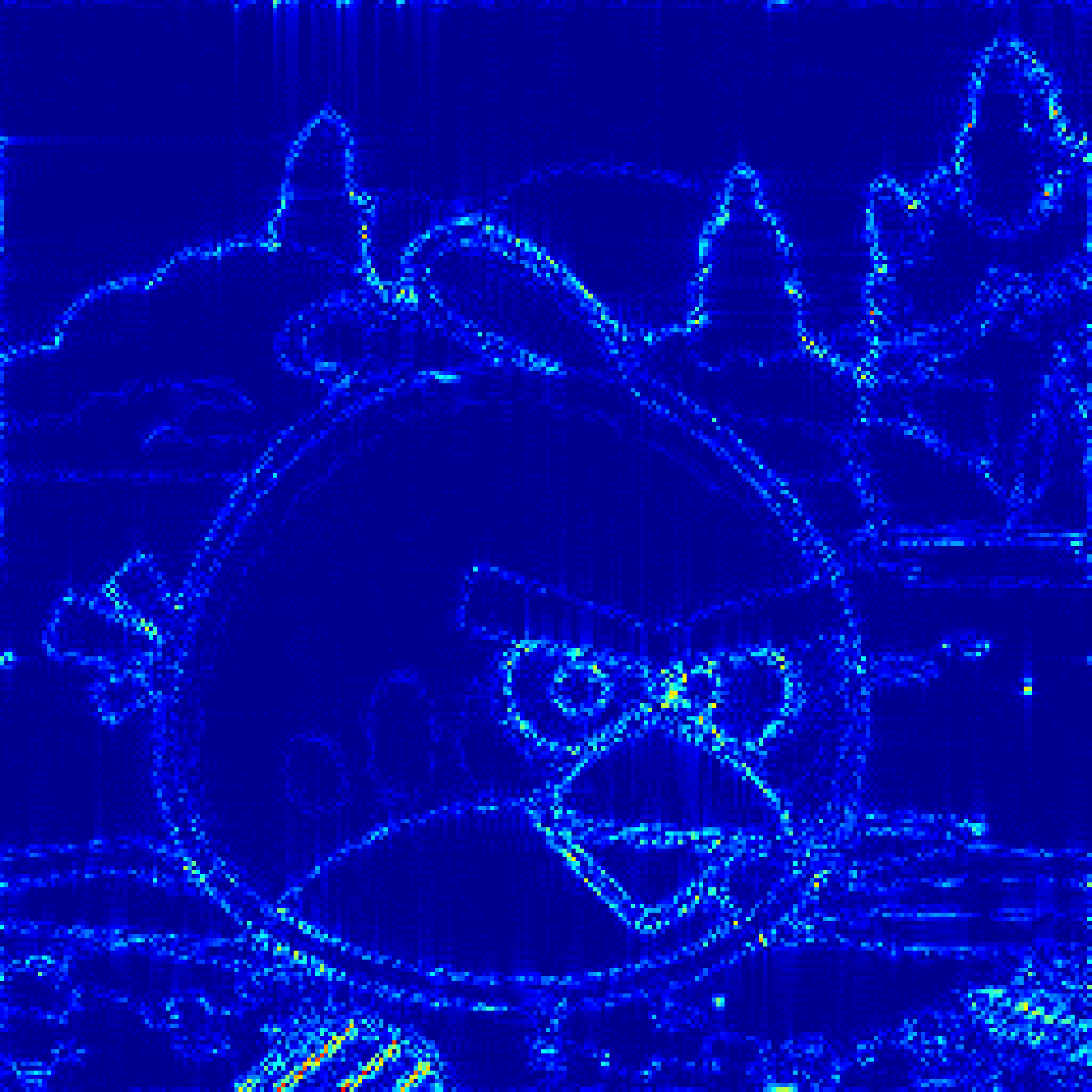}
\end{minipage}}\hspace{-0.05cm}
\subfloat[LRHDDTF]{\label{AngryBirdsRandomLRHDDTFError}\begin{minipage}{3cm}
\includegraphics[width=3cm]{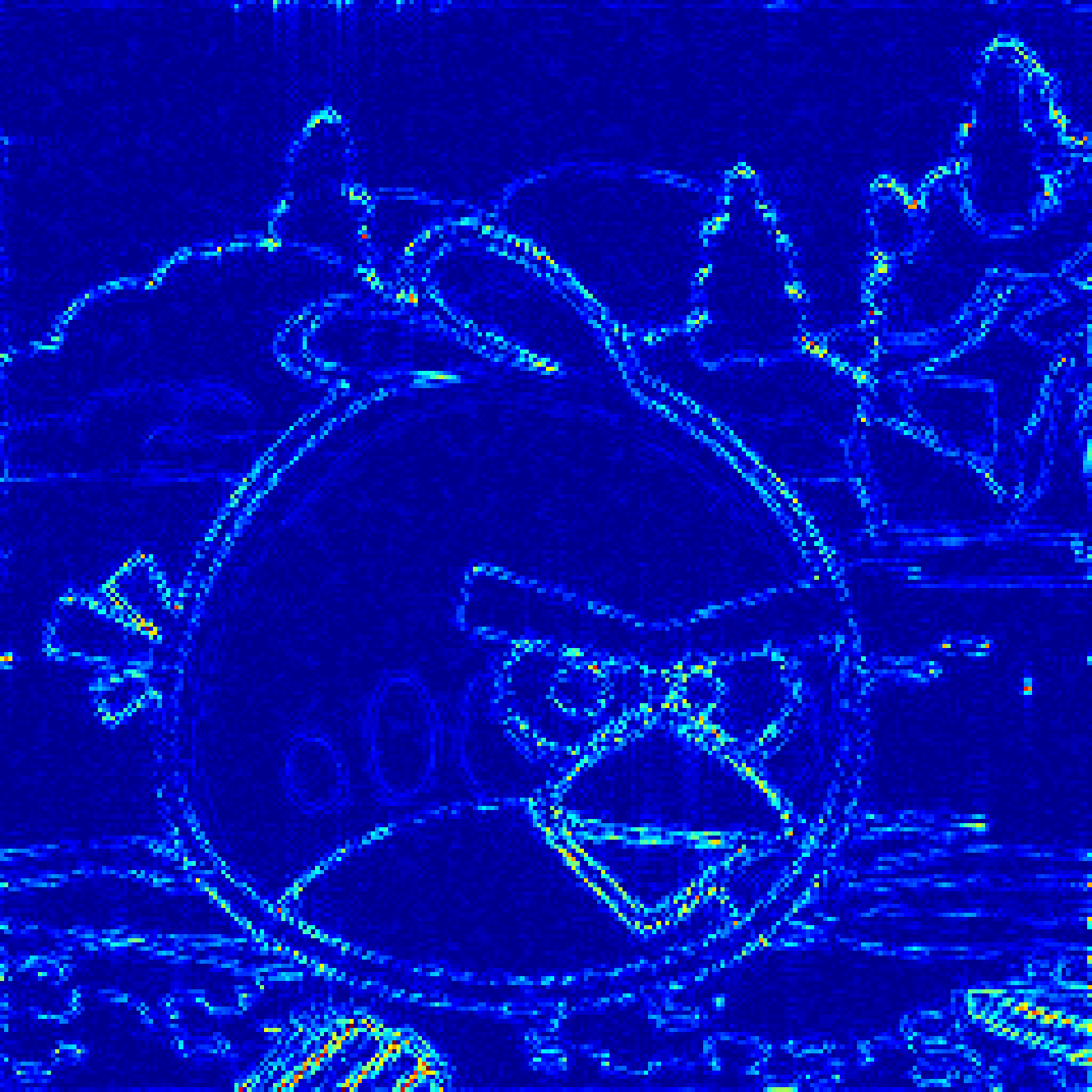}
\end{minipage}}\vspace{-0.25cm}
\subfloat[Schatten $0$]{\label{AngryBirdsRandomGIRAFError}\begin{minipage}{3cm}
\includegraphics[width=3cm]{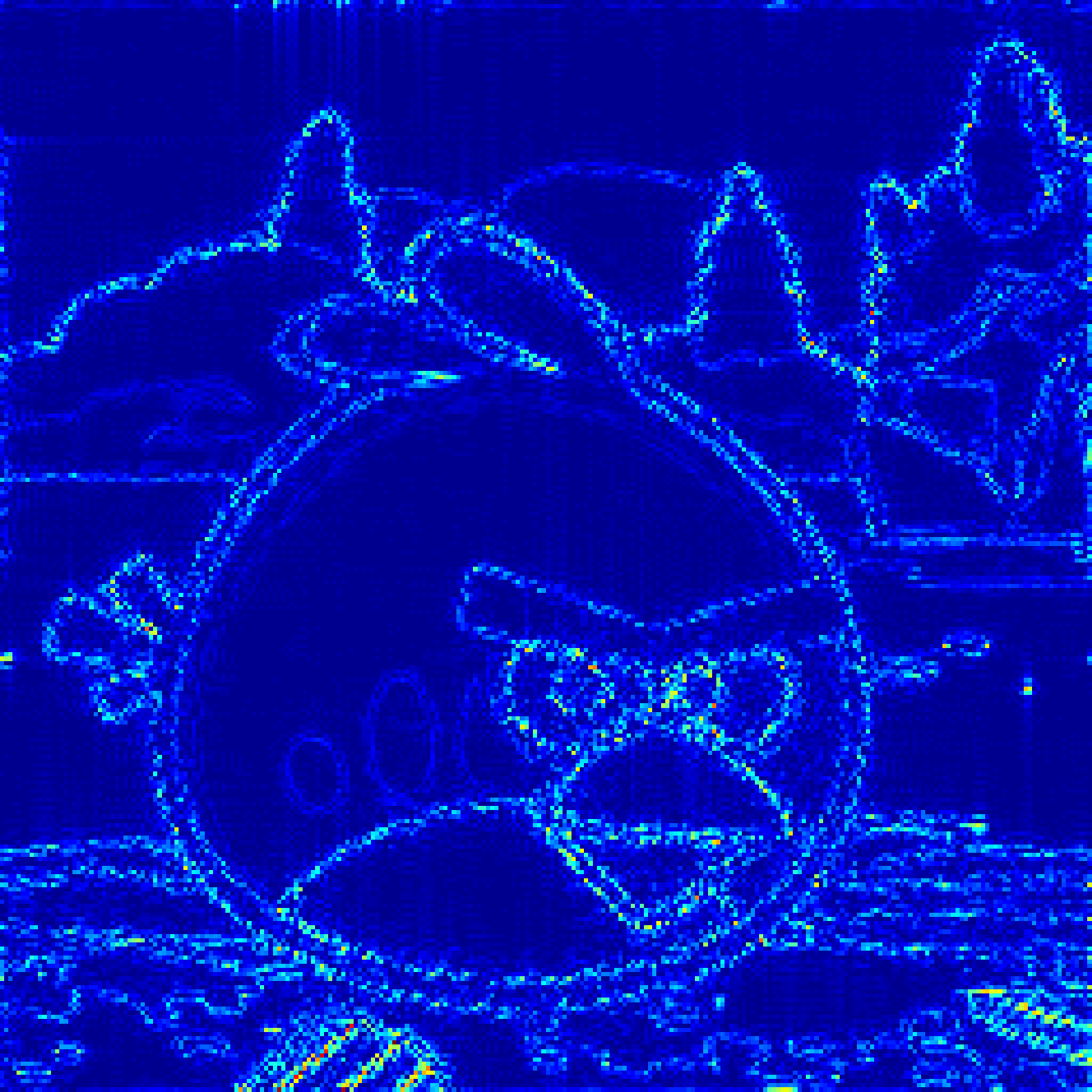}
\end{minipage}}\hspace{-0.05cm}
\subfloat[TV]{\label{AngryBirdsRandomTVError}\begin{minipage}{3cm}
\includegraphics[width=3cm]{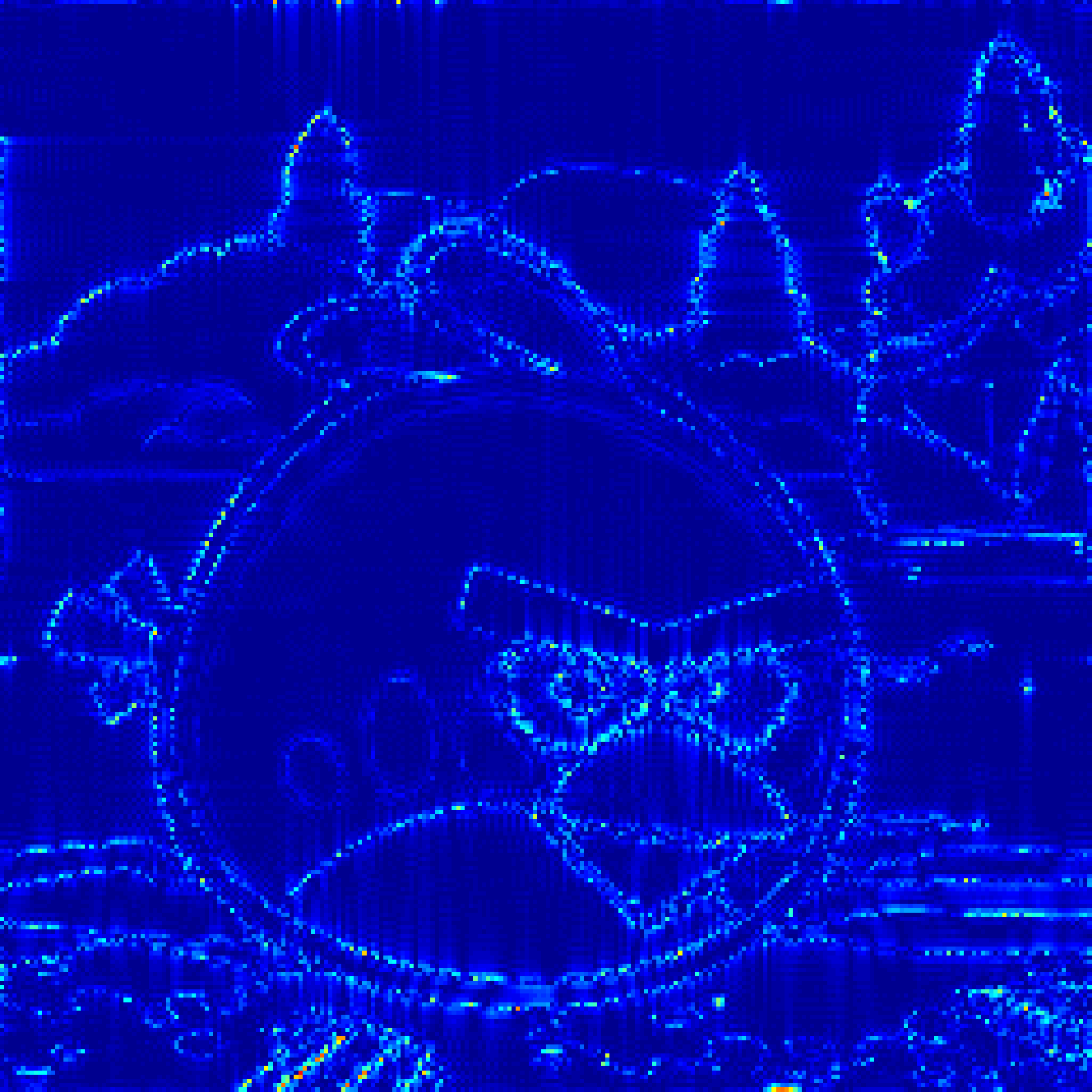}
\end{minipage}}\hspace{-0.05cm}
\subfloat[Haar]{\label{AngryBirdsRandomHaarError}\begin{minipage}{3cm}
\includegraphics[width=3cm]{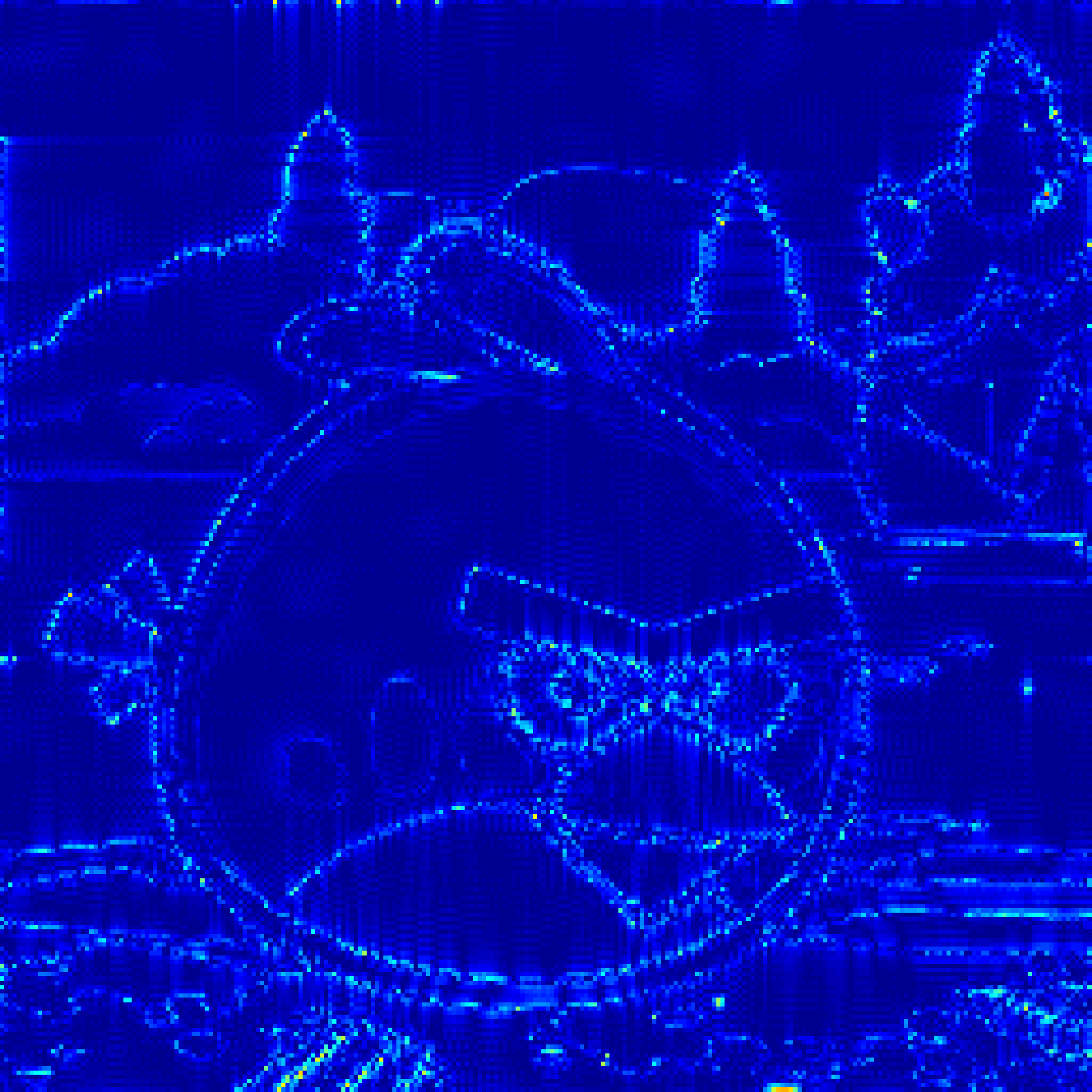}
\end{minipage}}\hspace{-0.05cm}
\subfloat[DDTF]{\label{AngryBirdsRandomDDTFError}\begin{minipage}{3cm}
\includegraphics[width=3cm]{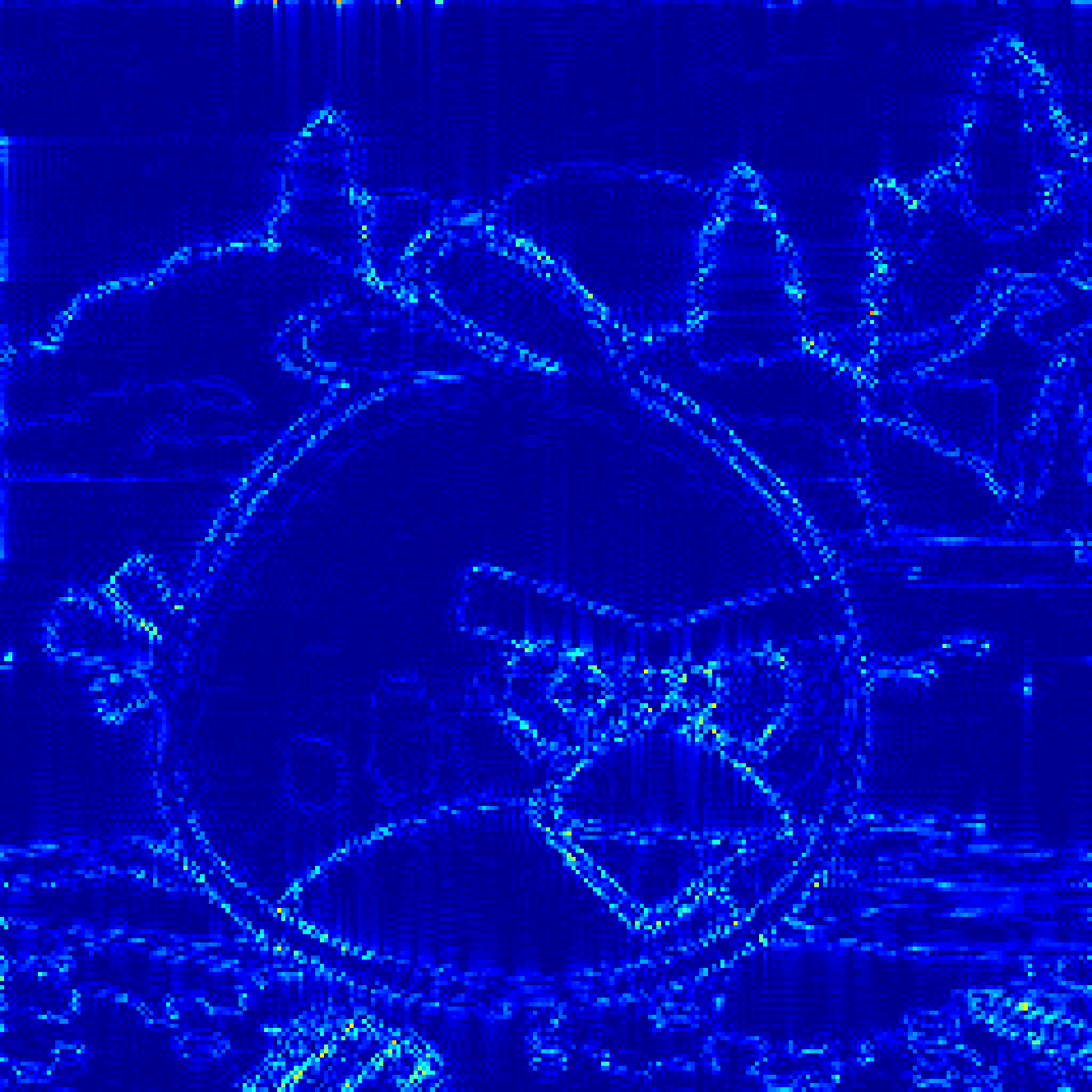}
\end{minipage}}
\caption{Error maps of \cref{AngryBirdsRandomResults}.}\label{AngryBirdsRandomError}
\end{figure}

\begin{figure}[t]
\centering
\subfloat[Ref.]{\label{AngryBirdsILFOriginal}\begin{minipage}{3cm}
\includegraphics[width=3cm]{AngryBirdsOriginal.pdf}
\end{minipage}}\hspace{-0.05cm}
\subfloat[Proposed]{\label{AngryBirdsILFProposed}\begin{minipage}{3cm}
\includegraphics[width=3cm]{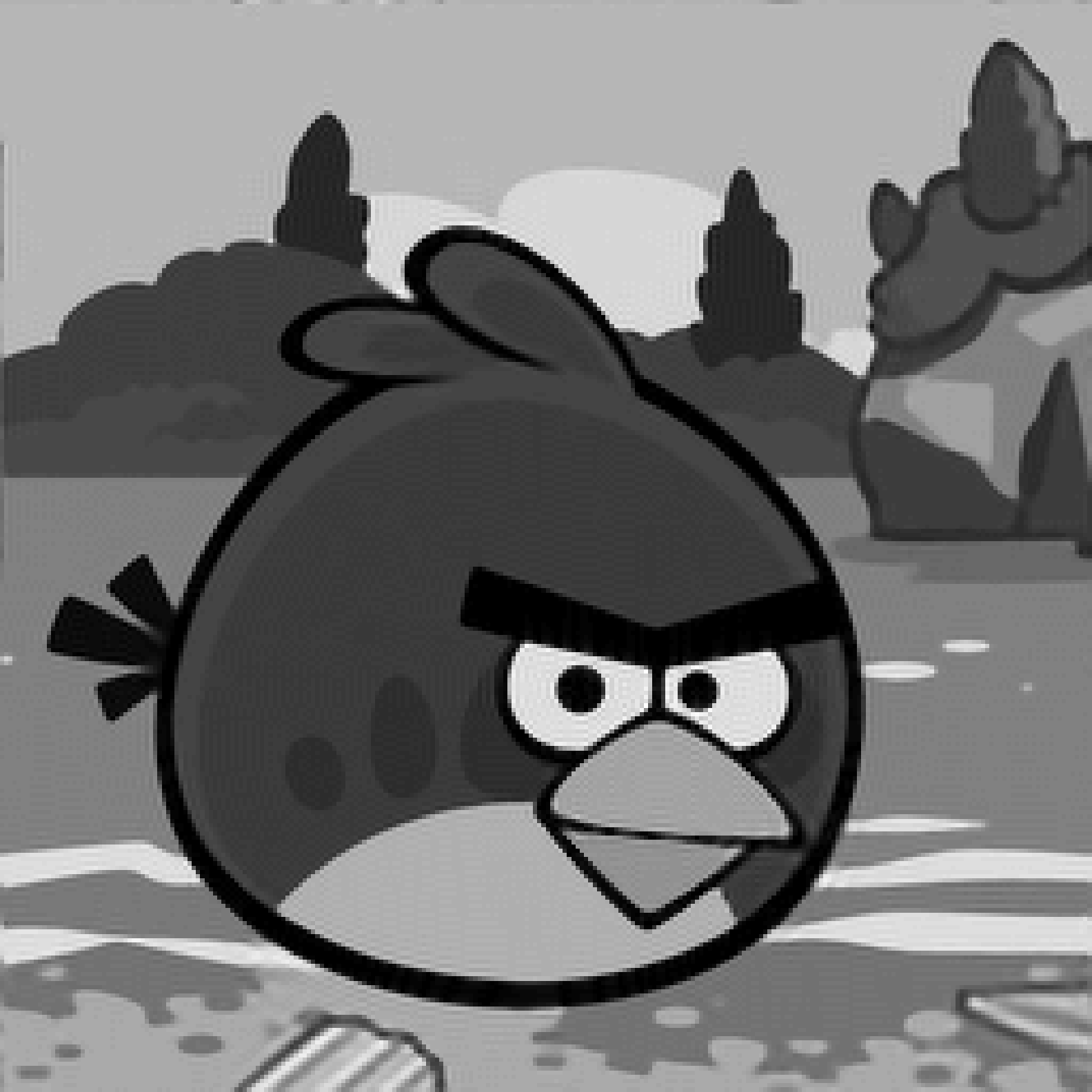}
\end{minipage}}\hspace{-0.05cm}
\subfloat[LSLP]{\label{AngryBirdsILFLSLP}\begin{minipage}{3cm}
\includegraphics[width=3cm]{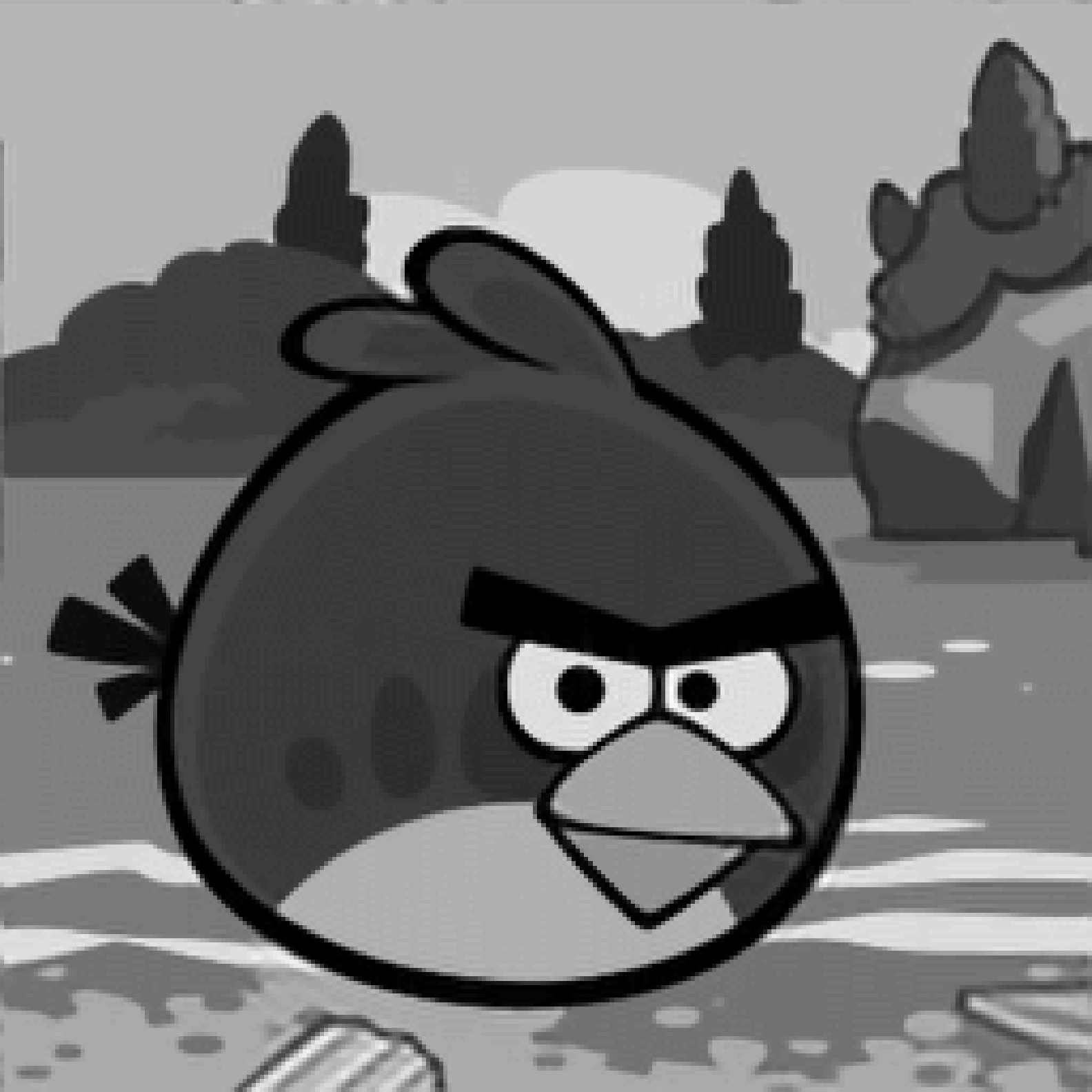}
\end{minipage}}\hspace{-0.05cm}
\subfloat[LRHDDTF]{\label{AngryBirdsILFLRHDDTF}\begin{minipage}{3cm}
\includegraphics[width=3cm]{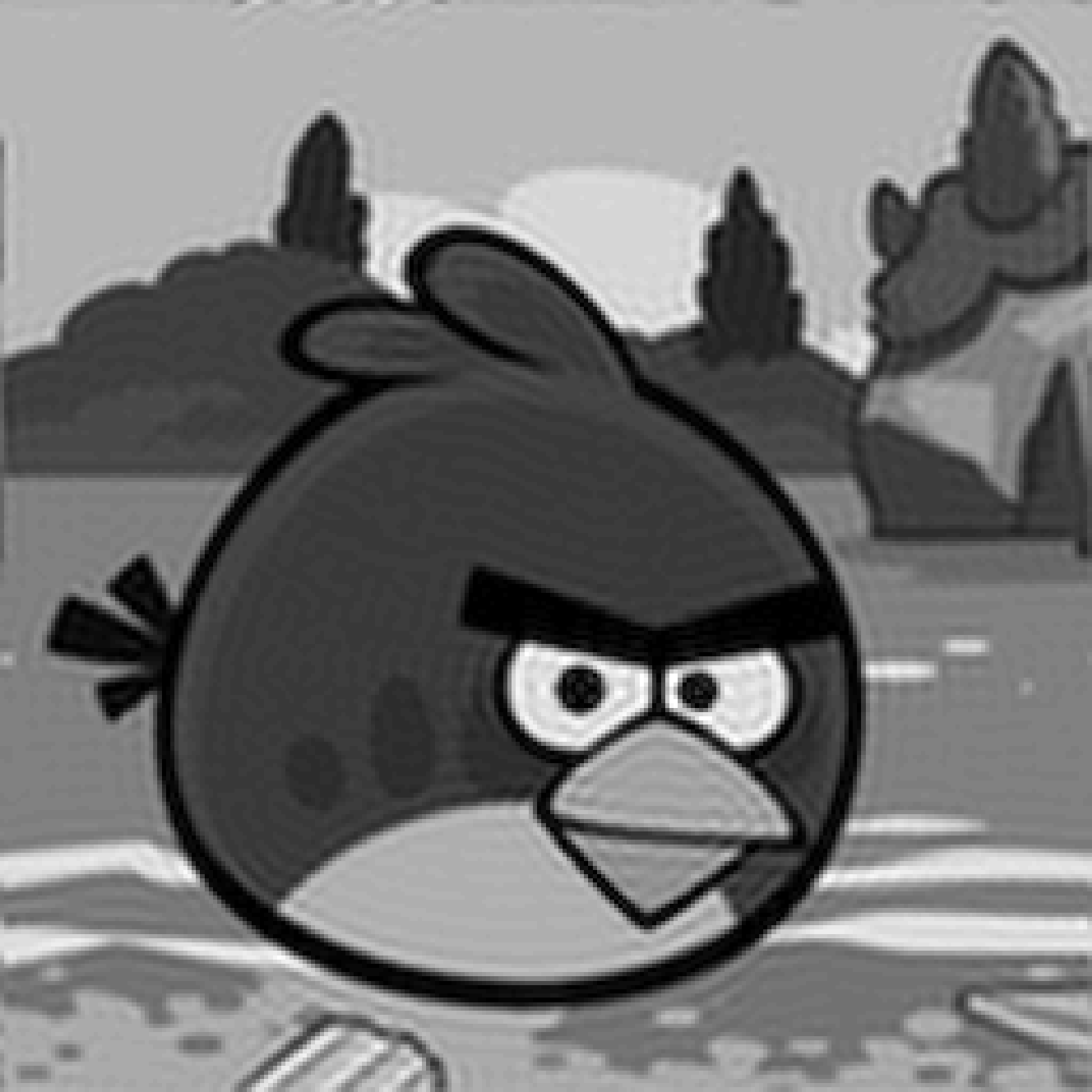}
\end{minipage}}\vspace{-0.25cm}
\subfloat[Schatten $0$]{\label{AngryBirdsILFGIRAF}\begin{minipage}{3cm}
\includegraphics[width=3cm]{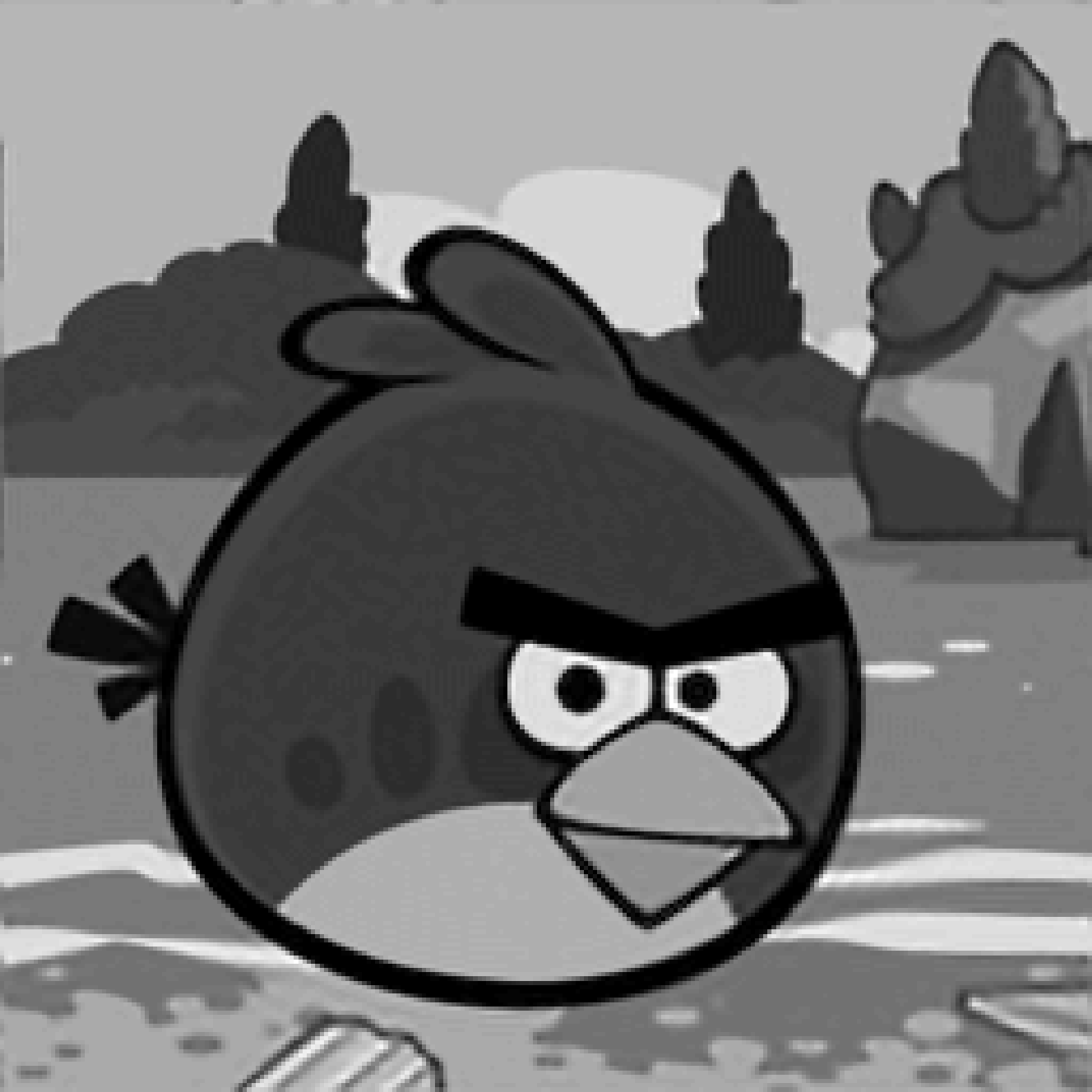}
\end{minipage}}\hspace{-0.05cm}
\subfloat[TV]{\label{AngryBirdsILFTV}\begin{minipage}{3cm}
\includegraphics[width=3cm]{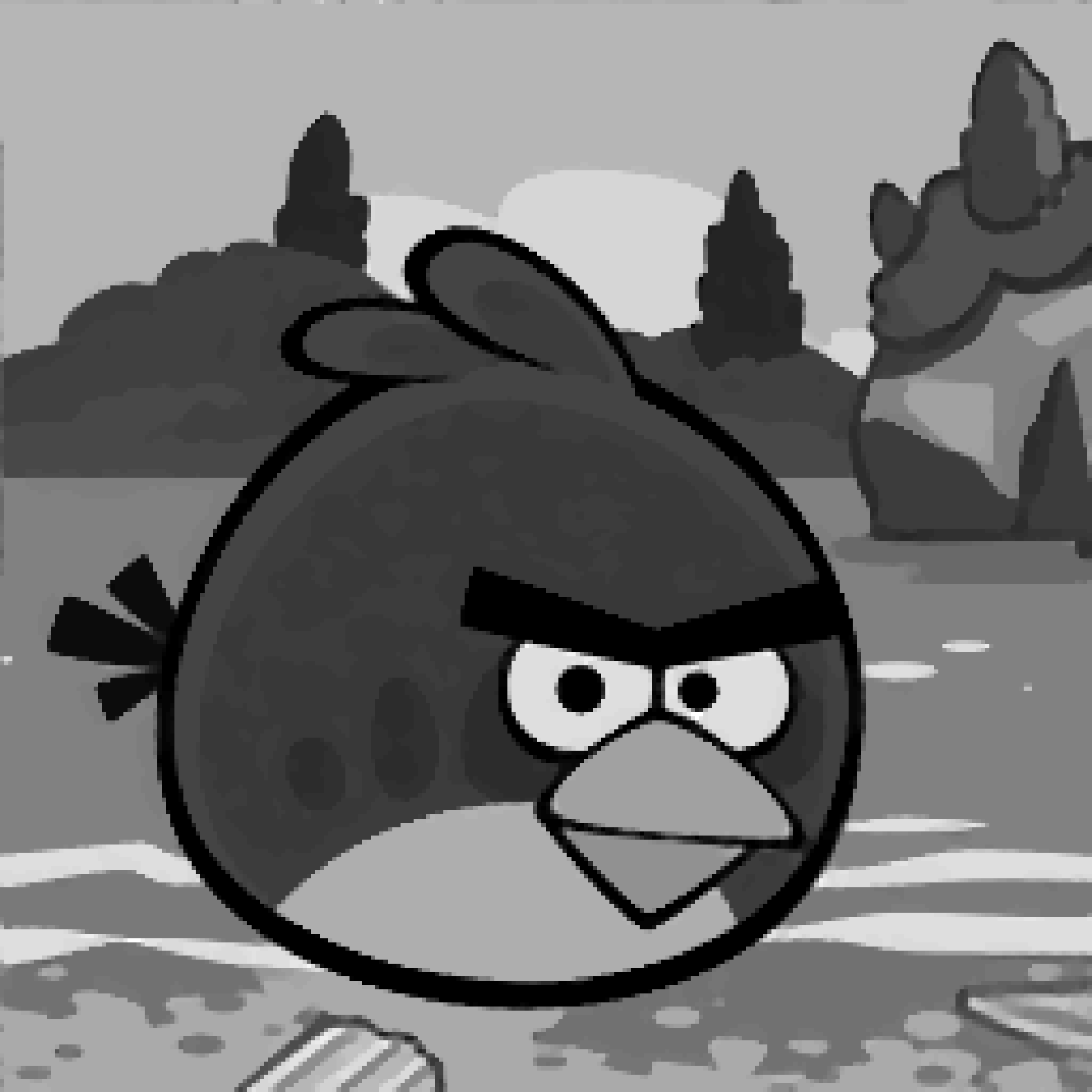}
\end{minipage}}\hspace{-0.05cm}
\subfloat[Haar]{\label{AngryBirdsILFHaar}\begin{minipage}{3cm}
\includegraphics[width=3cm]{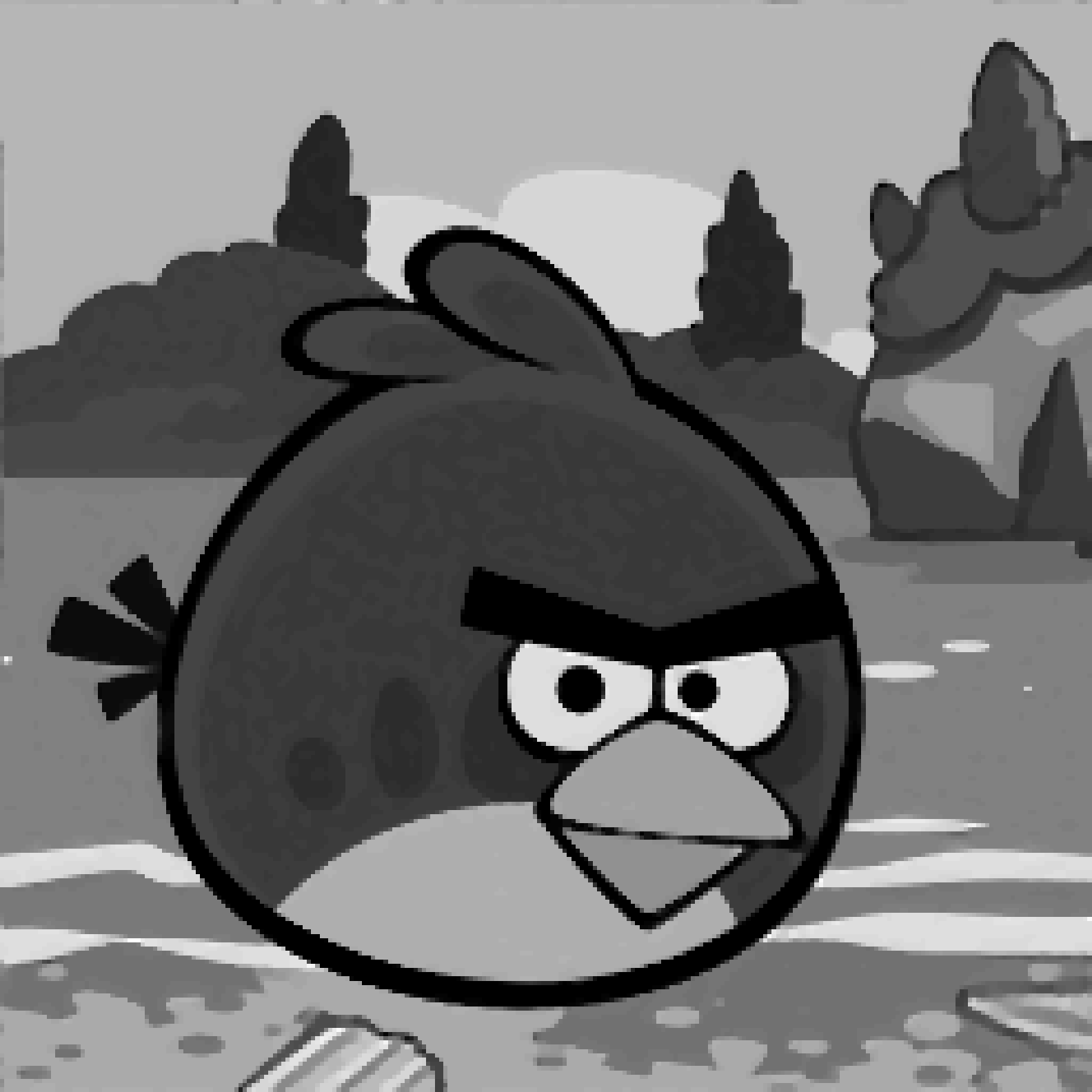}
\end{minipage}}\hspace{-0.05cm}
\subfloat[DDTF]{\label{AngryBirdsILFDDTF}\begin{minipage}{3cm}
\includegraphics[width=3cm]{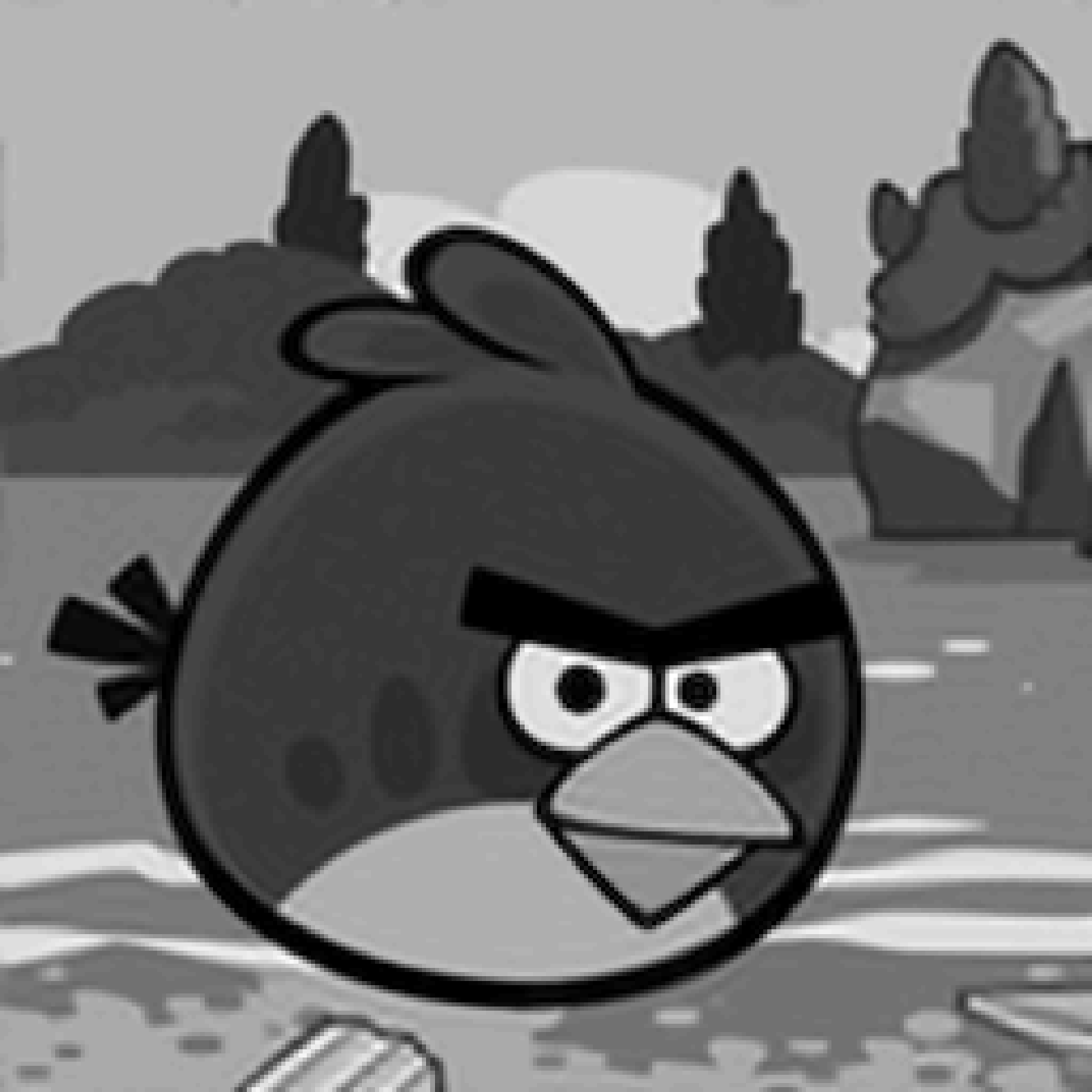}
\end{minipage}}
\caption{Visual comparison of ``Angry Birds'' for ideal low-pass filter deconvolution.}\label{AngryBirdsILFResults}
\end{figure}

\begin{figure}[t]
\centering
\subfloat[TF]{\label{AngryBirdsILFTF}\begin{minipage}{3cm}
\includegraphics[width=3cm]{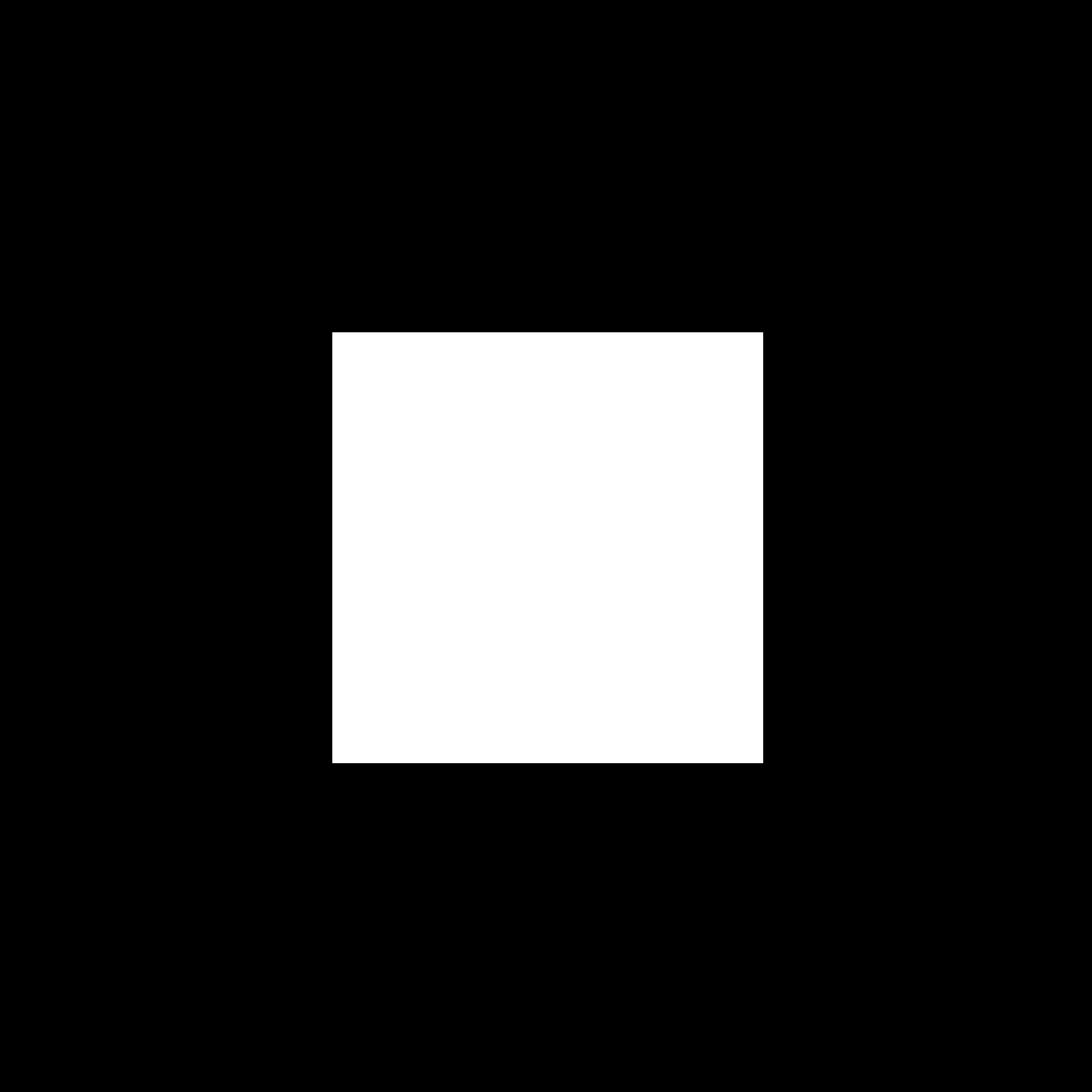}
\end{minipage}}\hspace{-0.05cm}
\subfloat[Proposed]{\label{AngryBirdsILFProposedError}\begin{minipage}{3cm}
\includegraphics[width=3cm]{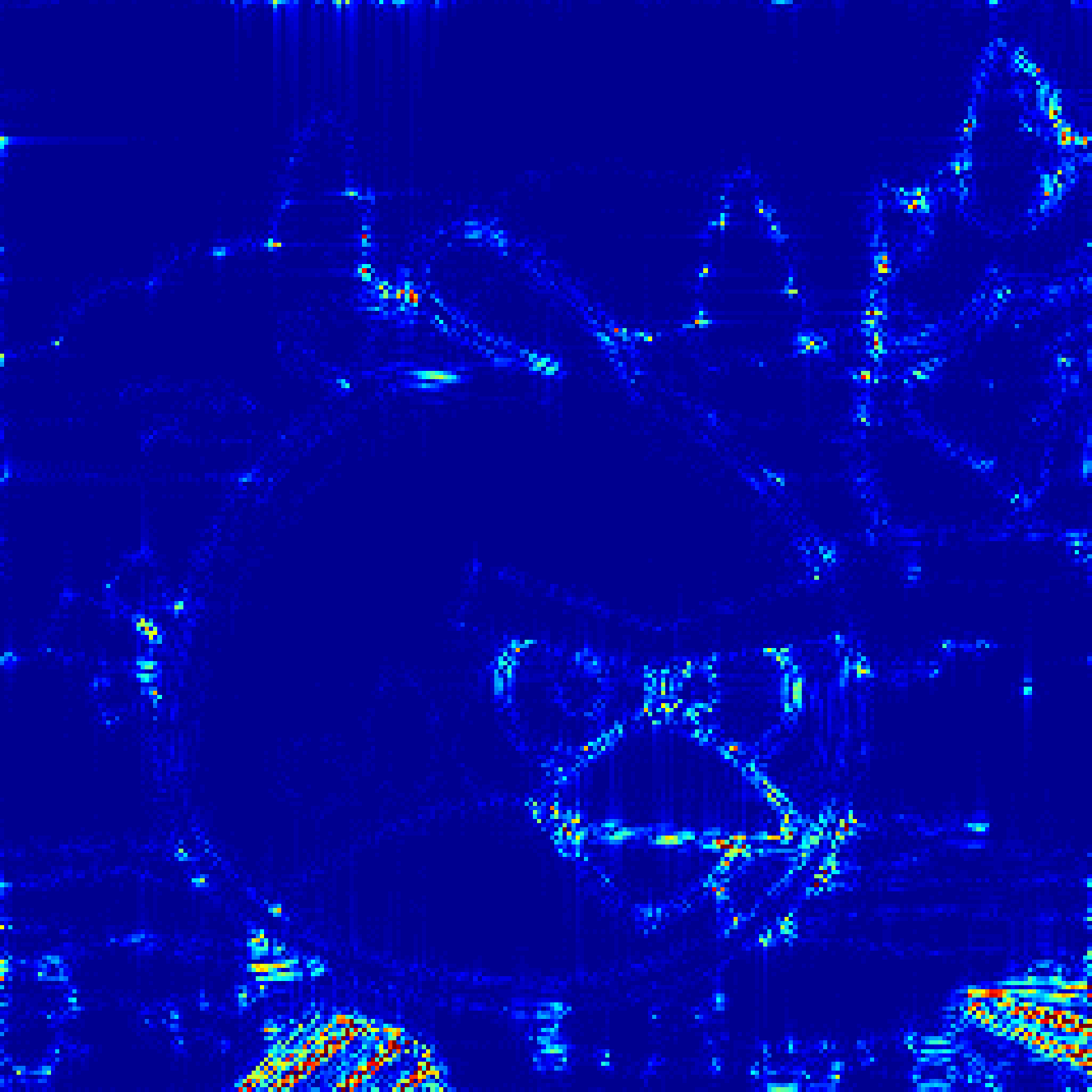}
\end{minipage}}\hspace{-0.05cm}
\subfloat[LSLP]{\label{AngryBirdsILFLSLPError}\begin{minipage}{3cm}
\includegraphics[width=3cm]{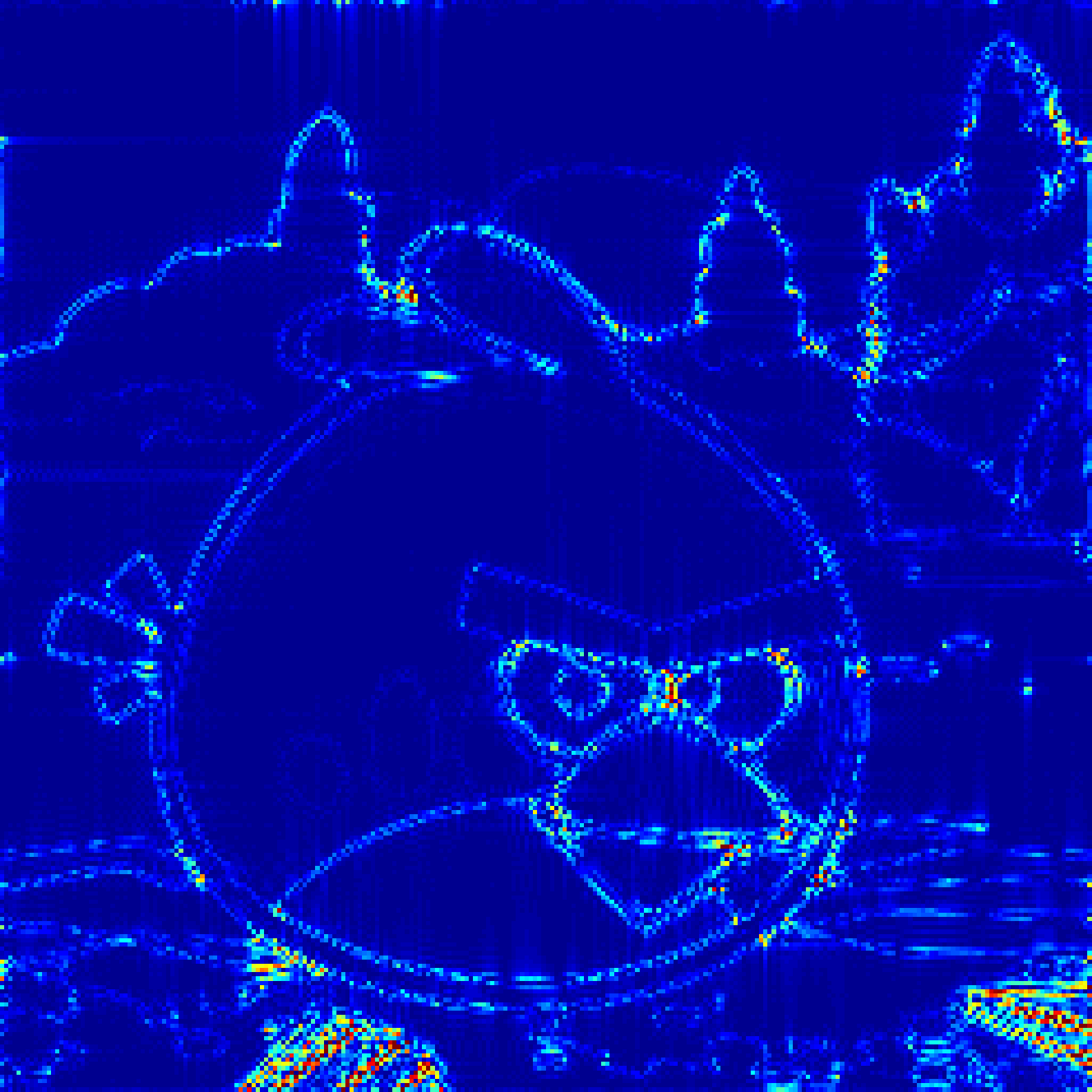}
\end{minipage}}\hspace{-0.05cm}
\subfloat[LRHDDTF]{\label{AngryBirdsILFLRHDDTFError}\begin{minipage}{3cm}
\includegraphics[width=3cm]{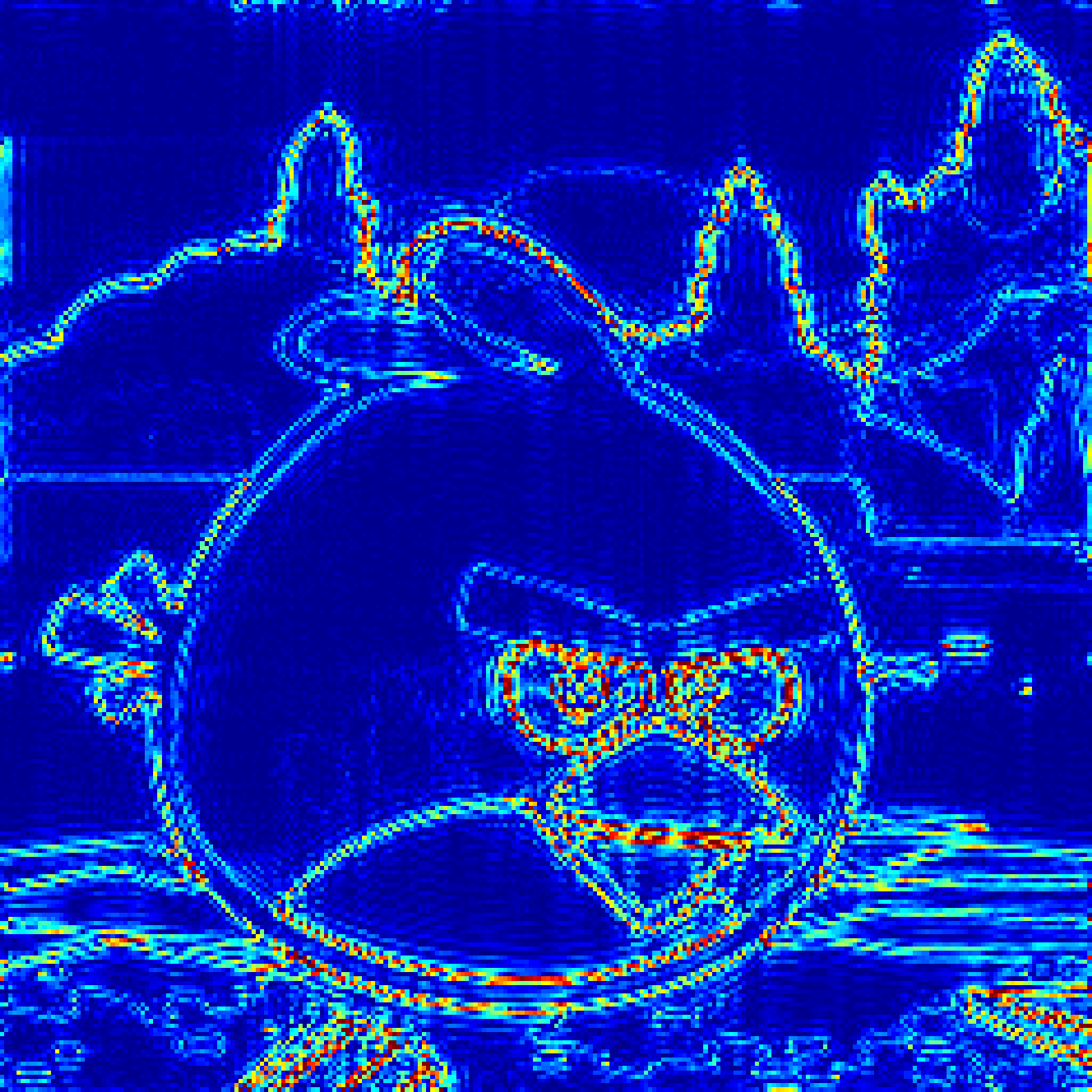}
\end{minipage}}\vspace{-0.25cm}
\subfloat[Schatten $0$]{\label{AngryBirdsILFGIRAFError}\begin{minipage}{3cm}
\includegraphics[width=3cm]{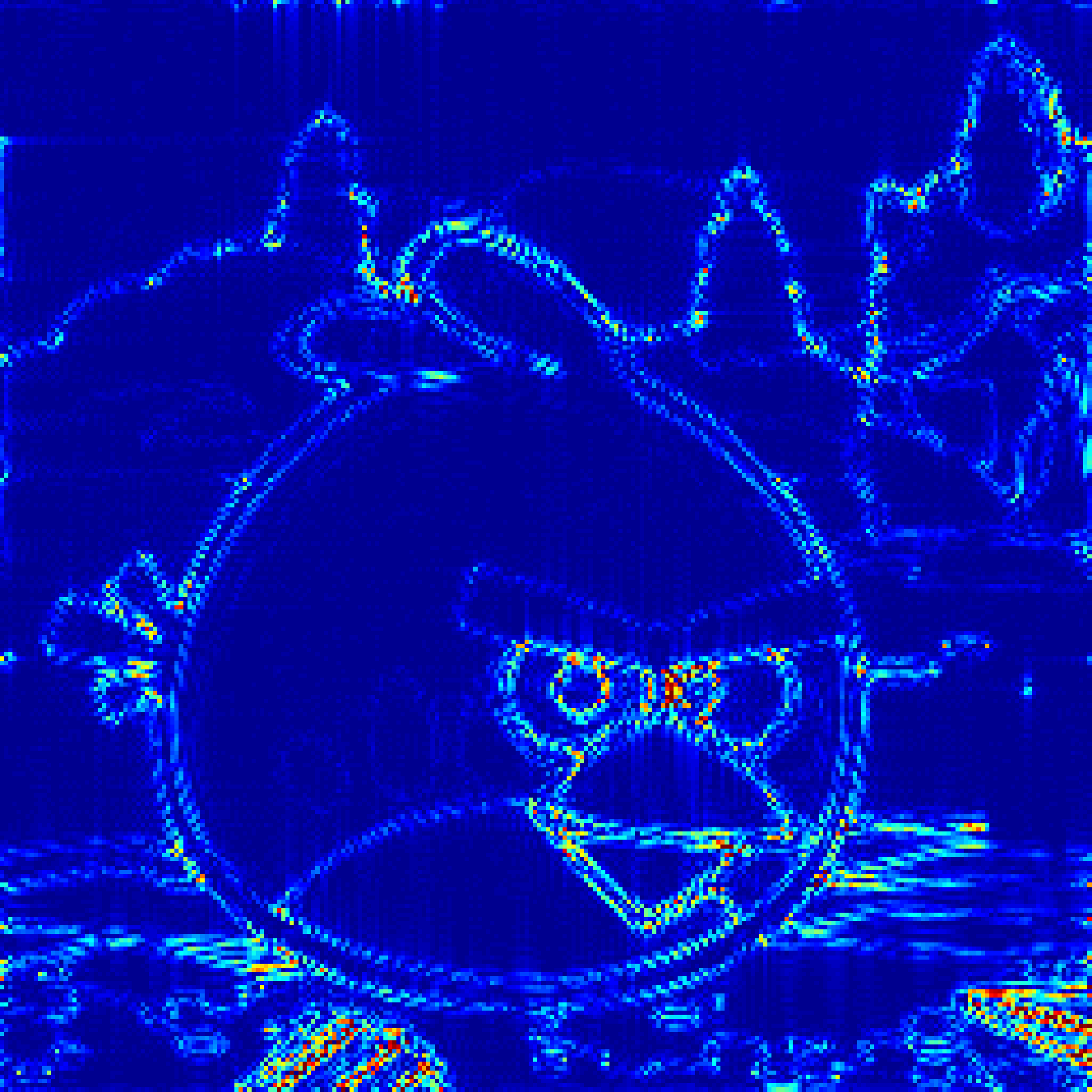}
\end{minipage}}\hspace{-0.05cm}
\subfloat[TV]{\label{AngryBirdsILFTVError}\begin{minipage}{3cm}
\includegraphics[width=3cm]{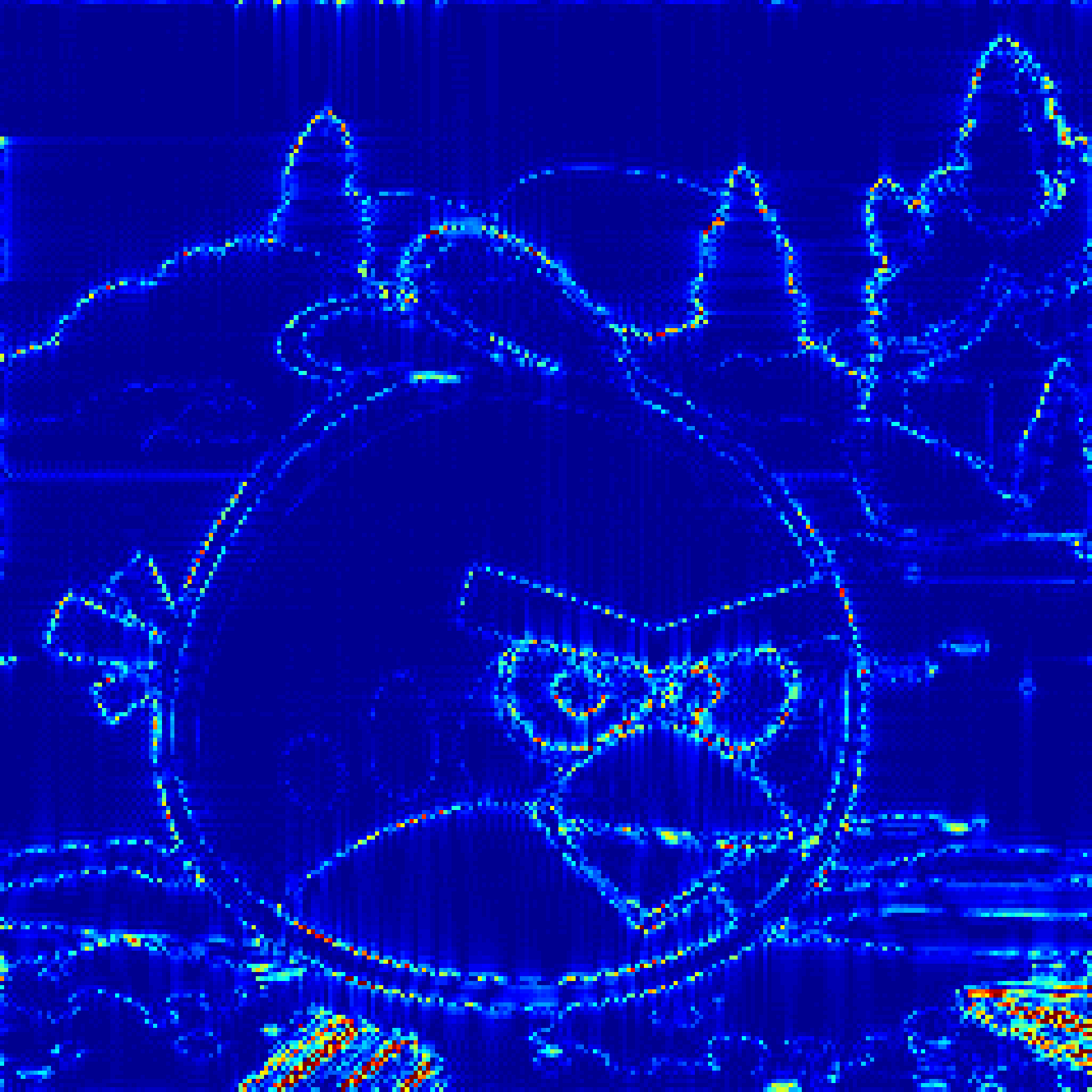}
\end{minipage}}\hspace{-0.05cm}
\subfloat[Haar]{\label{AngryBirdsILFHaarError}\begin{minipage}{3cm}
\includegraphics[width=3cm]{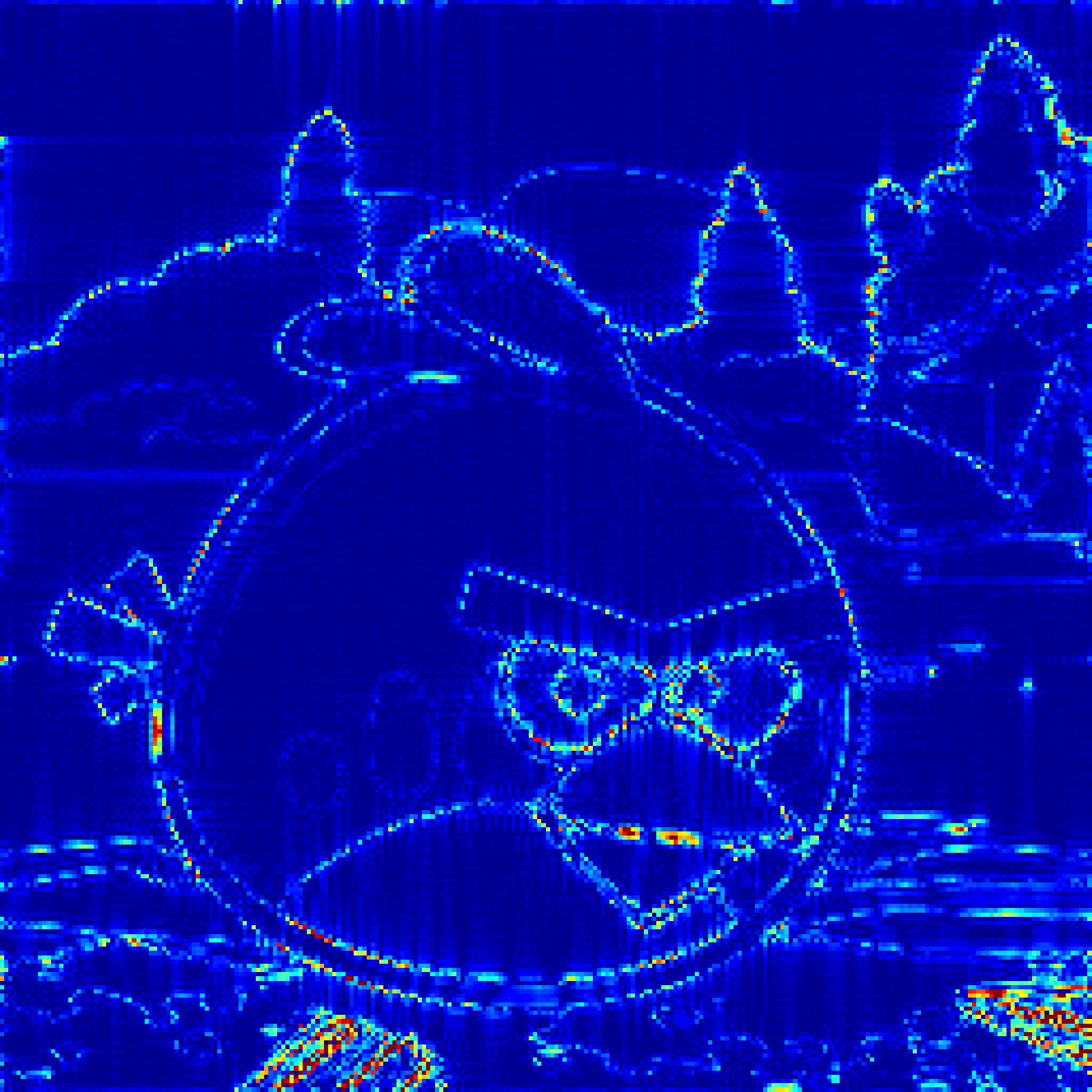}
\end{minipage}}\hspace{-0.05cm}
\subfloat[DDTF]{\label{AngryBirdsILFDDTFError}\begin{minipage}{3cm}
\includegraphics[width=3cm]{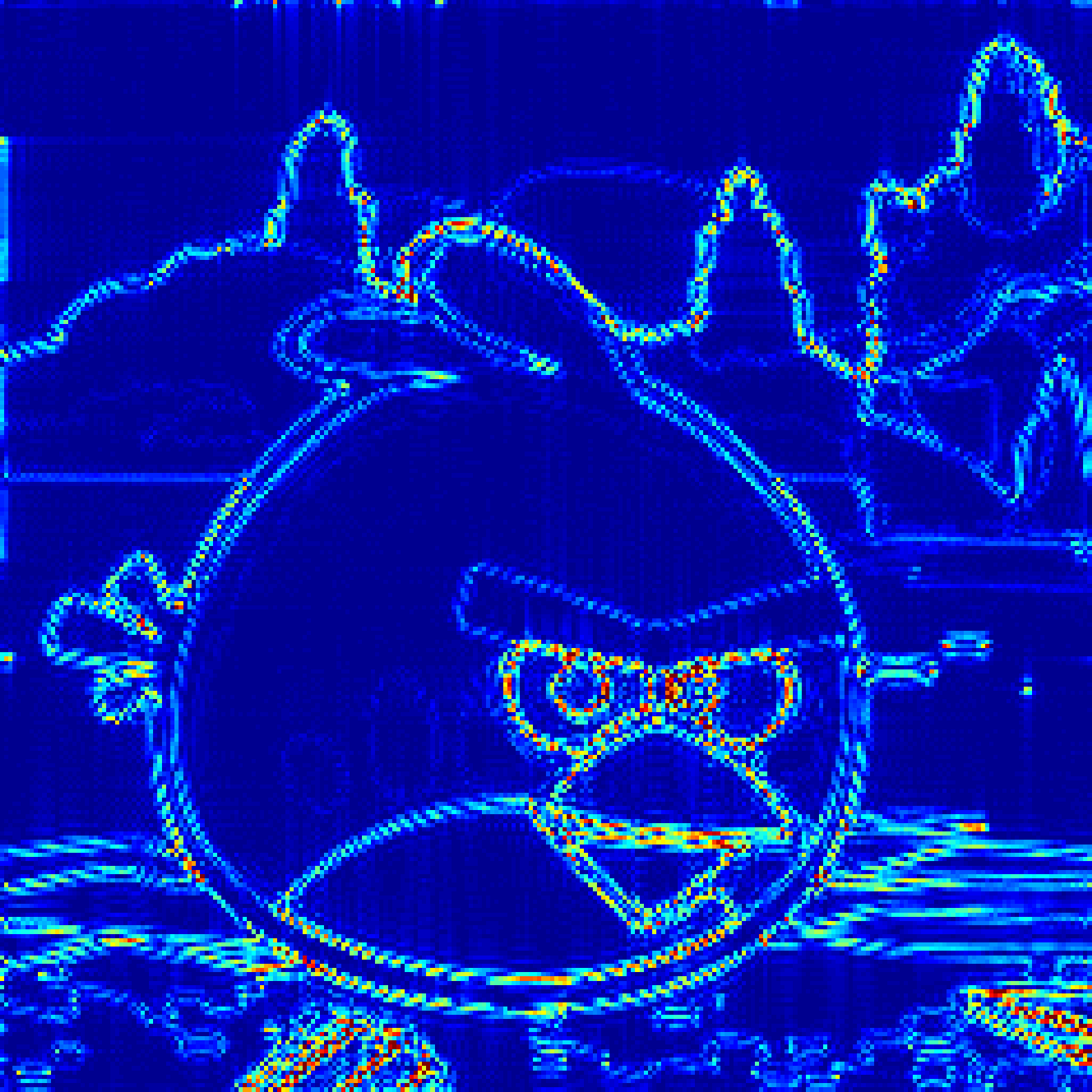}
\end{minipage}}
\caption{Error maps of \cref{AngryBirdsILFResults}.}\label{AngryBirdsILFError}
\end{figure}

\begin{figure}[t]
\centering
\subfloat[Ref.]{\label{SuperMarioRandomOriginal}\begin{minipage}{3cm}
\includegraphics[width=3cm]{SuperMarioOriginal.pdf}
\end{minipage}}\hspace{-0.05cm}
\subfloat[Proposed]{\label{SuperMarioRandomProposed}\begin{minipage}{3cm}
\includegraphics[width=3cm]{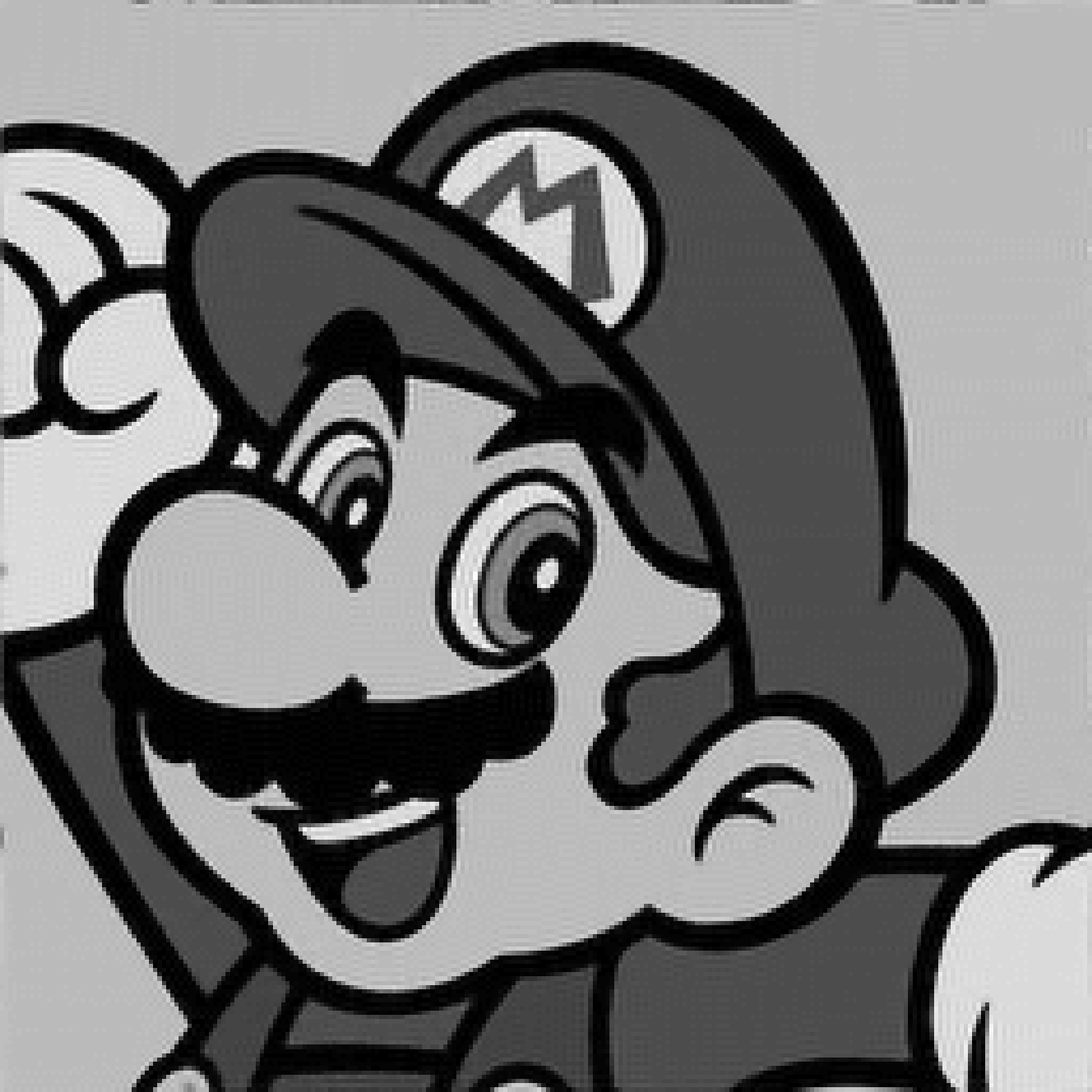}
\end{minipage}}\hspace{-0.05cm}
\subfloat[LSLP]{\label{SuperMarioRandomLSLP}\begin{minipage}{3cm}
\includegraphics[width=3cm]{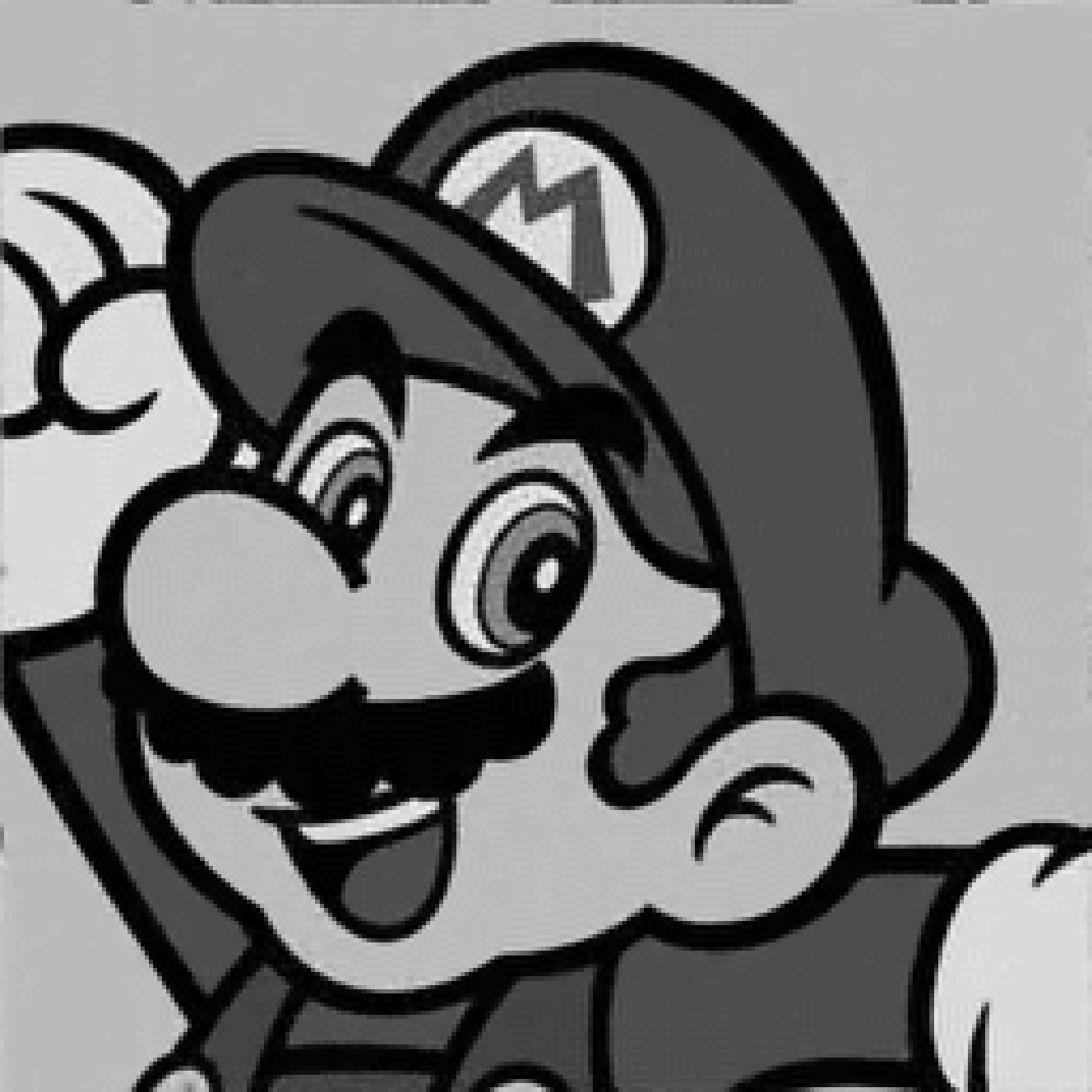}
\end{minipage}}\hspace{-0.05cm}
\subfloat[LRHDDTF]{\label{SuperMarioRandomLRHDDTF}\begin{minipage}{3cm}
\includegraphics[width=3cm]{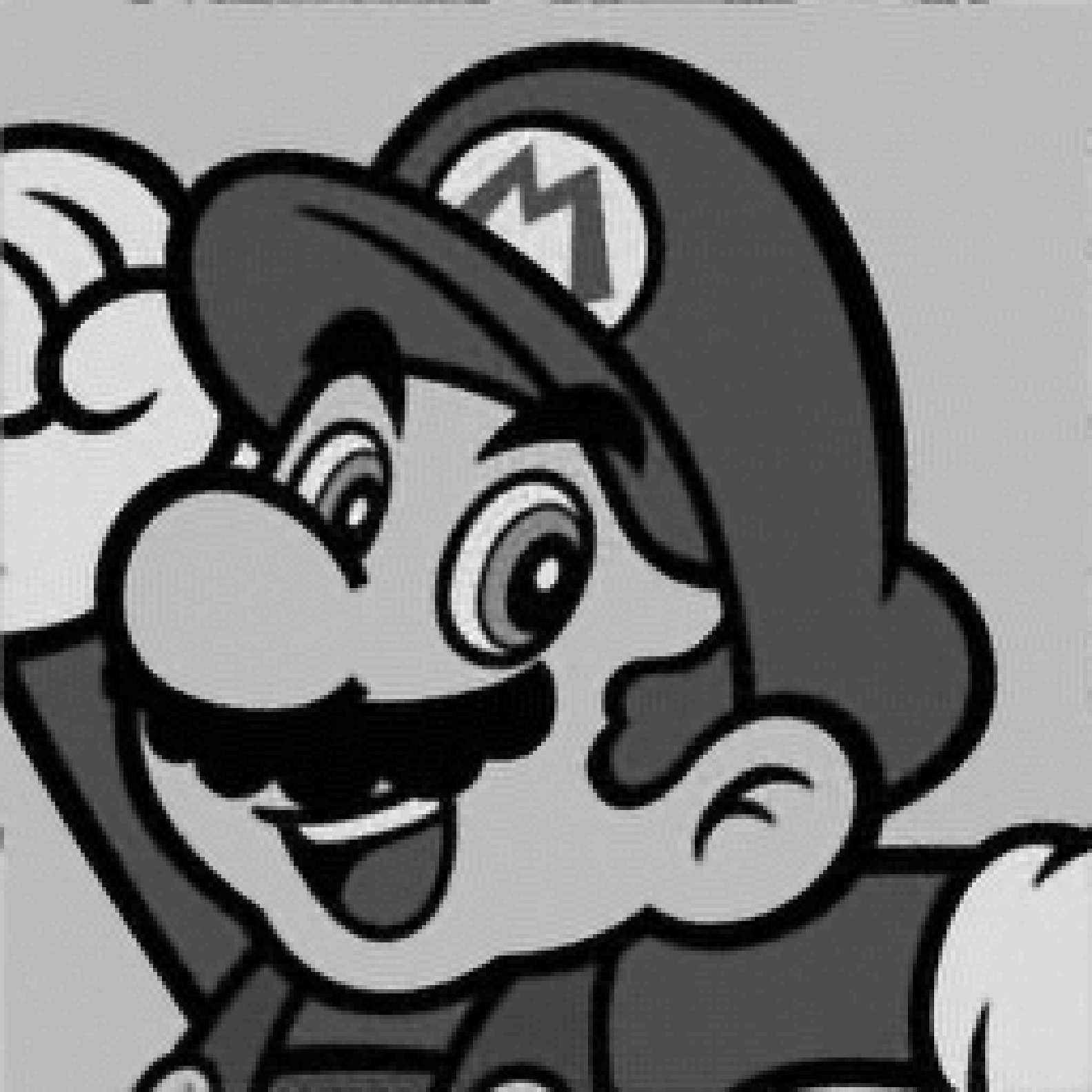}
\end{minipage}}\vspace{-0.25cm}
\subfloat[Schatten $0$]{\label{SuperMarioRandomGIRAF}\begin{minipage}{3cm}
\includegraphics[width=3cm]{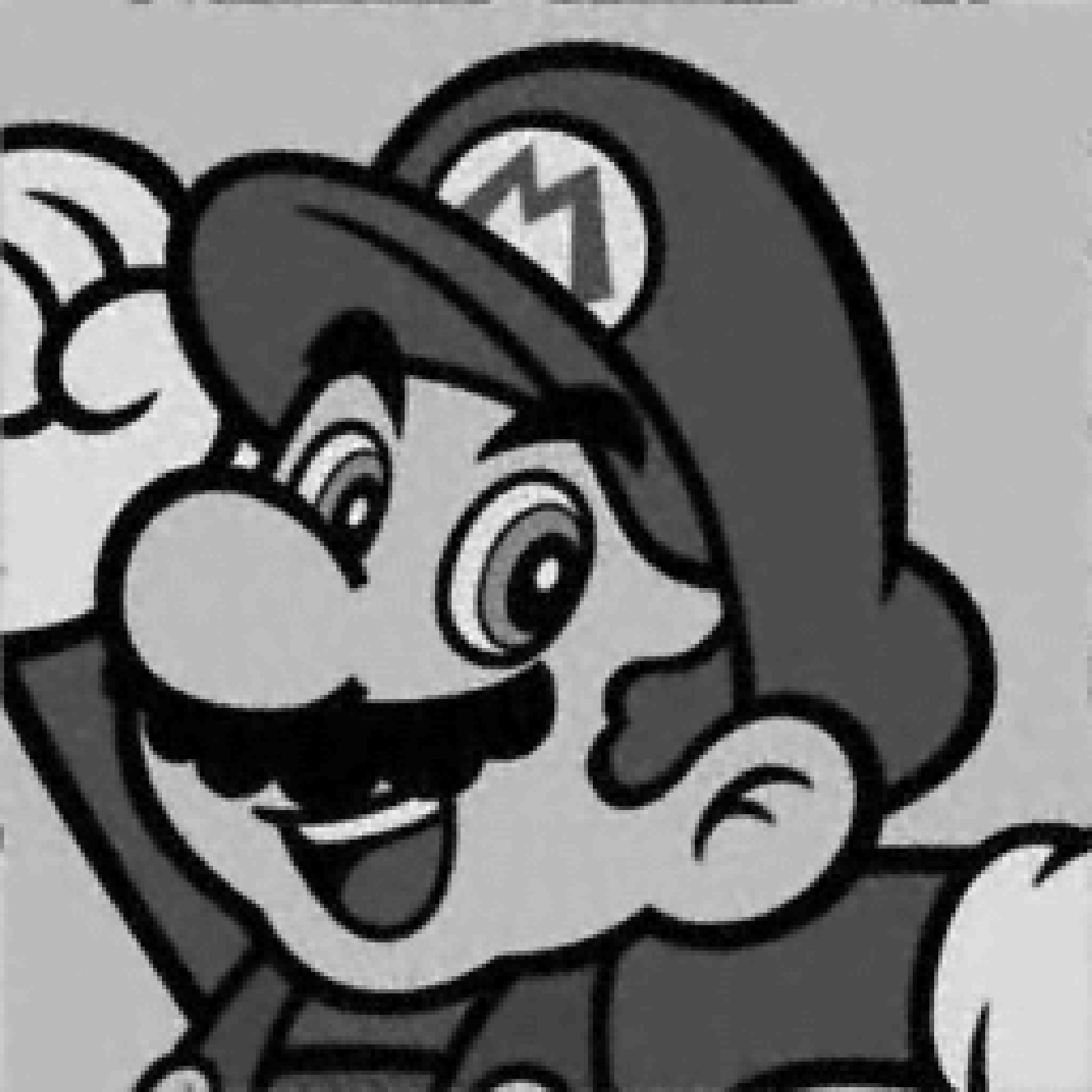}
\end{minipage}}\hspace{-0.05cm}
\subfloat[TV]{\label{SuperMarioRandomTV}\begin{minipage}{3cm}
\includegraphics[width=3cm]{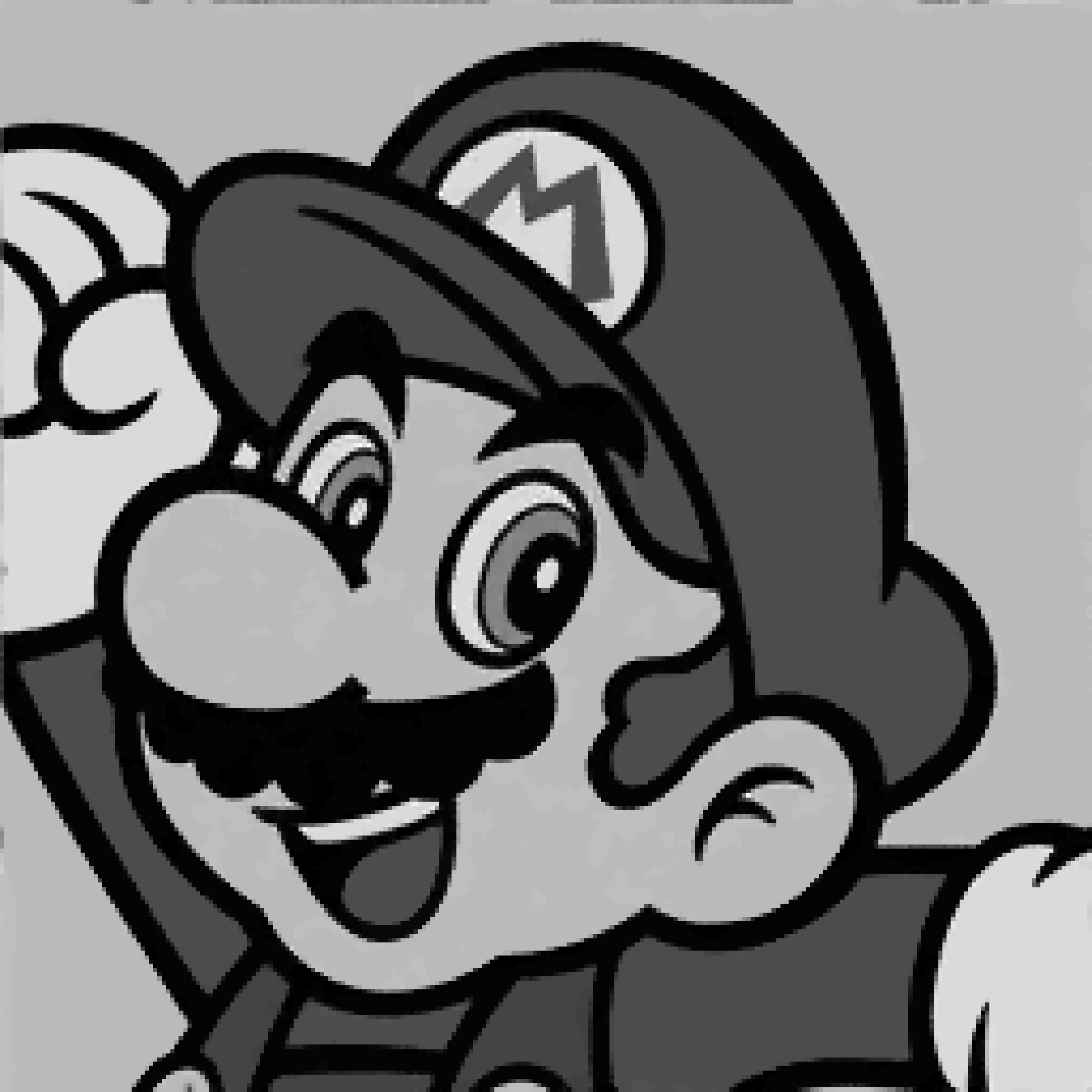}
\end{minipage}}\hspace{-0.05cm}
\subfloat[Haar]{\label{SuperMarioRandomHaar}\begin{minipage}{3cm}
\includegraphics[width=3cm]{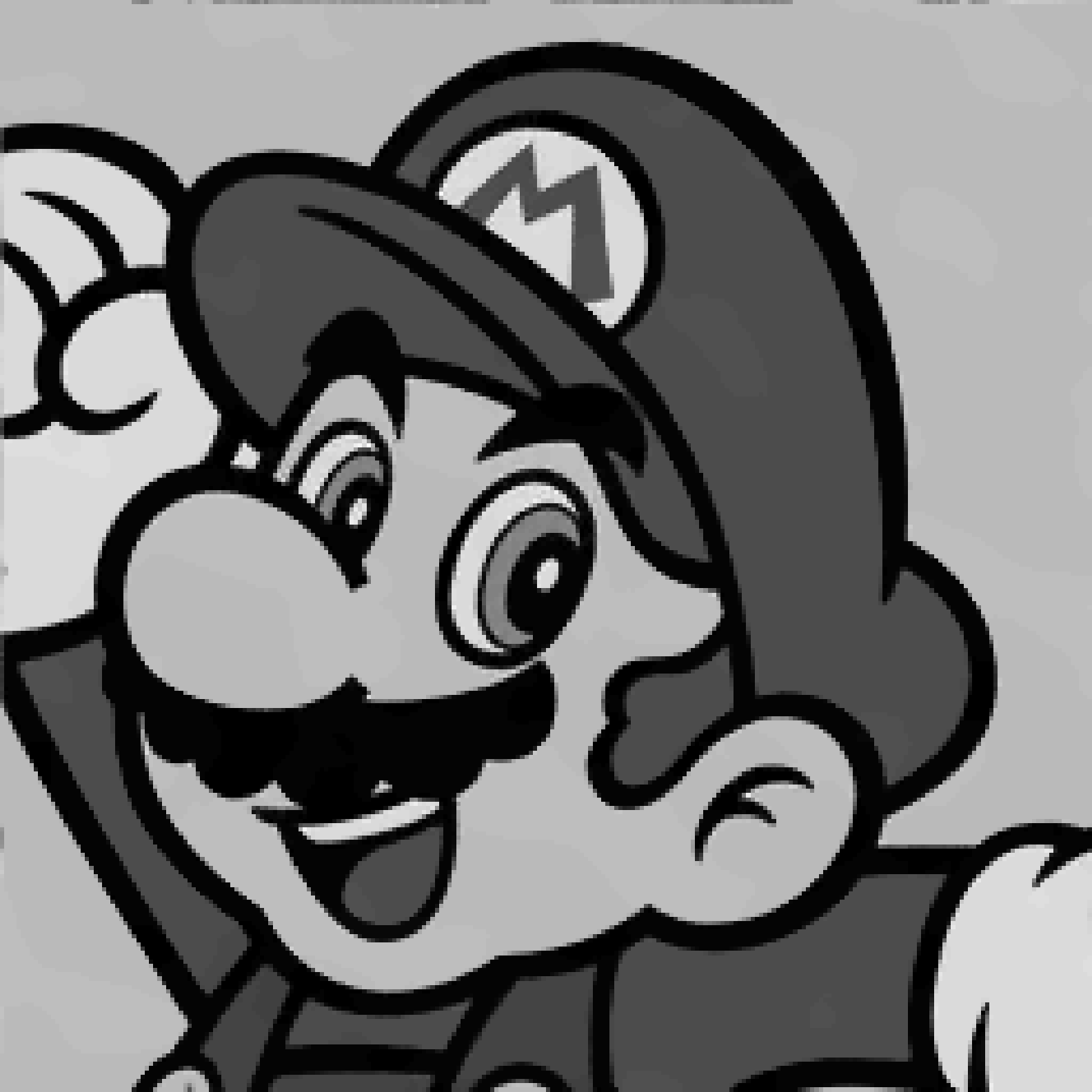}
\end{minipage}}\hspace{-0.05cm}
\subfloat[DDTF]{\label{SuperMarioRandomDDTF}\begin{minipage}{3cm}
\includegraphics[width=3cm]{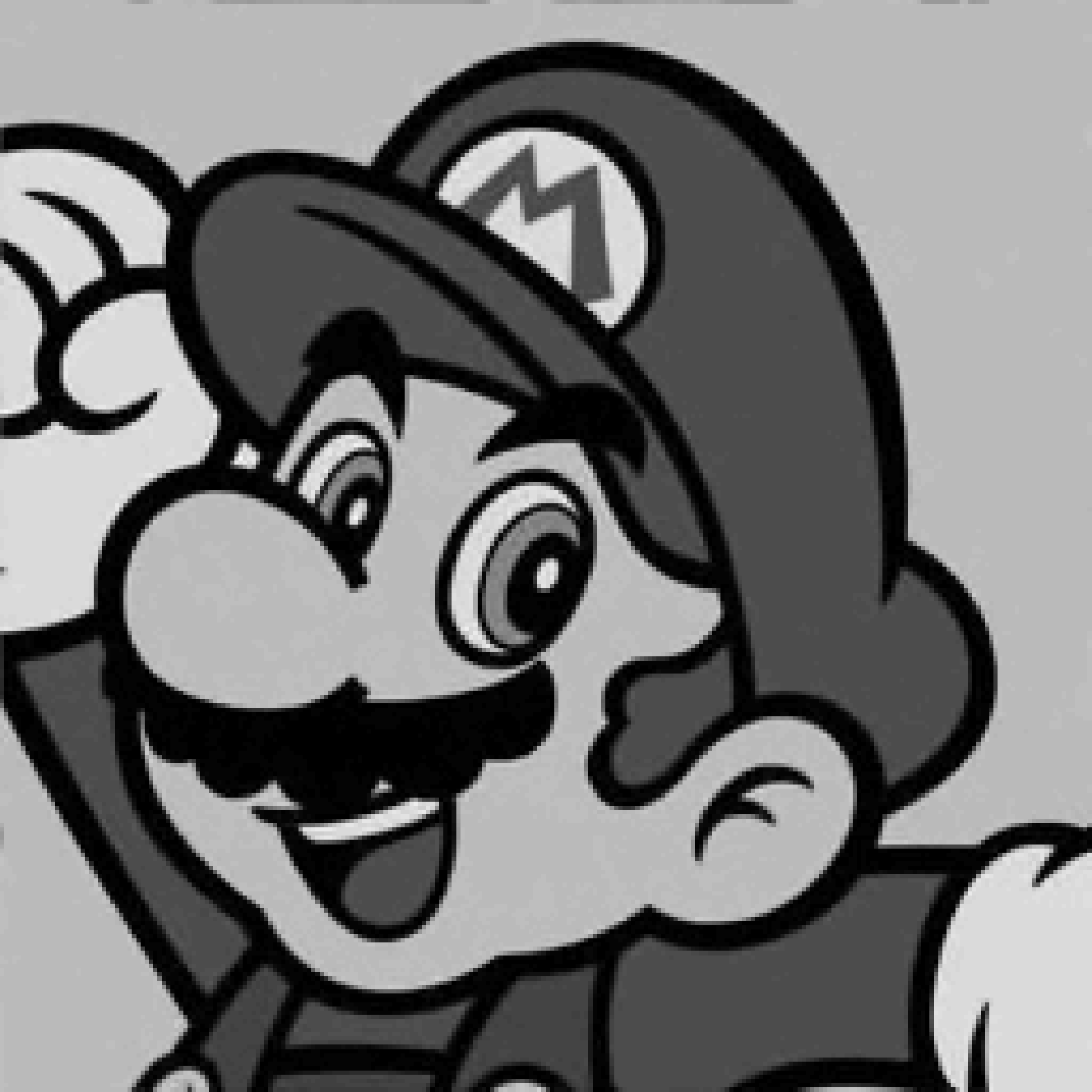}
\end{minipage}}
\caption{Visual comparison of ``Super Mario'' for random sampling.}\label{SuperMarioRandomResults}
\end{figure}

\begin{figure}[t]
\centering
\subfloat[Samples]{\label{SuperMarioRandomSample}\begin{minipage}{3cm}
\includegraphics[width=3cm]{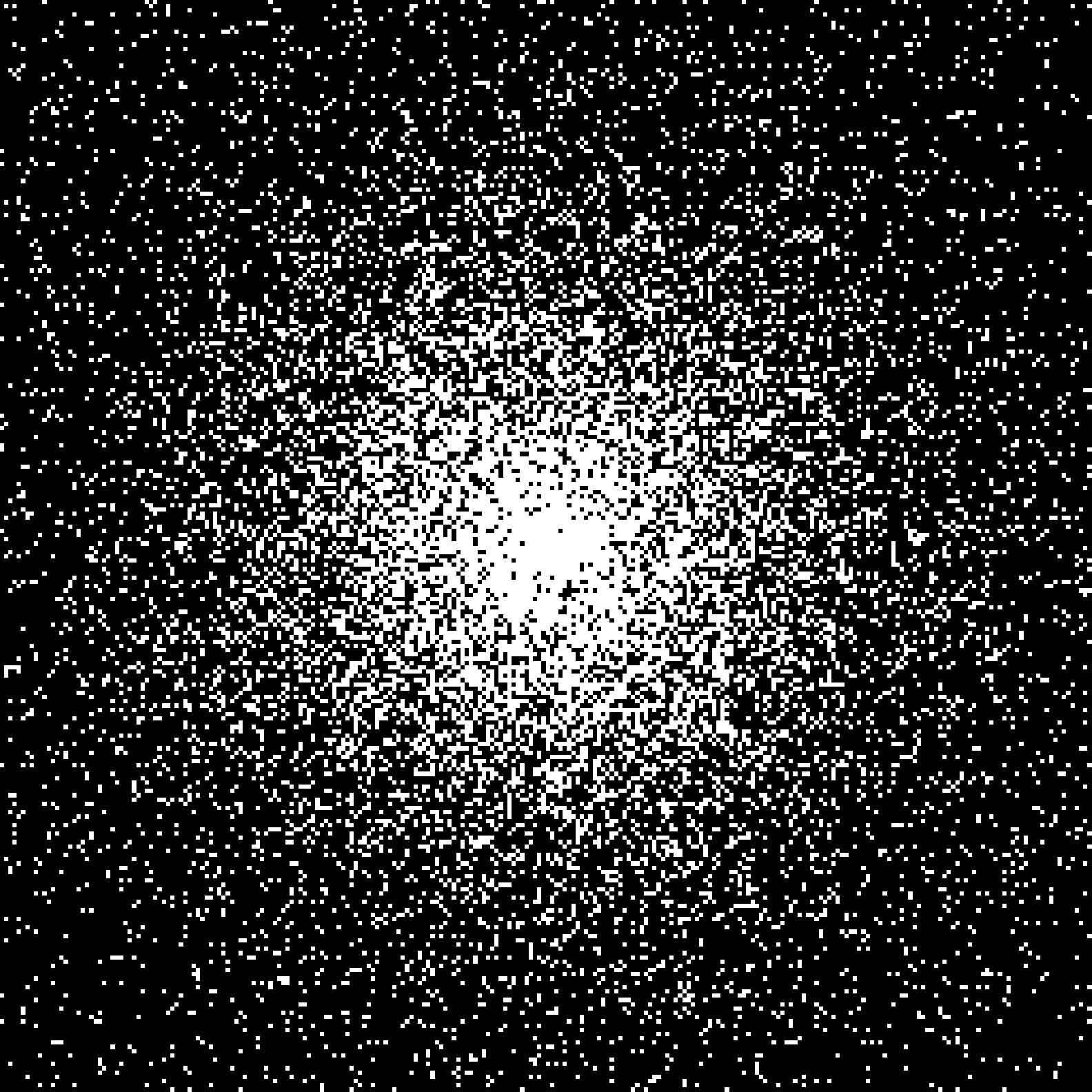}
\end{minipage}}\hspace{-0.05cm}
\subfloat[Proposed]{\label{SuperMarioRandomProposedError}\begin{minipage}{3cm}
\includegraphics[width=3cm]{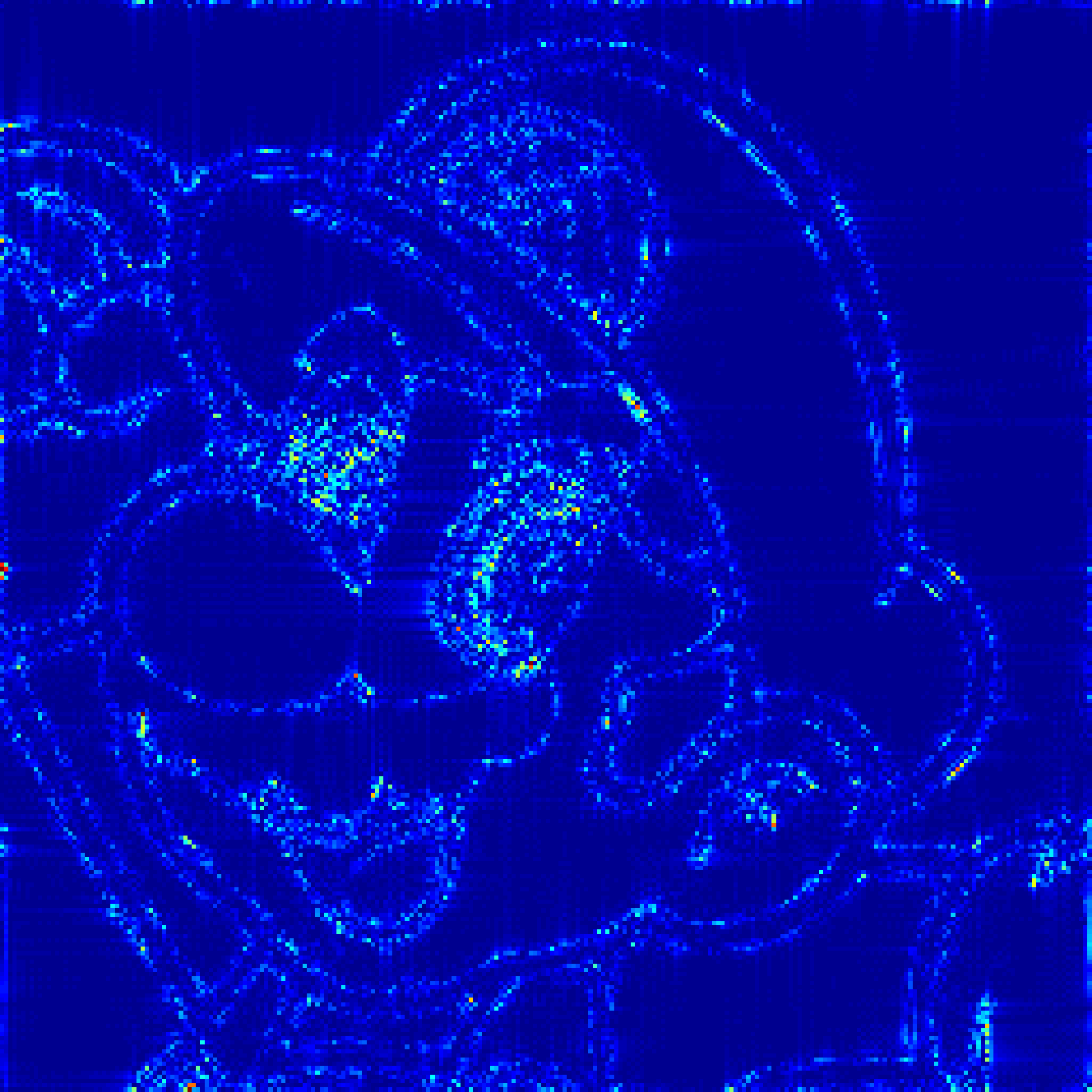}
\end{minipage}}\hspace{-0.05cm}
\subfloat[LSLP]{\label{SuperMarioRandomLSLPError}\begin{minipage}{3cm}
\includegraphics[width=3cm]{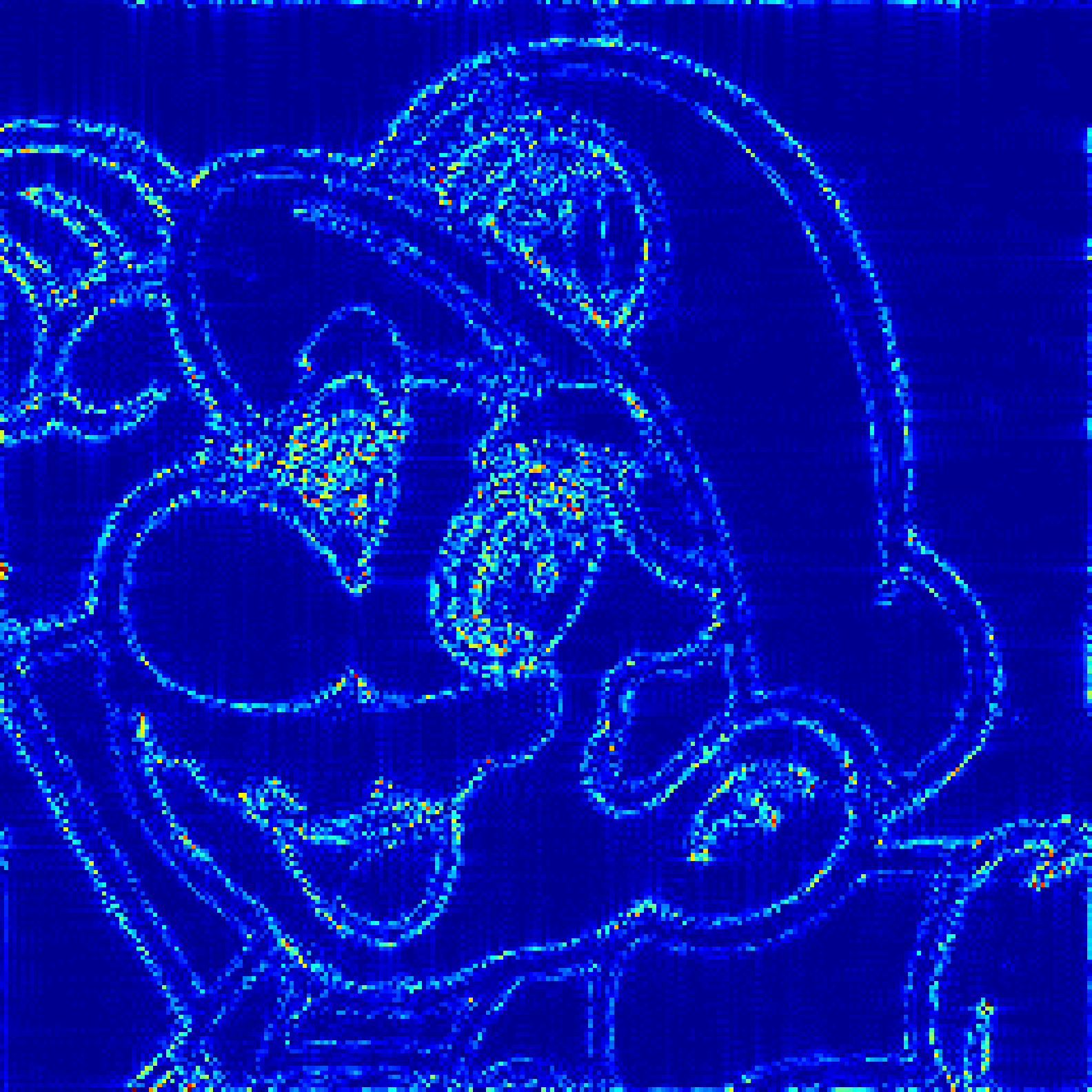}
\end{minipage}}\hspace{-0.05cm}
\subfloat[LRHDDTF]{\label{SuperMarioRandomLRHDDTFError}\begin{minipage}{3cm}
\includegraphics[width=3cm]{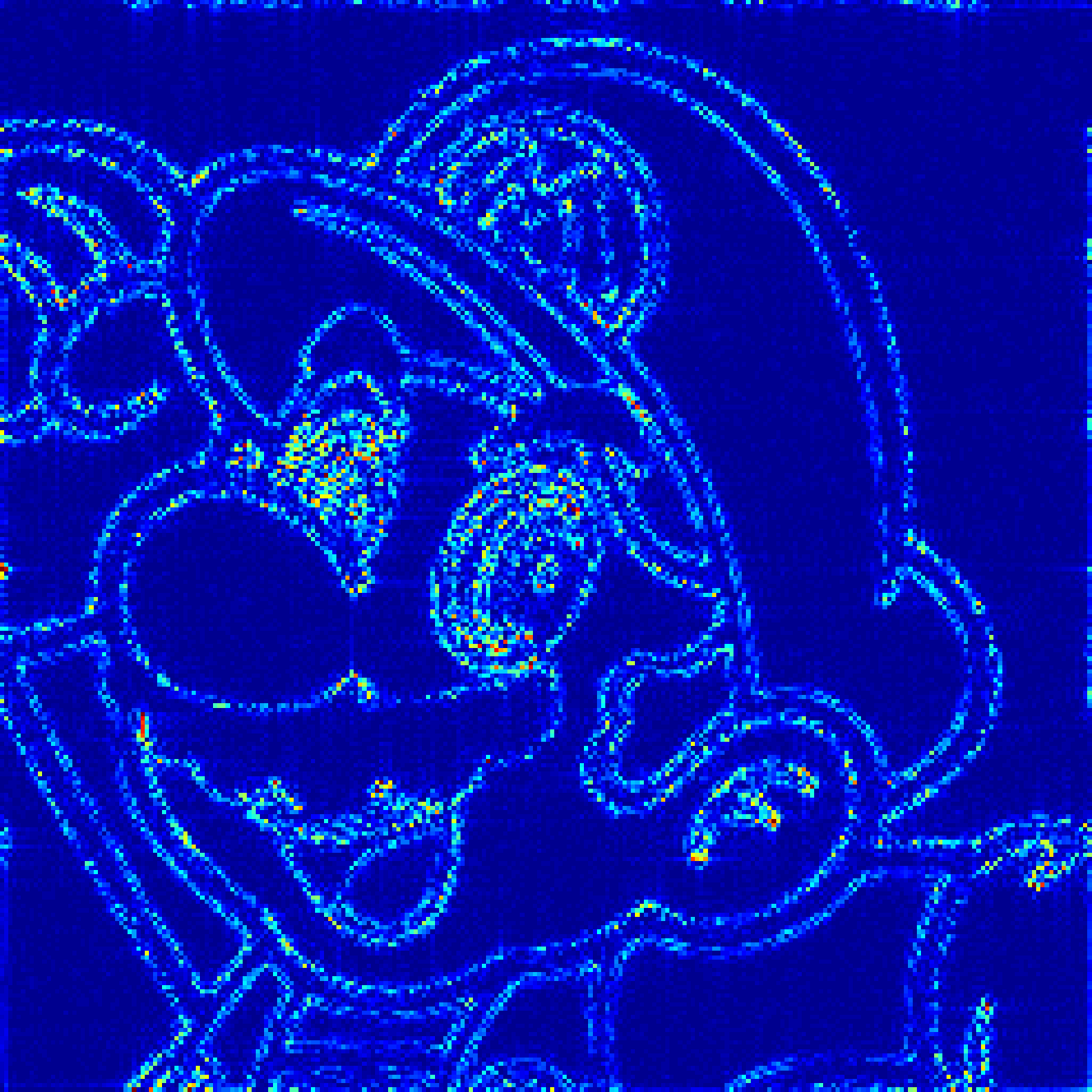}
\end{minipage}}\vspace{-0.25cm}
\subfloat[Schatten $0$]{\label{SuperMarioRandomGIRAFError}\begin{minipage}{3cm}
\includegraphics[width=3cm]{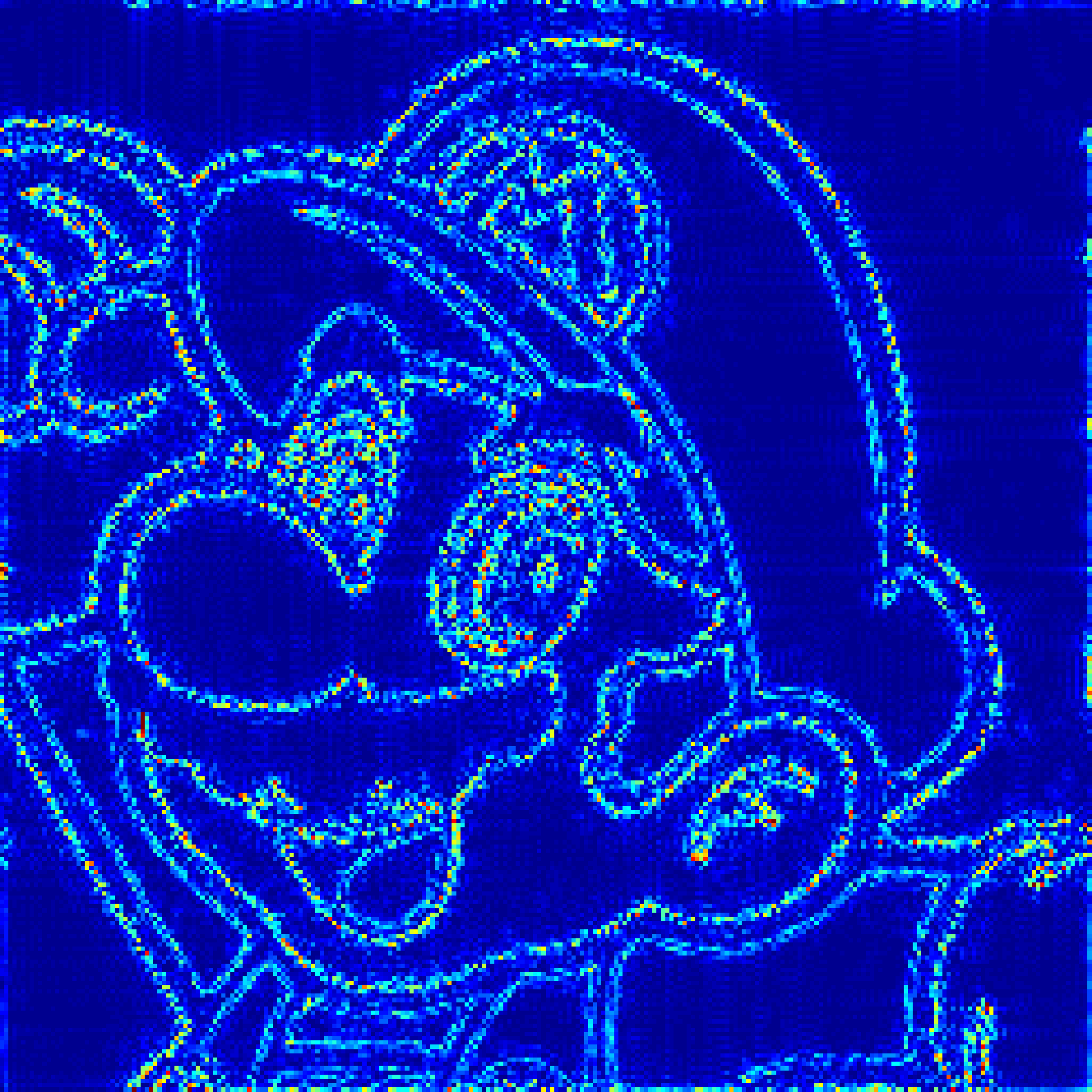}
\end{minipage}}\hspace{-0.05cm}
\subfloat[TV]{\label{SuperMarioRandomTVError}\begin{minipage}{3cm}
\includegraphics[width=3cm]{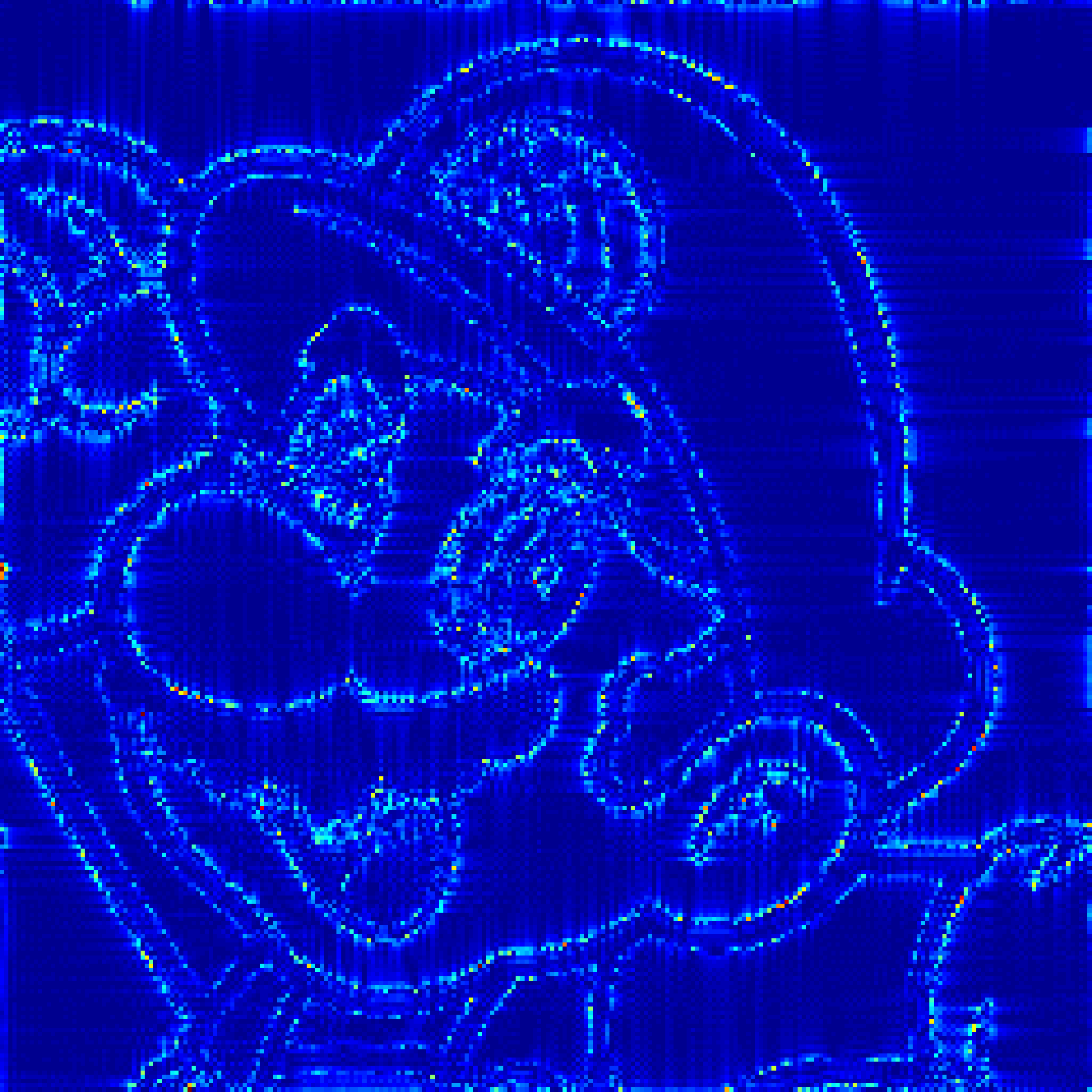}
\end{minipage}}\hspace{-0.05cm}
\subfloat[Haar]{\label{SuperMarioRandomHaarError}\begin{minipage}{3cm}
\includegraphics[width=3cm]{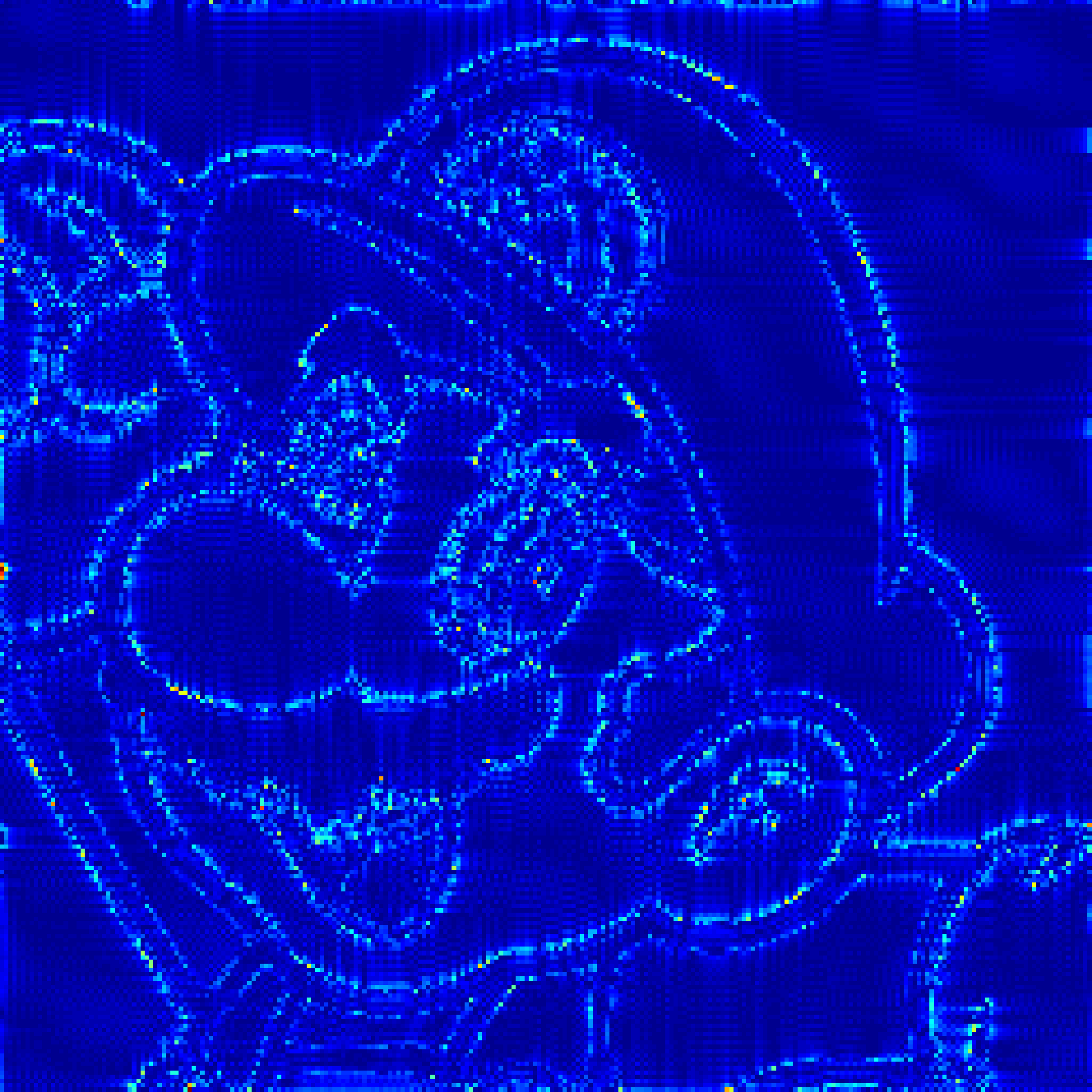}
\end{minipage}}\hspace{-0.05cm}
\subfloat[DDTF]{\label{SuperMarioRandomDDTFError}\begin{minipage}{3cm}
\includegraphics[width=3cm]{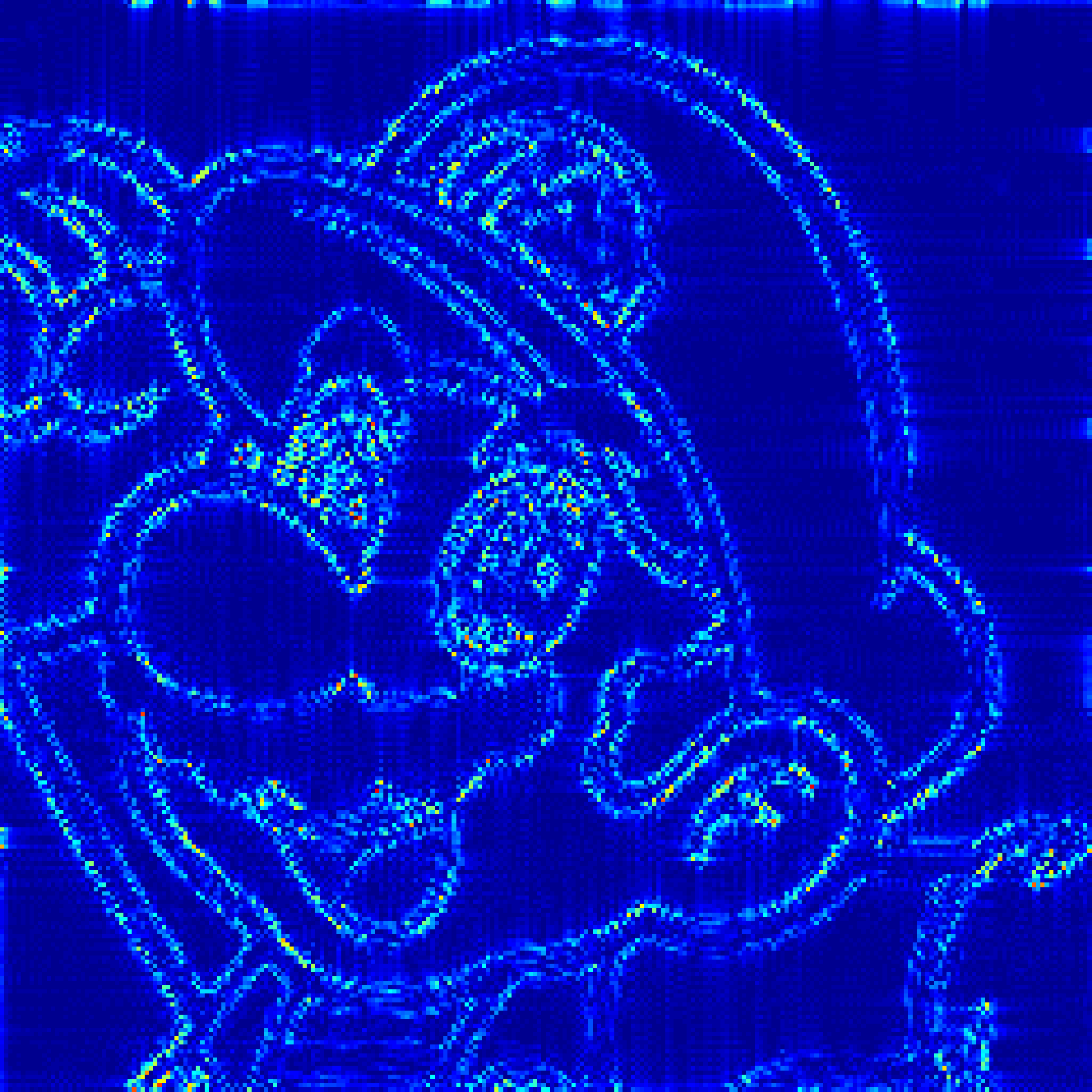}
\end{minipage}}
\caption{Error maps of \cref{SuperMarioRandomResults}.}\label{SuperMarioRandomError}
\end{figure}

\begin{figure}[t]
\centering
\subfloat[Ref.]{\label{SuperMarioILFOriginal}\begin{minipage}{3cm}
\includegraphics[width=3cm]{SuperMarioOriginal.pdf}
\end{minipage}}\hspace{-0.05cm}
\subfloat[Proposed]{\label{SuperMarioILFProposed}\begin{minipage}{3cm}
\includegraphics[width=3cm]{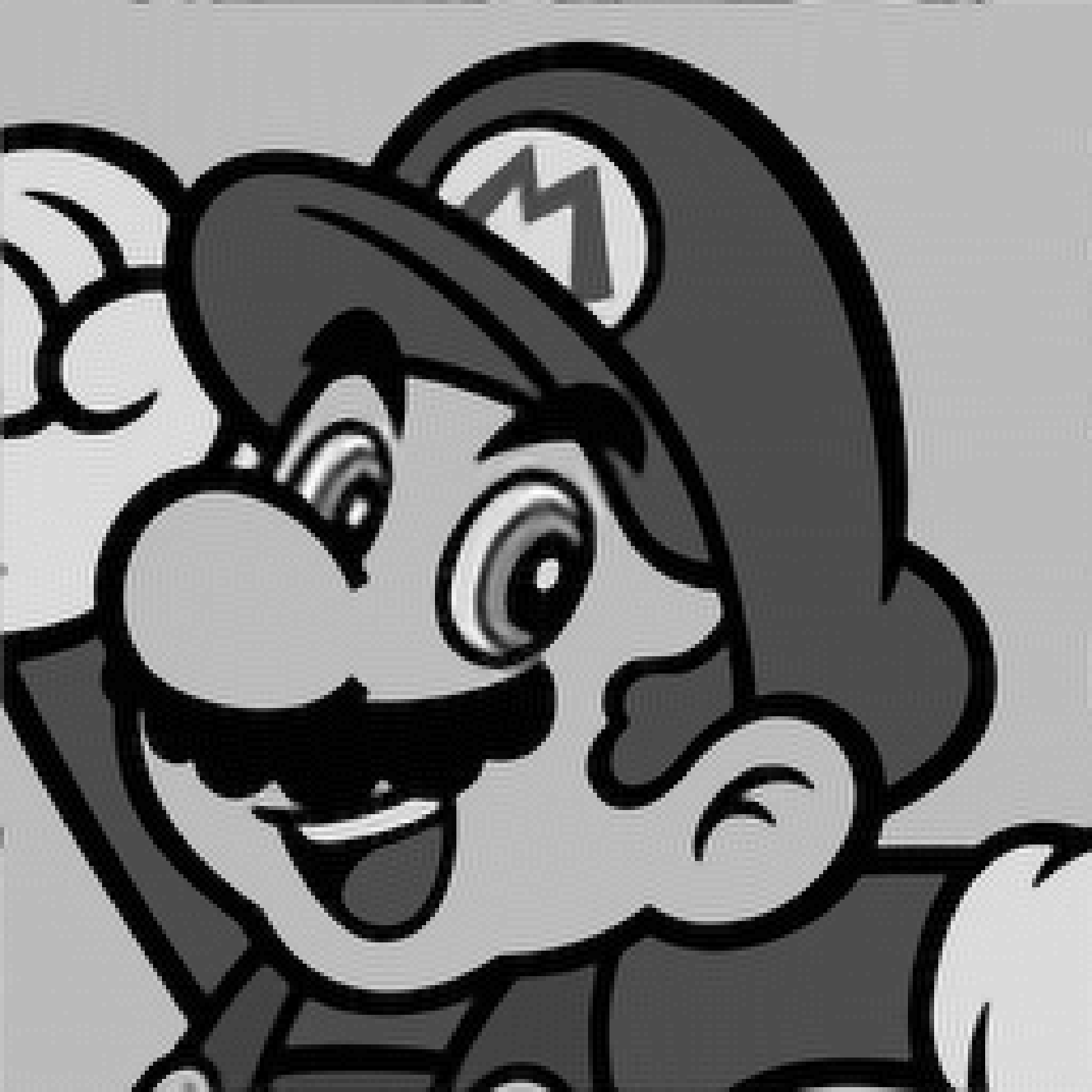}
\end{minipage}}\hspace{-0.05cm}
\subfloat[LSLP]{\label{SuperMarioILFLSLP}\begin{minipage}{3cm}
\includegraphics[width=3cm]{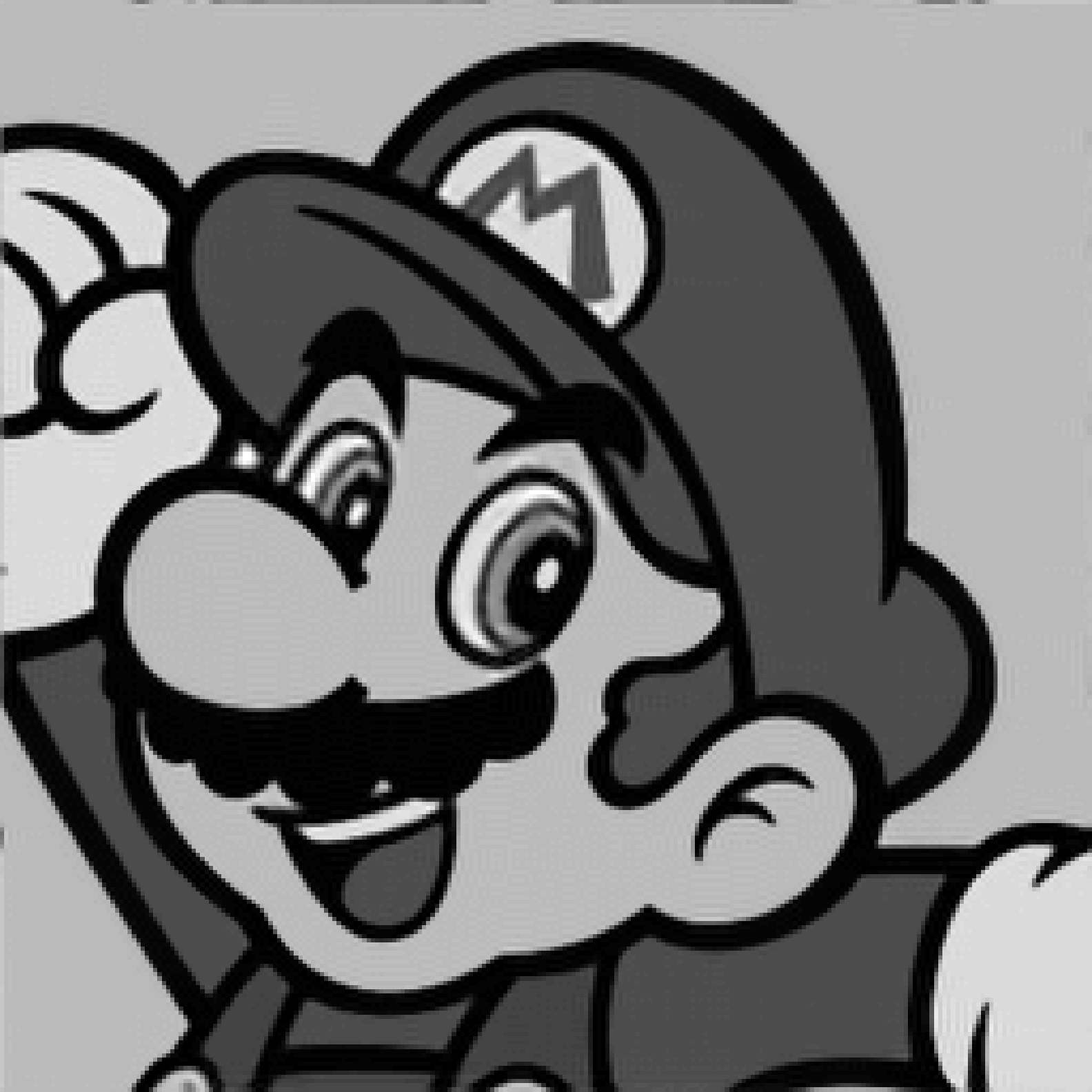}
\end{minipage}}\hspace{-0.05cm}
\subfloat[LRHDDTF]{\label{SuperMarioILFLRHDDTF}\begin{minipage}{3cm}
\includegraphics[width=3cm]{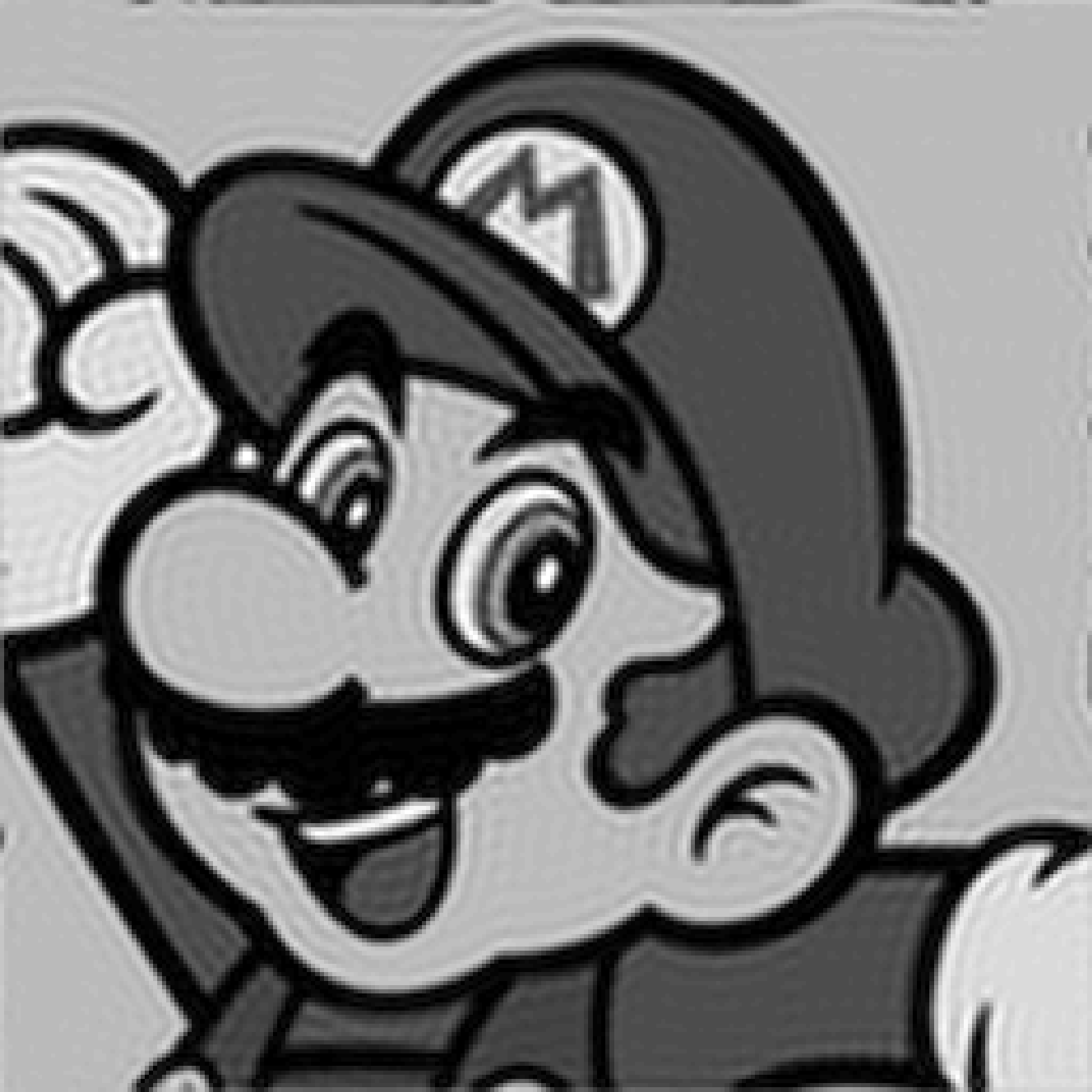}
\end{minipage}}\vspace{-0.25cm}
\subfloat[Schatten $0$]{\label{SuperMarioILFGIRAF}\begin{minipage}{3cm}
\includegraphics[width=3cm]{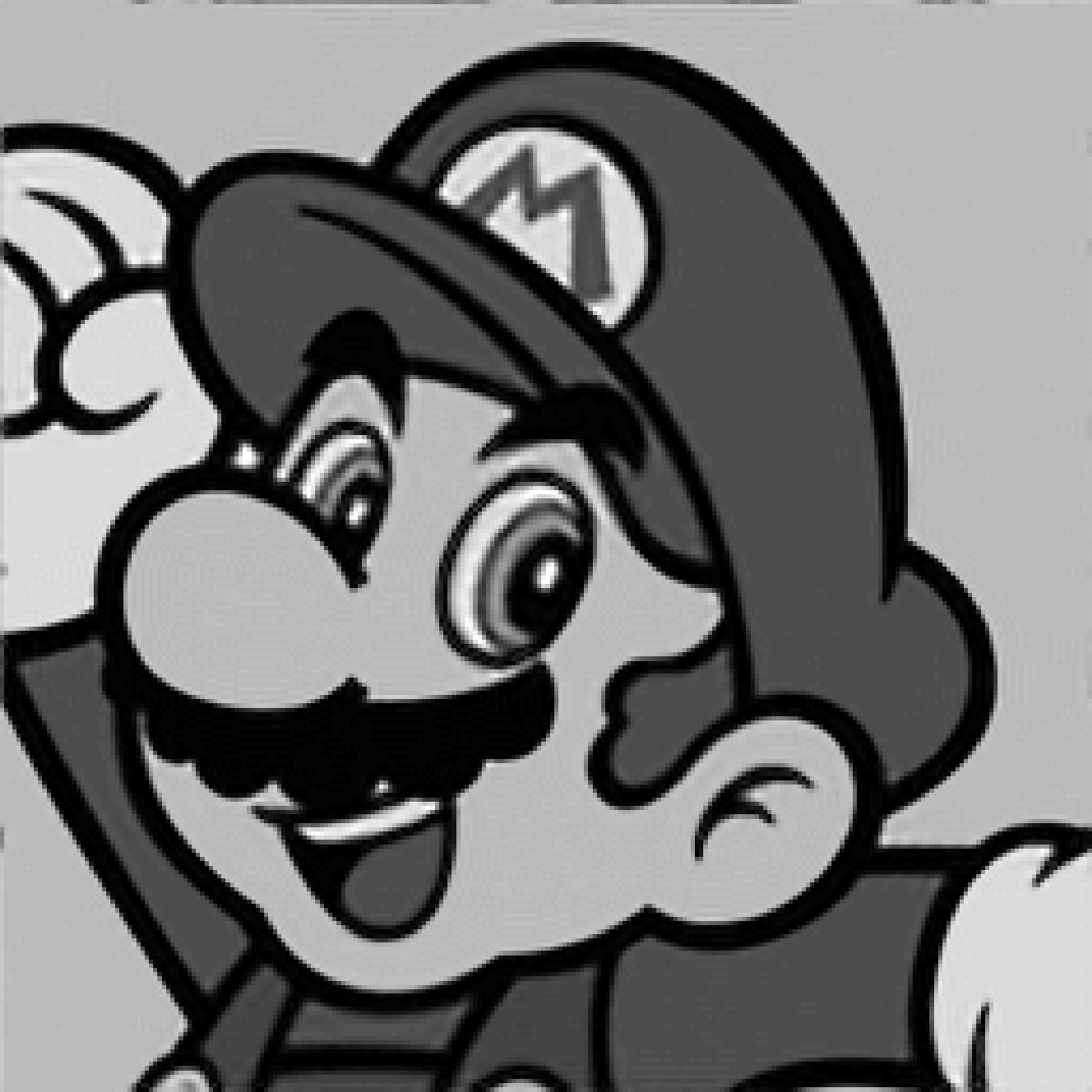}
\end{minipage}}\hspace{-0.05cm}
\subfloat[TV]{\label{SuperMarioILFTV}\begin{minipage}{3cm}
\includegraphics[width=3cm]{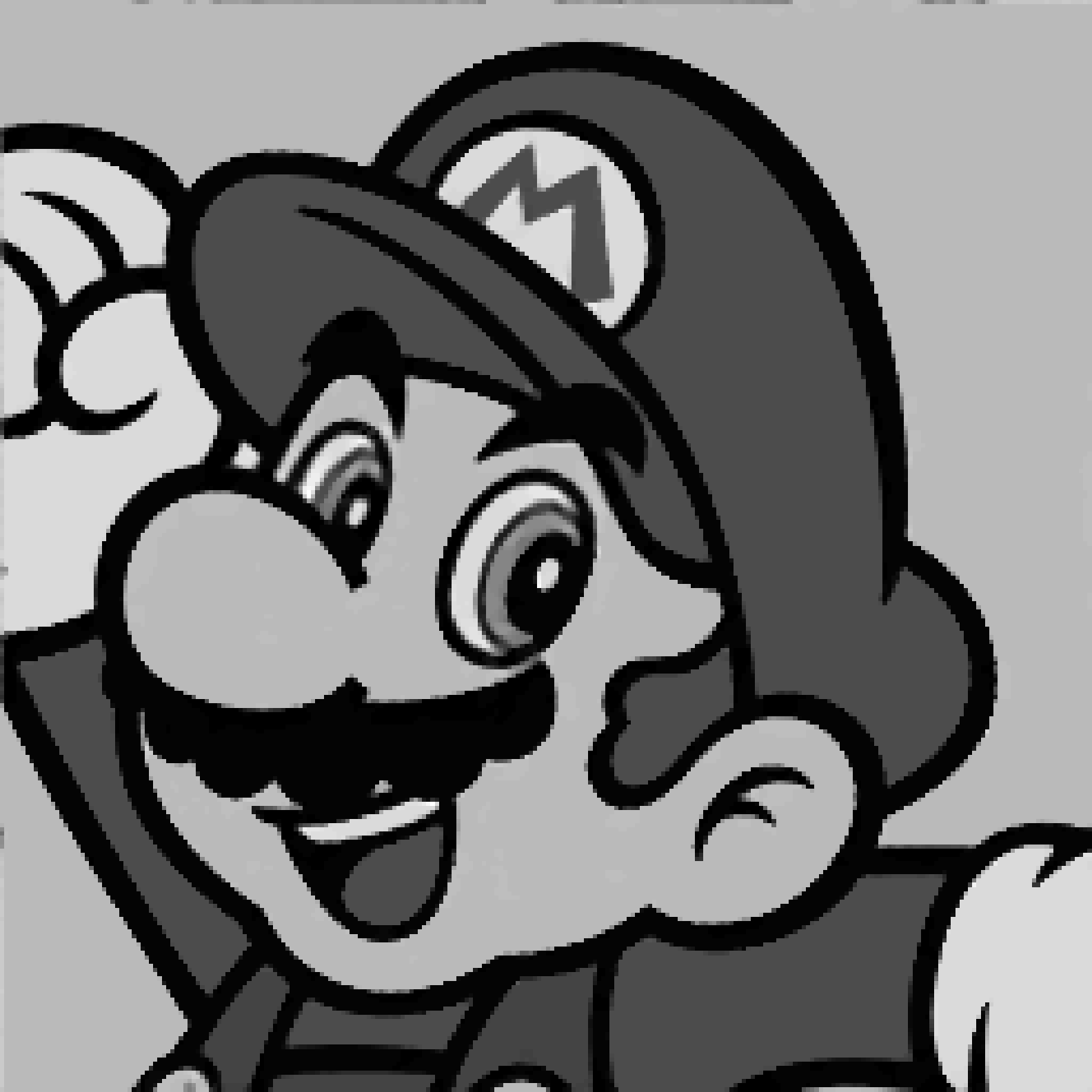}
\end{minipage}}\hspace{-0.05cm}
\subfloat[Haar]{\label{SuperMarioILFHaar}\begin{minipage}{3cm}
\includegraphics[width=3cm]{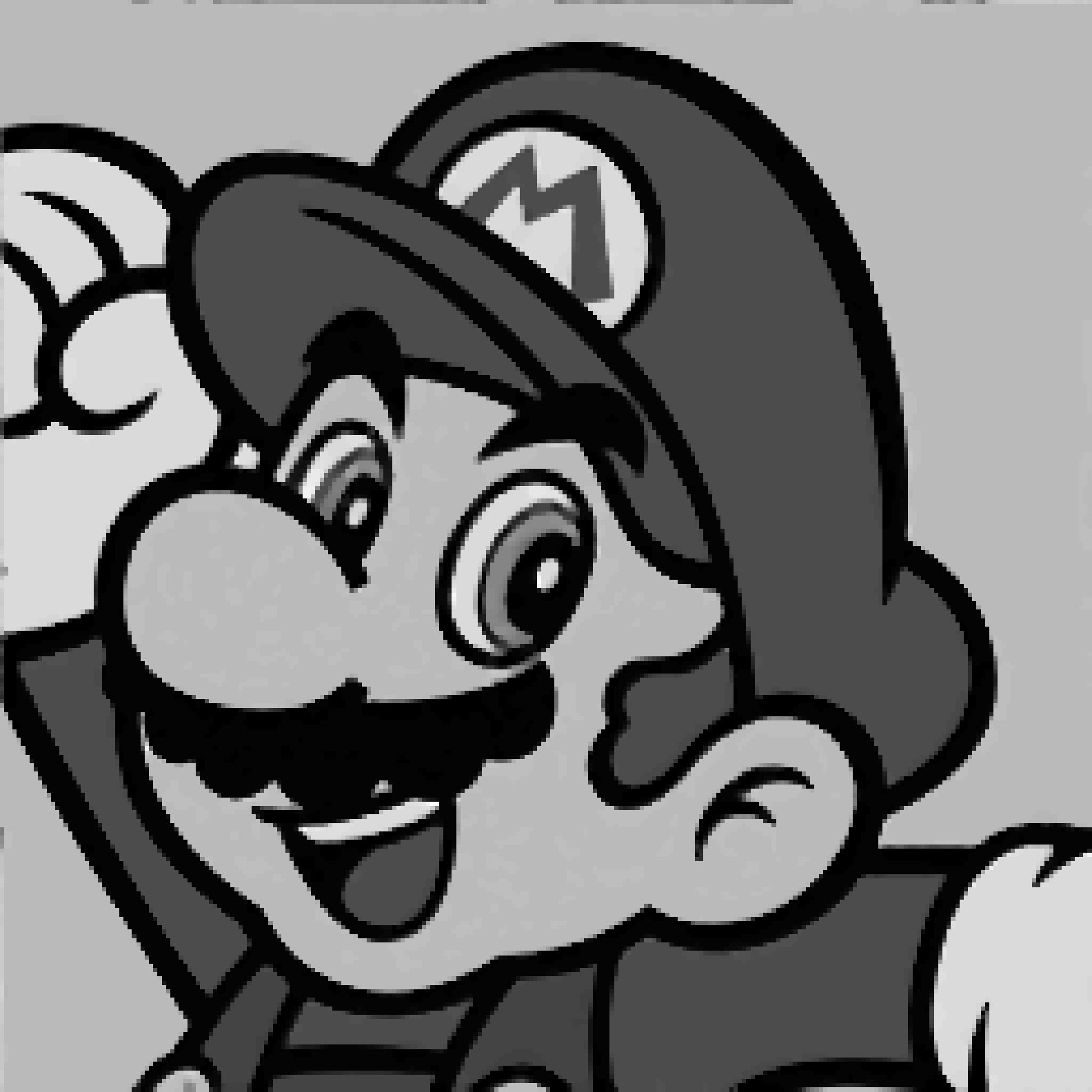}
\end{minipage}}\hspace{-0.05cm}
\subfloat[DDTF]{\label{SuperMarioILFDDTF}\begin{minipage}{3cm}
\includegraphics[width=3cm]{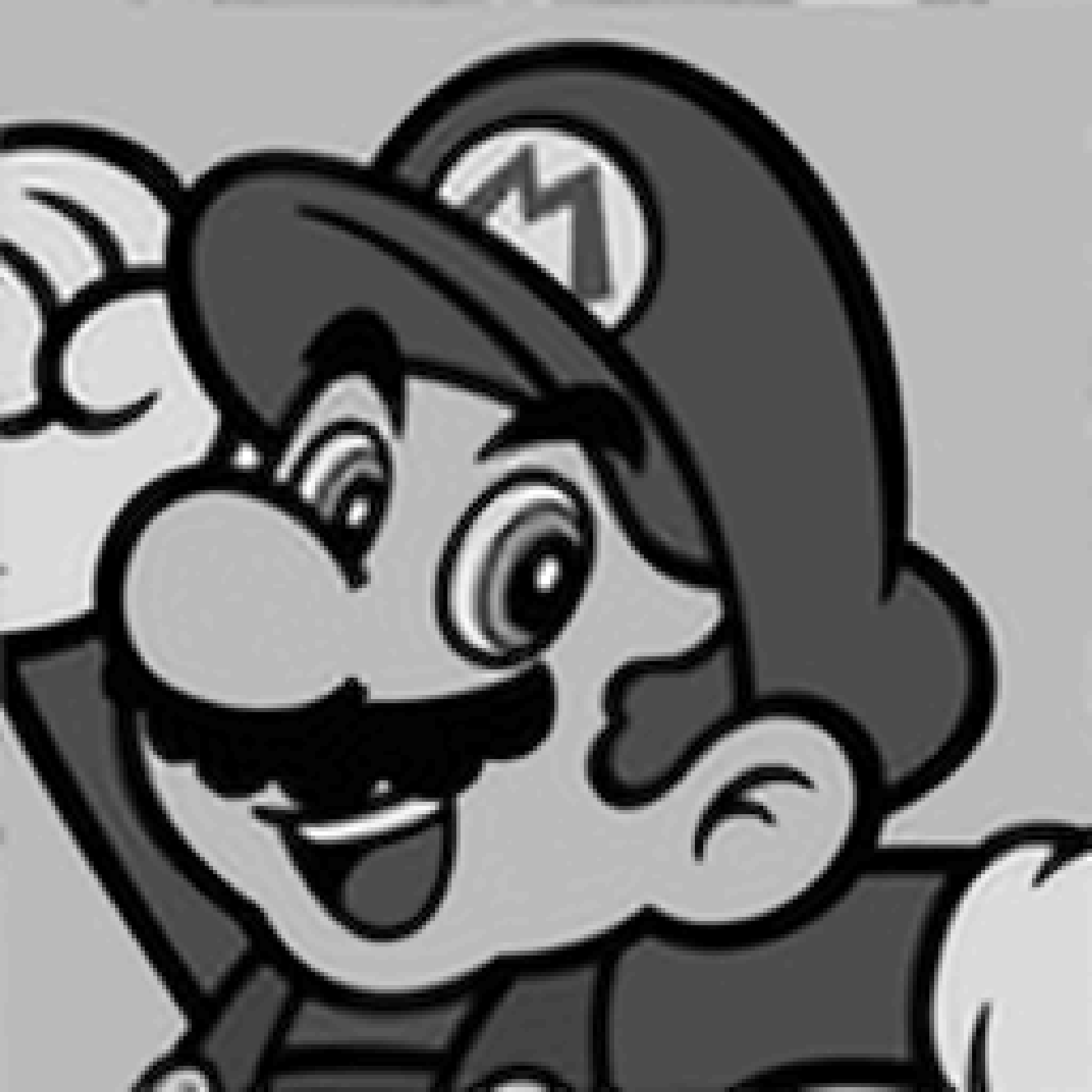}
\end{minipage}}
\caption{Visual comparison of ``Super Mario'' for ideal low-pass filter deconvolution.}\label{SuperMarioILFResults}
\end{figure}

\begin{figure}[t]
\centering
\subfloat[TF]{\label{SuperMarioILFTF}\begin{minipage}{3cm}
\includegraphics[width=3cm]{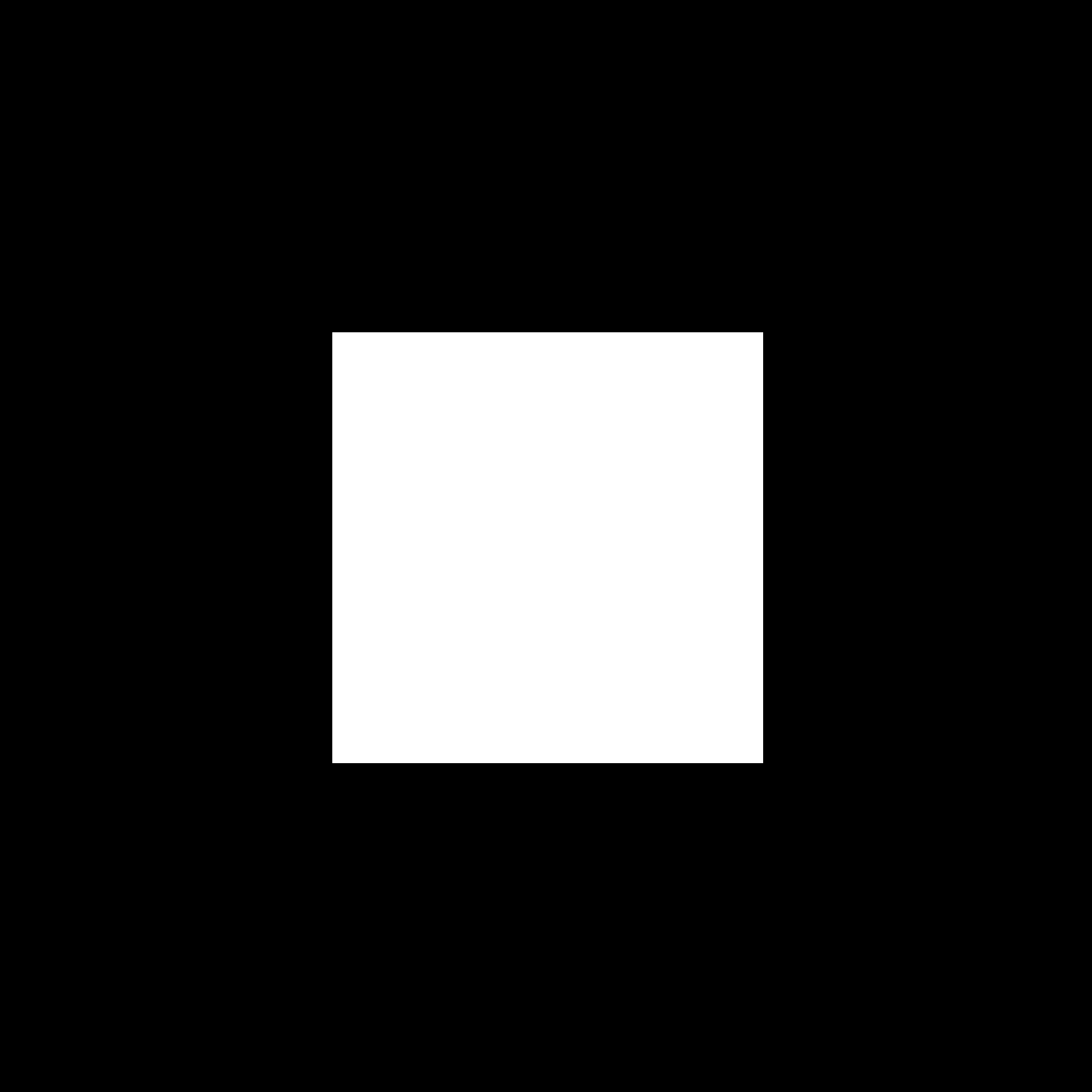}
\end{minipage}}\hspace{-0.05cm}
\subfloat[Proposed]{\label{SuperMarioILFProposedError}\begin{minipage}{3cm}
\includegraphics[width=3cm]{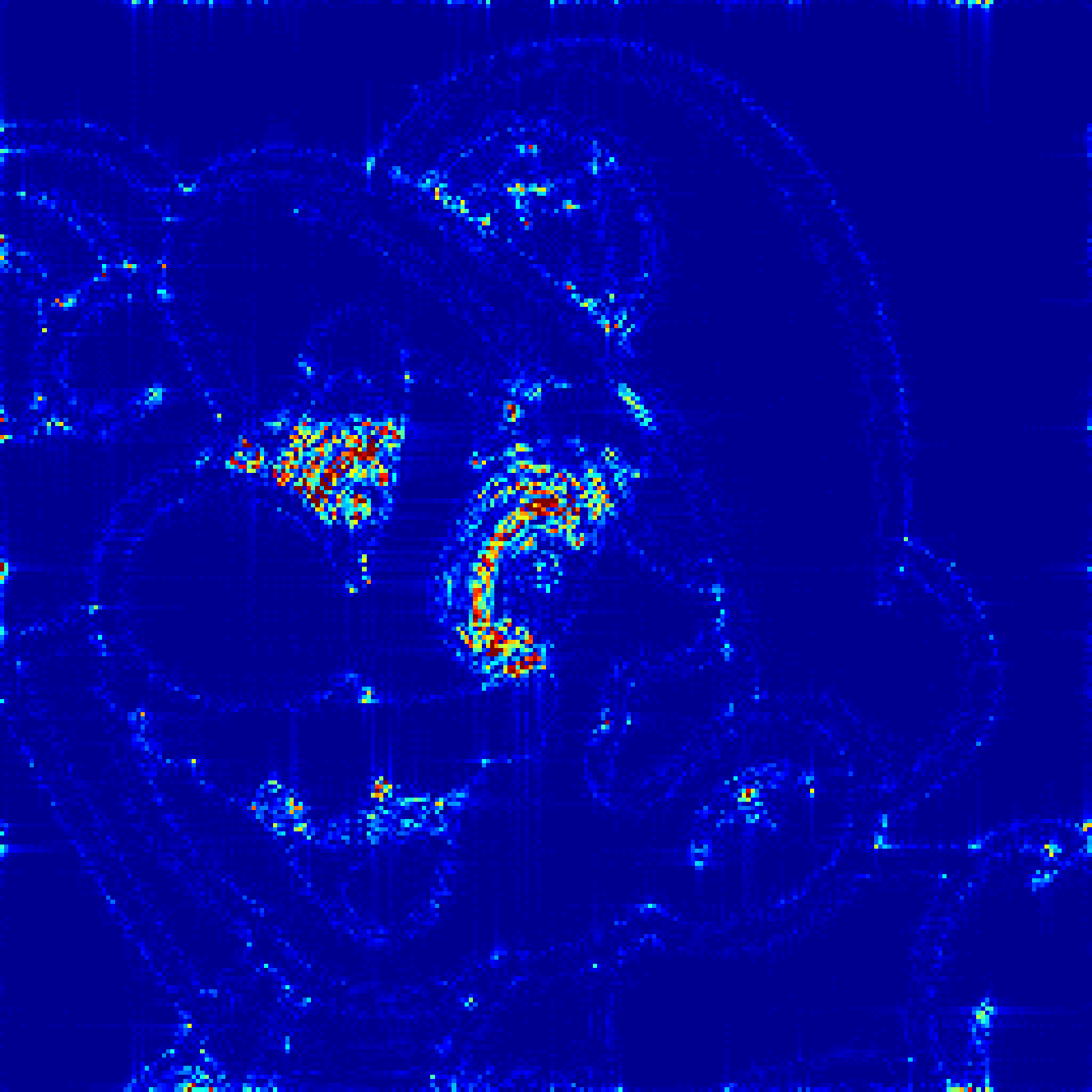}
\end{minipage}}\hspace{-0.05cm}
\subfloat[LSLP]{\label{SuperMarioILFLSLPError}\begin{minipage}{3cm}
\includegraphics[width=3cm]{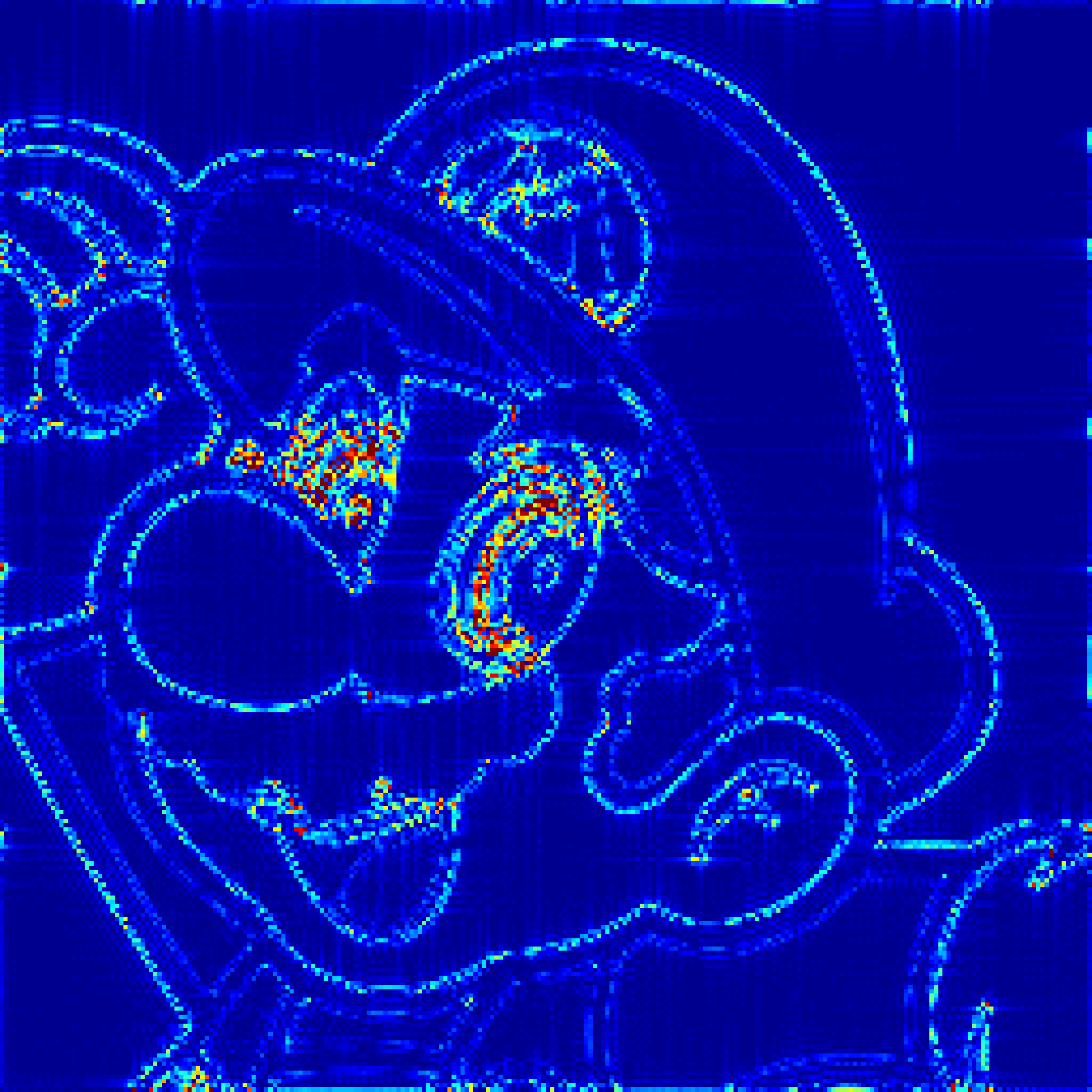}
\end{minipage}}\hspace{-0.05cm}
\subfloat[LRHDDTF]{\label{SuperMarioILFLRHDDTFError}\begin{minipage}{3cm}
\includegraphics[width=3cm]{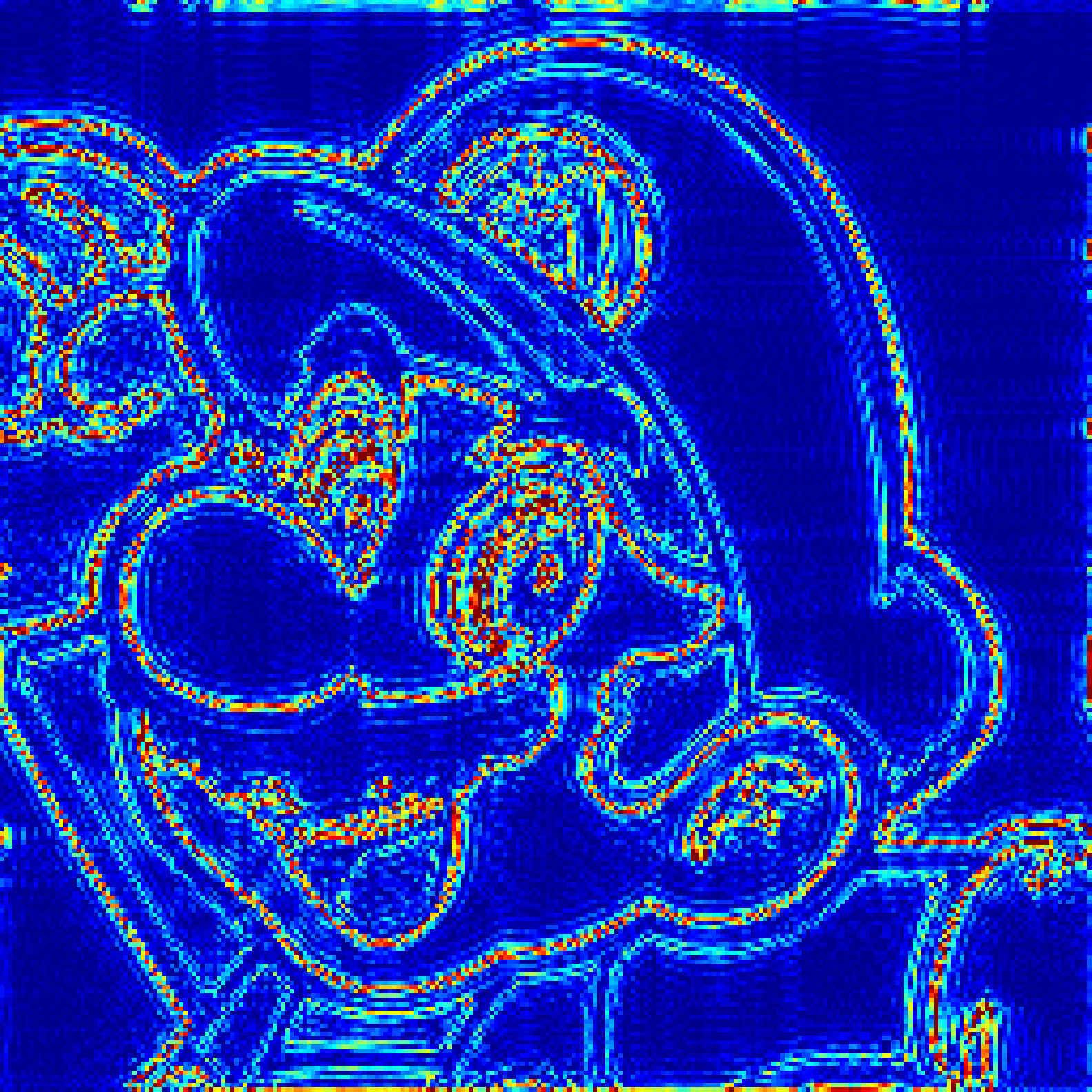}
\end{minipage}}\vspace{-0.25cm}
\subfloat[Schatten $0$]{\label{SuperMarioILFGIRAFError}\begin{minipage}{3cm}
\includegraphics[width=3cm]{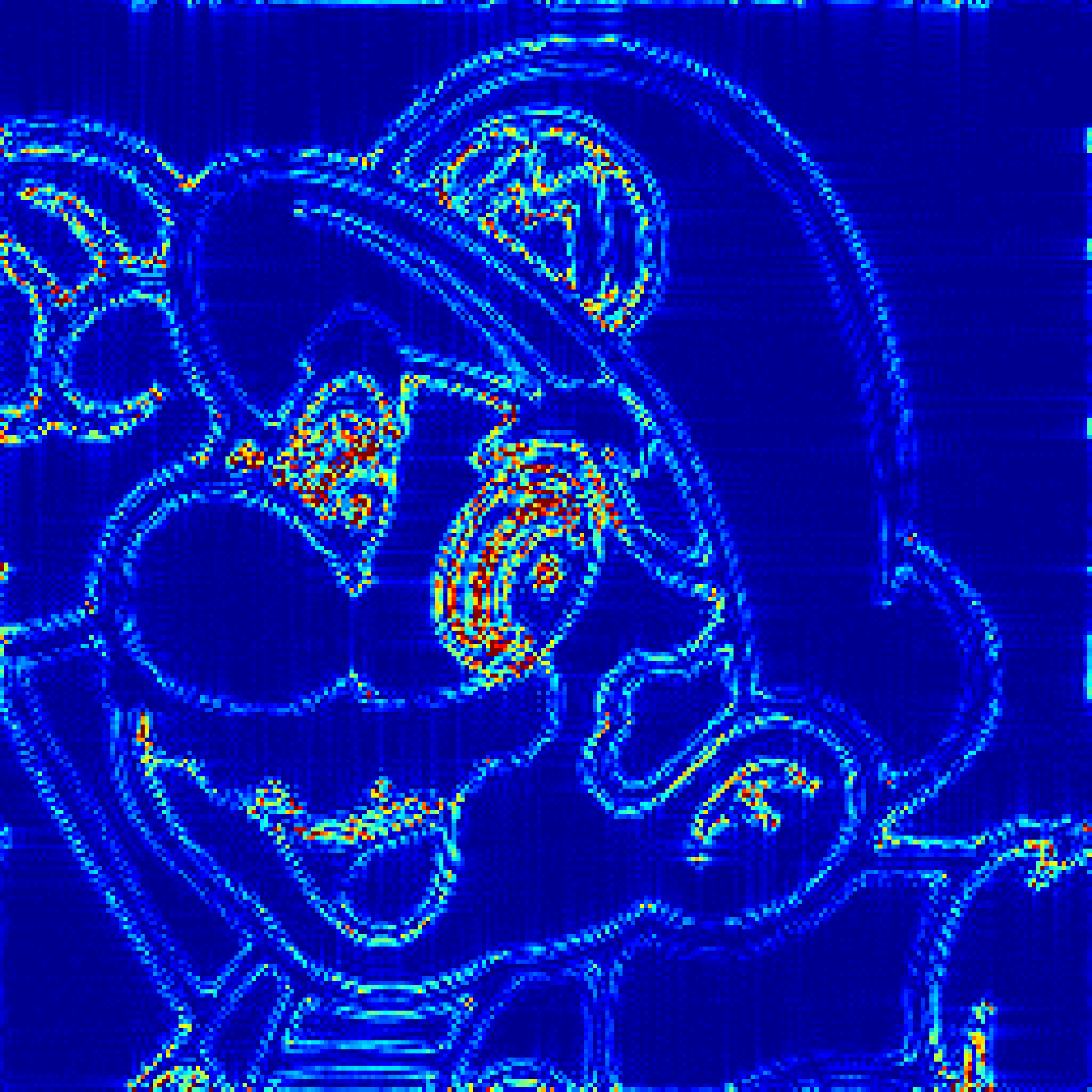}
\end{minipage}}\hspace{-0.05cm}
\subfloat[TV]{\label{SuperMarioILFTVError}\begin{minipage}{3cm}
\includegraphics[width=3cm]{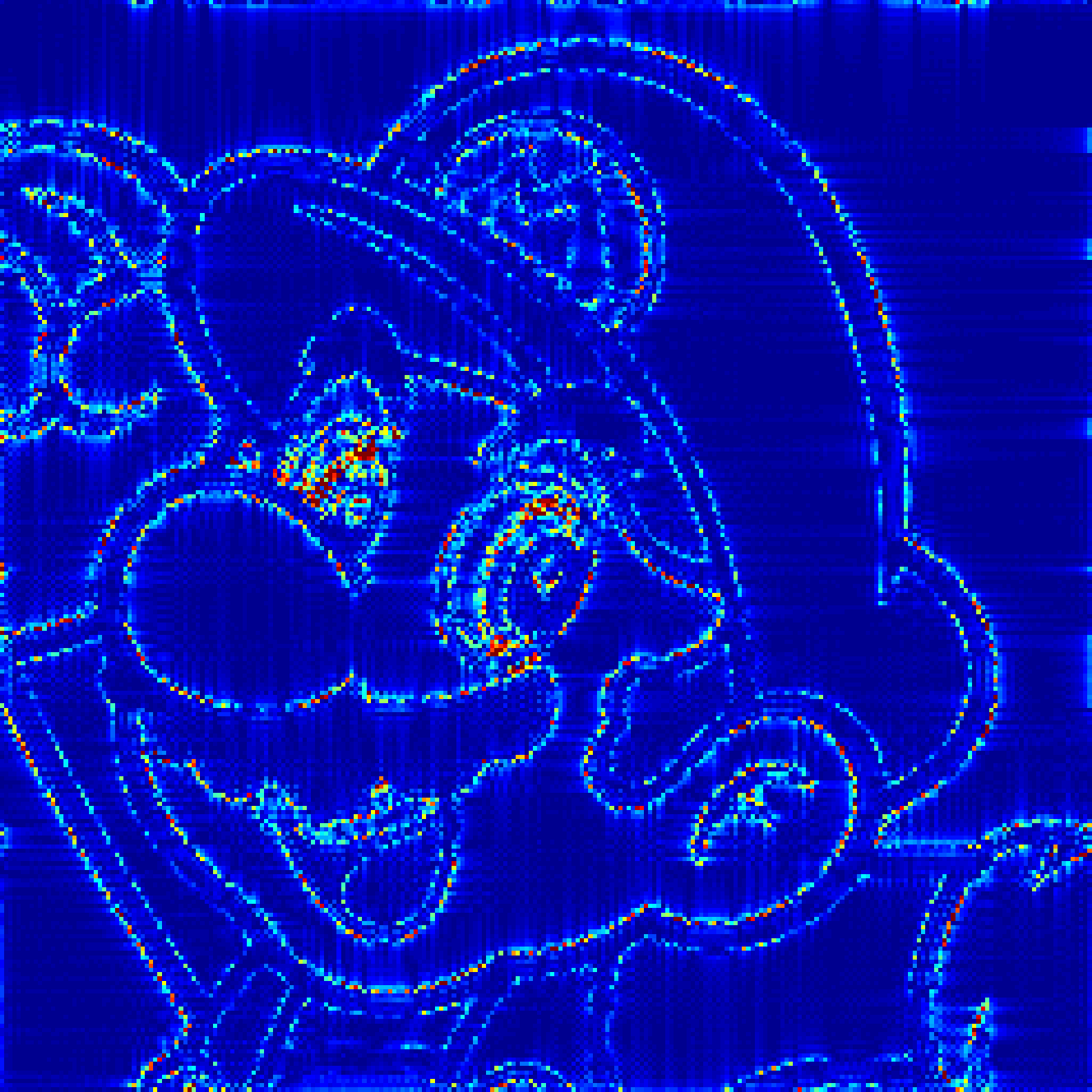}
\end{minipage}}\hspace{-0.05cm}
\subfloat[Haar]{\label{SuperMarioILFHaarError}\begin{minipage}{3cm}
\includegraphics[width=3cm]{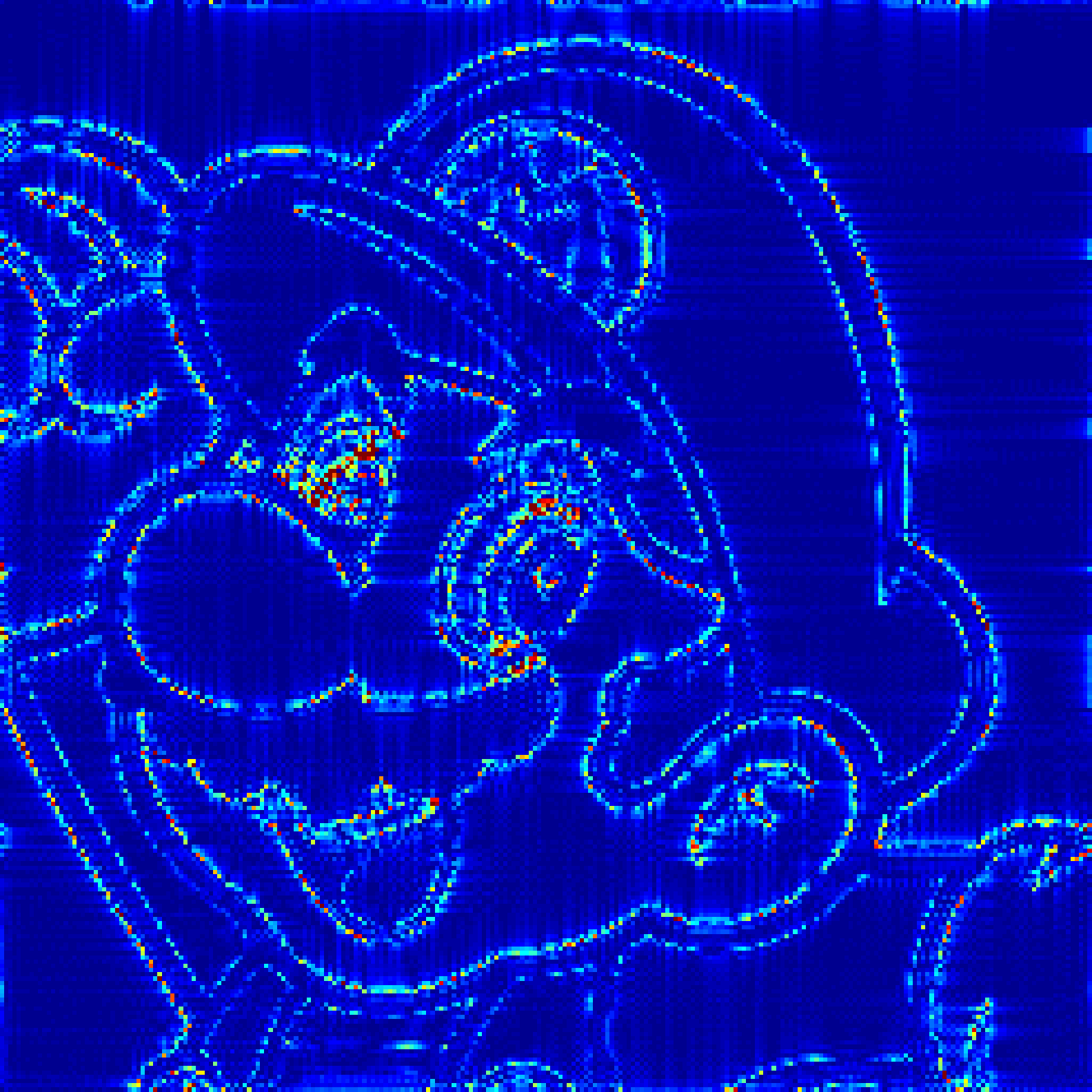}
\end{minipage}}\hspace{-0.05cm}
\subfloat[DDTF]{\label{SuperMarioILFDDTFError}\begin{minipage}{3cm}
\includegraphics[width=3cm]{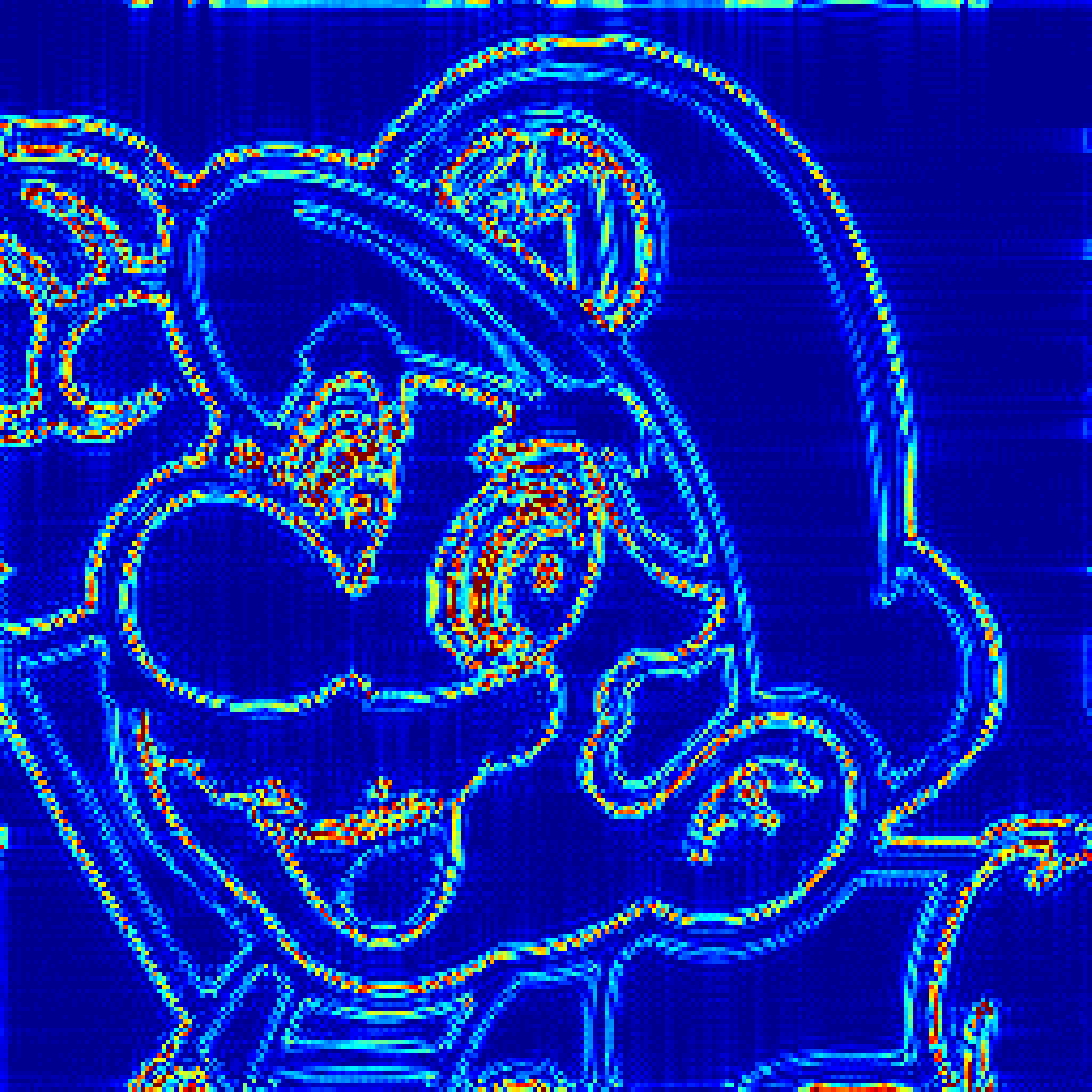}
\end{minipage}}
\caption{Error maps of \cref{SuperMarioILFResults}.}\label{SuperMarioILFError}
\end{figure}

\begin{figure}[t]
\centering
\subfloat[Ref.]{\label{MCRandomOriginal}\begin{minipage}{3cm}
\includegraphics[width=3cm]{MCOriginal.pdf}
\end{minipage}}\hspace{-0.05cm}
\subfloat[Proposed]{\label{MCRandomProposed}\begin{minipage}{3cm}
\includegraphics[width=3cm]{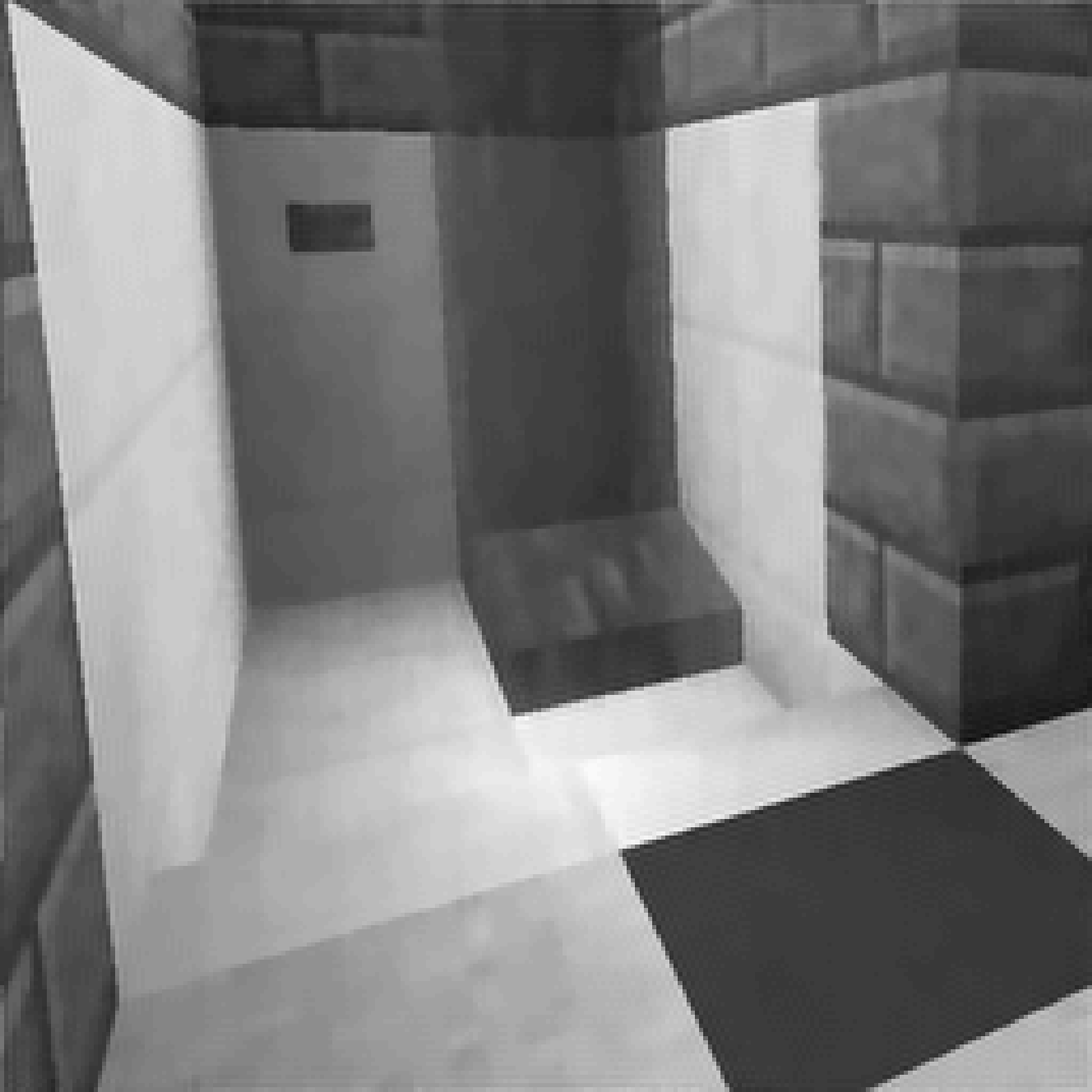}
\end{minipage}}\hspace{-0.05cm}
\subfloat[LSLP]{\label{MCRandomLSLP}\begin{minipage}{3cm}
\includegraphics[width=3cm]{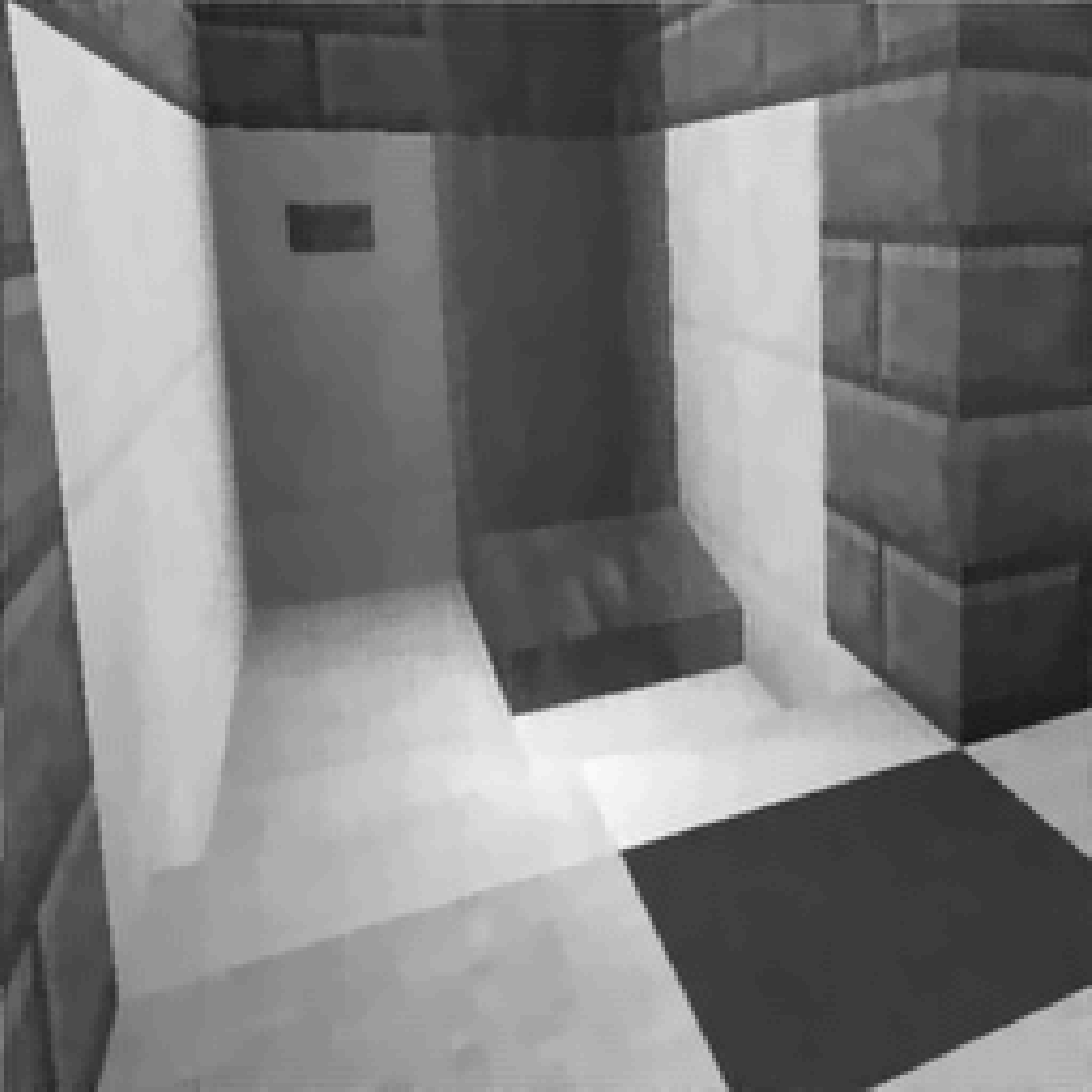}
\end{minipage}}\hspace{-0.05cm}
\subfloat[LRHDDTF]{\label{MCRandomLRHDDTF}\begin{minipage}{3cm}
\includegraphics[width=3cm]{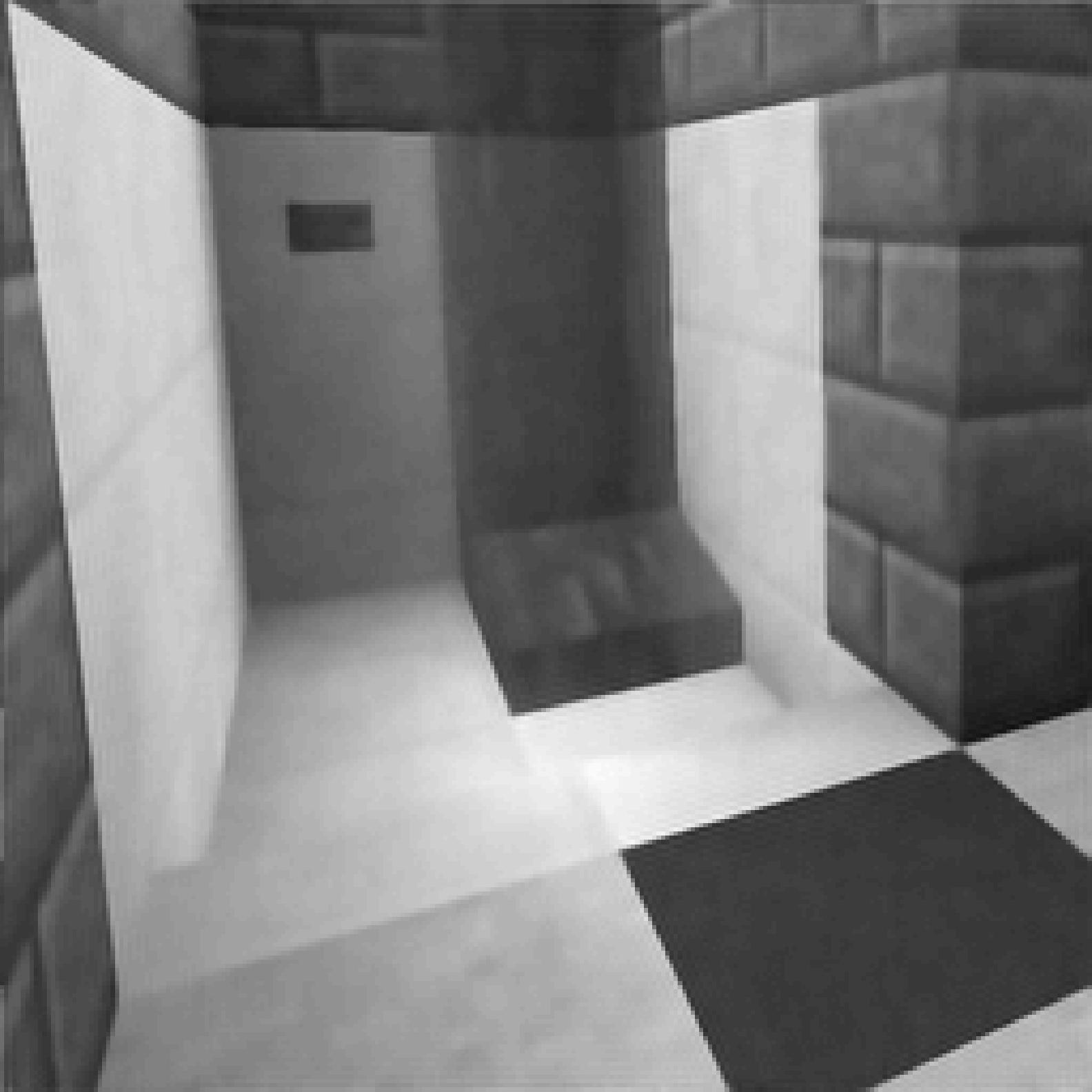}
\end{minipage}}\vspace{-0.25cm}
\subfloat[Schatten $0$]{\label{MCRandomGIRAF}\begin{minipage}{3cm}
\includegraphics[width=3cm]{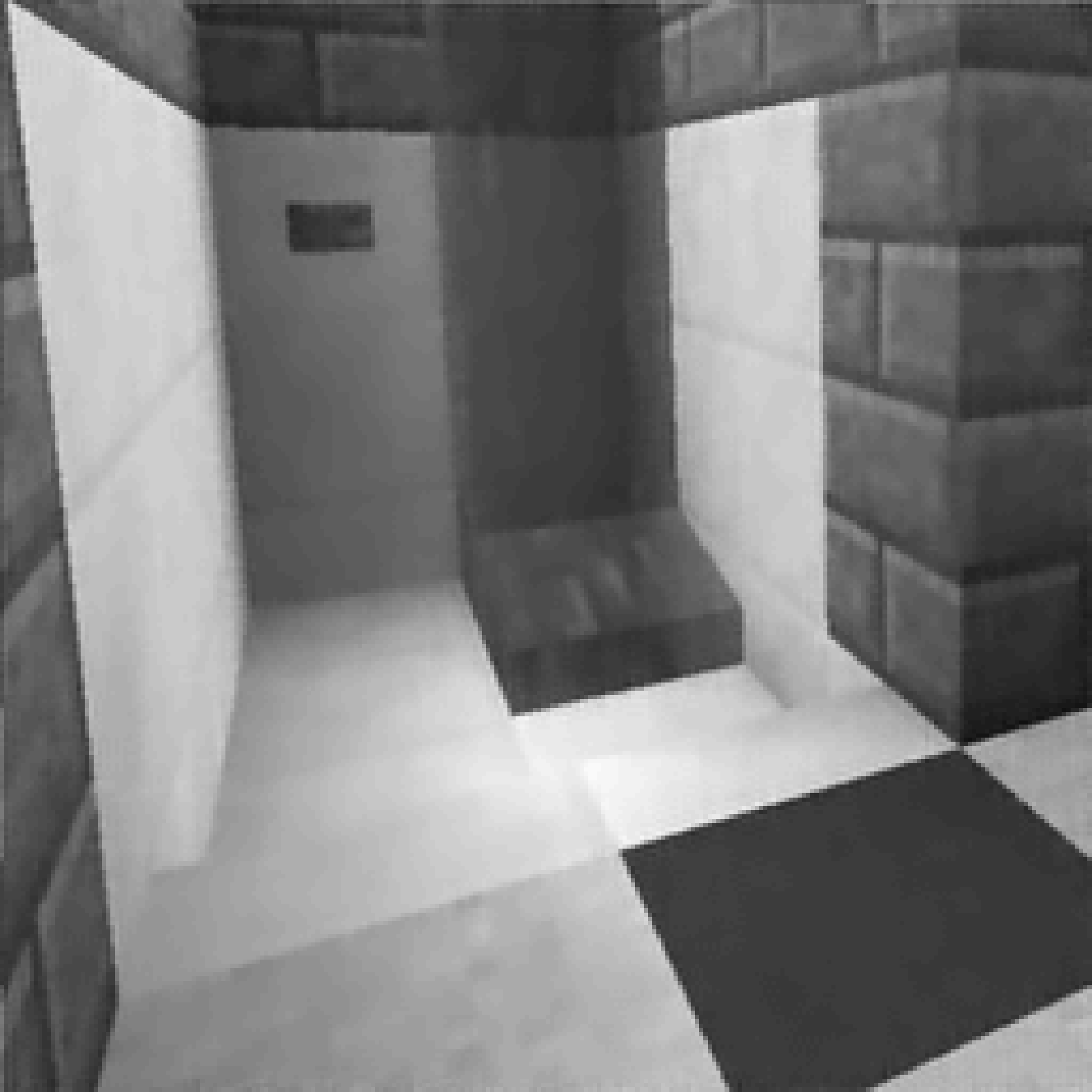}
\end{minipage}}\hspace{-0.05cm}
\subfloat[TV]{\label{MCRandomTV}\begin{minipage}{3cm}
\includegraphics[width=3cm]{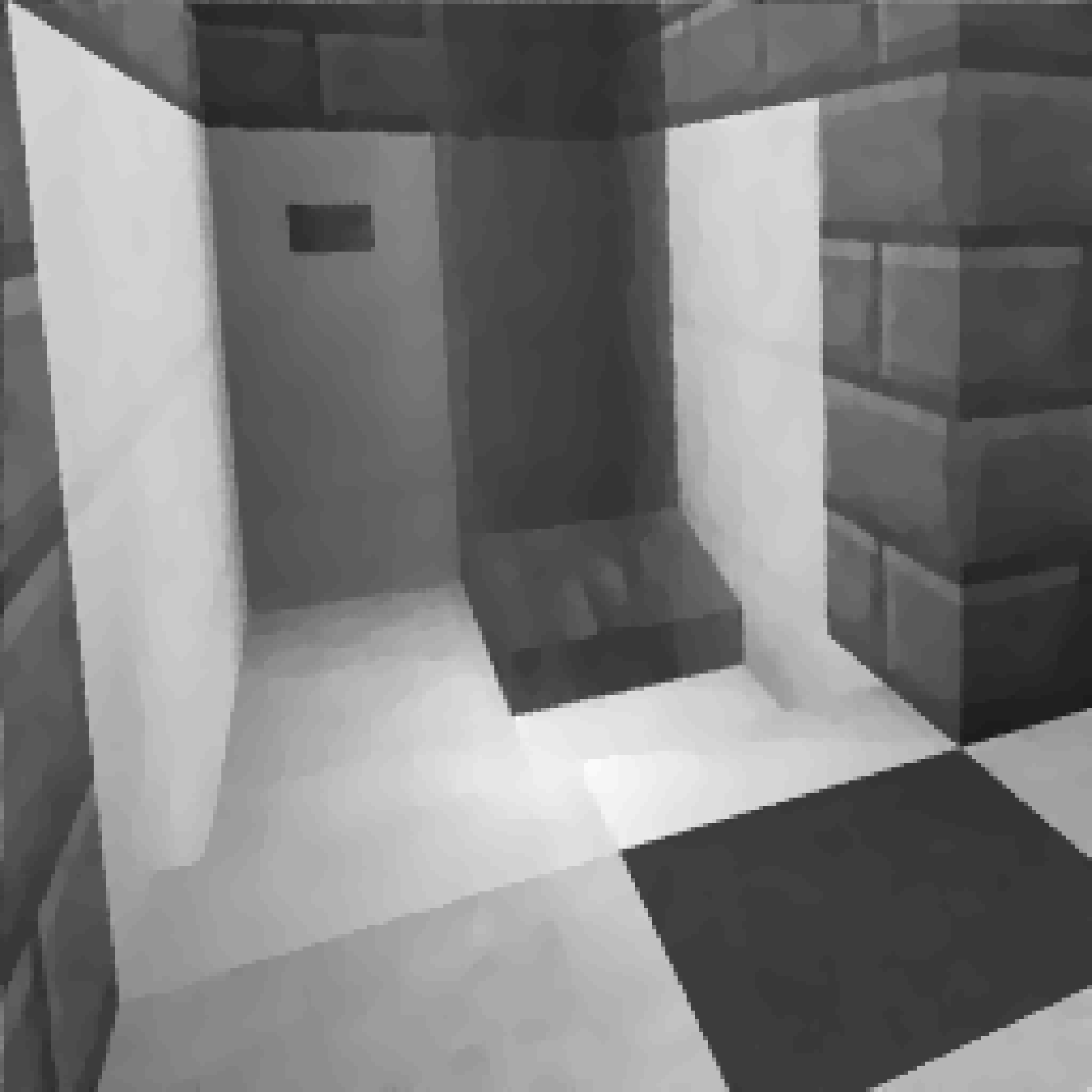}
\end{minipage}}\hspace{-0.05cm}
\subfloat[Haar]{\label{MCRandomHaar}\begin{minipage}{3cm}
\includegraphics[width=3cm]{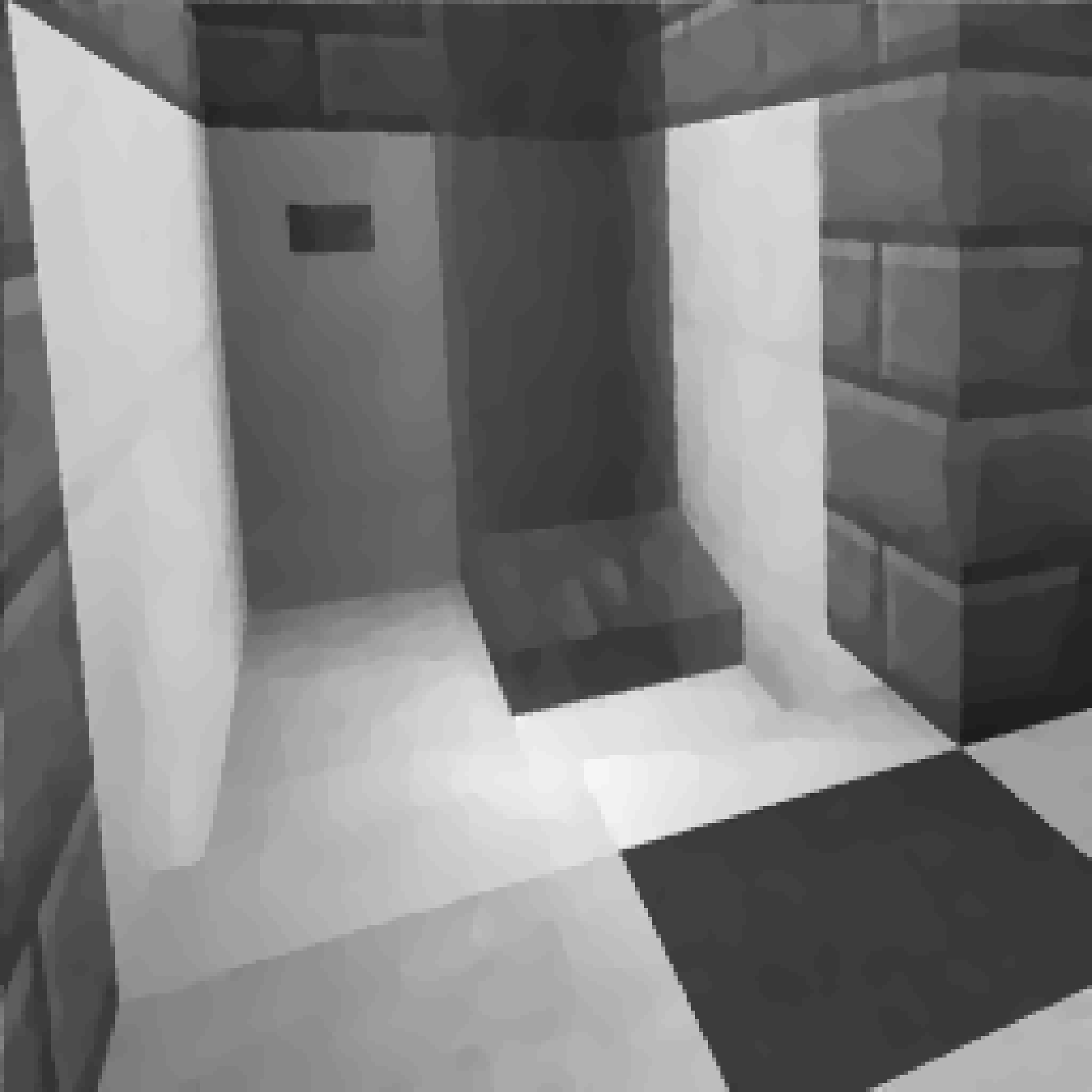}
\end{minipage}}\hspace{-0.05cm}
\subfloat[DDTF]{\label{MCRandomDDTF}\begin{minipage}{3cm}
\includegraphics[width=3cm]{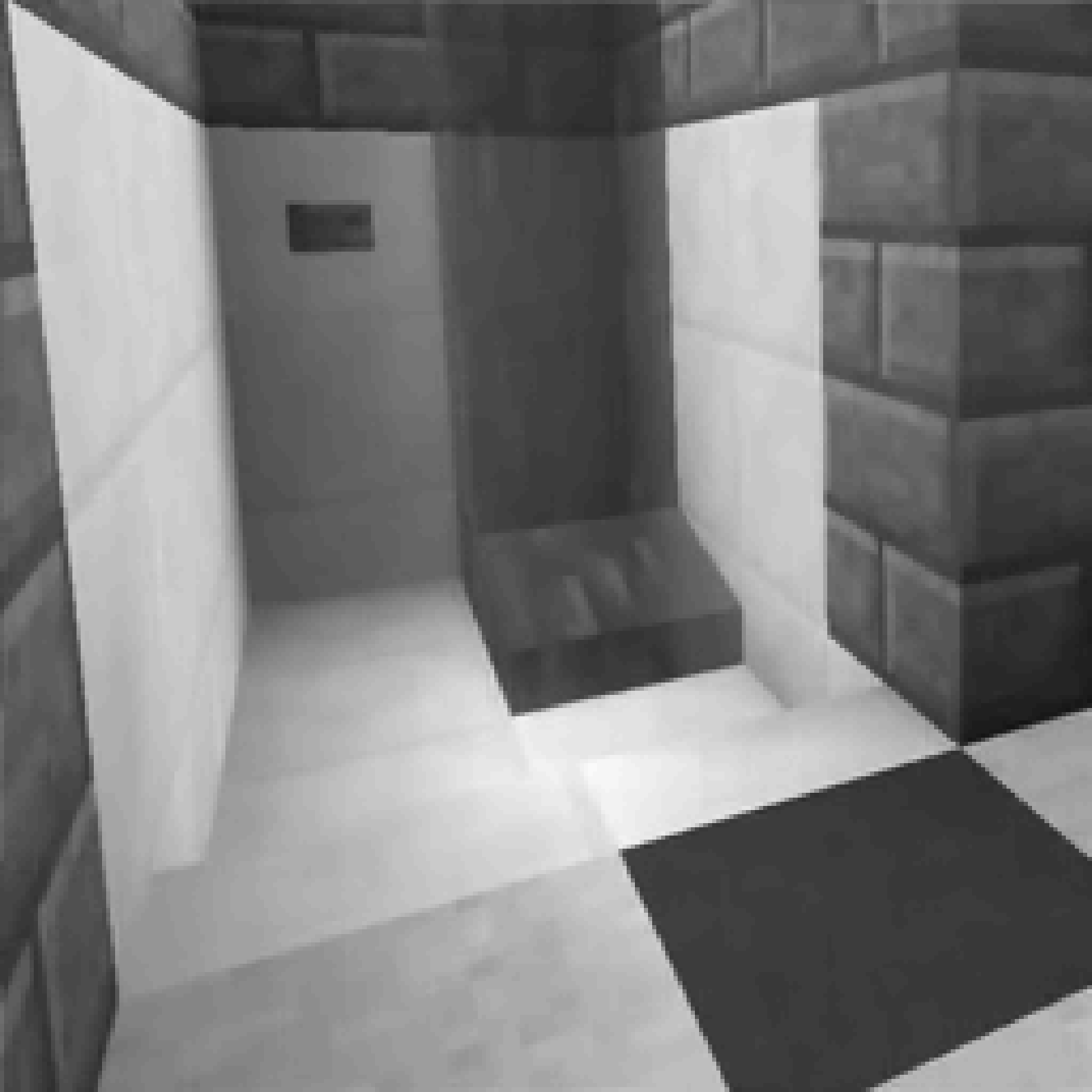}
\end{minipage}}
\caption{Visual comparison of ``Minecraft'' for random sampling.}\label{MCRandomResults}
\end{figure}

\begin{figure}[t]
\centering
\subfloat[Samples]{\label{MCRandomSample}\begin{minipage}{3cm}
\includegraphics[width=3cm]{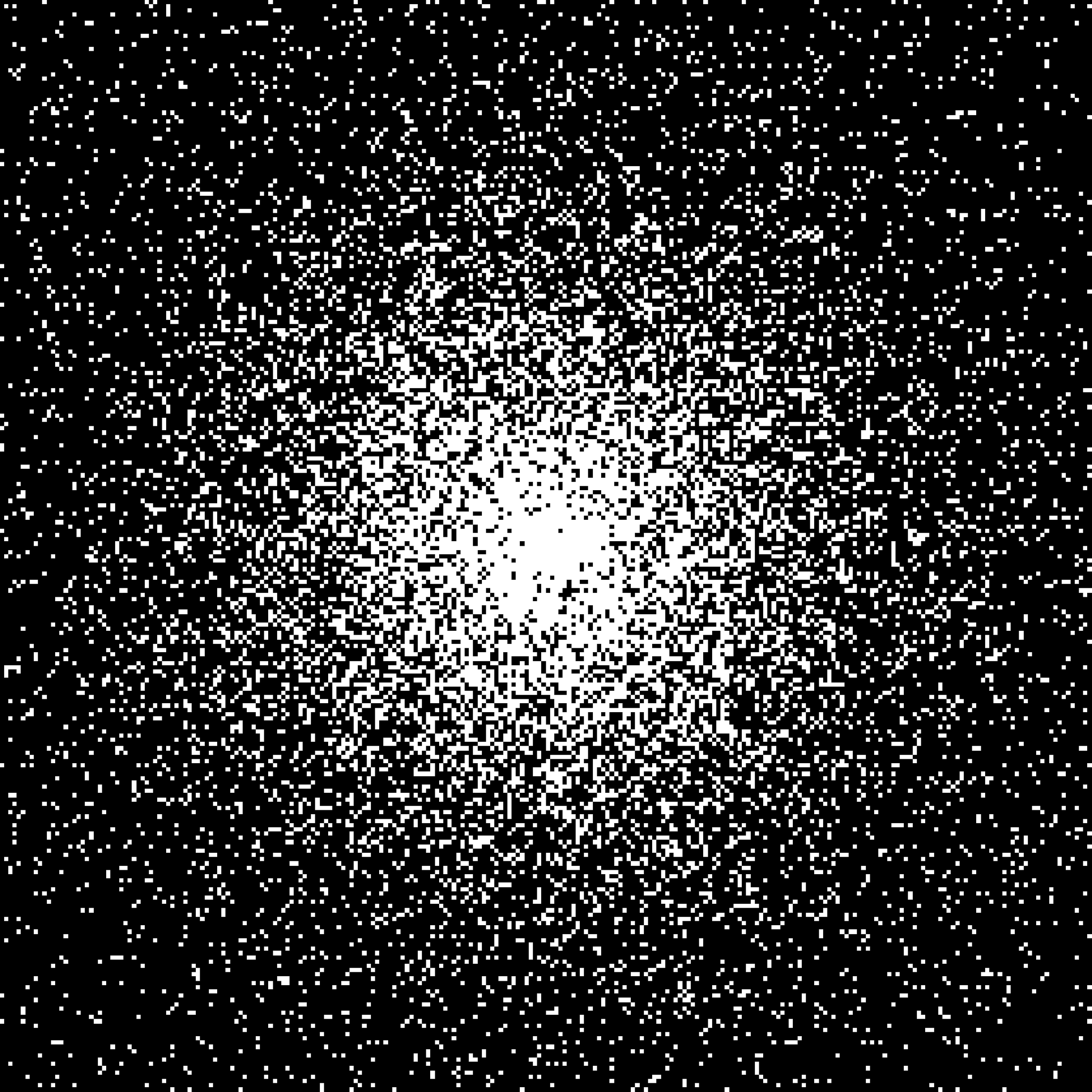}
\end{minipage}}\hspace{-0.05cm}
\subfloat[Proposed]{\label{MCRandomProposedError}\begin{minipage}{3cm}
\includegraphics[width=3cm]{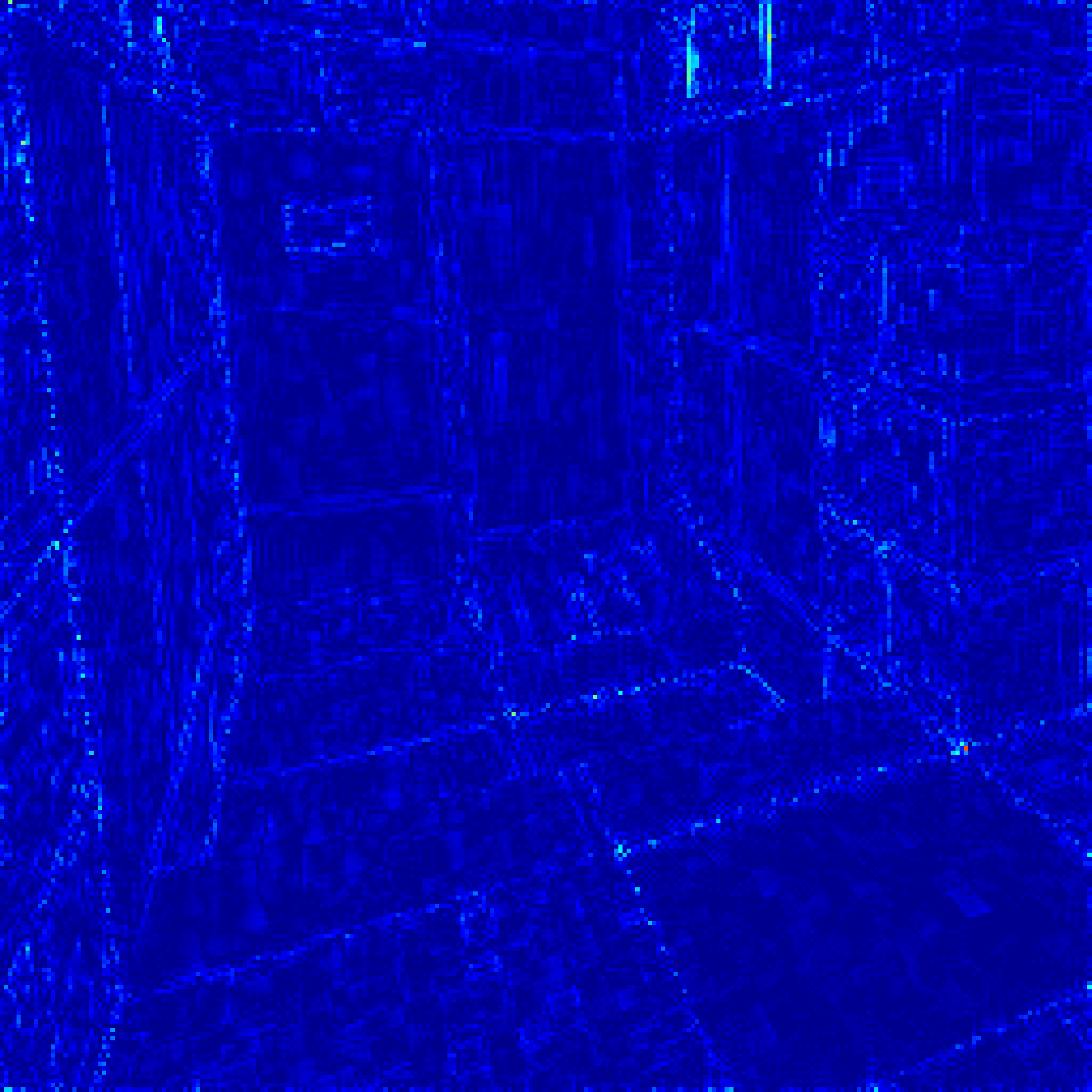}
\end{minipage}}\hspace{-0.05cm}
\subfloat[LSLP]{\label{MCRandomLSLPError}\begin{minipage}{3cm}
\includegraphics[width=3cm]{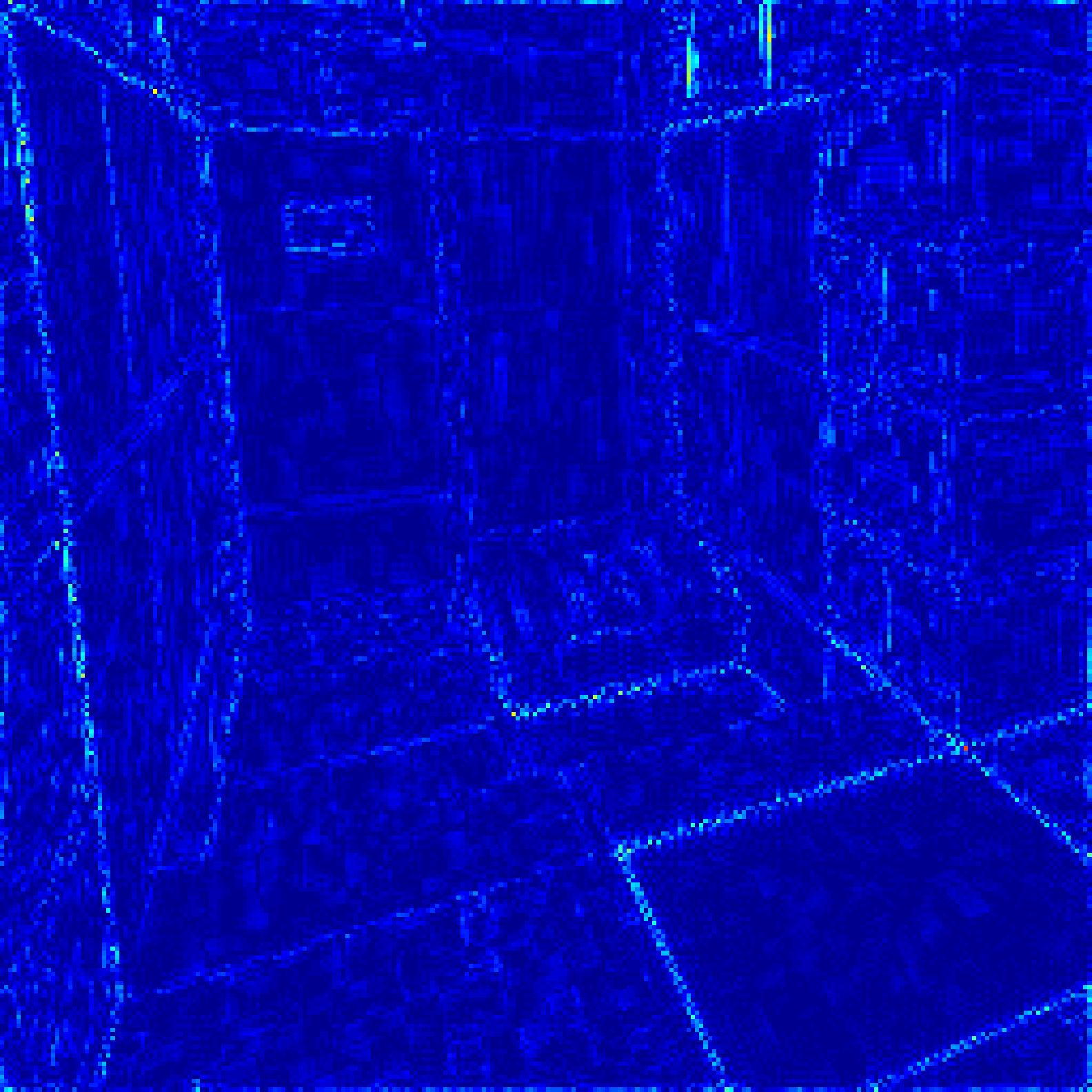}
\end{minipage}}\hspace{-0.05cm}
\subfloat[LRHDDTF]{\label{MCRandomLRHDDTFError}\begin{minipage}{3cm}
\includegraphics[width=3cm]{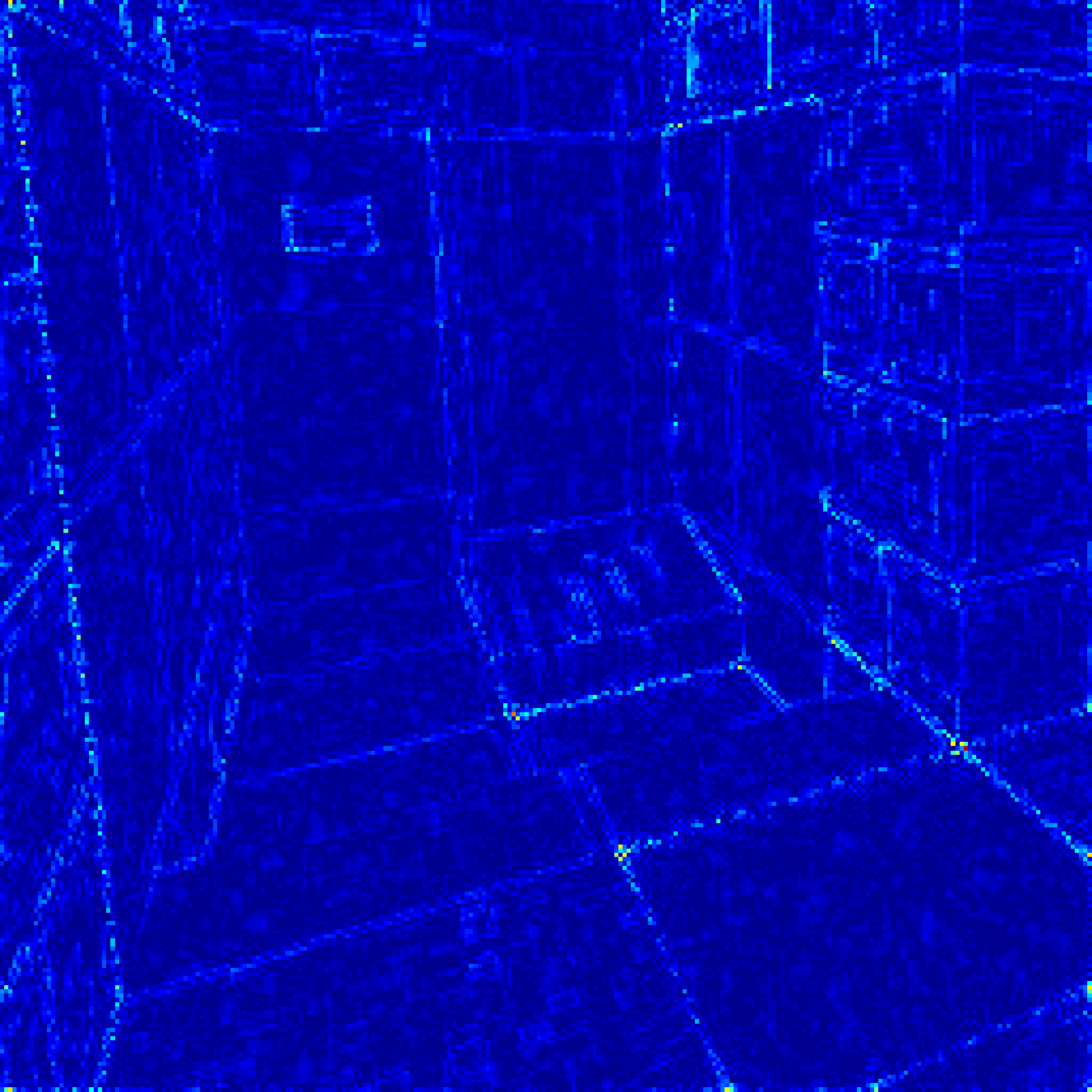}
\end{minipage}}\vspace{-0.25cm}
\subfloat[Schatten $0$]{\label{MCRandomGIRAFError}\begin{minipage}{3cm}
\includegraphics[width=3cm]{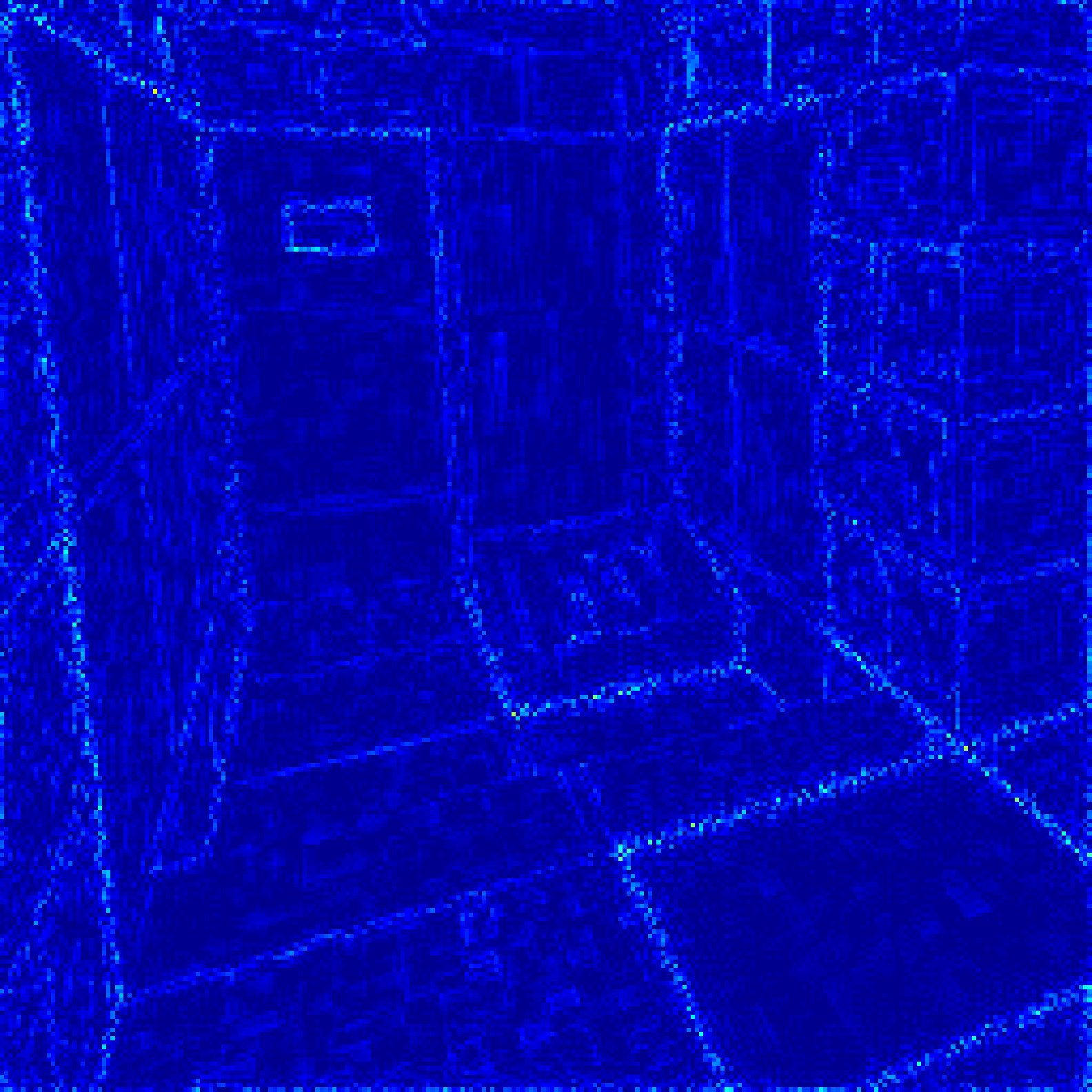}
\end{minipage}}\hspace{-0.05cm}
\subfloat[TV]{\label{MCRandomTVError}\begin{minipage}{3cm}
\includegraphics[width=3cm]{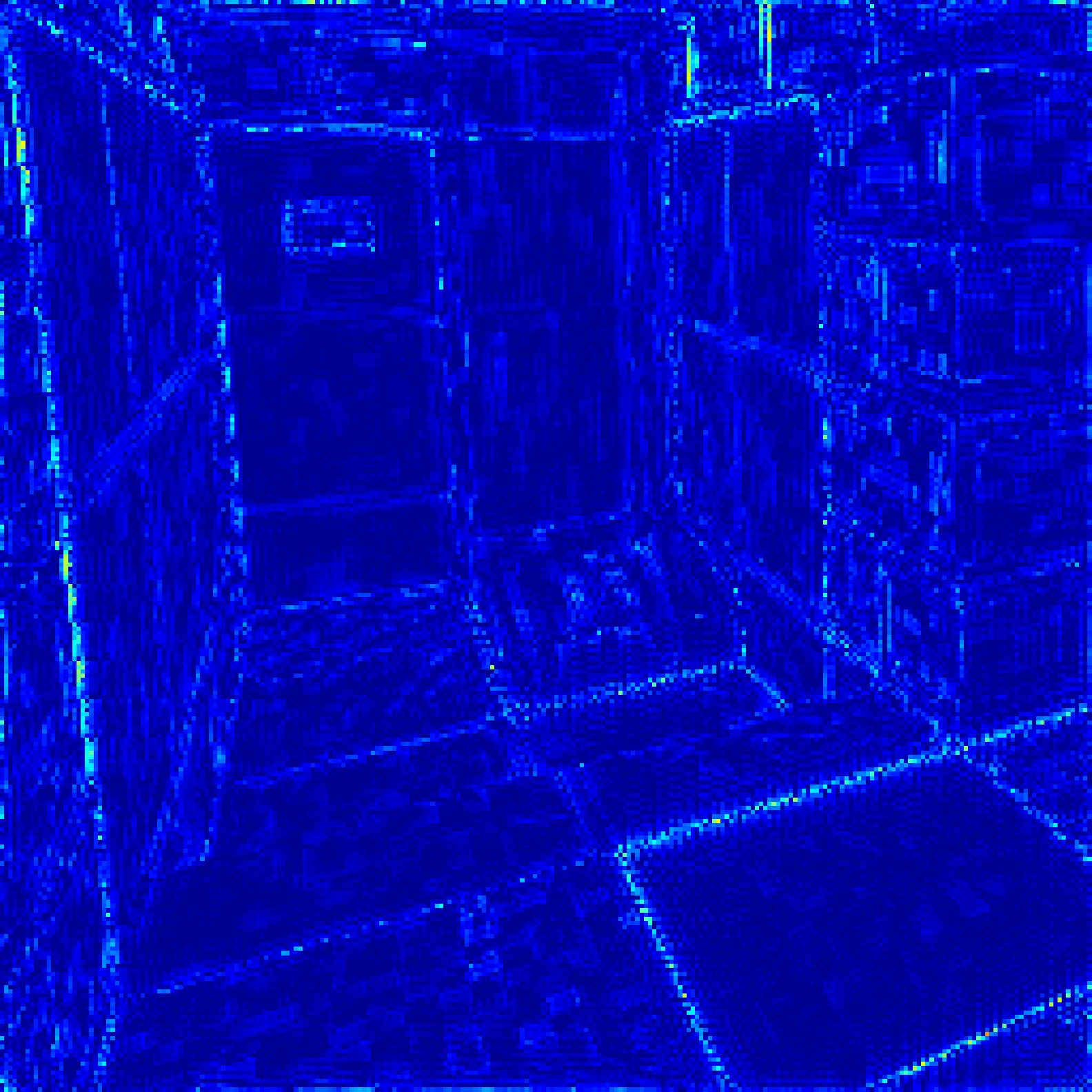}
\end{minipage}}\hspace{-0.05cm}
\subfloat[Haar]{\label{MCRandomHaarError}\begin{minipage}{3cm}
\includegraphics[width=3cm]{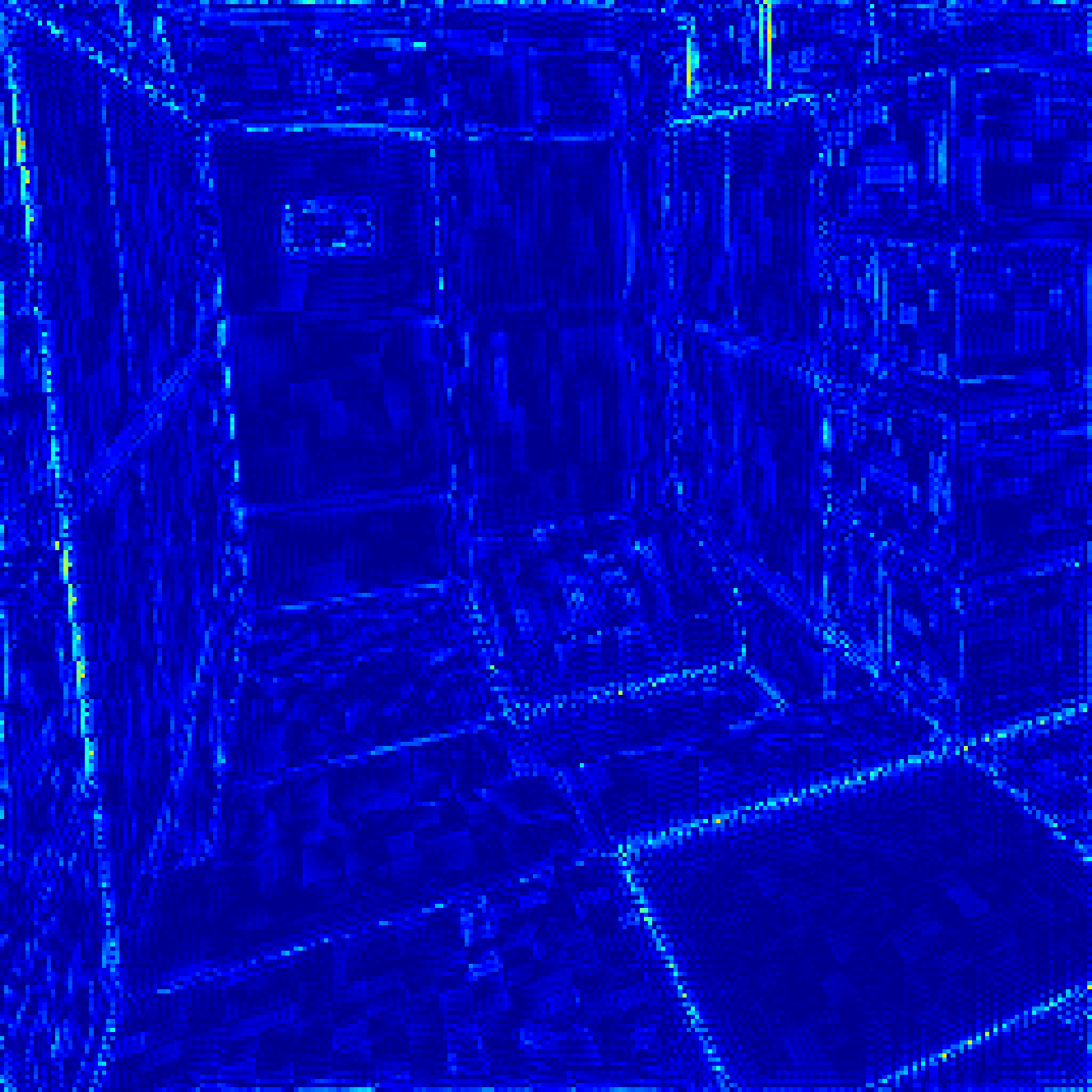}
\end{minipage}}\hspace{-0.05cm}
\subfloat[DDTF]{\label{MCRandomDDTFError}\begin{minipage}{3cm}
\includegraphics[width=3cm]{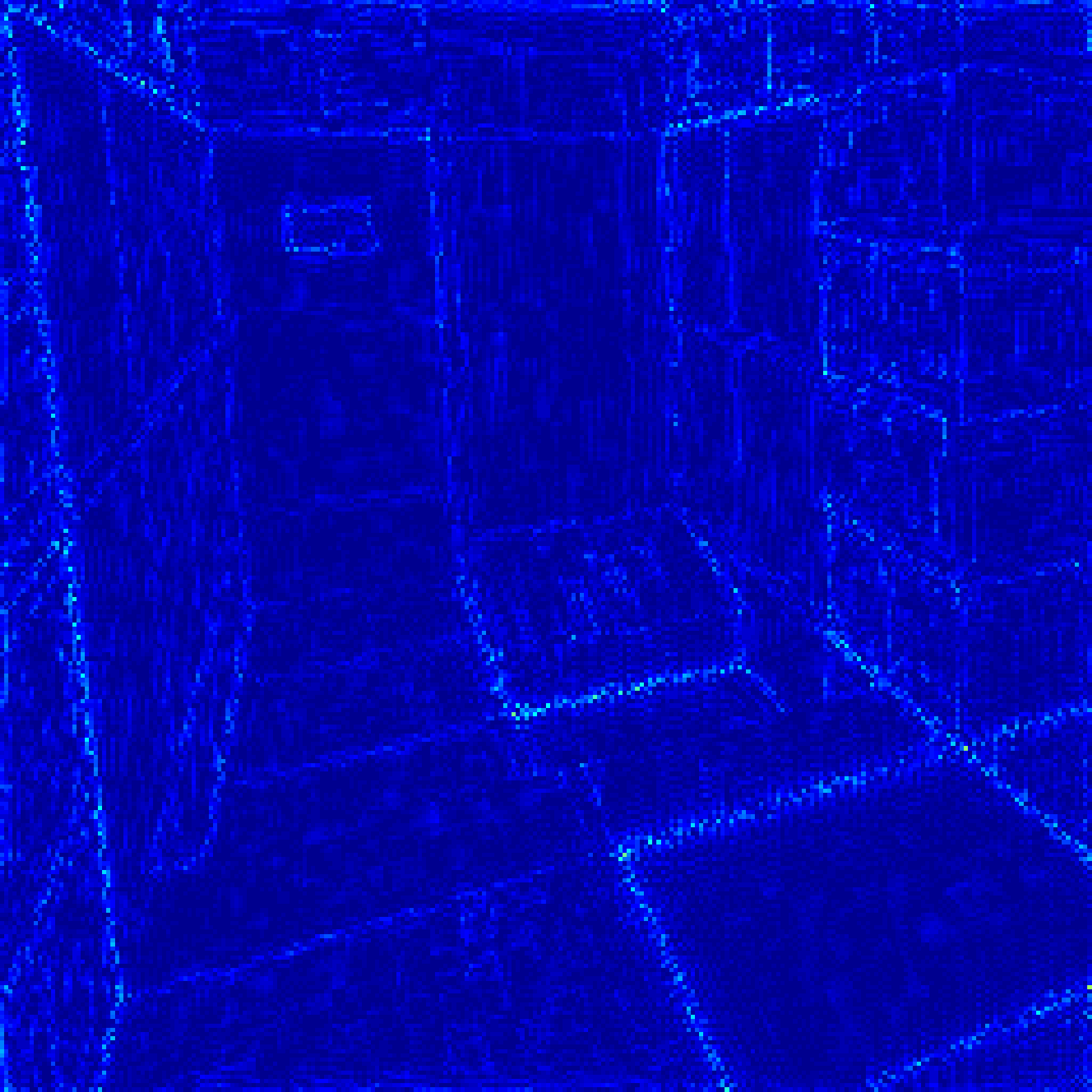}
\end{minipage}}
\caption{Error maps of \cref{MCRandomResults}.}\label{MCRandomError}
\end{figure}

\begin{figure}[t]
\centering
\subfloat[Ref.]{\label{MCILFOriginal}\begin{minipage}{3cm}
\includegraphics[width=3cm]{MCOriginal.pdf}
\end{minipage}}\hspace{-0.05cm}
\subfloat[Proposed]{\label{MCILFProposed}\begin{minipage}{3cm}
\includegraphics[width=3cm]{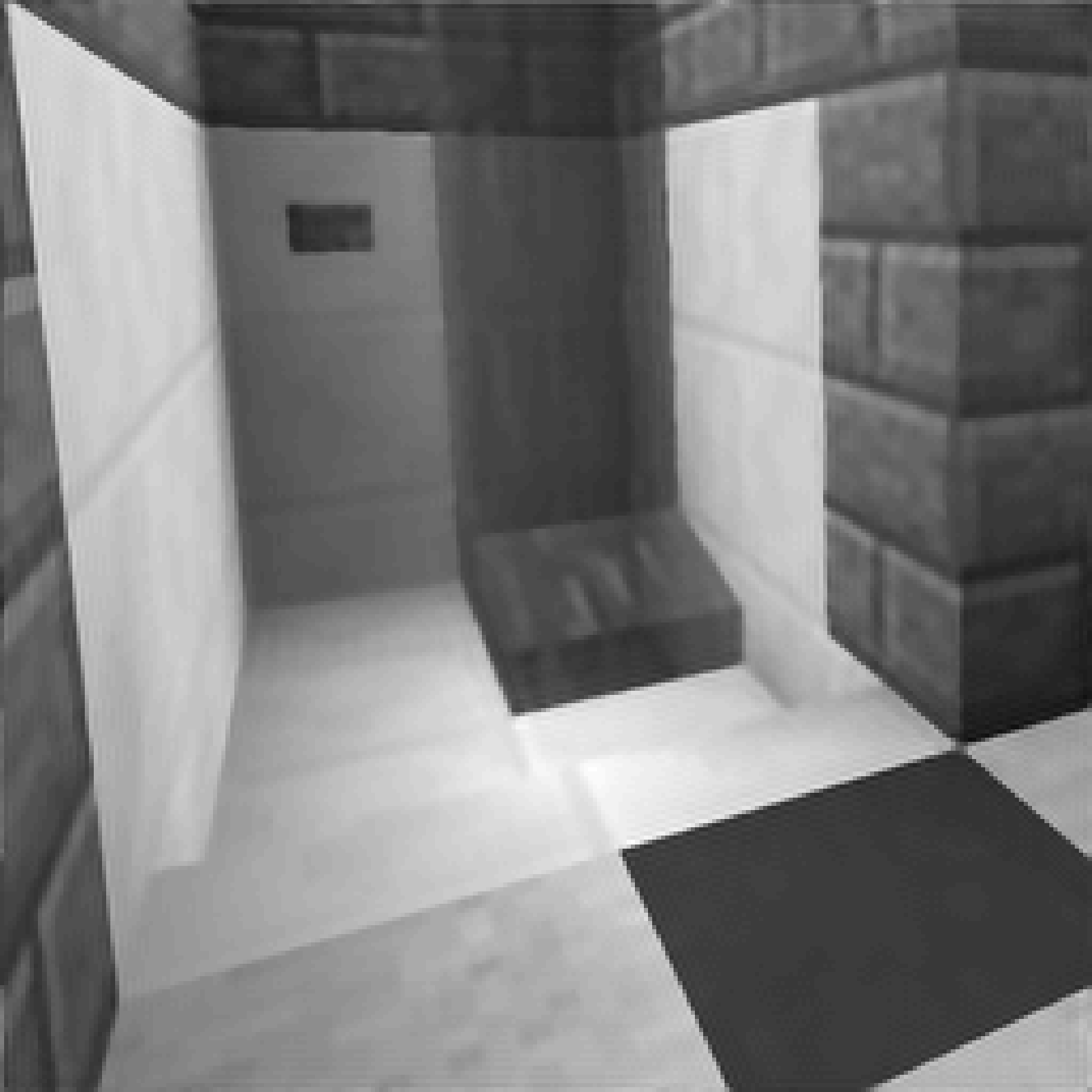}
\end{minipage}}\hspace{-0.05cm}
\subfloat[LSLP]{\label{MCILFLSLP}\begin{minipage}{3cm}
\includegraphics[width=3cm]{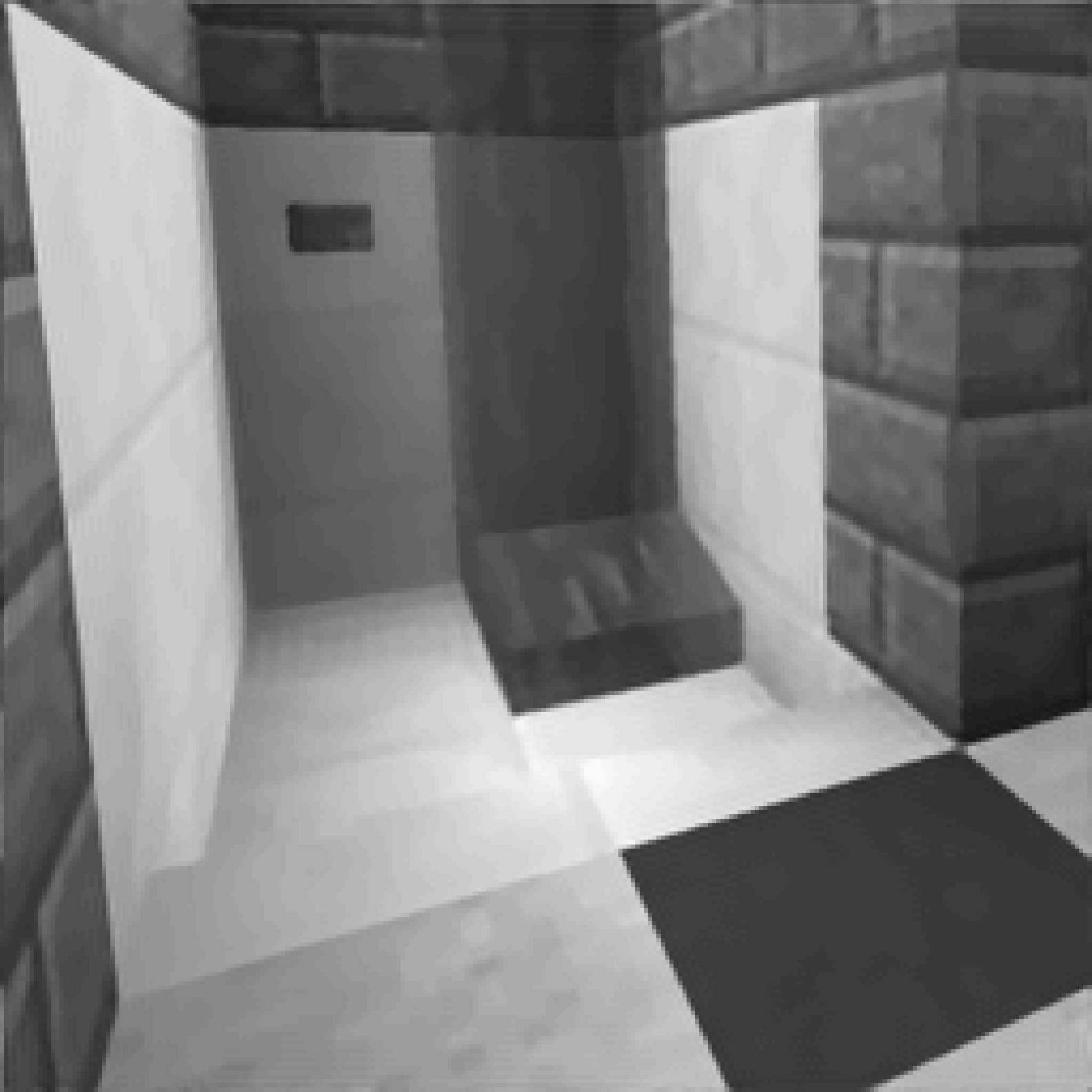}
\end{minipage}}\hspace{-0.05cm}
\subfloat[LRHDDTF]{\label{MCILFLRHDDTF}\begin{minipage}{3cm}
\includegraphics[width=3cm]{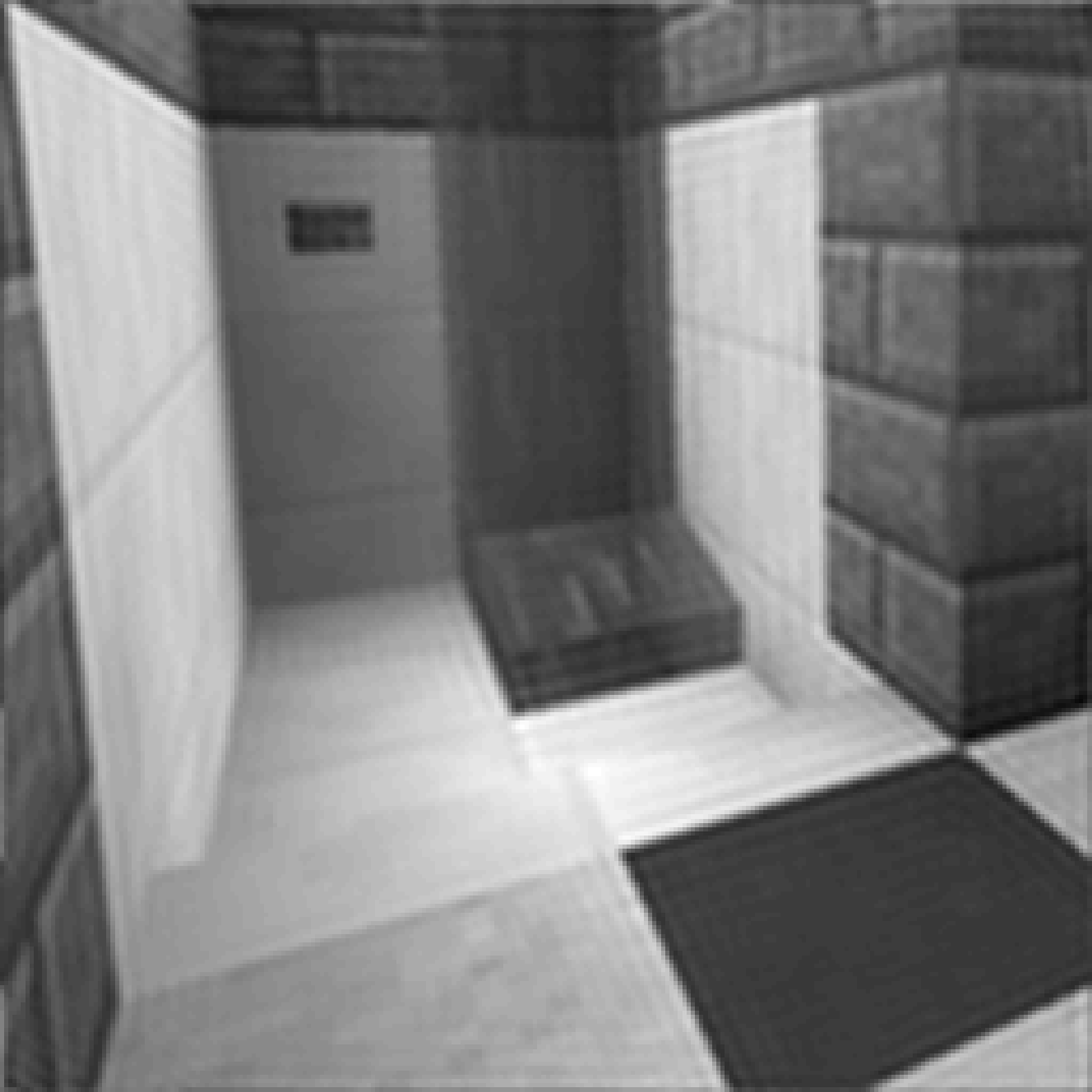}
\end{minipage}}\vspace{-0.25cm}
\subfloat[Schatten $0$]{\label{MCILFGIRAF}\begin{minipage}{3cm}
\includegraphics[width=3cm]{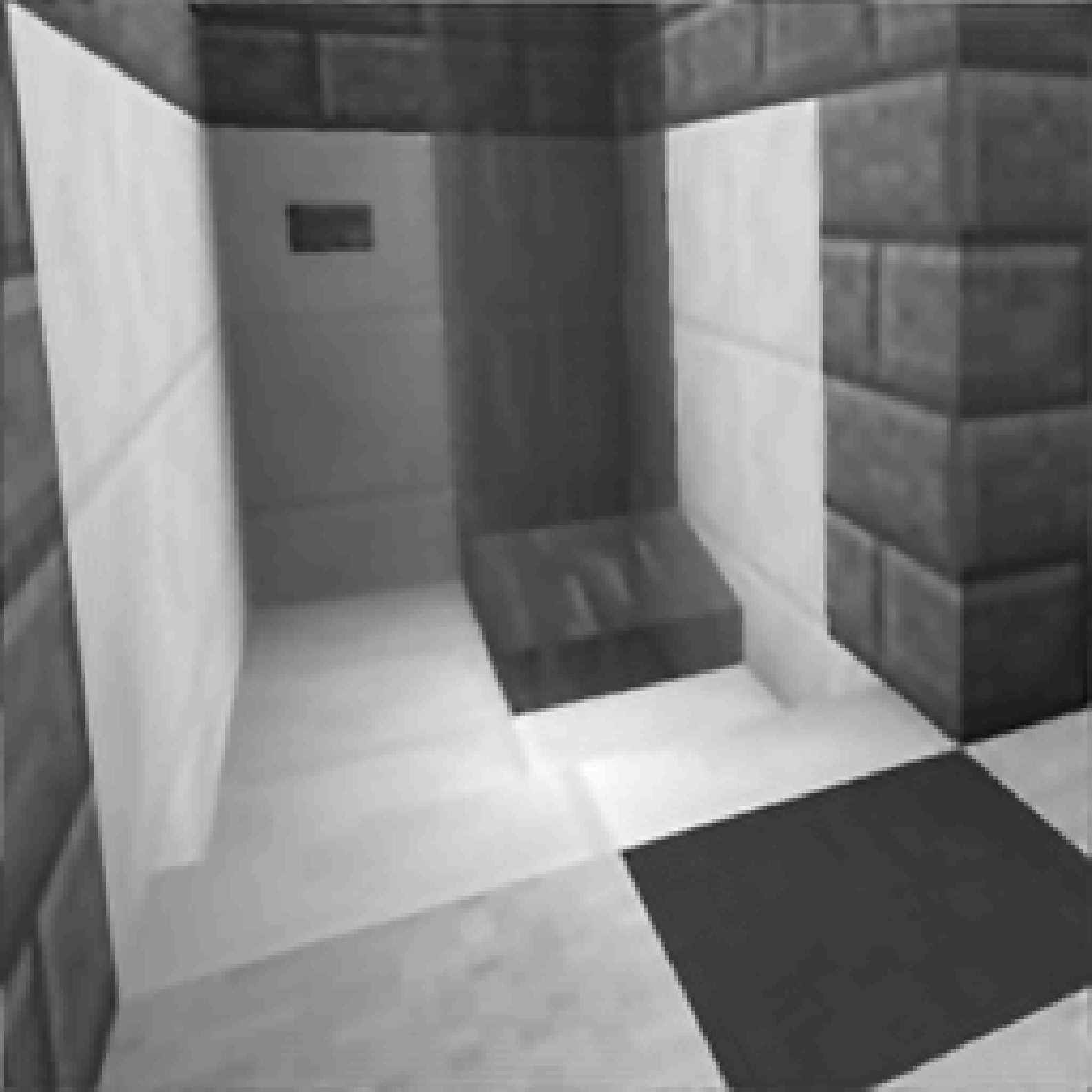}
\end{minipage}}\hspace{-0.05cm}
\subfloat[TV]{\label{MCILFTV}\begin{minipage}{3cm}
\includegraphics[width=3cm]{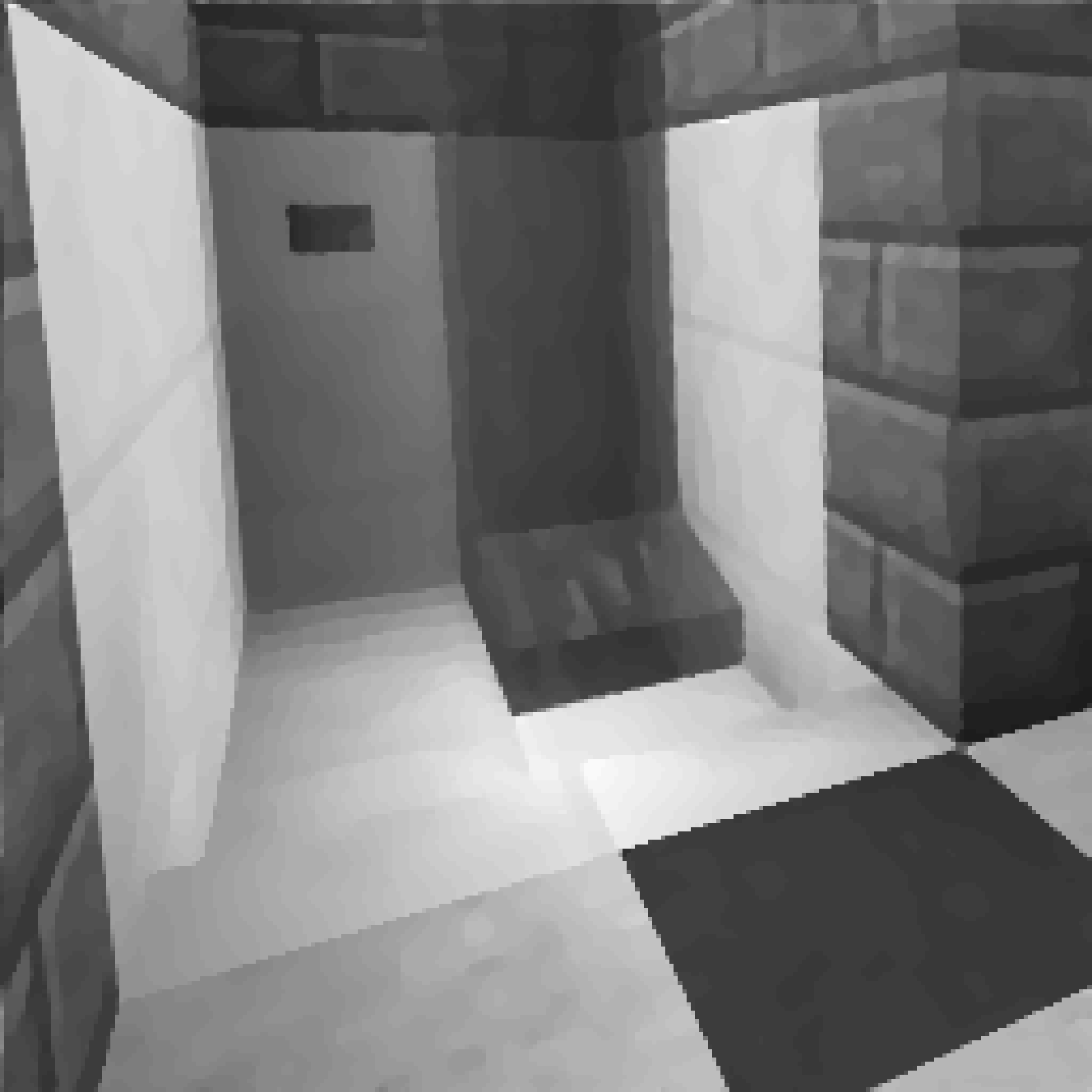}
\end{minipage}}\hspace{-0.05cm}
\subfloat[Haar]{\label{MCILFHaar}\begin{minipage}{3cm}
\includegraphics[width=3cm]{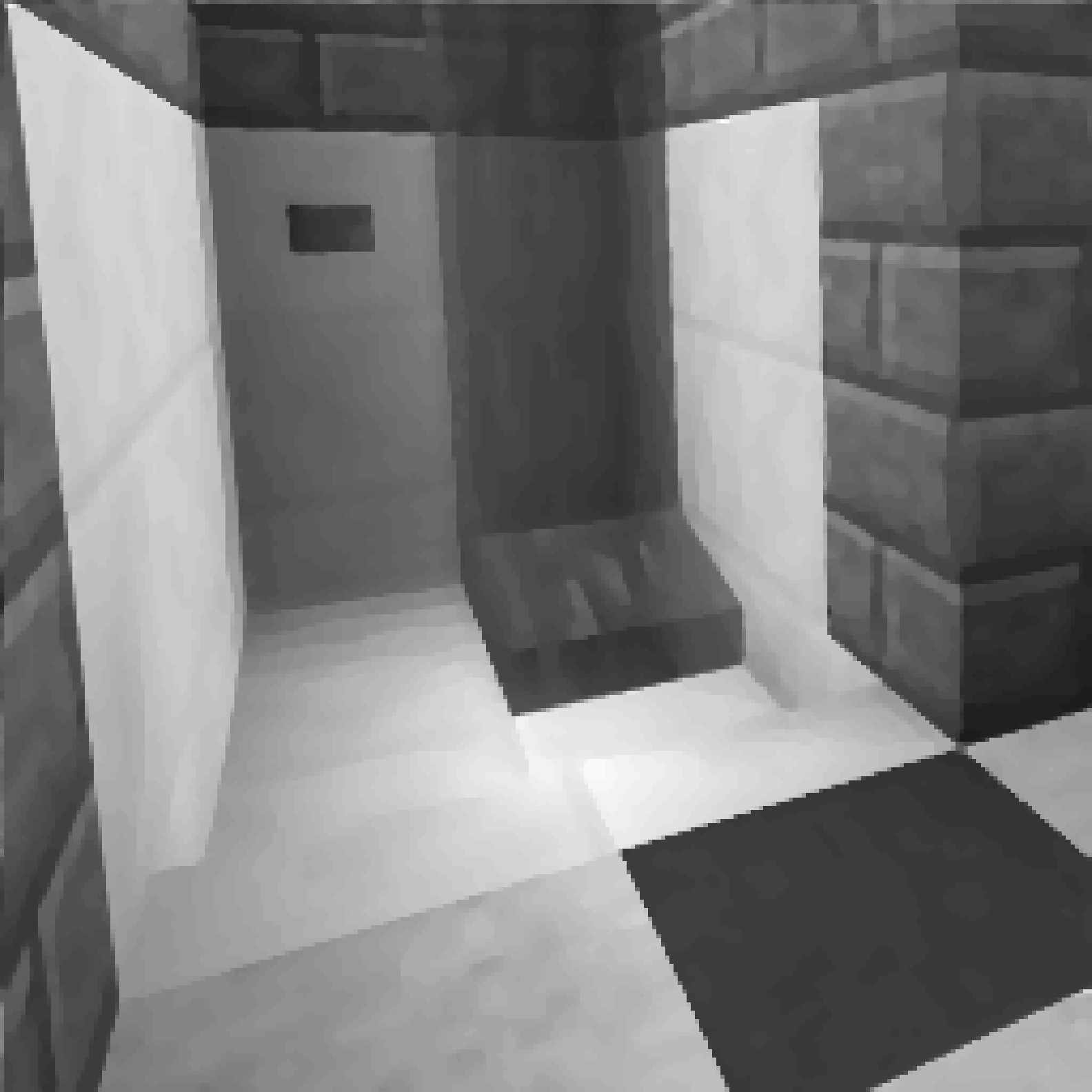}
\end{minipage}}\hspace{-0.05cm}
\subfloat[DDTF]{\label{MCILFDDTF}\begin{minipage}{3cm}
\includegraphics[width=3cm]{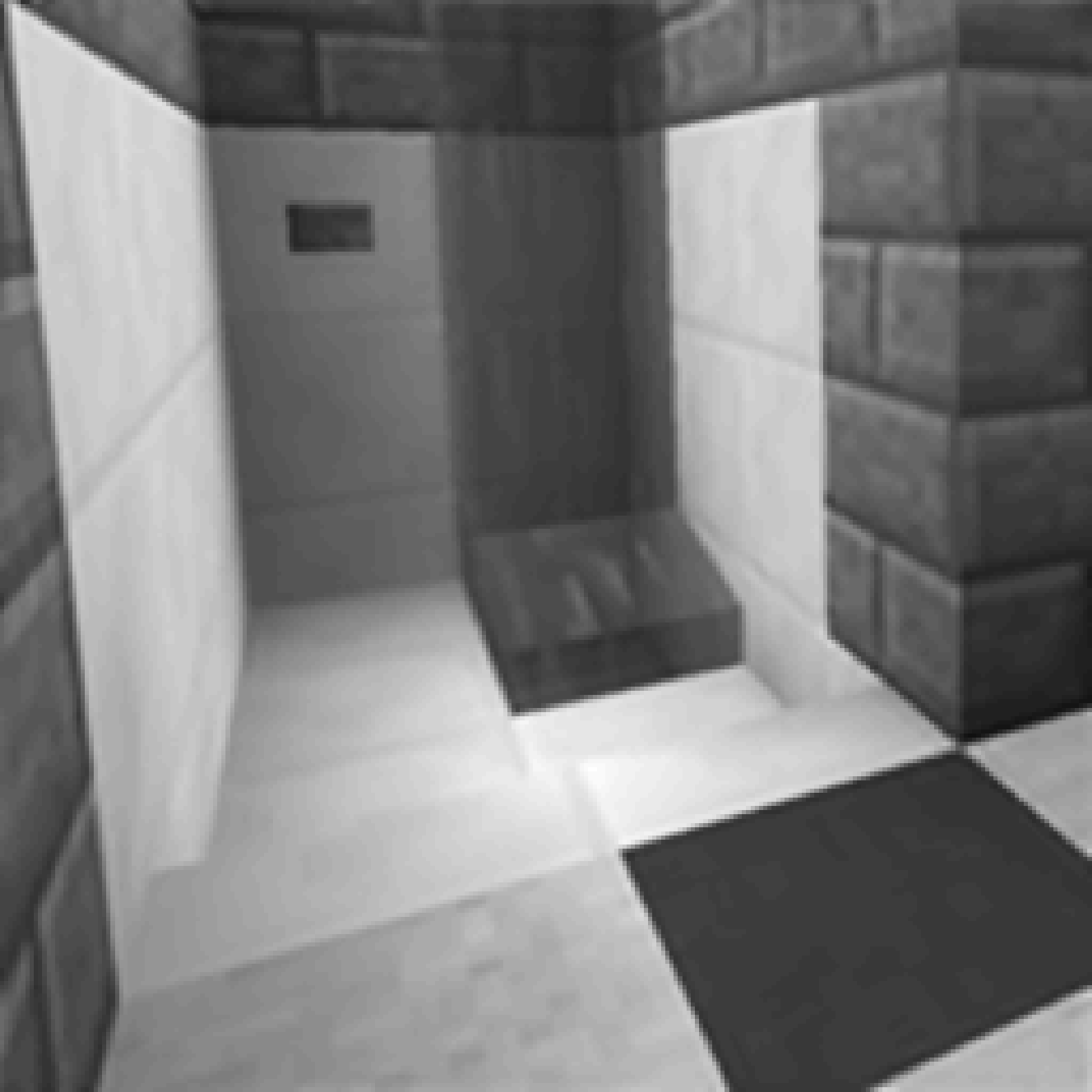}
\end{minipage}}
\caption{Visual comparison of ``Minecraft'' for ideal low-pass filter deconvolution.}\label{MCILFResults}
\end{figure}

\begin{figure}[t]
\centering
\subfloat[TF]{\label{MCILFTF}\begin{minipage}{3cm}
\includegraphics[width=3cm]{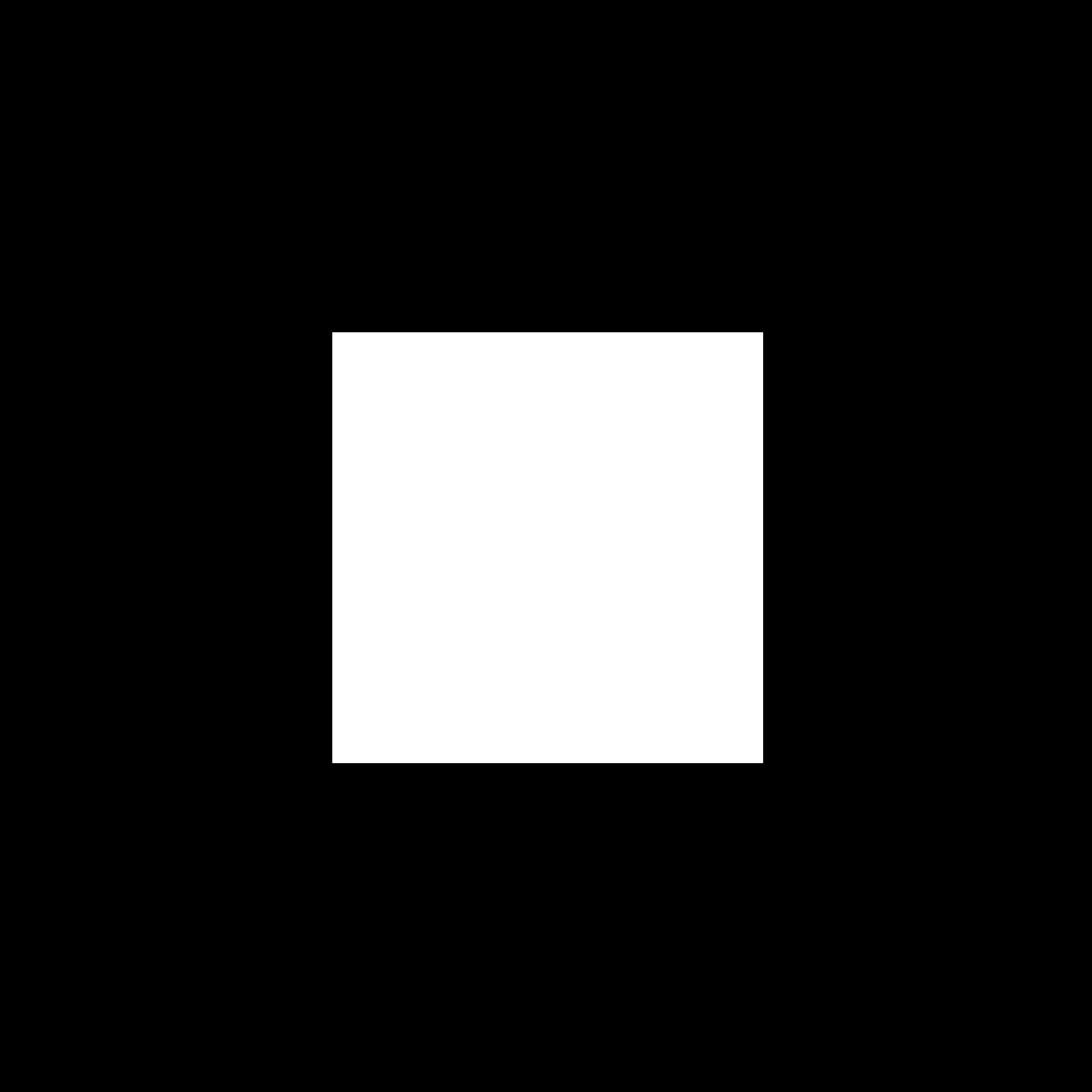}
\end{minipage}}\hspace{-0.05cm}
\subfloat[Proposed]{\label{MCILFProposedError}\begin{minipage}{3cm}
\includegraphics[width=3cm]{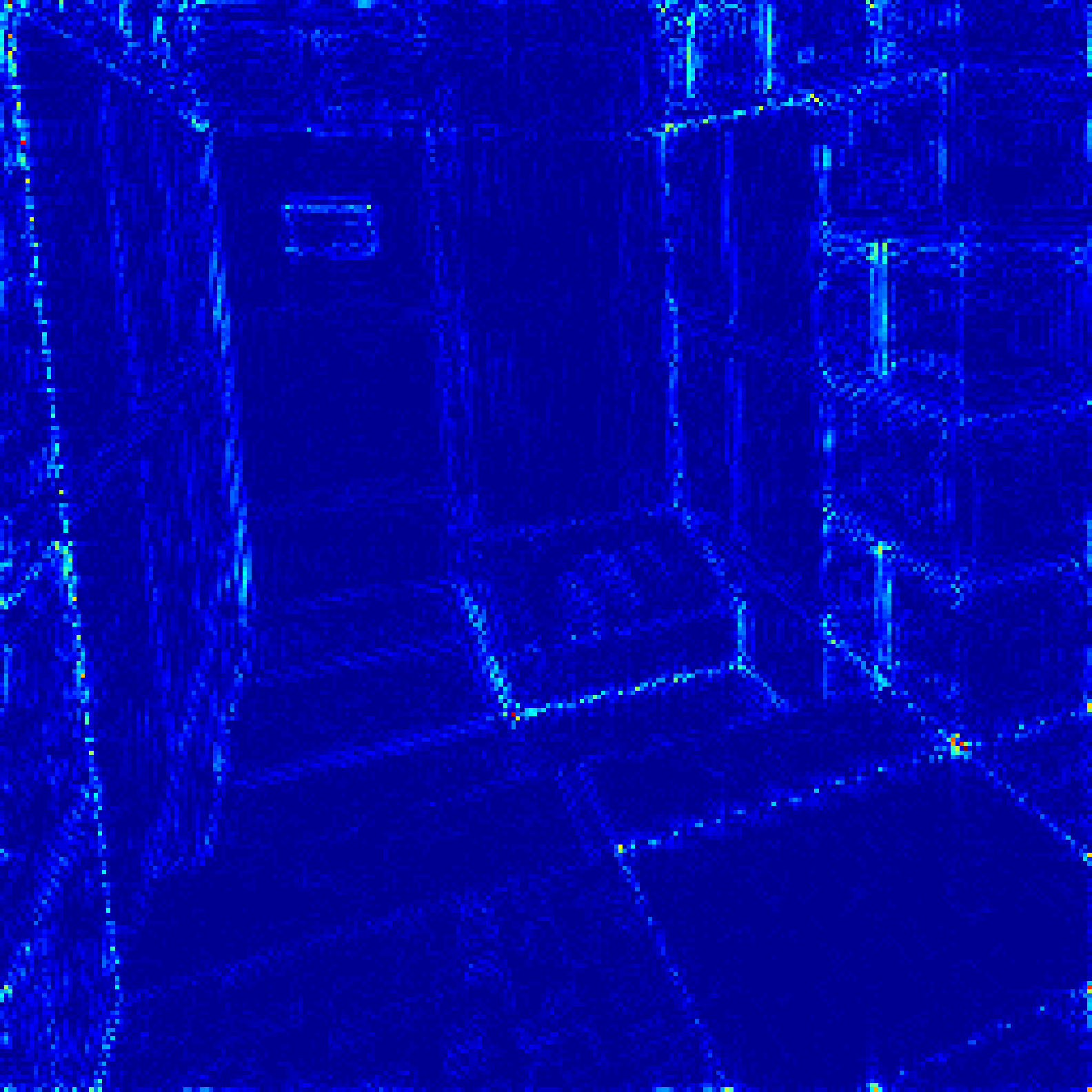}
\end{minipage}}\hspace{-0.05cm}
\subfloat[LSLP]{\label{MCILFLSLPError}\begin{minipage}{3cm}
\includegraphics[width=3cm]{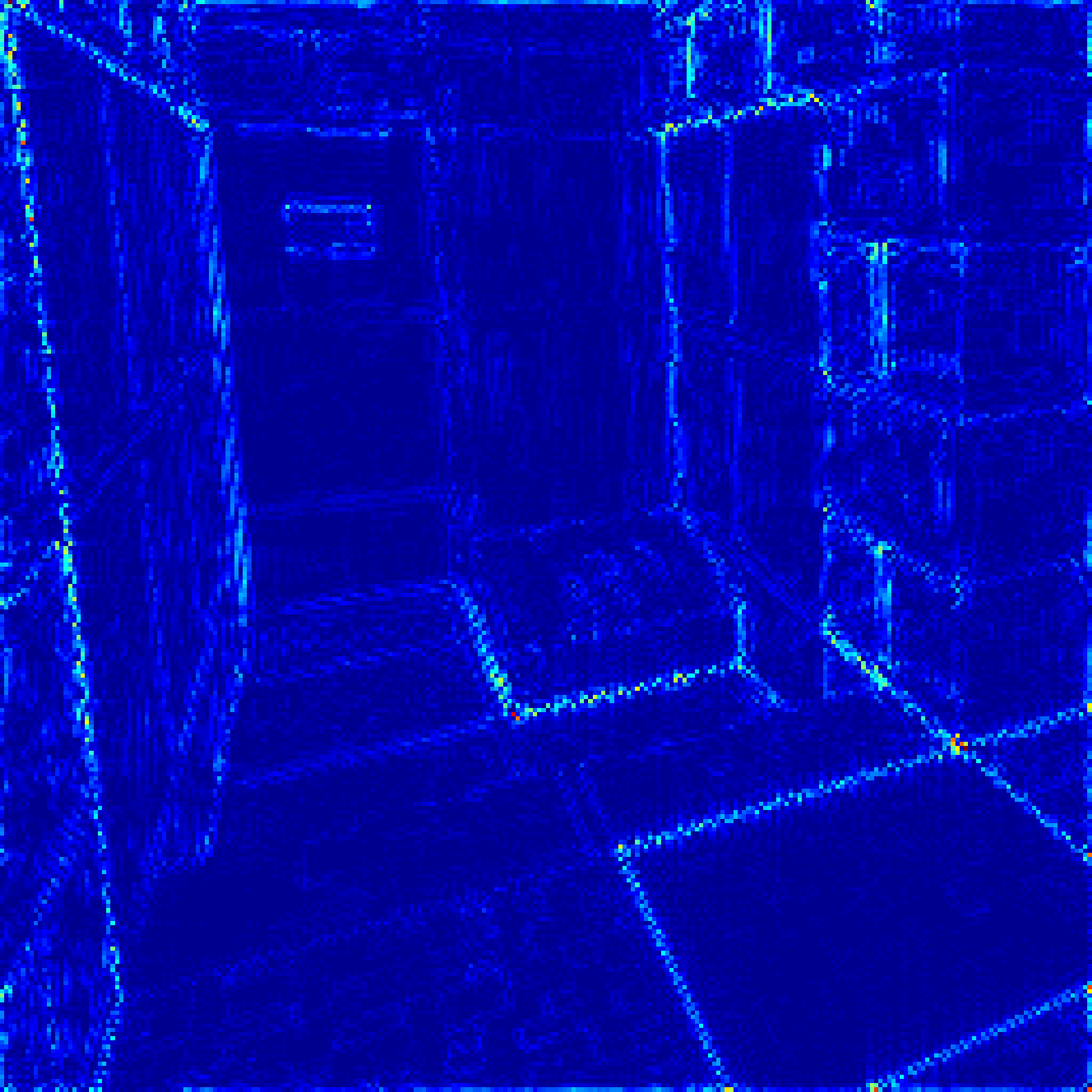}
\end{minipage}}\hspace{-0.05cm}
\subfloat[LRHDDTF]{\label{MCILFLRHDDTFError}\begin{minipage}{3cm}
\includegraphics[width=3cm]{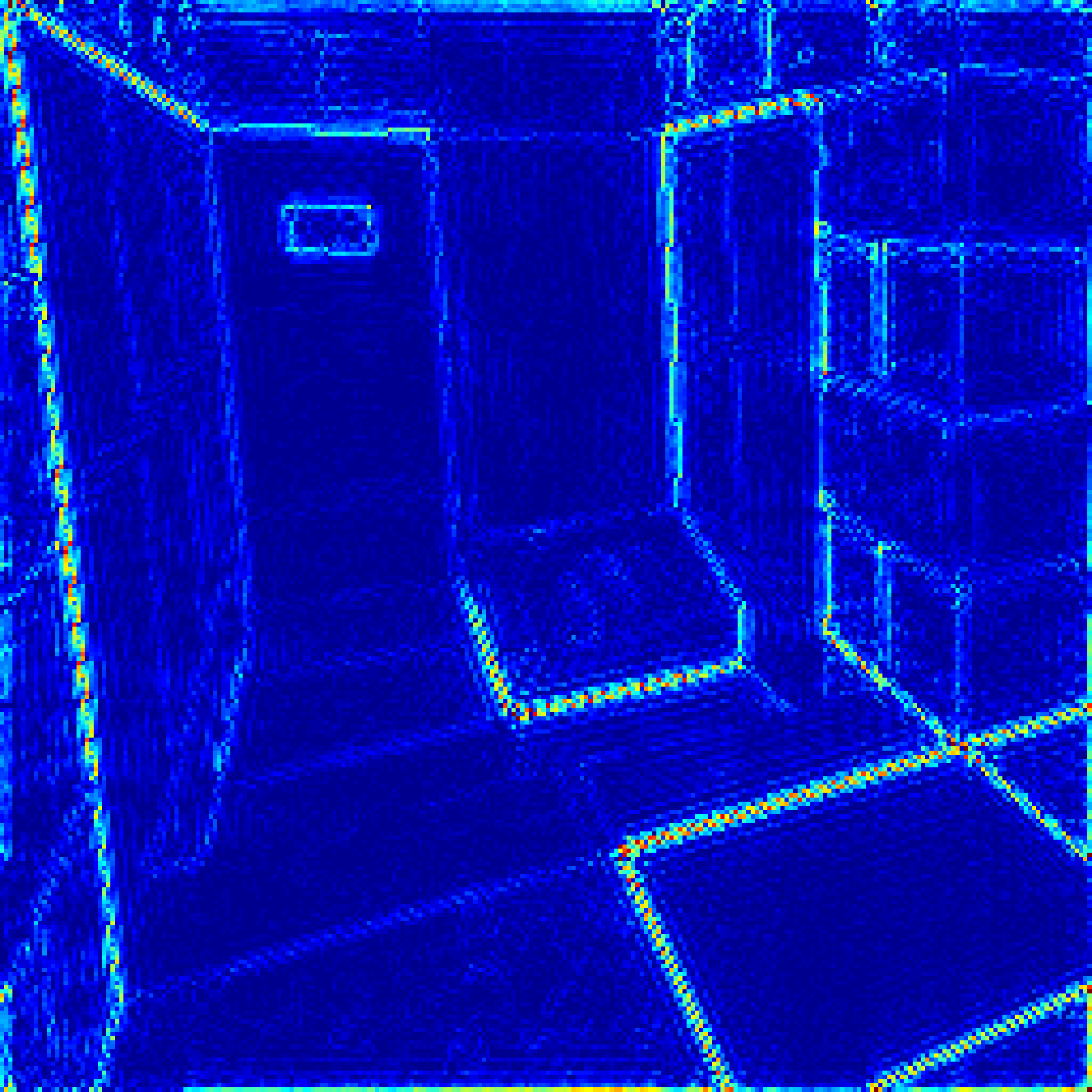}
\end{minipage}}\vspace{-0.25cm}
\subfloat[Schatten $0$]{\label{MCILFGIRAFError}\begin{minipage}{3cm}
\includegraphics[width=3cm]{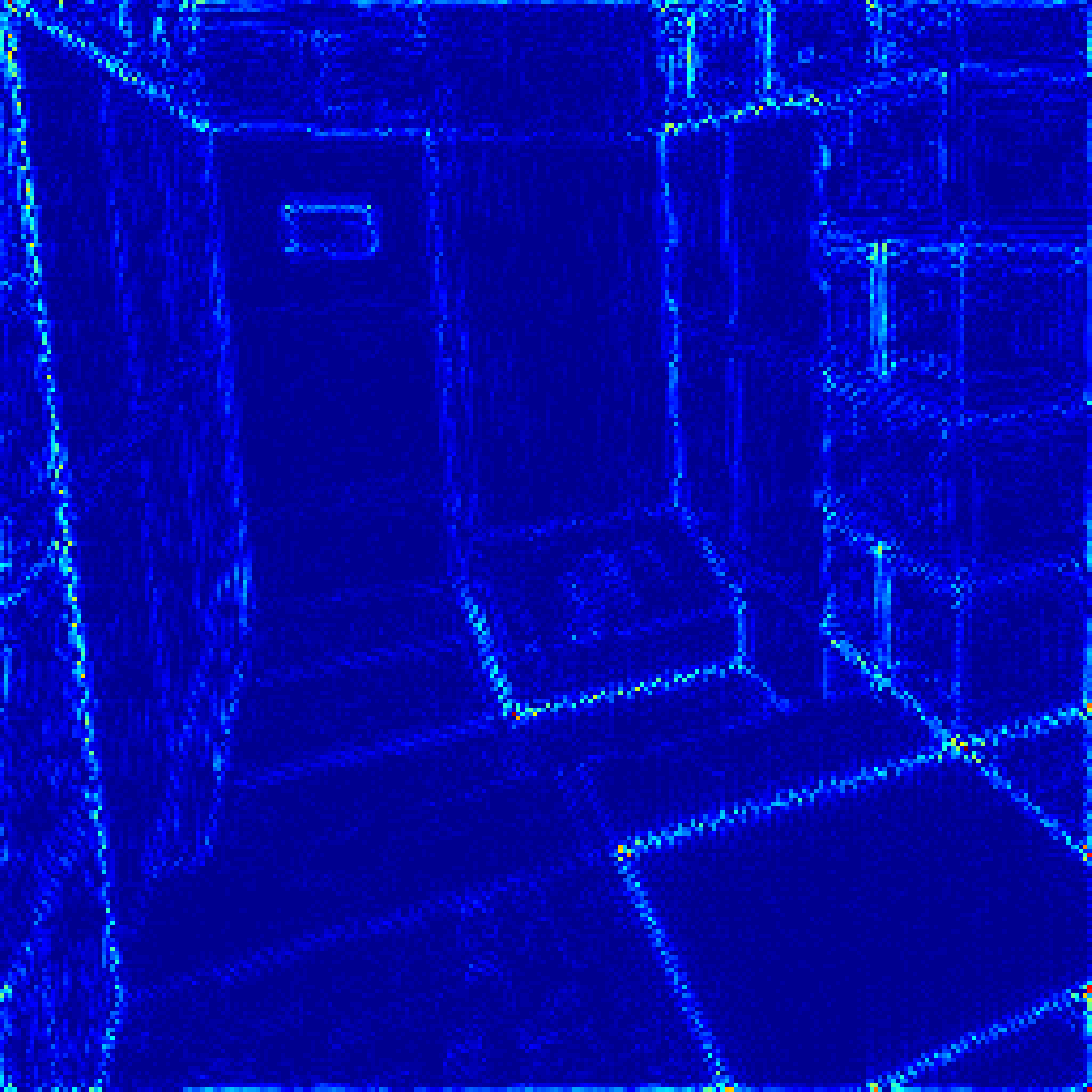}
\end{minipage}}\hspace{-0.05cm}
\subfloat[TV]{\label{MCILFTVError}\begin{minipage}{3cm}
\includegraphics[width=3cm]{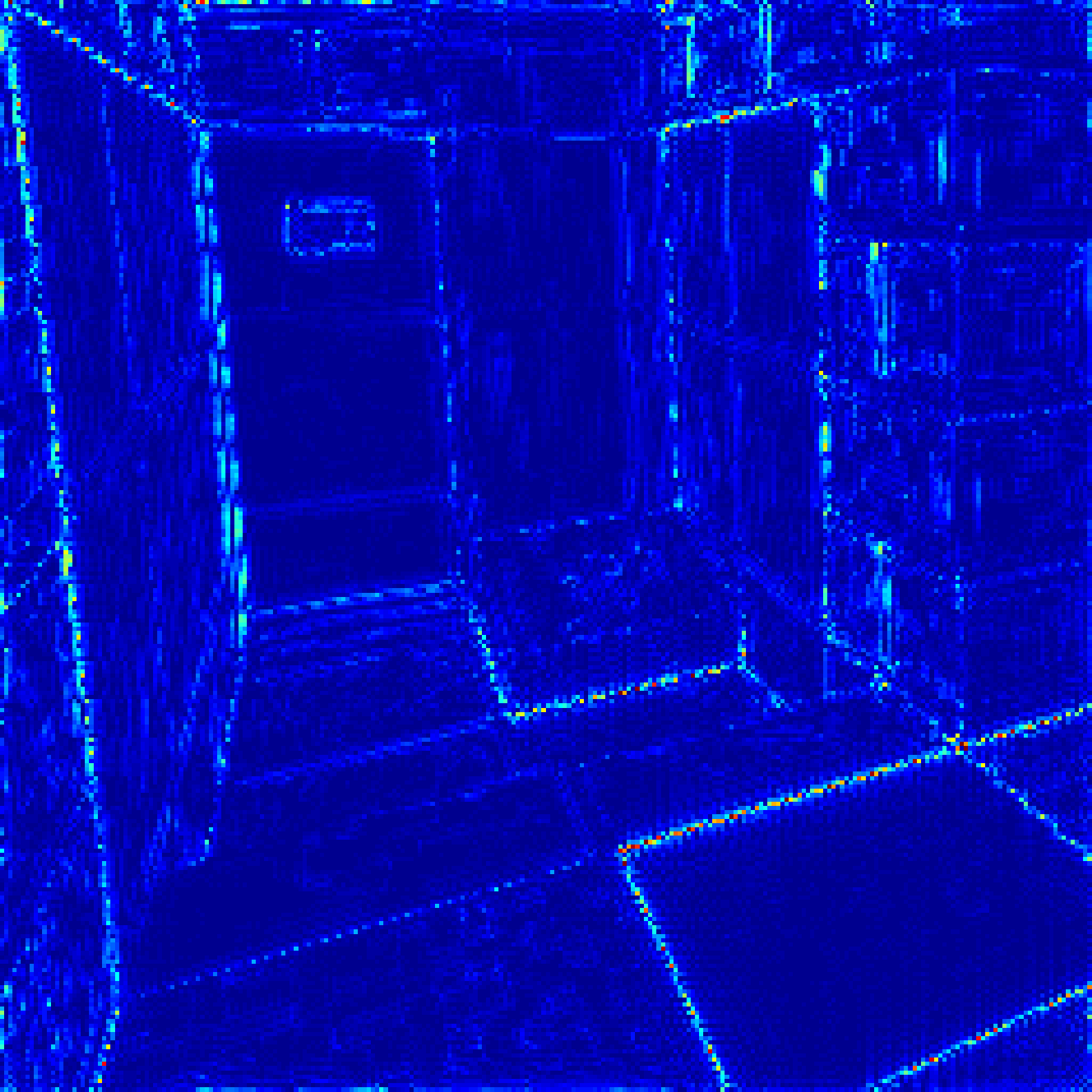}
\end{minipage}}\hspace{-0.05cm}
\subfloat[Haar]{\label{MCILFHaarError}\begin{minipage}{3cm}
\includegraphics[width=3cm]{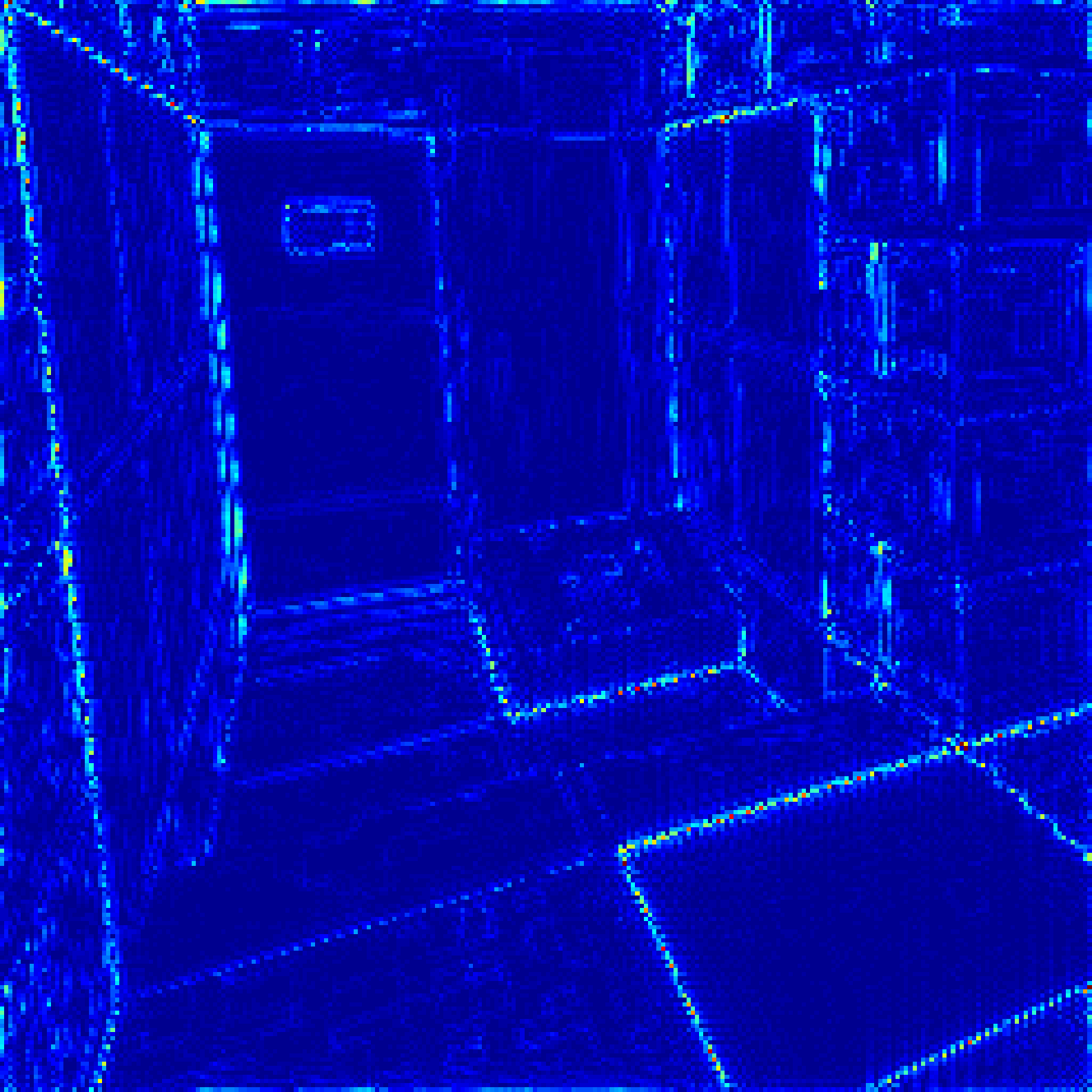}
\end{minipage}}\hspace{-0.05cm}
\subfloat[DDTF]{\label{MCILFDDTFError}\begin{minipage}{3cm}
\includegraphics[width=3cm]{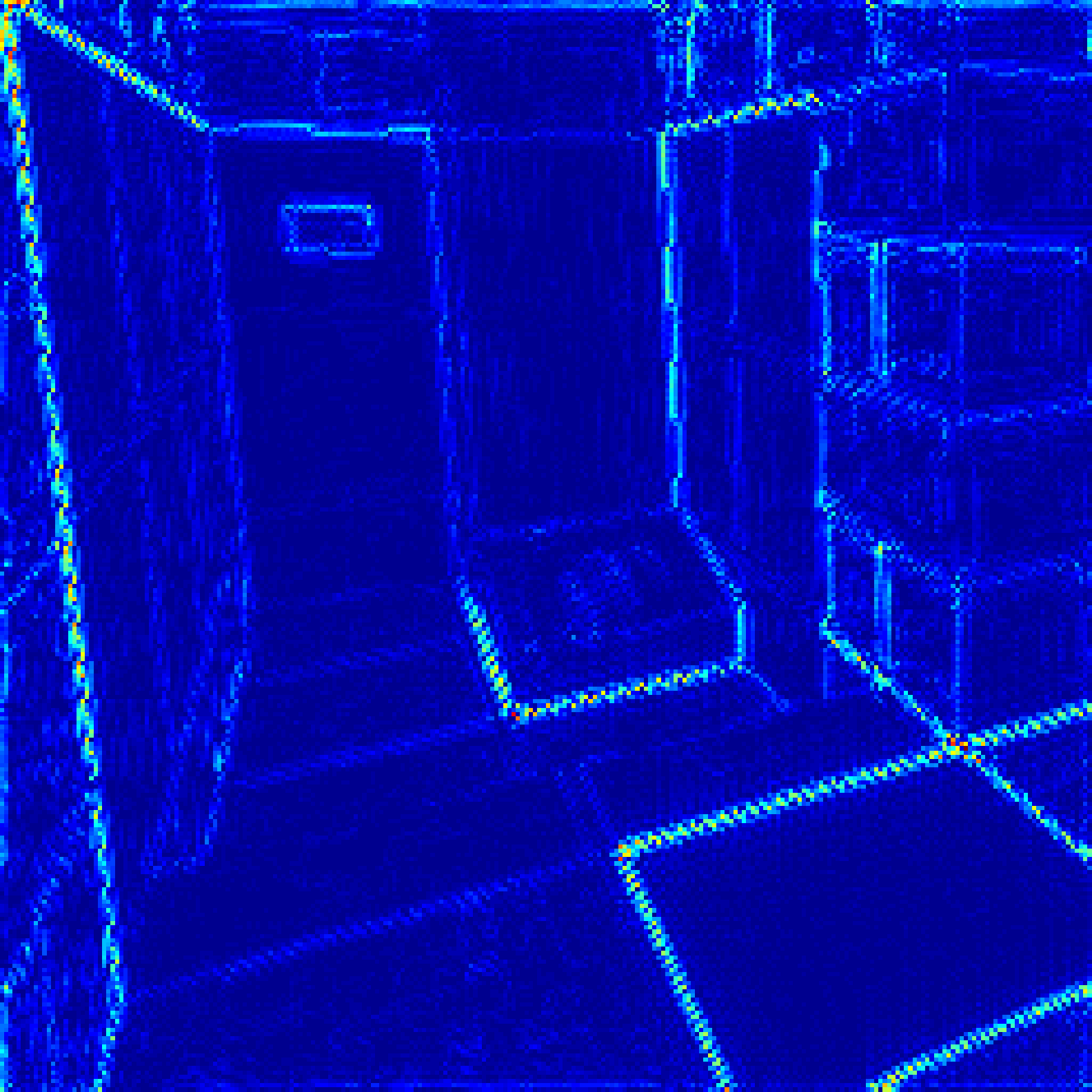}
\end{minipage}}
\caption{Error maps of \cref{MCILFResults}.}\label{MCILFError}
\end{figure}

\cref{NaturalImageResults} summarizes the SNR, HFEN, and SSIM for the non-synthetic images, and \cref{AngryBirdsRandomResults,AngryBirdsILFResults,SuperMarioRandomResults,SuperMarioILFResults,MCRandomResults,MCILFResults} display the visual comparisons with the corresponding error maps in \cref{AngryBirdsRandomError,AngryBirdsILFError,SuperMarioRandomError,SuperMarioILFError,MCRandomError,MCILFError}. Overall, we can see that the pros and the cons of the restoration methods are similar to the synthetic images, and our proposed model \cref{ProposedModel2} performs best in almost every index in all scenarios, with visually observable improvements. We first note that compared to the synthetic images (``Ellipse'' and ``Rectangle''), non-synthetic images include more complicated geometry. Nevertheless, the model \cref{ProposedModel2} is able to achieve a noticeable improvements over other image restoration methods, as we can see from ``Angry Birds'' and ``Super Mario''. Since the degree of a trigonometric curve puts a bound on the number of $\Omega_j$'s in \cref{uModel} \cite{G.Ongie2016}, we can handle the images with densely located edges by choosing a size of patches appropriately, as long as the target images can be well modeled by a piecewise constant function. Finally, it is also noting that our proposed analysis approach shows a performance gain even when the target image becomes far from a piecewise constant image by textures, as we can see from ``Minecraft''. Hence, the numerical results for non-synthetic images show that, even though textures are different from piecewise constant images, our approach is relatively robust to the image model \cref{uModel}.

\section{Conclusion and future directions}\label{Conclusion}

In this paper, we propose a two-stage continuous domain regularization approach for the piecewise constant image restoration. Our two-stage approach is based on the fact that the SVD of a Hankel matrix induces tight frame filter banks due to its underlying convolutional structure. More precisely, the first stage \cref{DDTFApproach} learns tight frame filter banks from low frequency samples of a degraded measurement. Using the tight frame filter banks learned from the first stage, the second stage restores a piecewise constant image via a wavelet frame based analysis approach \cref{ProposedModel}. Since we use the data-driven tight frame as an alternative to the structured low-rank matrix, the proposed \cref{DDTFApproach} achieves a faster computation than the conventional SVD based Cadzow denoising \cref{CadzowVS}, with a numerically equivalent filter learning result. In addition, since we adopts the analysis approach, our proposed image restoration model reflects the underlying structured low-rank matrix better than the previous balanced approach in \cite{J.F.Cai2020}. Finally, since we consider the sparsity in the entire transform domain rather than the group sparsity according to the index of filters, our proposed model \cref{ProposedModel} is more robust to the noise/the image modeling.

For the future work, we plan to extend our two-stage approach to the piecewise smooth image restoration. Indeed, since the wavelet frame based sparse regularization model for the piecewise smooth image restoration via the continuous domain regularization is readily available in \cite{J.F.Cai2022} as a byproduct of the structured low-rank matrix framework for the piecewise smooth functions, it suffices to develop a method to learn tight frame filter banks from a given degraded measurement. To do this, we need to establish a sufficient condition for the number of low frequency samples to guarantee the estimation of image singularities from annihilating relations. Unfortunately, since the annihilating relations for piecewise smooth functions in \cite{J.F.Cai2022} involves blind deconvolutions, it is not clear for us at this moment under what conditions we can restore the singularity set from uniform Fourier samples. Nevertheless, this is definitely a future direction we would like to work on, as it is more appropriate to model natural images as piecewise smooth functions (e.g. \cite{J.F.Cai2016}).

To broaden the scope of applications, we also plan to apply the idea in this paper to other image restoration tasks. First of all, since the ideal low-pass filter deconvolution can be viewed as a Shannon's wavelet superresolution problem (e.g. \cite{Mallat2008}), it will not be difficult to apply our approach to the generic wavelet superresolution problem. In addition, even though a careful design on the resampling strategy on the frequency domain \cite{A.Kiperwas2017} is required, the application to the sparse angle CT is another important application of this approach. Finally, motivated by the recent works on the deep learning techniques combined with the structured low-rank matrix approaches \cite{M.Jacob2020,T.H.Kim2019,A.Pramanik2019,A.Pramanik2020}, we can also consider developing a deep learning/a convolutional neural network framework based on the structured low-rank matrix frameworks. In particular, we will focus on developing a convolutional neural network framework from the cascade structure of the structured low-rank matrix framework in \cite{J.F.Cai2020}, designed for the piecewise smooth image restoration via continuous domain regularization.

{\appendix

\section{Preliminaries on tight wavelet frames}\label{PreliminariesTightFrame}

In this section, we provide a brief introduction on tight wavelet frames and the wavelet frame based image restoration. Interested readers may consult \cite{J.F.Cai2008,B.Dong2013,B.Dong2015,A.Ron1997,Shen2010} for detailed surveys on tight wavelet frames, and \cite{J.F.Cai2012} for detailed reviews on the wavelet frame based image restoration. For the sake of simplicity, we only discuss the real valued wavelet tight frame systems, but note that it is not difficult to extend the idea to the complex case.

Let $\msH$ be a Hilbert space equipped with an inner product $\la\cdot,\cdot\ra$. A (Bessel) sequence $\{\bvphi_n:n\in\Z\}\subseteq\msH$ is called a tight frame on $\msH$ if
\begin{align}\label{TightFrame}
\|\bu\|^2=\sum_{n\in\Z}|\la\bu,\bvphi_n\ra|^2~~~~~\text{for all}~~~\bu\in\msH.
\end{align}
Given $\{\bvphi_n:n\in\Z\}\subseteq\msH$, we define the analysis operator $\bmW:\msH\to\ell_2(\Z)$ as
\begin{align*}
\bu\in\msH\mapsto \bmW\bu=\{\la\bu,\bvphi_n\ra:n\in\Z\}\in\ell_2(\Z).
\end{align*}
The synthesis operator $\bmW^T:\ell_2(\Z)\to\msH$ is defined as the adjoint of $\bmW$:
\begin{align*}
\bc\in\ell_2(\Z)\mapsto\bmW^T\bc=\sum_{n\in\Z}\bc(n)\bvphi_n\in\msH.
\end{align*}
Then $\{\bvphi_n:n\in\Z\}$ is a tight frame on $\msH$ if and only if $\bmW^T\bmW=\bmI$. It follows that, for a given tight frame $\{\bvphi_n:n\in\Z\}$, we have the following canonical expression:
\begin{align*}
\bu=\sum_{n\in\Z}\la\bu,\bvphi_n\ra\bvphi_n,
\end{align*}
with $\bmW\bu=\{\la\bu,\bvphi_n\ra:n\in\Z\}$ being called the canonical tight frame coefficients. Hence, the tight frames are extensions of orthonormal bases to the redundant systems. In fact, a tight frame is an orthonormal basis if and only if $\|\bvphi_n\|=1$ for all $n\in\Z$.

One of the most widely used class of tight frames is the discrete wavelet frame generated by a set of finitely supported filters $\{\a_1,\cdots,\a_M\}$. Throughout this paper, we only discuss the two dimensional undecimated wavelet frames on $\ell_2(\Z^2)$ for our purpose. Given $\a\in\ell_1(\Z^2)$, we define a convolution operator $\bmS_{\a}:\ell_2(\Z^2)\to\ell_2(\Z^2)$ by
\begin{align}\label{DiscreteConv}
(\bmS_{\a}\bu)(\bk)=(\a\ast\bu)(\bk)=\sum_{\bsl\in\Z^2}\a(\bk-\bsl)\bu(\bsl)~~~~~\text{for}~~~\bu\in\ell_2(\Z^2).
\end{align}
Given a set of finitely supported filters $\{\a_1,\cdots,\a_M\}$, define the analysis operator $\bmW$ and the synthesis operator $\bmW^T$ by
\begin{align}
\bmW&=\big[\bmS_{\a_1(-\cdot)}^T,\bmS_{\a_2(-\cdot)}^T,\cdots,\bmS_{\a_M(-\cdot)}^T\big]^T,\label{AnalConv2}\\
\bmW^T&=\big[\bmS_{\a_1},\bmS_{\a_2},\cdots,\bmS_{\a_M}\big],\label{SyntConv2}
\end{align}
respectively. Then, the direct computation can show that the rows of $\bmW$ form a tight frame on $\ell_2(\Z^2)$ (i.e. $\bmW^T\bmW=\bmI$) if and only if the filters $\{\a_1,\cdots,\a_m\}$ satisfy
\begin{align}\label{SomeUEP}
\sum_{m=1}^M\sum_{\bsl\in\Z^2}\a_m(\bk+\bsl)\a_m(\bsl)=\dde(\bk)=\left\{\begin{array}{cl}
1~&\text{if}~\bk=\0,\\
0~&\text{if}~\bk\neq\0,
\end{array}\right.
\end{align}
called the {\emph{unitary extension principle}} (UEP) condition \cite{B.Han2011}. Finally, for a two dimensional discrete image on the finite grid, $\bmS_{\a}$ in \cref{DiscreteConv} with a finitely supported $\a$ denotes the discrete convolution under the periodic boundary condition throughout this paper.

For the image restoration, we assume that there exists a tight wavelet frame $\bmW$ defined as \cref{AnalConv2} with filters $\left\{\a_1,\cdots,\a_M\right\}$ such that $\bu$ is sparse under $\bmW$. This leads us to solve
\begin{align}\label{FrameModelGeneric}
\min_{\bc}~\f{1}{2}\|\bmA\bmW^T\bc-\bsf\|_2^2+\f{1}{2}\|(\bmI-\bmW\bmW^T)\bc\|_2^2+\|\gga\cdot\bc\|_1
\end{align}
and restore $\bu=\bmW^T\bc$. In \cref{FrameModelGeneric}, the first term is a data fidelity term, the third $\ell_1$ norm with the parameters $\gga=\left[\gamma_1,\ldots,\gamma_M\right]^T$ and $\gamma_j>0$ promotes the sparsity of $\bc$. The second term penalizes the distance between $\bc$ and its projection onto $\msR(\bmW)$. Since we obtain the image via $\bu=\bmW^T\bc$, the second term in fact forces the sparse coefficient $\bc$ close to the canonical coefficient of the restored image. Since the magnitude/decay of the canonical coefficient reflects the regularity of the image under some mild conditions on the tight frame system $\bmW$ \cite{L.Borup2004}, we can see that the second and third terms promote a balance between the sparsity of coefficient and the regularity of image.

We can further impose the flexibility on \cref{FrameModelGeneric} by introducing the weight on the second term:
\begin{align}\label{Balanced}
\min_{\bc}~\f{1}{2}\|\bmA\bmW^T\bc-\bsf\|_2^2+\f{\mu}{2}\|(\bmI-\bmW\bmW^T)\bc\|_2^2+\|\gga\cdot\bc\|_1,
\end{align}
called the balanced approach \cite{J.F.Cai2008,R.H.Chan2003}. If $\mu=0$, then \cref{Balanced} becomes the synthesis approach \cite{I.Daubechies2007,M.J.Fadili2007,M.A.T.Figueiredo2003}
\begin{align}\label{Synthesis}
\min_{\bc}~\f{1}{2}\|\bmA\bmW^T\bc-\bsf\|_2^2+\|\gga\cdot\bc\|_1,
\end{align}
as we aim to find the sparsest coefficient that synthesizes the image. On the other hand, if $\mu=\infty$, then it must be $(\bmI-\bmW\bmW^T)\bc=\0$. Since we then have $\bc\in\msR(\bmW)$, it follows that $\bc=\bmW\bu$ for some $\bu$. Hence, \cref{FrameModelGeneric} becomes the analysis approach \cite{J.F.Cai2009/10}
\begin{align}\label{Analysis}
\min_{\bu}~\f{1}{2}\|\bmA\bu-\bsf\|_2^2+\|\gga\cdot\bmW\bu\|_1.
\end{align}
Since the analysis approach promotes the sparse canonical coefficient, it emphasizes the regularity of the restored image.

These three approaches are equivalent if and only if $\bmW^T\bmW=\bmW\bmW^T=\bmI$, i.e. $\bmW$ forms an orthonormal basis. However, since we consider the redundant system, each approach will lead to the different restoration results. Notice that the coefficient vector $\bc$ obtained from the synthesis approach \cref{Synthesis} is much sparser than the canonical coefficient $\bmW\bu$ obtained from the analysis approach \cref{Analysis}. It is empirically observed that the synthesis approach \cref{Synthesis} tends to yield some artifacts in the restored image. In contrast, the analysis approach \cref{Analysis} usually have less artifacts as it promotes the regularity of the image along the singularities. Theoretically, it has been proved in \cite{J.F.Cai2012} that, under a suitable choice of parameters, the analysis approach \cref{Analysis} can be seen as sophisticated discretization of a variational model including the Rudin-Osher-Fatemi model \cite{L.I.Rudin1992}, which enables a geometric interpretation on \cref{Analysis}.

\section{Data driven tight frames}\label{DataDrivenTightFrame}

As mentioned in the previous section, wavelet frames are widely used for the sparse approximation of an image. To explore a better sparse approximation, the authors in \cite{J.F.Cai2014} proposed a \emph{data-driven tight frame} approach, inspired by \cite{M.Aharon2006,S.Ravishankar2013}. Specifically, given a two dimensional image $\bu\in\mV$, a tight frame system $\bmW$ defined as in \cref{AnalConv2}, which is generated by a set of finitely supported $p\times p$ filters $\{\a_1,\cdots,\a_{p^2}\}$ supported on a $p\times p$ grid $\KK$ and satisfying \cref{SomeUEP}, is constructed via the following minimization
\begin{align}\label{DDTFModel}
\min_{\bc,\bmW}~\|\bc-\bmW\bu\|_2^2+\gamma^2\|\bc\|_0~~~~~\text{subject to}~~~~\bmW^T\bmW=\bmI,
\end{align}
with the $\ell_0$ norm $\|\bc\|_0$ encoding the number of nonzero entries in the coefficient vector $\bc$.

To solve \cref{DDTFModel}, the authors in \cite{J.F.Cai2014} presented that if $\bsA\in\R^{p^2\times p^2}$ satisfies $\bsA\bsA^T=p^{-2}\bsI$, after reformulating into $p\times p$ filters supported on $\KK$, its column vectors generate filters satisfying the UEP condition \cref{SomeUEP}. Hence, under such an assumption on the filters $\a_1,\cdots,\a_{p^2}$, we can rewrite \cref{DDTFModel} as
\begin{align}\label{DDTFModelModi}
\min_{\bC,\bsA}~\|\bC-\left(\bmH\bu\right)\bsA\|_F^2+\gamma^2\|\bC\|_0~~~~~\text{subject to}~~~~\bsA\bsA^T=p^{-2}\bsI,
\end{align}
to construct filter banks. In \cref{DDTFModelModi}, $\bmH\bu\in\R^{(N-p+1)^2\times p^2}$ is a Hankel matrix generated by the $p\times p$ patches of $\bu$, $\bC$ is the corresponding matrix formulation of $\bc$, and $\left\|\cdot\right\|_F$ is the Frobenius norm of a matrix. To solve \cref{DDTFModelModi}, the alternating minimization method with closed form solutions for each stage is presented in \cite{J.F.Cai2014}. More explicitly, for $n=0$, $1$, $2$, $\cdots$, we have
\begin{align}
\bC_{(n+1)}&=\argmin_{\bC}~\lambda^2\|\bC\|_0+\|\bC-\left(\bmH\bu\right)\bsA_{(n)}\|_F^2=\bmT_{\lambda}\left(\left(\bmH\bu\right)\bsA_{(n)}\right)\label{Csub}\\
\bsA_{(n+1)}&=\argmin_{\bsA\bsA^T=p^{-2}\bsI}~\left\|\left(\bmH\bu\right)\bsA-\bC_{(n+1)}\right\|_F^2=p^{-1}\bU\bV^T\label{Dsub}
\end{align}
where $\bmT_{\lambda}$ is a hard thresholding defined as
\begin{align}\label{HardThresholding}
\bmT_{\lambda}\left(\bC\right)^{(j,k)}=\left\{\begin{array}{cl}
\bC^{(j,k)}~&~\text{if}~|\bC^{(j,k)}|>\lambda\vspace{0.25em}\\
0~&~\text{otherwise},
\end{array}\right.
\end{align}
and $\bU$ and $\bV$ are singular vectors of $\left(\bmH\bu\right)^T\bC_{(n+1)}$; $\left(\bmH\bu\right)^T\bC_{(n+1)}=\bU\bsS\bV^T$. In addition, the following proximal alternating minimization (PAM) scheme
\begin{align}
\bC_{(n+1)}&=\argmin_{\bC}~\lambda^2\|\bC\|_0+\|\bC-\left(\bmH\bu\right)\bsA_{(n)}\|_F^2+\beta_1\left\|\bC-\bC_{(n)}\right\|_F^2\label{CsubPAM}\\
\bsA_{(n+1)}&=\argmin_{\bsA\bsA^T=p^{-2}\bsI}~\left\|\left(\bmH\bu\right)\bsA-\bC_{(n+1)}\right\|_F^2+\beta_2\left\|\bsA-\bsA_{(n)}\right\|_F^2\label{DsubPAM}
\end{align}
is proposed in \cite{C.Bao2015}, together with global convergence to stationary points.

\section*{Acknowledgments} The authors thank Prof. Ke Wei in the School of Data Sciences at Fudan University, au author of \cite{J.F.Cai2022,J.F.Cai2020} for the valuable discussions on this research. The authors also thank Dr. Greg Ongie in the Department of Mathematical and Statistical Sciences at Marquette University, an author of \cite{G.Ongie2018,G.Ongie2015,G.Ongie2016,G.Ongie2017} for making the data sets as well as the MATLAB toolbox available so that the experiments can be implemented.}

\bibliographystyle{aims}
\providecommand{\href}[2]{#2}
\providecommand{\arxiv}[1]{\href{http://arxiv.org/abs/#1}{arXiv:#1}}
\providecommand{\url}[1]{\texttt{#1}}
\providecommand{\urlprefix}{URL }

\end{document}